\documentclass[11pt]{article}
\usepackage{bm}
\usepackage[font=small]{caption}
\RequirePackage{amsthm,amsmath,amsfonts,amssymb}
\usepackage{omega}
\usepackage[a4paper,left=2.2cm,right=2.2cm,top=2.5cm,bottom=2.5cm]{geometry}
\usepackage{booktabs}

\usepackage{bm}
\usepackage{natbib}
\usepackage[ruled]{algorithm2e}
\usepackage[dvipsnames]{xcolor}
\usepackage[colorlinks = true,
            linkcolor = NavyBlue,
            urlcolor  = NavyBlue,
            citecolor = NavyBlue,
            anchorcolor = NavyBlue]{hyperref}

\usepackage{url}
\usepackage[graphicx]{realboxes}
\usepackage{enumerate}
\usepackage{comment}
\usepackage{subfigure}
\usepackage[title]{appendix}
\usepackage{chngcntr}
\RequirePackage{subfigure}
\RequirePackage{algorithmic}

\author{}

\numberwithin{equation}{section}

\theoremstyle{plain}
\newtheorem{remark1}{Remark}
\newtheorem{theorem1}{Theorem}
\newtheorem{lemma1}{Lemma}
\newtheorem{definition1}{Definition}
\newtheorem{proposition1}{Proposition}

\usepackage{setspace}

\newcommand{\ee}{\end{aligned} \end{equation}}
\newcommand{\eq}{\end{quote}}

\newcommand{\ep}{\end{parts}}

\newcommand{\bqp}{\begin{quote}\begin{parts}}

\newcommand{\epq}{\end{parts}\end{quote}}

\newcommand{\Rom}[1]{\text{\uppercase\expandafter{\romannumeral #1\relax}}}
\newcommand{\bee}{\begin{equation}\begin{aligned}}

\newcommand{\emm}{\end{bmatrix}}

\numberwithin{equation}{section}
\newcommand{\vertiii}[1]{{\vert\kern-0.25ex\vert\kern-0.25ex\vert #1 
    \vert\kern-0.25ex\vert\kern-0.25ex\vert}}

\newcommand{\bq}{\begin{quote}}
\newcommand{\bp}{\begin{parts}}

\newcommand\piu{\overline{\bm{\Pi}}_{\mathbf{U}}}
\newcommand\ma{\mathbf{A}}
\newcommand\mb{\mathbf{B}}
\newcommand\mc{\mathbf{C}}
\newcommand\md{\mathbf{D}}
\newcommand\me{\mathbf{E}}
\newcommand\mh{\mathbf{H}}
\newcommand\mi{\mathbf{I}}
\newcommand\mq{\mathbf{Q}}

\newcommand\mr{\mathbf{R}}
\newcommand\ms{\mathbf{S}}

\newcommand\mv{\mathbf{V}}
\newcommand\mx{\mathbf{X}}
\newcommand\mz{\mathbf{Z}}
\newcommand\mw{\mathbf{W}}
\newcommand\mm{\mathbf{M}}
\newcommand\ml{\mathbf{L}}
\newcommand\mt{\mathbf{T}}

\newcommand\mpp{\mathbf{P}}
\newcommand\mo{\mathbf{O}}
\newcommand\my{\mathbf{Y}}
\newcommand\muu{\mathbf{U}}
\newcommand\mSigma{\bm{\Sigma}}
\newcommand\mOmega{\bm{\Sigma}}
\newcommand\mUpsilon{\bm{\Sigma}}
\newcommand\mLambda{\bm{\Lambda}}

\newcommand{\tabincell}[2]{\begin{tabular}{@{}#1@{}}#2\end{tabular}}

\newcommand{\bea}{\begin{eqnarray}}
	\newcommand{\eea}{\end{eqnarray}}
\newcommand{\beas}{\begin{eqnarray*}}
	\newcommand{\eeas}{\end{eqnarray*}}

\title{Limit results for distributed estimation of invariant subspaces in multiple networks inference and PCA}
\begin{document}

\author{Runbing Zheng\\{Department of Applied Mathematics and Statistics, Johns
  Hopkins University}\\
  Minh Tang \\{Department of Statistics,
North Carolina State University}
  }

\date{}

\maketitle

\begin{abstract}
Several statistical problems, such as multiple
heterogeneous graph analysis, distributed PCA, integrative data
analysis, and simultaneous dimension reduction of images, can involve
a collection of $m$ matrices whose leading subspaces $\muu^{(i)}$
consist of a shared subspace $\muu_c$ and individual subspaces
$\muu_s^{(i)}$. We consider a distributed estimation procedure that
first obtains $\hat\muu^{(i)}$ as the leading singular vectors for
each observed noisy matrix, then computes the leading left singular
vectors of the concatenated matrix
$[\hat\muu^{(1)}|\hat\muu^{(2)}|\dots|\hat\muu^{(m)}]$ as
$\hat\muu_c$, and finally computes the leading singular vectors of the
projection of each $\hat\muu^{(i)}$ onto the orthogonal complement of
$\hat\muu_c$ as $\hat\muu_s^{(i)}$.  In this paper, we provide a
framework for deriving limit results for such distributed estimation
procedures, including expansions of estimation errors in both common
and individual subspaces and their asymptotically normal
approximations.  We apply this framework specifically to (1) parameter
estimation for multiple heterogeneous random graphs
with shared subspaces, and (2) distributed PCA for independent
sub-Gaussian random vectors with spiked covariance
structures. Leveraging these results, we also consider a two-sample
test for the null hypothesis that a pair of random graphs have the
same edge probabilities, and present a test statistic whose limiting
distribution converges to a central (resp., non-central) $\chi^2$
distribution under the null (resp., local alternative) hypothesis.

\end{abstract}

\noindent%
\vfill

\vfill

\begin{sloppypar}

\section{Introduction}
\label{sec:intro}
Distributed estimation, also known as divide-and-conquer or aggregated
inference, is used in numerous methodological applications including
regression \citep{aggregated_inference_huo,ridge_regression_dobriban},
integrative data analysis \citep{jive,feng2018angle,hector}, multiple
network inference \citep{arroyo2019inference}, signal processing \cite{lin_scaglione2,jiao_gu,roy_vetterli,bertrand_moonen}, distributed PCA and
image population analysis
\citep{population_svd,sagonas2017robust,tang2021integrated,fan2019distributed,chen2021distributed,distributedPCA_review},
and is also a key component underlying federated learning
\citep{federated_survey}. Such procedures are particularly
important for analyzing large-scale datasets that are scattered across
multiple organizations and/or computing nodes where both the
computational complexities and communication costs (as well as
possibly privacy constraints) prevent the transfer of all the raw data
to a single location.

In this paper, we focus on distributed estimation for a collection of matrices with a shared subspace $\muu_c$ and potentially distinct individual subspaces $\muu_s^{(i)}$. We consider an algorithm that first obtains $\hat\muu^{(i)}$ as the leading singular vectors for each matrix, then integrates $\hat\muu^{(i)}$ across all matrices to obtain the estimated common subspace $\hat\muu_c$, and finally projects each $\hat\muu^{(i)}$ onto the orthogonal complement of $\hat\muu_c$ and computes its leading subspace as $\hat\muu_s^{(i)}$.


One widely studied example of such a problem is distributed PCA, in
which there are $N$ independent $D$-dimensional sub-Gaussian random
vectors $\{X_j\}_{j=1}^{N}$ with common covariance matrix $\bm{\Sigma}$ scattered across $m$ computing nodes, and the goal is to
find the leading eigenspace $\muu$ of $\bm{\Sigma}$. 
Letting $\mx^{(i)}$ be the $D \times
n_i$ matrix whose columns are the subsample of $\{X_j\}_{j=1}^{N}$
stored in node $i$, \cite{fan2019distributed} analyzes a procedure
where each node $i$ first computes the $D \times d$ matrix
$\hat{\muu}^{(i)}$ whose columns are the leading left singular vectors
of $\mx^{(i)}$. These $\hat{\muu}^{(i)}$ are then sent to a central
computing node which outputs the leading left singular vectors of
$[\hat\muu^{(1)}|\hat\muu^{(2)}|\dots|\hat\muu^{(m)}]$ as
$\hat\muu$. This algorithm, which is also widely used in signal
processing (see e.g., \cite{bertrand_moonen,distributedPCA_review} and
the references therein), is essentially the version of our
aforementioned algorithm when $\muu^{(i)}\equiv\muu_c=\muu$.
Another example of multiple matrices with common subspaces is simultaneous
dimension reduction of high-dimensional images
$\{\mathbf{Y}_i\}_{i=1}^{m}$, namely each $\mathbf{Y}_i$ is an $F
\times T$ matrix whose entries are measurements recorded for various
frequencies and various times, and the goal is to find a ``population
value decomposition'' of each $\mathbf{Y}_i$ as $\mathbf{Y}_i \approx
\mathbf{P} \mathbf{V}_i \mathbf{D}$. Here $\mathbf{P}$ and
$\mathbf{D}$ are $F \times A$ and $A \times T$ matrices (with $A \ll
\min\{F,T\}$) representing \emph{population} frames of reference, and
$\{\mathbf{V}_i\}$ are the \emph{subject-level} features; see
\cite{population_svd} for more details.  An example that includes both
common subspaces and individual subspaces is heterogeneous multiple
directed networks with probability matrices
$\mpp^{(i)}=\muu^{(i)}\mr^{(i)}\mv^{(i)\top}$, where
$\muu^{(i)}=[\muu_c|\muu_s^{(i)}]$ and $\mv^{(i)}=[\mv_c|\mv_s^{(i)}]$
contain common and possibly distinct individual left and right
subspaces for the networks, and $\mr^{(i)}$ are low-dimensional
matrices that are heterogeneous across networks. This setup
includes the widely-used COSIE model
\citep{arroyo2019inference} for multiple networks where $\muu^{(i)}\equiv\muu_c$ and
$\mv^{(i)}\equiv\mv_c$, and the estimation procedure proposed in
\cite{arroyo2019inference} is also a version of our aforementioned
algorithm.  As a final example, a typical setting for integrative data
analysis assumes that there is a collection of data matrices
$\{\mathbf{X}^{(i)}\}$ from multiple disparate sources and the goal is
to decompose each $\mathbf{X}^{(i)}$ as $\mathbf{X}^{(i)} =
\mathbf{J}^{(i)} + \mathbf{I}^{(i)} + \mathbf{N}^{(i)}$, where
$\{\mathbf{J}^{(i)}\}$ share a common row space $\mathbf{J}_*$ which
captures the \emph{joint} structure among all $\{\mathbf{X}^{(i)}\}$,
$\mathbf{I}^{(i)}$ represent the \emph{individual} structure in each
$\mathbf{X}^{(i)}$, and $\mathbf{N}^{(i)}$ are noise matrices. Several
algorithms, such as aJIVE and robust aJIVE
\citep{feng2018angle,ponzi2021rajive}, compute the estimate
$\hat{\mathbf{J}}_*$ by aggregating the leading (right) singular
vectors $\hat{\muu}^{(i)}$ of $\mathbf{X}^{(i)}$ and then estimate
each individual $\mathbf{I}^{(i)}$ by projecting $\mathbf{X}^{(i)}$
onto the orthogonal complement of $\hat{\mathbf{J}}_*$, and are 
thus equivalent to our aforementioned algorithm.


Despite the wide applicability of distributed estimators for
matrices with common subspaces such as those described above, their theoretical
results are still somewhat limited. For example, the papers that
proposed the aJIVE/rAJIVE procedures
\citep{feng2018angle,ponzi2021rajive} and the PVD
\citep{population_svd} do not consider any specific noise models and
thus do not present explicit error bounds for the
estimates. Similarly, in the context of the COSIE model and
distributed PCA, \cite{arroyo2019inference} and
\cite{population_svd,tang2021integrated,fan2019distributed,chen2021distributed}
only provide Frobenius norm upper bounds between $\hat{\mathbf{U}}$
and $\mathbf{U}$.

In this paper, we provide a general framework for analyzing these types of estimators, with special emphasis on uniform $\ell_{2 \to \infty}$ error bounds and normal approximations for the row-wise fluctuations of $\hat{\mathbf{U}}_c$ and $\hat\muu_s^{(i)}$ around $\mathbf{U}_c$ and $\muu_s^{(i)}$, respectively. This framework is based on the following result (see Section~\ref{sec:notations} for a description of the notation used here), which is also a key contribution of our paper.
\begin{theorem1}
  \label{thm3_general}
  Let $\{\muu^{(i)} = [\muu_c \,|\, \muu_s^{(i)}]\}_{i=1}^{m}$ be a
collection of $n \times d_i$ orthonormal matrices, where $\muu_c$
represents the set of $d_0$ columns shared across all $\muu^{(i)}$,
and $\muu_s^{(i)}$ denotes the set of $(d_i - d_0)$ columns
specific to each $\mathbf{U}^{(i)}$.
Denote $\bm{\Pi}_{s} = \frac{1}{m} \sum_{i=1}^m \muu_s^{(i)} \muu_s^{(i)\top}$.
For each $i \in [m]$, suppose that we have an estimate $\hat{\muu}^{(i)}$ of $\muu^{(i)}$ such that
\[
\hat{\muu}^{(i)} \mw_\muu^{(i)} - \muu^{(i)} = \mt_0^{(i)} + \mt^{(i)}
\]
for some orthogonal matrix $\mw_\muu^{(i)}$, where $\mt_0^{(i)}$ and $\mt^{(i)}$ satisfy
\begin{gather}
  \label{eq:cond_pis}
  \max_{i \in [m]} \Bigl(2 \|\mt_0^{(i)}\| + 2 \|\mt^{(i)}\| + 
   \|\mt_0^{(i)} + \mt^{(i)}\|^2\Bigr)
  \leq c(1 - \|\bm{\Pi}_{s}\|)
\end{gather}
for some constant $c < \frac{1}{2}$.
Define the quantities
\begin{equation}
  \label{eq:quantities_main_thm3}
  \begin{aligned}
  	 &\zeta_{\muu} = \max_{i\in[m]} \|\muu^{(i)}\|_{2 \to \infty}, \quad \epsilon_{\star} = \max_{i,j\in[m]}
 \|\muu^{(i)\top} \mt_0^{(j)}\|, \\
& \epsilon_{\mt_0} = \max_{i\in[m]} \|\mt_0^{(i)}\|, \quad \zeta_{\mt_0} =
    \max_{i\in[m]} \|\mt_0^{(i)}\|_{2 \to \infty}, \quad
    \epsilon_{\mt} = \max_{i\in[m]} \|\mt^{(i)}\|, \quad \zeta_{\mt} =
    \max_{i\in[m]} \|\mt^{(i)}\|_{2 \to \infty}.
  \end{aligned}
\end{equation}

Now let $\hat{\muu}_c$ denote the matrix whose columns are the $d_{0}$
leading eigenvectors of $m^{-1} \sum_{i=1}^m \hat{\muu}^{(i)}
\hat{\muu}^{(i)\top}$. Let $\mathbf{W}_{\mathbf{U}_c}$ be the minimizer of $\|\hat{\muu}_c \mathbf{O} -
\muu_c\|_{F}$ over all orthogonal
matrices $\mathbf{O}$. We then have
\begin{equation}
  \label{eq:expansion_hatU_seperate}
  \hat{\mathbf{U}}_c \mathbf{W}_{\mathbf{U}_c} -\mathbf{U}_c =
  \frac{1}{m}\sum_{i=1}^m \mt_0^{(i)} \muu^{(i)\top} \muu_c+\mq_{\muu_c},
\end{equation}
where $\mq_{\muu_c}$ is a matrix satisfying
\begin{equation}
\begin{aligned}
  \label{eq:bound_mq_uc_thm3}
	&\|\mq_{\muu_c}\|\lesssim \epsilon_{\star} + \epsilon_{\mt_0}^2
+ \epsilon_{\mt}, \\
&\|\mq_{\muu_c}\|_{2\to\infty}\lesssim \zeta_\muu(\epsilon_{\star}+\epsilon_{\mt_0}^2+\epsilon_{\mt})
       +\zeta_{\mt_0}(\epsilon_{\star}+\epsilon_{\mt_0}+\epsilon_{\mt})
       +\zeta_{\mt}. 
\end{aligned}
\end{equation}

Given $\hat{\muu}_c$, let $\hat{\muu}_s^{(i)}$ be the matrix whose columns are the $(d_i - d_{0})$ leading left singular vectors of $(\mi - \hat{\muu}_c \hat{\muu}_c^{\top}) \hat{\muu}^{(i)}$. For any $i \in [m]$, let $\mw^{(i)}_{\muu_s}$ be the minimizer of $\|\hat{\muu}_s^{(i)} \mo - \muu_s^{(i)}\|_{F}$ over all orthogonal matrices $\mo$. We then have
\begin{equation}
  \label{eq:expansion_Us}
  \hat{\muu}_s^{(i)} \mw^{(i)}_{\muu_s} - \muu_s^{(i)} = \mt_0^{(i)} \muu^{(i)\top} \muu_s^{(i)} + \mq_{\muu_s}^{(i)},
\end{equation}
where $\mq_{\muu_s}^{(i)}$ is a matrix satisfying the same upper bounds as those for $\mq_{\muu_c}$.
\end{theorem1} Theorem~\ref{thm3_general} is a deterministic matrix
perturbation bound and provides expansions for $\hat{\muu}_c$ and
$\hat{\muu}_s^{(i)}$ in terms of the expansions for the individual
$\hat{\muu}^{(i)}$. The upper bounds for $\mq_{\muu_c}$ and
$\mq_{\muu_s}^{(i)}$ in Theorem~\ref{thm3_general} depend only on
$\epsilon_\star$, $\epsilon_{\mt_0}$, $\zeta_{\mt_0}$,
$\epsilon_{\mt}$, $\zeta_{\mt}$, and $\zeta_{\muu}$, and the bounds
for these quantities can be derived in various settings.

{\color{black}
To illustrate the significance of these theoretical results with a potential concrete real-world example, consider the integrative analysis of gene expression data collected from $m$ different studies, where each study $i$ measures the expression levels of the same large collection of genes across $n_i$ samples. The leading singular vectors $\hat{\muu}^{(i)}$ of each dataset capture the dominant patterns of co-expression among the genes in study $i$. Some of these patterns are shared across all studies, reflecting biological pathways active in every cohort, and span the common subspace $\muu_c$; others are specific to a single study, reflecting study-specific effects, and span the individual subspace $\muu_s^{(i)}$. Here each row of $\muu_c$ (and of $\muu_s^{(i)}$) corresponds to a gene, and its entries characterize the gene's participation in the shared (respectively, study-specific) co-expression patterns. Consider the following estimation procedure: (i) each study locally computes $\hat{\muu}^{(i)}$ and shares only this low-dimensional matrix rather than the raw expression data (which is especially useful when the data cannot be pooled due to privacy or ownership constraints); (ii) the shared subspace $\hat{\muu}_c$ is estimated as the leading eigenvectors of $m^{-1}\sum_{i=1}^m \hat{\muu}^{(i)}\hat{\muu}^{(i)\top}$; and (iii) each study recovers its individual subspace $\hat{\muu}_s^{(i)}$ as the leading left singular vectors of $(\mi - \hat{\muu}_c \hat{\muu}_c^{\top}) \hat{\muu}^{(i)}$.
The decomposition of $\hat{\muu}^{(i)}$ required by Theorem~\ref{thm3_general} is available in a wide range of settings, and Theorem~\ref{thm3_general} then yields the following guarantees. The $\ell_{2\to\infty}$ bounds in Theorem~\ref{thm3_general} first guarantee that the estimated common subspace $\hat{\muu}_c$ is accurate for every gene, rather than only on average as controlled by the Frobenius error, which is what one needs, for example, to reliably identify which genes contribute to a shared pathway; the analogous bounds for $\hat{\muu}_s^{(i)}$ similarly guarantee accurate recovery of the study-specific patterns at every gene. Moreover, the row-wise normal approximations based on the $\ell_{2\to\infty}$ bounds (as we establish in the specific settings considered later in this paper) provide a limiting distribution for each gene, allowing one, for example, to test whether a particular gene is significantly involved in the shared structure, or to detect genes that behave anomalously in a specific study. We note that such gene-level inference cannot be obtained from Frobenius-norm error bounds alone. 
}

In this paper, based on the proposed Theorem~\ref{thm3_general}, we
specifically analyze two problems: inference for heterogeneous multiple networks and
distributed PCA, as these problems have been widely studied and yet our
results are still novel. Specifically, our model for
heterogeneous multiple networks is a natural extension
to the COSIE model in \cite{arroyo2019inference} and also encompasses other existing models such as the
MultiNeSS model \citep{macdonald2022latent} and multilayer SBMs
\citep{holland1983stochastic}. Furthermore, while existing results for these models
\citep{arroyo2019inference,macdonald2022latent,paul2020spectral,jing2021community,lei2020bias}
primarily focus on spectral or Frobenius norm error bounds (with \cite{arroyo2019inference} also
providing row-wise upper error bounds), we provide 
limiting distributions for the row-wise fluctuations of $\hat{\muu}_c$ and $\hat{\muu}_s^{(i)}$, as well as normal approximations
for $\mr^{(i)}$. 
Similarly, for distributed PCA, existing works
\citep{chen2021distributed,charisopoulos2021communication,fan2019distributed,liang2014improved}
also focus on spectral or Frobenius norm error bounds for $\hat{\muu}_c$ instead of the more refined row-wise fluctuations presented here. 
A detailed comparison between our results and existing works is provided
in Sections~\ref{sec:related_works} and \ref{sec:related_works_pca}.

The structure of our paper is as follows. In Section~\ref{sec:COISIE},
we study the heterogeneous multiple networks model with probability
matrices $\mpp^{(i)}=\muu^{(i)}\mr^{(i)}\mv^{(i)\top}$, where
$\muu^{(i)}=[\muu_c|\muu_s^{(i)}]$ and
$\mv^{(i)}=[\mv_c|\mv_s^{(i)}]$. We show that the rows of the
estimates $\hat{\mathbf{U}}_c$, $\hat{\mathbf{U}}_s^{(i)}$,
$\hat{\mathbf{V}}_c$, $\hat{\mathbf{V}}_s^{(i)}$ obtained from the
observed adjacency matrices $\{\ma^{(i)}\}$ are normally distributed
around the rows of their true counterparts. Furthermore, we consider
the COSIE model in \cite{arroyo2019inference} as a special case with $\muu^{(i)}\equiv \muu_c=\muu$,
$\mv^{(i)}\equiv \mv_c=\mv$, and prove that
$\hat{\mathbf{R}}^{(i)}=\hat\muu^\top\ma^{(i)}\hat\mv$ also converges
to a multivariate normal distribution centered around
$\mathbf{R}^{(i)}$ for any $i \in [m]$. We then consider two-sample (and multi-sample)
testing for the null hypothesis that some networks from the COSIE
model have the same probability matrix. Leveraging the theoretical
results for $\{\hat{\mr}^{(i)}\}$, we derive a test statistic whose
limiting distribution converges to a central $\chi^2$
(resp. non-central $\chi^2$) under the null (resp. local alternative)
hypothesis.  In Section~\ref{sec:distributed_PCA}, we study the
distributed PCA setting and derive normal approximations for the rows
of the leading principal components when the data exhibit a spiked
covariance structure.  Numerical simulations and experiments on real
data are presented in Section~\ref{sec:simulation and real data}.
Detailed proofs of all stated results are presented in the
supplementary material.

\subsection{Notation}
\label{sec:notations}
We summarize some notation used in this paper. 
We denote by $\mathcal{O}_d$ the set of $d\times d$ orthogonal matrices, and by $\mathcal{O}_{n\times d}$ the set of $n\times d$ matrices with orthonormal columns.
For a positive integer $p$, we denote by $[p]$ the set $\{1,\dots,p\}$. For two non-negative sequences $\{a_n\}_{n \geq 1}$ and $\{b_n\}_{n \geq 1}$, we write $a_n \lesssim b_n$ (resp. $a_n \gtrsim b_n$) if there exists some constant $C>0$ such that $a_n \leq C b_n$ (resp. $a_n \geq C b_n$) for all $n \geq 1$, and we write $a_n \asymp b_n$ if $a_n\lesssim b_n$ and $a_n\gtrsim b_n$.
The notation $a_n \ll b_n$ (resp. $a_n \gg b_n$) means that there exists some sufficiently small (resp. large) constant $C>0$ such that $a_n \leq Cb_n$ (resp. $a_n \geq Cb_n$).
If $a_n/b_n$ stays bounded away from $+\infty$, we write $a_n=O(b_n)$ and $b_n=\Omega(a_n)$, and we use the notation $a_n=\Theta(b_n)$ to indicate that $a_n=O(b_n)$ and $a_n=\Omega(b_n)$.
If $a_n/b_n\to 0$, we write $a_n=o(b_n)$ and $b_n=\omega(a_n)$.
We say a sequence of events $\mathcal{A}_n$ holds with high probability if for any $c > 0$, there exists a finite constant $n_0$ depending only on $c$ such that $\mathbb{P}(\mathcal{A}_n)\geq 1-n^{-c}$ for all $n \geq n_0$.
We write $a_n = O_p(b_n)$ (resp. $a_n = o_p(b_n)$) to denote that $a_n = O(b_n)$ (resp. $a_n = o(b_n)$) holds with high probability.
Given a matrix $\mathbf{M}$, we denote its spectral, Frobenius, and infinity norms by $\|\mathbf{M}\|$, $\|\mathbf{M}\|_{F}$, and $\|\mathbf{M}\|_{\infty}$, respectively.
We also denote the maximum entry (in modulus) of $\mathbf{M}$ by $\|\mathbf{M}\|_{\max}$ and the $2 \to \infty$ norm of $\mathbf{M}$ by
$$\|\mathbf{M}\|_{2 \to \infty} = \max_{\|\bm{x}\| = 1} \|\mathbf{M} \bm{x}\|_{\infty} = \max_{i} \|m_i\|,$$
where $m_i$ denotes the $i$-th row of $\mathbf{M}$, i.e., $\|\mathbf{M}\|_{2 \to \infty}$ is the maximum of the $\ell_2$ norms of the rows of $\mathbf{M}$. We note that the $2 \to \infty$ norm is \emph{not} sub-multiplicative. However, for any matrices $\mathbf{M}$ and $\mathbf{N}$ of conformable dimensions, we have
\begin{equation*}
 \|\mathbf{M} \mathbf{N} \|_{2 \to \infty} \leq
\min\{\|\mathbf{M}\|_{2 \to \infty} \times \|\mathbf{N}\|,
\|\mathbf{M}\|_{\infty} \times \|\mathbf{N}\|_{2\to\infty}\};
\end{equation*}
see Proposition~6.5 in \cite{cape2019two}. Perturbation bounds using
the $2 \to \infty$ norm for the eigenvectors and/or singular vectors
of a noisily observed matrix have recently attracted significant interest from the
statistics community; see
\cite{chen2021spectral,cape2019two,lei2019unified,damle,fan2018eigenvector,abbe2020entrywise}
and the references therein.

\section{Multiple Heterogeneous Networks with Common and Individual Subspaces}
\label{sec:COISIE}

Inference for multiple networks is an important and nascent research area with applications across diverse scientific fields, including neuroscience \citep{bullmore2009complex, battiston2017multilayer, de2017multilayer, kong2021multiplex}, economics \citep{schweitzer2009economic, lee2016strength}, and social sciences \citep{papalexakis2013more, greene2013producing}. 
Multiple networks with shared vertices typically assume that the networks share a common structure.
One prominent example is the multilayer stochastic block model (SBM) \citep{holland1983stochastic, han2015consistent, paul2020spectral, lei2023bias, lei2024computational}, which assumes that vertices share common community assignments across different layers while allowing for layer-specific block probabilities.

Other examples include multilayer eigenscaling models
\citep{nielsen2018multiple, wang2019joint, draves2020bias,
weylandt2022multivariate} and the common subspace independent edge
(COSIE) model \citep{arroyo2019inference}. In particular, the COSIE
model for directed networks $\{\mathcal{G}_i\}_{i=1}^{m}$ assumes that
each $\mathcal{G}_i$ is an edge-independent random graph on the same
set of $n$ vertices where the edge probabilities are given by
$\mathbf{P}^{(i)} = \mathbf{U} \mathbf{R}^{(i)}
\mathbf{V}^{\top}$. Here, $\mathbf{U}, \mathbf{V} \in
\mathcal{O}_{n\times d}$ represent the common subspaces, and the $d
\times d$ matrices $\{\mathbf{R}^{(i)}\}$ capture the heterogeneity
across networks. The COSIE model is quite flexible and encompasses
many popular multiple network models, including the multilayer SBM and
multilayer eigenscaling models mentioned above.

In this paper, we consider the following extension of the COSIE model in which the $\{\mathbf{P}^{(i)}\}$
share some common invariant subspaces $\mathbf{U}_c$ and $\mathbf{V}_c$,
while also allowing for distinct subspaces $\{\mathbf{U}_{s}^{(i)}, \mathbf{V}_{s}^{(i)}\}$ that are specific to each network.
\begin{definition1}[Common and individual subspaces independent edge graphs (COISIE)] 
\label{COISIEgraph} 
For each $i \in [m]$, let $\mathbf{R}^{(i)}$ be a $d_i \times d_i$ matrix, and let $\mathbf{U}^{(i)} = [\mathbf{U}_c \mid \mathbf{U}_s^{(i)}]$ and $\mathbf{V}^{(i)} = [\mathbf{V}_c \mid \mathbf{V}_s^{(i)}]$ be $n \times d_i$ orthonormal matrices. Here, $\mathbf{U}_c \in \mathcal{O}_{n \times d_{0,\muu}}$ and $\mathbf{V}_c \in \mathcal{O}_{n \times d_{0,\mv}}$ represent the shared subspaces across all $i$, while $\mathbf{U}_s^{(i)} \in \mathcal{O}_{n \times (d_i - d_{0,\muu})}$ and $\mathbf{V}_s^{(i)} \in \mathcal{O}_{n \times (d_i - d_{0,\mv})}$ are possibly different between $i$.
Suppose that $u_k^{(i)\top} \mathbf{R}^{(i)} v_\ell^{(i)} \in [0,1]$ for all $k,\ell \in [n]$ and $i \in [m]$,  
where $u_k^{(i)}$ and $v_{\ell}^{(i)}$ denote the $k$th and $\ell$th rows of $\mathbf{U}^{(i)}$ and $\mathbf{V}^{(i)}$, respectively.
We say that the random adjacency matrices $\{\mathbf{A}^{(i)}\}_{i=1}^m$ are jointly distributed according to the common and individual subspaces independent edge graphs model with $\mathbf{U}_c$, $\mathbf{V}_c$, $\{\mathbf{U}_s^{(i)}, \mathbf{V}_s^{(i)}, \mathbf{R}^{(i)}\}_{i=1}^m$, if, for each $i \in [m]$, $\mathbf{A}^{(i)}$ is an $n \times n$ random matrix whose entries $\{\mathbf{A}_{k\ell}^{(i)}\}$ are independent Bernoulli random variables with  
$
\mathbb{P}[\mathbf{A}_{k \ell}^{(i)} = 1] = u_k^{(i)\top} \mathbf{R}^{(i)} v_{\ell}^{(i)}.
$
In other words,  
$$
\mathbb{P}\big(\mathbf{A}^{(i)} \mid \mathbf{U}_c, \mathbf{V}_c, \mathbf{U}_s^{(i)}, \mathbf{V}_s^{(i)}, \mathbf{R}^{(i)}\big) 
= \prod_{k\in[n]} \prod_{\ell \in[n]} \big(u_k^{(i)\top} \mathbf{R}^{(i)} v_{\ell}^{(i)}\big)^{\mathbf{A}_{k \ell}^{(i)}} \big(1 - u_k^{(i)\top} \mathbf{R}^{(i)} v_{\ell}^{(i)}\big)^{1-\mathbf{A}_{k \ell}^{(i)}}.
$$ 
We denote the multiple networks by $\left(\mathbf{A}^{(1)}, \ldots, \mathbf{A}^{(m)}\right) \sim \operatorname{COISIE}(\mathbf{U}_c,\mathbf{V}_c,\{\mathbf{U}_s^{(i)}, \mathbf{V}_s^{(i)}, \mathbf{R}^{(i)}\}_{i=1}^m)$, and write  
\begin{equation}
    \label{eq:shared_new}
    \mathbf{P}^{(i)} = \mathbf{U}^{(i)} \mathbf{R}^{(i)} \mathbf{V}^{(i)\top} = [\mathbf{U}_c \mid \mathbf{U}_{s}^{(i)}] \mathbf{R}^{(i)} [\mathbf{V}_c \mid \mathbf{V}_{s}^{(i)}]^\top
  \end{equation}
to represent the (unobserved) edge probabilities matrix for each $\mathbf{A}^{(i)}$.
\end{definition1}

Note that the dimensions $d_i$ can vary between networks, and the
number of columns in $\mathbf{V}_c$ (denoted by $d_{0,\mathbf{V}}$)
can differ from that in $\mathbf{U}_c$ (denoted by
$d_{0,\mathbf{U}}$).  Moreover, $d_{0,\mathbf{U}}$ or
$d_{0,\mathbf{V}}$ (or both) can be zero, allowing networks to share a common
left subspace $\mathbf{U}_c$ while maintaining distinct subspaces
$\{\mathbf{V}^{(i)}\}$, or vice versa.

The definition presented here is written for directed networks. For undirected networks, we simply require $\mathbf{U}_c = \mathbf{V}_c$, $\mathbf{U}_s^{(i)} = \mathbf{V}_s^{(i)}$, and enforce $\mathbf{R}^{(i)}$ and $\mathbf{A}^{(i)}$ to be symmetric. 
Our subsequent theoretical results, although stated for directed graphs, remain valid for the undirected COISIE model after accounting for the symmetry; see Remark~\ref{rm:undirect} and Remark~\ref{rm:undirect2} for further details.
The COISIE model is also equivalent to a version of the MultiNeSS model \citep{macdonald2022latent} which assumes
\[
\mathbf{P}^{(i)} = \mathbf{X}_c \mathbf{I}_{p_0,q_0} \mathbf{X}_c^{\top} + \mathbf{X}_s^{(i)} \mathbf{I}_{p_i,q_i} \mathbf{X}_s^{(i)\top}.
\]
Here, $\mathbf{I}_{r_+,r_-} = \text{diag}(\mathbf{I}_{r_+}, -\mathbf{I}_{r_-})$ is a diagonal matrix with $r_+$ entries of $+1$ and $r_-$ entries of $-1$ on the diagonal.

We emphasize that $\{\mathbf{A}^{(i)}\}$ are not necessarily independent in the statement of Definition~\ref{COISIEgraph}. While the assumption that $\{\mathbf{A}^{(i)}\}$ are mutually independent appears extensively in the literature (see, for example, the COSIE model \citep{arroyo2019inference}, the multilayer random dot product graph model \citep{jones2020multilayer}, multilayer SBMs \citep{han2015consistent, tang2009clustering, paul2016consistent, lei2020bias, paul2020spectral}, and the MultiNeSS model \citep{macdonald2022latent}), this assumption is either unnecessary or can be relaxed for the theoretical results presented in this paper. See Remark~\ref{rem:independence} for further details.

Given $\left(\mathbf{A}^{(1)}, \ldots, \mathbf{A}^{(m)}\right) \sim \operatorname{COISIE}(\mathbf{U}_c, \mathbf{V}_c, \{\mathbf{U}_s^{(i)}, \mathbf{V}_s^{(i)}, \mathbf{R}^{(i)}\}_{i=1}^m)$, we estimate the parameters using Algorithm~\ref{Alg_COISIE} below. 
\begin{algorithm}[htbp!]
\caption{Estimation of COISIE parameters}	
\label{Alg_COISIE}
\begin{algorithmic}
\REQUIRE Adjacency matrices $\mathbf{A}^{(1)}, \dots, \mathbf{A}^{(m)}$, embedding dimensions $d_1, \dots, d_m$, and common dimensions $d_{0,\mathbf{U}}, d_{0,\mathbf{V}}$. 
\begin{enumerate}
	\item For each $i \in [m]$, obtain $\hat{\mathbf{U}}^{(i)}$ and $\hat{\mathbf{V}}^{(i)}$ as the $n \times d_i$ matrices whose columns are the $d_i$ leading left and right singular vectors of $\mathbf{A}^{(i)}$, respectively.
	
	\item Compute $\hat{\mathbf{U}}_c$ as the $n \times d_{0,\mathbf{U}}$ matrix whose columns are the leading left singular vectors of $[\hat{\mathbf{U}}^{(1)} \mid \cdots \mid \hat{\mathbf{U}}^{(m)}]$, and compute $\hat{\mathbf{V}}_c$ as the $n \times d_{0,\mathbf{V}}$ matrix whose columns are the leading left singular vectors of $[\hat{\mathbf{V}}^{(1)} \mid \cdots \mid \hat{\mathbf{V}}^{(m)}]$.

    \item For each $i \in [m]$, compute $\hat{\mathbf{U}}_s^{(i)}$ as the $n \times (d_i - d_{0,\mathbf{U}})$ matrix whose columns are the leading left singular vectors of $(\mathbf{I} - \hat{\mathbf{U}}_c \hat{\mathbf{U}}_c^\top) \hat{\mathbf{U}}^{(i)}$, and compute $\hat{\mathbf{V}}_s^{(i)}$ as the $n \times (d_i - d_{0,\mathbf{V}})$ matrix whose columns are the leading left singular vectors of $(\mathbf{I} - \hat{\mathbf{V}}_c \hat{\mathbf{V}}_c^\top) \hat{\mathbf{V}}^{(i)}$.
	
	\item For each $i \in [m]$, compute $\hat{\mathbf{R}}^{(i)} = \tilde{\mathbf{U}}^{(i)\top} \mathbf{A}^{(i)} \tilde{\mathbf{V}}^{(i)}$, where $\tilde{\mathbf{U}}^{(i)} = [\hat{\mathbf{U}}_c \mid \hat{\mathbf{U}}_s^{(i)}]$ and $\tilde{\mathbf{V}}^{(i)} = [\hat{\mathbf{V}}_c \mid \hat{\mathbf{V}}_s^{(i)}]$.
\end{enumerate} 
\ENSURE $\hat{\mathbf{U}}_c, \hat{\mathbf{V}}_c, \{\hat{\mathbf{U}}_s^{(i)}, \hat{\mathbf{V}}_s^{(i)}, \hat{\mathbf{R}}^{(i)}\}_{i=1}^m$.
\end{algorithmic}
\end{algorithm}

\subsection{Theoretical results}
\label{sec:theoretical_COISIE}
We shall make the following assumptions on the edge probability matrices $\mathbf{P}^{(i)}$ for $1 \leq i \leq m$. We emphasize that, because our theoretical results address either large-sample approximations or limiting distributions, these assumptions should be interpreted in the regime where $n$ is arbitrarily large and/or $n \rightarrow \infty$. 
We also assume, unless stated otherwise, that the number of graphs $m$ is bounded as (1) in many applications, we only observe a bounded number of networks even when the number of vertices $n$ per graph is large, and (2) if the graphs are not too sparse, allowing $m \rightarrow \infty$ leads to more accurate
estimation of $\mathbf{U}_c$ and $\mathbf{V}_c$, while having no detrimental effect on the estimation of $\{\mathbf{U}_s^{(i)}$, $\mathbf{V}_s^{(i)},\mathbf{R}^{(i)}\}_{i=1}^m$.
\begin{assumption}
  \label{ass:main_part2}
  The following conditions hold for sufficiently large $n$.
  \begin{itemize}
    \item The matrices $\mathbf{U}^{(i)}$ and $\mathbf{V}^{(i)}$ are $n \times d_i$ matrices with bounded coherence, i.e.,
    \[
    \|\mathbf{U}^{(i)}\|_{2 \to \infty} \lesssim d_i^{1/2}n^{-1/2} \quad \text{and} \quad \|\mathbf{V}^{(i)}\|_{2 \to \infty} \lesssim d_i^{1/2}n^{-1/2}.
    \]
    \item There exists a factor $\rho_n \in [0,1]$ depending on $n$ such that for each $i \in [m]$, $\mathbf{R}^{(i)}$ is a $d_i \times d_i$ matrix with $\|\mathbf{R}^{(i)}\| = \Theta(n \rho_n)$, where $n\rho_n \geq C\log n$ for some sufficiently large but finite constant $C > 0$. We interpret $n \rho_n$ as the growth rate for the average degree of the network $\mathbf{A}^{(i)}$ generated from $\mathbf{P}^{(i)}$. 
    \item The matrices $\{\mathbf{R}^{(i)}\}_{i=1}^{m}$ have bounded condition numbers, i.e., there exists a finite constant $M$ such that
    \[
    \max_{i \in [m]} \frac{\sigma_1(\mathbf{R}^{(i)})}{\sigma_{d_i}(\mathbf{R}^{(i)})} \leq M,
    \]
    where $\sigma_{1}(\mathbf{R}^{(i)})$ and $\sigma_{d_i}(\mathbf{R}^{(i)})$ denote the largest and smallest singular values of $\mathbf{R}^{(i)}$, respectively.
    \item There exists
a constant $c_s > 0$ not depending on $n$ such that
$$\max\Bigl\{\Bigl\| \frac{1}{m} \sum_{i=1}^m \muu_s^{(i)} \muu_s^{(i)\top}\Bigr\|,
\Bigl\| \frac{1}{m} \sum_{i=1}^m \mv_s^{(i)} \mv_s^{(i)\top}\Bigr\|\Bigr\} \leq 1 - c_s.$$
  \end{itemize}
\end{assumption}

\begin{remark1}
  \label{rem:assumptions}
  We provide some brief discussions surrounding Assumption~\ref{ass:main_part2}. The first condition on
  bounded coherence of $\mathbf{U}^{(i)}$ and $\mathbf{V}^{(i)}$ is a widely used and typically mild assumption in random graphs and other high-dimensional statistical inference problems, including matrix completion, covariance estimation, and subspace estimation; see, e.g., \cite{candes_recht, fan2018eigenvector, lei2019unified, abbe2020entrywise, cape2019two, cai2021subspace}. 
Bounded coherence together with the second condition $\|\mathbf{R}^{(i)}\| \asymp n \rho_n = \Omega(\log n)$ implies that the average degree of each graph $\mathbf{A}^{(i)}$ grows poly-logarithmically in $n$. This semisparse regime $n \rho_n = \Omega(\log n)$ is generally necessary for spectral methods to work, i.e., if $n \rho_n = o(\log n)$, then the singular vectors of any individual $\mathbf{A}^{(i)}$ are no longer consistent estimates of $\mathbf{U}^{(i)}$ and $\mathbf{V}^{(i)}$. 
The third condition of bounded condition number ensures that each $\mathbf{R}^{(i)}$ is full-rank and hence the column space (resp. row space) of each $\mathbf{P}^{(i)}$ is identical to that of $\mathbf{U}^{(i)}$ (resp. $\mathbf{V}^{(i)}$).
{\color{black}The last condition ensures that the individual subspaces $\{\mathbf{U}^{(i)}_s\}_i$ and $\{\mathbf{V}^{(i)}_s\}_i$ are sufficiently diverse and thus neither of them is part of the common subspaces $\muu_c$ and $\mv_c$, respectively. We now provide some intuition in terms of $\{\muu_s^{(i)}\}_i$ and $\muu_c$ (the discussion for $\{\mv_s^{(i)}\}_i$ and $\mv_c$ is analogous). Since each $\muu_s^{(i)} \muu_s^{(i)\top}$ is a projection matrix, the average $\bm\Pi_s := \frac{1}{m}\sum_{i=1}^m \muu_s^{(i)}\muu_s^{(i)\top}$ satisfies $\|\bm\Pi_s\| \leq 1$, with equality if and only if there exists a direction $\mathbf{w}$ that is contained in the column space of {every} $\muu_s^{(i)}$ (such a direction gives $\mathbf{w}^\top \bm\Pi_s \mathbf{w} = \frac{1}{m}\sum_{i=1}^m \|\muu_s^{(i)\top}\mathbf{w}\|^2 = 1$). Any such shared direction reflects structure common to all layers, and should therefore be attributed to the common subspace $\muu_c$ rather than to the individual subspaces. The condition $\|\bm\Pi_s\| \leq 1 - c_s$ rules out such shared directions and thereby guarantees that $\muu_c$ captures {all} of the shared structure, where the gap $c_s$ is required to be a constant independent of $n$ so that the individual subspaces remain sufficiently diverse as $n \to \infty$.
}
  \end{remark1}
  
  We now present uniform error bounds and normal approximations for the row-wise fluctuations of $\hat{\mathbf{U}}_c$ and $\hat{\mathbf{U}}_s^{(i)}$ (resp. $\hat{\mathbf{V}}_c$ and $\hat{\mathbf{V}}_s^{(i)}$) around $\mathbf{U}_c$ and $\mathbf{U}_s^{(i)}$ (resp. $\mathbf{V}_c$ and $\mathbf{V}_s^{(i)}$).
  These results offer significantly stronger theoretical guarantees compared to the Frobenius norm error bounds commonly
  encountered in the literature; see Section~\ref{sec:related_works} for further discussion.

\begin{theorem1}
  \label{thm:UhatW-U=EVR^{-1}+xxx_part2} 
  Consider $\left(\mathbf{A}^{(1)}, \ldots, \mathbf{A}^{(m)}\right) \sim \operatorname{COISIE}(\mathbf{U}_c, \mathbf{V}_c, \{\mathbf{U}_s^{(i)}, \mathbf{V}_s^{(i)}, \mathbf{R}^{(i)}\}_{i=1}^m)$ under the conditions in Assumption~\ref{ass:main_part2}. Let $\hat{\mathbf{U}}_c$ be the estimate of $\mathbf{U}_c$ obtained by Algorithm~\ref{Alg_COISIE}, and let $\mathbf{W}_{\mathbf{U}_c}$ be the minimizer of $\|\hat{\mathbf{U}}_c \mathbf{O} - \mathbf{U}_c\|_{F}$ over all $d_{0,\mathbf{U}} \times d_{0,\mathbf{U}}$ orthogonal matrices $\mathbf{O}$.
Then
\begin{equation}
  \label{eq:expansion_hatU_seperate_3}
  \hat{\mathbf{U}}_c \mathbf{W}_{\mathbf{U}_c} -\mathbf{U}_c =
  \frac{1}{m}\sum_{i=1}^m \me^{(i)} \mv^{(i)}(\mr^{(i)})^{-1}\muu^{(i)\top} \muu_c+\mq_{\muu_c},
\end{equation}
where $\me^{(i)}=\ma^{(i)}-\mpp^{(i)}$ and $\mq_{\muu_c}$ is a random matrix satisfying
\begin{equation*}
  \begin{aligned}
 &\|\mq_{\muu_c}\| \lesssim (n \rho_n)^{-1}\max\{1,d_{\max}^{1/2}\rho_n^{1/2}\log^{1/2} n\}, \\
&\|\mq_{\muu_c}\|_{2\to\infty} \lesssim d_{\max}^{1/2}n^{-1/2}(n\rho_n)^{-1}\log n
  \end{aligned}
\end{equation*}
with high probability, where $d_{\max} = \max_{i\in[m]} d_i$. 
Also, for any $k \in [n]$, the $k$th row $q_{\muu_c,k}$ of $\mq_{\muu_c}$ satisfies $$\|q_{\muu_c,k}\| \lesssim d_{\max}^{1/2} n^{-1/2} (n \rho_n)^{-1} t$$ with probability at least $1 - n^{-c}- O(me^{-t})$ for any $c>0$. 

For each $i\in[m]$, let $\hat{\mathbf{U}}_s^{(i)}$ be the estimate of
$\muu_s^{(i)}$
obtained by Algorithm~\ref{Alg_COISIE}, and let $\mw^{(i)}_{\muu_s}$ be the minimizer of $\|\hat{\muu}_s^{(i)} \mo -
\muu_s^{(i)}\|_{F}$ over all $(d_i-d_{0,\muu})\times (d_i-d_{0,\muu})$ orthogonal matrices $\mo$. Then
\begin{equation*}
  \hat{\muu}_s^{(i)} \mw^{(i)}_{\muu_s} - \muu_s^{(i)} = \me^{(i)}
  \mv^{(i)} (\mr^{(i)})^{-1} \muu^{(i)\top} \muu_s^{(i)} + \mq_{\muu_s}^{(i)},
\end{equation*}
where the random matrix 
$\mathbf{Q}_{\mathbf{U}_s}^{(i)}$ and its $k$th row $q_{\mathbf{U}_s,k}^{(i)}$ satisfy the same upper bounds as those for $\mathbf{Q}_{\mathbf{U}_c}$ and $q_{\mathbf{U}_c,k}$.

The estimates $\hat{\mathbf{V}}_c$, $\hat{\mathbf{V}}_s^{(i)}$ have similar expansions and analogous bounds, with $\mathbf{E}^{(i)}$, $\mathbf{R}^{(i)}$, $\mq_{\muu_c}$, and $\mq_{\muu_s}^{(i)}$ replaced by $\mathbf{E}^{(i)\top}$, $\mathbf{R}^{(i)\top}$, $\mq_{\mv_c}$, and $\mq_{\mv_s}^{(i)}$, respectively, and the roles of $\mathbf{V}^{(i)}$, $\mathbf{U}_c$, and $\mathbf{U}_s^{(i)}$ swapped with $\mathbf{U}^{(i)}$, $\mathbf{V}_c$, and $\mathbf{V}_s^{(i)}$.
\end{theorem1}

For ease of exposition, we assume that $\{d_i\}_{i=1}^m$, $d_{0,\mathbf{U}}$, and $d_{0,\mathbf{V}}$ are known in the statement of Theorem~\ref{thm:UhatW-U=EVR^{-1}+xxx_part2}. If $\{d_i\}$ are unknown, they can be estimated using the following approach: for each $i \in [m]$, let $\hat{d}_i$ be the number of eigenvalues of $\mathbf{A}^{(i)}$ exceeding $4 \sqrt{\delta(\mathbf{A}^{(i)})}$ in modulus, where $\delta(\mathbf{A}^{(i)})$ denotes the maximum degree of $\mathbf{A}^{(i)}$. 
Under the conditions in Assumption~\ref{ass:main_part2}, we can show that $\hat{d}_i$ is a consistent estimate of $d_i$ by combining tail bounds for $\|\mathbf{A}^{(i)} - \mathbf{P}^{(i)}\|$ (such as those in \cite{rinaldo_2013, oliveira2009concentration}) with Weyl's inequality; the details are omitted here.
If $d_{0,\mathbf{U}}$ (resp. $d_{0,\mathbf{V}}$) is unknown, it can be consistently estimated by selecting the number of eigenvalues of  
$
\hat{\bm{\Pi}}_\mathbf{U} := m^{-1} \sum_{i=1}^m \hat{\mathbf{U}}^{(i)} \hat{\mathbf{U}}^{(i)\top}$ (resp. $\hat{\bm{\Pi}}_\mathbf{V} := m^{-1} \sum_{i=1}^m \hat{\mathbf{V}}^{(i)} \hat{\mathbf{V}}^{(i)\top}
$)
that are approximately $1$. For example, let $\{\lambda_k(\hat{\bm{\Pi}}_\mathbf{U})\}_{k \geq 1}$ denote the eigenvalues of $\hat{\bm{\Pi}}_{\mathbf{U}}$
and define
$
\hat{d}_{0,\mathbf{U}} = |\{k \colon \lambda_k(\hat{\bm{\Pi}}_\mathbf{U}) \geq 1 - (n \rho_n)^{-1/2} \log n\}|.
$
Then under Assumption~\ref{ass:main_part2} we have $\hat{d}_{0,\mathbf{U}} \rightarrow d_{0,\mathbf{U}}$ (resp. $\hat{d}_{0,\mathbf{V}} \rightarrow d_{0,\mathbf{V}}$) almost surely.

\begin{remark1}
\label{rmk:akin to sample average}
If we fix an $i \in [m]$ and let $\hat{\muu}^{(i)}$ denote the leading left singular vectors of $\ma^{(i)}$, then there exists an orthogonal matrix $\mw_{\muu}^{(i)}$ such that  
$$
\hat{\muu}^{(i)} \mw_{\muu}^{(i)} - \muu = \me^{(i)} \mv^{(i)} (\mr^{(i)})^{-1} + \mq^{(i)}_\muu,
$$
where $\mq^{(i)}_\muu$ satisfies the same bounds as those for
$\mq_{\muu_c}$ and $\mq_{\muu_s}^{(i)}$ in
Theorem~\ref{thm:UhatW-U=EVR^{-1}+xxx_part2}.  This type of expansion
for the leading eigenvectors of a single $\ma^{(i)}$ is well-known in
the literature; see, e.g., \cite{cape2019signal, xie2021entrywise,
abbe2020entrywise}.  The primary conceptual and technical contribution
of Theorem~\ref{thm:UhatW-U=EVR^{-1}+xxx_part2} is in showing that,
while $\hat{\muu}_c$ is a \emph{nonlinear} function of
$\{\hat{\muu}^{(i)}\}_{i=1}^{m}$, the expansion for $\hat{\muu}_c$ can still be written as a linear
combination of the expansions for $\{\hat{\muu}^{(i)}\}$.
\end{remark1}

\begin{remark1}\label{rmk:independent1}
As mentioned previously, 
Theorem~\ref{thm:UhatW-U=EVR^{-1}+xxx_part2} does not require 
$\{\ma^{(i)}\}_{i=1}^{m}$ to be mutually independent. As a simple example, 
let $m = 2$ and suppose $\ma^{(1)}$ is an edge-independent random graph with edge 
probabilities $\mpp = \muu \mr \mv^{\top}$, while $\ma^{(2)}$ is a partially 
observed copy of $\ma^{(1)}$ where the entries are set to $0$ with 
probability $1 - p_*$ completely at random for some $p_* > 0$. 
Note that $\ma^{(2)}$ is dependent on $\ma^{(1)}$ but is also \emph{marginally} 
an edge-independent random graph with edge probabilities $p_* \mpp$. 
Hence, by Theorem~\ref{thm:UhatW-U=EVR^{-1}+xxx_part2} with $\muu_c = \muu$, $\mv_c = \mv$, 
$\muu_s^{(i)} = \mv_s^{(i)} = \mathbf{0}$, and $\mr^{(i)} = \mr$, we have
$$
\hat{\muu} \mw - \muu = (\ma^{(1)} - \mpp) \mv \mr^{-1} +
\frac{1}{2} \left(p_*^{-1} \ma^{(2)} - \ma^{(1)} \right) \mv \mr^{-1} +
\mq_{\muu},
$$
where $\mq_{\muu}$ satisfies the bounds as stated for $\mq_{\muu_c}$ 
in Theorem~\ref{thm:UhatW-U=EVR^{-1}+xxx_part2} with high probability. 
The difference between $\hat{\muu}^{(1)}$ (which depends only on $\ma^{(1)}$) and 
$\hat{\muu}$ thus corresponds to $p_*^{-1} \ma^{(2)} - \ma^{(1)}$.
\end{remark1}

\begin{remark1}\label{rmk:general error} Note that, for the COISIE
model, the entries of the noise $\mathbf{E}^{(i)} = \mathbf{A}^{(i)} -
\mathbf{P}^{(i)}$ are (centered) Bernoulli random variables. Our
theoretical results, however, can be easily adapted to a more general
setting where each $\mathbf{E}^{(i)}$ can be decomposed as the sum of
two mean-zero random matrices, $\mathbf{E}^{(i,1)}$ and
$\mathbf{E}^{(i,2)}$, where $\{\mathbf{E}^{(i,1)}\}$ have independent
bounded entries satisfying
$\max_{i,s,t}\mathbb{E}[(\mathbf{E}_{st}^{(i,1)})^2] \lesssim \rho_n$,
and $\{\mathbf{E}^{(i,2)}\}$ have independent sub-Gaussian entries
satisfying $\max_{i,s,t}\|\mathbf{E}_{st}^{(i,2)}\|_{\psi_2} \lesssim
\rho_n^{1 / 2}$.  In particular, the proofs in
Section~\ref{Appendix:thm:UhatW-U=EVR^{-1}+xxx_part2} and
Section~\ref{sec:proof_normality1_part2} of the supplementary material are written for this more
general noise model. The reason for presenting only (centered)
Bernoulli noise in this section is purely for simplicity of
exposition, as the COISIE model aligns well with many existing random
graph models.  For more general settings, we have the same theoretical
results with the caveat that the variance of $\me^{(i)}$ may have
different expressions under different settings. For example, the quantity
$\mathbf{\Xi}^{(i,k)}_{\ell\ell}$ in
Theorem~\ref{thm:WhatU_k-U_k->norm_part2} is actually the variance of
$\me^{(i)}_{k,\ell}$ and may need to be adjusted in different
settings, and similarly for $\tilde\md^{(i)}, \breve\md^{(i)},
\md^{(i)}$ in Theorem~\ref{thm:What_U Rhat What_v^T-R->norm}.
\end{remark1}

\begin{remark1}
\label{rk:generalized model}
Theorem~\ref{thm:UhatW-U=EVR^{-1}+xxx_part2} can be applied to the MultiNeSS model for multiplex networks in \cite{macdonald2022latent}.
More specifically, the MultiNeSS model assumes that we have a collection of symmetric matrices 
$$
\mpp^{(i)} = \mx_c \mi_{p_0,q_0} \mx_c^\top + \mx_s^{(i)} \mi_{p_i,q_i} \mx_s^{(i)\top},
$$
where $\mi_{p,q} = \text{diag}(\mi_p, -\mi_q)$ is a diagonal matrix with $r$ entries of ``1'' and $s$ entries of ``-1'' on the diagonal. Given a collection of noisily observed matrices $\ma^{(i)} = \mpp^{(i)} + \me^{(i)}$, where the upper triangular entries of $\me^{(i)}$ are independent mean-zero random variables, \cite{macdonald2022latent} proposes estimating $\mathbf{F} = \mx_c \mi_{p_0,q_0} \mx_c^\top$ and $\mathbf{G}^{(i)} = \mx_s^{(i)} \mi_{p_i,q_i} \mx_s^{(i)\top}$ by solving a convex optimization problem of the form
\begin{equation}
\label{eq:2}
\min_{\mathbf{F}, \{\mathbf{G}^{(i)}\}_{i=1}^{m}} \ell(\mathbf{F}, \{\mathbf{G}^{(i)}\}_{i=1}^{m} \mid \{\mathbf{A}^{(i)}\}_{i=1}^{m}) + \lambda \|\mathbf{F}\|_{\ast} + \sum_{i=1}^{m} \lambda \alpha_i \|\mathbf{G}^{(i)}\|_{\ast},
\end{equation}
where the minimization is over the set of $n \times n$ matrices $\{\mathbf{F}, \mathbf{G}^{(1)}, \dots, \mathbf{G}^{(m)}\}$. Here, $\ell(\cdot)$ is a loss function (e.g., the negative log-likelihood of $\ma^{(i)}$ assuming some parametric distribution for the entries of $\me^{(i)}$), $\|\cdot\|_{\ast}$ is the nuclear norm, and $\lambda, \alpha_1, \dots, \alpha_m$ are tuning parameters.
Denoting the minimizers of Eq.~\eqref{eq:2} by $\{\hat{\mathbf{F}}, \hat{\mathbf{G}}^{(1)}, \dots, \hat{\mathbf{G}}^{(m)}\}$, \cite{macdonald2022latent} provides upper bounds for $\|\mathbf{F} - \hat{\mathbf{F}}\|_{F}$ and $\|\hat{\mathbf{G}}^{(i)} - \mathbf{G}^{(i)}\|_{F}$. Letting $\hat{\mx}_c$ (resp. $\hat{\mx}_s^{(i)}$) be the minimizer of $\|\mz \mi_{p_0,q_0} \mz^\top - \hat{\mathbf{F}}\|_{F}$ (resp. $\|\mz \mi_{p_i,q_i} \mz^\top - \hat{\mathbf{G}}^{(i)}\|_{F}$) over all $\mz$ with the same dimensions as $\mx_c$ (resp. $\mx_s^{(i)}$), \cite{macdonald2022latent} also provides upper bounds for $\min_{\mathbf{W}} \|\hat\mx_c \mw - {\mx}_c\|_{F}$ and $\min_{\mathbf{W}^{(i)}} \|\hat{\mx}_s^{(i)}\mw^{(i)}-\mx_s^{(i)}\|_{F}$, where the minimization is over all (indefinite) orthogonal matrices $\mw, \mw^{(1)}, \dots, \mw^{(m)}$ of appropriate dimensions. See Theorem~2 and Proposition~2 in \cite{macdonald2022latent} for more details.

Instead of solving the optimization in Eq.~\eqref{eq:2}, one could
also estimate $\hat{\mx}_c$ and $\{\hat{\mx}_s^{(i)}\}$ using
Algorithm~\ref{Alg_COISIE}. Furthermore, by applying
Theorem~\ref{thm:UhatW-U=EVR^{-1}+xxx_part2}, one could obtain $2
\to \infty$ norm error bounds for these estimates, which would yield
uniform entrywise bounds for $\|\hat{\mathbf{F}} -
\mathbf{F}\|_{\max}$ and $\|\hat{\mathbf{G}}^{(i)} -
\mathbf{G}^{(i)}\|_{\max}$ for all $i \in [m]$. These $2 \to \infty$
error bounds and uniform entrywise bounds can be viewed as
refinements of the Frobenius norm upper bounds in
\cite{macdonald2022latent}.  Due to space constraints, we leave the
precise statement of these theoretical results to the interested reader and
instead present, in Section~\ref{sec:MultiNeSS}, some numerical
results comparing the estimates obtained from
Algorithm~\ref{Alg_COISIE} with those from \cite{macdonald2022latent}.
  \end{remark1}

  We now note several results that can be directly obtained from the expansions in Theorem~\ref{thm:UhatW-U=EVR^{-1}+xxx_part2}. The first result 
  provides a collection of $2 \to \infty$ and Frobenius norm bounds for $\hat{\muu}_c$ and $\hat{\muu}_s^{(i)}$.  
   \begin{proposition1}
  \label{thm:Vhat-VW_part2}
Consider the setting in Theorem~\ref{thm:UhatW-U=EVR^{-1}+xxx_part2} and
furthermore assume that $\{\ma^{(i)}\}_{i=1}^{m}$ are mutually independent.
Then
\begin{equation}
  \label{eq:expansion_hatU_2inf2_part2}
  \begin{aligned}
  	  &\|\hat{\mathbf{U}}_c\mathbf{W}_{\mathbf{U}_c} -\mathbf{U}_c\|_{2 \to
    \infty}  \lesssim  d_{\max}^{1/2} (mn)^{-1/2} (n \rho_n)^{-1/2} \log^{1/2}n+d_{\max}^{1/2}n^{-1/2}(n\rho_n)^{-1}\log n,\\
   & \|\hat{\mathbf{U}}_s^{(i)}\mathbf{W}_{\mathbf{U}_s}^{(i)} -\mathbf{U}_s^{(i)}\|_{2 \to
    \infty}  \lesssim  d_{\max}^{1/2} n^{-1/2} (n \rho_n)^{-1/2} \log^{1/2}n,
  \end{aligned}
\end{equation}
\begin{equation}
  \label{eq:expansion_hatU_F_part2}
  \begin{aligned}
&\|\hat{\mathbf{U}}_c\mw_{\muu_c}-\mathbf{U}_c \|_{F} 
\lesssim d_{\max}^{1/2}m^{-1/2}(n\rho_n)^{-1/2}
 +d_{0,\muu}^{1/2}(n\rho_n)^{-1}\max\{1, (d_{\max}\rho_n\log n)^{1/2}\},\\
 &\|\hat{\mathbf{U}}_s^{(i)}\mw_{\muu_s}^{(i)}-\mathbf{U}_{s}^{(i)} \|_{F} 
\lesssim d_{\max}^{1/2}(n\rho_n)^{-1/2}
  \end{aligned}
\end{equation}
with high probability.
Similar results hold for $\hat{\mathbf{V}}_c$ and $\hat{\mathbf{V}}_s^{(i)}$.
\end{proposition1}

\begin{remark1}
  \label{rem:condition m}
Note that, while we had generally assumed that $m$ is bounded (see the beginning of this subsection),
Eq.~\eqref{eq:expansion_hatU_2inf2_part2} holds as long as $m = O(n^{c})$ for some 
finite constant $c > 0$. Indeed, for any $c' \geq c$ we can choose a sufficiently large $C$ depending only on $c'$ such that
$\mt^{(i)} \lesssim C d_i^{1/2} n^{-1/2} (n \rho_n)^{-1} \log n$ with probability at least $1 - n^{-c'}$ (see Lemma~\ref{lemma:Uhat-UW} in the supplementary material) and thus, by taking a union bound over all $i \in [m]$ with $m = O(n^{c})$ we can still have Eq.~\eqref{eq:expansion_hatU_seperate_3} with the same bounds.
For $\|\hat{\mathbf{U}}_c \mathbf{W}_{\mathbf{U}_c} - \mathbf{U}_c\|_{2 \to \infty}$,
if $m=O(n\rho_n)$, 
then the first term $d_{\max}^{1/2} (mn)^{-1/2} (n \rho_n)^{-1/2} \log^{1/2} n$ dominates, and the error decreases as $m$ increases (assuming $n$ and $\rho_n$ are fixed). In contrast, if $m = \omega(n \rho_n)$ then the second term
dominates, i.e., increasing $m$ with $n$ and $\rho_n$ fixed does not guarantee smaller
errors. 
The Frobenius norm bound in Eq.~\eqref{eq:expansion_hatU_F_part2} exhibits similar behavior; see Theorem~3 in \cite{arroyo2019inference} for a similar result.
These results indicate that, for the estimation of the shared subspaces in the COISIE model to achieve the ``optimal" error rate, we need $m$ not to be too large compared to $n \rho_n$.
\end{remark1}

The next result provides normal approximations for the rows of $\hat\muu_c$ and $\hat\muu_s^{(i)}$. 
\begin{theorem1}
  \label{thm:WhatU_k-U_k->norm_part2}
Consider the setting in Theorem~\ref{thm:UhatW-U=EVR^{-1}+xxx_part2} and
further assume that $\{\ma^{(i)}\}_{i=1}^{m}$ are mutually independent.
For any $i\in[m]$ and $k\in[n]$, let $\mathbf{\Xi}^{(i,k)}$ be a $n\times n$ diagonal matrix whose diagonal elements are of the form
$
\mathbf{\Xi}^{(i,k)}_{\ell\ell}=\mpp_{k  \ell}^{(i)}(1-\mpp_{k    \ell}^{(i)}).
$
Define $\bm{\Upsilon}^{(k)}_{\muu_c}$ as the $d_{0,\muu} \times d_{0,\muu}$ symmetric matrix
$$
\begin{aligned}
	&\bm{\Upsilon}_{\muu_c}^{(k)}=\frac{1}{m^2}\sum_{i=1}^{m}\muu_c^\top\muu^{(i)}
(\mr^{(i)\top})^{-1}\mv^{(i)\top}\mathbf{\Xi}^{(k,i)}\mv^{(i)}(\mr^{(i)})^{-1}\muu^{(i)\top}
\muu_c.
\end{aligned}
$$
Note that $\|\bm{\Upsilon}^{(k)}_{\muu_c}\|\lesssim (mn^2\rho_n)^{-1}$. 
Further suppose $\sigma_{\min}(\bm{\Upsilon}^{(k)}_{\muu_c})\gtrsim
(mn^2\rho_n)^{-1}$. 
Then for the $k$th rows $\hat{u}_{c,k}$ and ${u}_{c,k}$ of $\hat\muu_c$ and $\muu_c$, we have
\begin{equation}
  \label{eq:clt_COSIE_part2}
  \begin{aligned}
  	(\bm{\Upsilon}^{(k)}_{\muu_c})^{-1/2} \big(\mw_{\muu_c}^\top \hat{u}_{c,k}-{u}_{c,k}\big)
\rightsquigarrow \mathcal{N}\big(\mathbf{0},\mi_{d_{0,\muu}}\big)
  \end{aligned}
\end{equation}
as $n \rightarrow \infty$.

For each $i\in[m]$, define $\bm{\Upsilon}^{(i,k)}_{\muu_s}$ as the $(d_i-d_{0,\muu}) \times (d_i-d_{0,\muu})$ symmetric matrix
$$
\bm{\Upsilon}_{\muu_s}^{(i,k)}=\muu_s^{(i)\top}\muu^{(i)}
(\mr^{(i)\top})^{-1}\mv^{(i)\top}\mathbf{\Xi}^{(k,i)}\mv^{(i)}(\mr^{(i)})^{-1}\muu^{(i)\top}
\muu_s^{(i)}.
$$
Note that $\|\bm{\Upsilon}^{(i,k)}_{\muu_s}\|\lesssim (n^2\rho_n)^{-1}$. Further suppose $\sigma_{\min}(\bm{\Upsilon}^{(i,k)}_{\muu_s})\gtrsim
(n^2\rho_n)^{-1}$.
Then for the $k$th rows $\hat{u}_{s,k}^{(i)}$ and ${u}_{s,k}^{(i)}$ of $\hat\muu_s^{(i)}$ and $\muu_s^{(i)}$, we have
$$
(\bm{\Upsilon}^{(i,k)}_{\muu_s})^{-1/2} \big(\mw_{\muu_s}^{(i)\top} \hat{u}_{s,k}^{(i)}-{u}_{s,k}^{(i)}\big)
\rightsquigarrow \mathcal{N}\big(\mathbf{0},\mi_{(d_i-d_{0,\muu})}\big)
$$
as $n \rightarrow \infty$.

Similar results hold for $\hat{\mathbf{V}}_c$, $\hat{\mathbf{V}}_s^{(i)}$ and their rows $\hat v_{c,k}$, $\hat v_{s,k}^{(i)}$ with $\mpp^{(i)}$ and $\mathbf{R}^{(i)}$ replaced by $\mpp^{(i)\top}$ and $\mathbf{R}^{(i)\top}$, respectively, and the roles of $\mathbf{V}^{(i)}$, $\mathbf{U}_c$, and $\mathbf{U}_s^{(i)}$ swapped with $\mathbf{U}^{(i)}$, $\mathbf{V}_c$, and $\mathbf{V}_s^{(i)}$.
\end{theorem1}

\begin{remark1}\label{rem:nondegeneracy}
The row-wise normal approximations in
Eq.~\eqref{eq:clt_COSIE_part2} 
assumes
that the minimum eigenvalue of $\bm{\Upsilon}^{(k)}_{\muu_c}$ grows at rate $(mn^{2}\rho_n)^{-1}$,
and this condition holds whenever the entries of
$\mpp^{(i)}$ are {\em homogeneous}, e.g., 
suppose $\min_{k\ell} \mpp^{(i)}_{k \ell} \asymp
\max_{k\ell} \mpp^{(i)}_{k \ell} \asymp \rho_n$, then for any $i\in[m]$ 
we have $\min_{k,i,\ell} \mathbf{\Xi}^{(k,i)}_{\ell\ell} \gtrsim \rho_n$ and hence 
$$
\begin{aligned}
	&\sigma_{\min}\big(\muu_c^\top\muu^{(i)}(\mr^{(i)\top})^{-1}\mv^{(i)\top}\mathbf{\Xi}^{(k,i)}\mv^{(i)}(\mr^{(i)})^{-1}\muu^{(i)\top}
\muu_c\big)\\
&\geq \min_{\ell\in[n]}(\mathbf{\Xi}^{(k,i)}_{\ell\ell})\cdot  \sigma_{\min}\big(\muu_c^\top\muu^{(i)}(\mr^{(i)\top})^{-1}\mv^\top\mv(\mr^{(i)})^{-1}\muu^{(i)\top}
\muu_c\big)\\
&\geq \min_{\ell\in[n]}(\mathbf{\Xi}^{(k,i)}_{\ell\ell})\cdot \sigma_{\min}^2\big((\mr^{(i)})^{-1}\big)
 \gtrsim (n^{2}\rho_n)^{-1}.
\end{aligned}
$$
Weyl's inequality then implies
$$
\sigma_{\min}(\bm{\Upsilon}^{(k)}_{\muu_c})
\geq\frac{1}{m^2}\sum_{i=1}^{m}\sigma_{\min}\big(\muu_c^\top\muu^{(i)}(\mr^{(i)\top})^{-1}\mv^{(i)\top}\mathbf{\Xi}^{(k,i)}\mv^{(i)}(\mr^{(i)})^{-1}\muu^{(i)\top}
\muu_c\big)
 \gtrsim (mn^{2}\rho_n)^{-1}.
$$
The main reason for requiring a lower bound for the eigenvalues of $\bm{\Upsilon}^{(k)}_{\muu_c}$ is that we do not require $\bm{\Upsilon}^{(k)}_{\muu_c}$ to converge to any fixed matrix as $n \to \infty$, and thus  we cannot directly use $\bm{\Upsilon}^{(k)}_{\muu_c}$ in our limiting normal approximation. Rather, we need to scale $\mw_{\muu_c}^\top \hat{u}_{c,k} - u_{c,k}$ by $(\bm{\Upsilon}^{(k)}_{\muu_c})^{-1/2}$, and to ensure that this scaling is well-behaved, we need to control the smallest eigenvalue of $\bm{\Upsilon}^{(k)}_{\muu_c}$. 
{\color{black}
We note that this lower-bound condition can in fact be removed, albeit at a slight cost in the form of the limiting distribution. More specifically, without assuming a lower bound on $\sigma_{\min}(\bm{\Upsilon}^{(k)}_{\muu_c})$, we still have $m^{1/2} n \rho_n^{1/2} \big(\mw_{\muu_c}^\top \hat{u}_{c,k}-{u}_{c,k}\big) \rightsquigarrow \mathcal{N}\big(\mathbf{0}, m n^2 \rho_n \bm{\Upsilon}^{(k)}_{\muu_c}\big)$. However, since the covariance matrix $m n^2 \rho_n \bm{\Upsilon}^{(k)}_{\muu_c}$ depends on $\{\muu^{(i)},\mv^{(i)},\mr^{(i)}\}_i$ and thus varies with $n$, it is not a bona fide limit result (even though the covariance has bounded Frobenius norm) and is therefore less elegant than Eq.~\eqref{eq:clt_COSIE_part2}. Alternatively, one may stabilize the scaling by replacing $(\bm{\Upsilon}^{(k)}_{\muu_c})^{-1/2}$ with $(\bm{\Upsilon}^{(k)}_{\muu_c} + \gamma \mathbf{I})^{-1/2}$ for some regularization parameter $\gamma = O((m n^2 \rho_n)^{-1})$, in which case $(\bm{\Upsilon}^{(k)}_{\muu_c} + \gamma \mathbf{I})^{-1/2} \big(\mw_{\muu_c}^\top \hat{u}_{c,k}-{u}_{c,k}\big) \rightsquigarrow \mathcal{N}\big(\mathbf{0}, \mathbf{I} - \gamma (\bm{\Upsilon}^{(k)}_{\muu_c} + \gamma \mathbf{I})^{-1}\big)$. Such a stabilized scaling can still be useful for subsequent inference; for example, $\|(\bm{\Upsilon}^{(k)}_{\muu_c} + \gamma \mathbf{I})^{-1/2}(\mw_{\muu_c}^\top \hat{u}_{c,k}-{u}_{c,k})\|^2$ converges to a $\chi^2$ distribution for a suitable choice of $\gamma$.
}
A similar analysis applies to the condition on $\bm{\Upsilon}^{(i,k)}_{\muu_s}$.
Finally, if we allow $m$ to grow, then Eq.~\eqref{eq:clt_COSIE_part2} also holds for $m \log^2 m = o(n \rho_n)$, as we still have $(mn^2 \rho_n)^{1/2} q_{\muu_c,k} \to 0$ in probability, where $q_{\muu_c,k}$ is the term appearing in Eq.~\eqref{eq:expansion_hatU_seperate_3}.
\end{remark1}

\begin{remark1}\label{rem:independence}
For simplicity of presentation we assume in Theorem~\ref{thm:WhatU_k-U_k->norm_part2} that 
$\{\ma^{(i)}\}_{i=1}^{m}$ are mutually independent, but our result also holds 
under weaker conditions. More specifically, the normal approximation of $\hat{u}_{c,k}$ 
in Theorem~\ref{thm:WhatU_k-U_k->norm_part2} is based on
Eq.~\eqref{eq:thm2_proof1}, 
where $q_{\muu_c,k}$ is negligible in the limit. If $\{\ma^{(i)}\}$ are 
mutually independent, then the right-hand side of 
Eq.~\eqref{eq:thm2_proof1}
(ignoring $q_{\muu_c,k}$) is a sum of independent, mean $\mathbf{0}$ random vectors. 
In this case, we can apply the Lindeberg-Feller central limit theorem to show that 
$\mw_{\muu_c}^\top \hat{u}_{c,k} - u_{c,k}$ is approximately multivariate normal. 
Now suppose we make the weaker assumption that, for a fixed index $k\in[n]$, 
$\xi_{k1}, \xi_{k2}, \dots, \xi_{kn}$ are mutually independent random vectors, 
where $\xi_{k \ell} = (\me^{(1)}_{k \ell}, \dots, \me^{(m)}_{k \ell})$ 
for each $\ell \in [n]$. Then, under certain mild conditions on the covariance 
matrix for each $\xi_{k \ell}$, we have
$
(\tilde{\bm{\Upsilon}}^{(k)}_{\muu_c})^{-1/2} \bigl(\mw_{\muu_c}^\top \hat{u}_{c,k} - {u}_{c,k}\bigr)
\rightsquigarrow \mathcal{N}(\mathbf{0}, \mathbf{I})
$
as $n \to \infty$, where $\tilde{\bm{\Upsilon}}^{(k)}_{\muu_c}$ is a 
$(d_i-d_{0,\muu}) \times (d_i-d_{0,\muu})$ covariance matrix of the form
$$
\tilde{\bm{\Upsilon}}^{(k)}_{\muu_c} = \frac{1}{m^2} \sum_{\ell=1}^{n} \sum_{i=1}^{m} \sum_{j=1}^{m}
\mathrm{Cov}\bigl(\me^{(i)}_{k \ell}, \me^{(j)}_{k \ell}\bigr)
\cdot \muu_c^\top \muu^{(i)} (\mr^{(i)\top})^{-1} v_{\ell}^{(i)} v_{\ell}^{(i)\top} 
(\mr^{(j)})^{-1}\muu^{(i)\top} \muu_c,
$$
where $v_{\ell}^{(i)}$ denotes the $\ell$th row of $\mv^{(i)}$.
For example, suppose that the entries of $\xi_{k \ell}$ are pairwise uncorrelated, 
i.e., $\mathbb{E}[\me^{(i)}_{k \ell} \me^{(j)}_{k \ell}] = 0$ for all $i \neq j$ and all 
$\ell \in [n]$. Then $\mathrm{Var}[\xi_{k \ell}]$ is a diagonal matrix for all $\ell$, 
in which case $\tilde{\bm{\Upsilon}}^{(k)}_{\muu_c}$ coincides with 
$\bm{\Upsilon}^{(k)}_{\muu_c}$ as given in Theorem~\ref{thm:WhatU_k-U_k->norm_part2}. 
As another example, suppose $\ma^{(i)}$ and $\ma^{(j)}$ are pairwise $\rho$-correlated 
random graphs \citep{zheng2022vertex} for all $i \neq j$. Then
\[
\begin{aligned}
\tilde{\bm{\Upsilon}}^{(k)}_{\muu_c} = \frac{1}{m^2} \sum_{\ell=1}^{n} \sum_{i=1}^{m} \sum_{j=1}^{m}
&
\Bigl(\mathrm{Var}[\ma^{(i)}_{k \ell}] \mathrm{Var}[\ma^{(j)}_{k \ell}]\Bigr)^{1/2}
\Bigl(\rho \bm{1}\{i \neq j\} + \bm{1}\{i = j\}\Bigr) 
\\& \muu_c^\top \muu^{(i)} (\mr^{(i)\top})^{-1} v_{\ell}^{(i)} v_{\ell}^{(i)\top} 
(\mr^{(j)})^{-1} \muu^{(i)\top} \muu_c.
\end{aligned}
\] 
Similar remarks also hold for the normal approximations of $\hat{u}_{s,k}^{(i)}$.
\end{remark1}

\begin{remark1}\label{rm:compare_single}
  We now compare our inference results for multiple networks
  against existing results for the spectral embedding of a \textit{single} network. In particular, the COISIE model with $m = 1$ is equivalent to the GRDPG model \citep{grdpg1}, and thus our limiting results for $m = 1$ are the same as those for the adjacency spectral decomposition of a single GRDPG; e.g., Theorem~3.1 in \cite{xie2021entrywise} and Theorem~3 in \cite{tang2018eigenvalues} correspond to special cases of Theorem~\ref{thm:WhatU_k-U_k->norm_part2} and the following Theorem~\ref{thm:What_U Rhat What_v^T-R->norm} in this paper. 
If $m > 1$ and $\mpp^{(i)} = \mpp^{(1)}$ for all $i$, then for any $k \in [n]$, we have
$
\bm{\Upsilon}^{(k)}_{\muu_c} = m^{-1} \bm{\Upsilon}^{(1,k)}_{\muu_c}
$ and $
\bm{\Upsilon}^{(i,k)}_{\muu_s} = \bm{\Upsilon}^{(1,k)}_{\muu_s},
$
where $\bm{\Upsilon}^{(1,k)}_{\muu_c}$ and $\bm{\Upsilon}^{(1,k)}_{\muu_s}$ are the asymptotic covariance matrices for the corresponding entries in the adjacency spectral decomposition of a single GRDPG with edge probability matrix $\mpp^{(1)}$ (as given in Theorem~3.1 of \cite{xie2021entrywise}). 
If $\{\mpp^{(i)}\}$ are heterogeneous, then $\bm{\Upsilon}^{(k)}_{\muu_c}$ has a more complicated form (as it depends on the full collection $\{\mpp^{(i)}\}_{i=1}^{m}$), but nevertheless we still have 
$
\|\bm{\Upsilon}^{(k)}_{\muu_c}\| \lesssim (m n^2 \rho_n)^{-1},
$
while $\bm{\Upsilon}^{(i,k)}_{\muu_s}$ depends only on $\mpp^{(i)}$.
In summary, having $m > 1$ graphs with a common subspace leads to better estimation accuracy for $\muu_c$ and $\mv_c$ compared to that of a single GRDPG, as we can leverage information across multiple graphs.
In contrast, the estimation accuracy for $\muu_s^{(i)}$ and $\mv_s^{(i)}$ is not improved even when we have $m > 1$ graphs (see Theorem~\ref{thm:WhatU_k-U_k->norm_part2} and Proposition~\ref{thm:Vhat-VW_part2}), and the same holds for the estimation accuracy of $\mr^{(i)}$ (see Theorem~\ref{thm:What_U Rhat What_v^T-R->norm}). This is because $\muu_s^{(i)}$, $\mv_s^{(i)}$, and $\mr^{(i)}$ may be heterogeneous across different $i$, and thus each is estimated using only the corresponding $\ma^{(i)}$.
\end{remark1}

\begin{remark1}
\label{rm:undirect}
Theorem~\ref{thm:UhatW-U=EVR^{-1}+xxx_part2}, Theorem~\ref{thm:WhatU_k-U_k->norm_part2}, and Proposition~\ref{thm:Vhat-VW_part2},
with minimal changes, also hold when the $\ma^{(i)}$ are adjacency matrices for {\em undirected} graphs.
In particular, the expansion in Eq.~\eqref{eq:expansion_hatU_seperate_3} still 
holds for undirected graphs with $\mv^{(i)} = \muu^{(i)}$. Given this expansion, the bounds in Proposition~\ref{thm:Vhat-VW_part2} and
the normal approximations in 
Theorem~\ref{thm:WhatU_k-U_k->norm_part2} can be derived using the same arguments as those 
presented in the supplementary material.
  \end{remark1}

\subsection{Application to the COSIE model and two-sample hypothesis testing}
\label{sec:COSIE}
We now present our theoretical results for the COSIE model as a special case
of the COISIE model in which $\muu^{(i)} \equiv \muu_c$ and $\mv^{(i)} \equiv \mv$ for all $i$, so that there are no individual subspaces.  
In particular, we will consider the two-sample hypothesis testing problem for detecting similarities or
differences between multiple networks, which is of both theoretical and practical interest;
e.g., this type of problem arises naturally in neuroscience \citep{mheich2020brain,zalesky2012connectivity} and social networks \citep{fan2015similarity} applications. 

Recall that the edge probabilities matrices for the COSIE model are of
the form $\mpp^{(i)}=\muu\mr^{(i)}\mv^\top$ for all $i$.
See Section~\ref{sec:formal COSIE} in the supplementary material for
a more formal definition. 
We will
denote a collection of networks from the COSIE model as
$\left(\mathbf{A}^{(1)}, \ldots, \mathbf{A}^{(m)}\right) \sim
\operatorname{COSIE}(\mathbf{U}, \mathbf{V},
\{\mathbf{R}^{(i)}\}_{i=1}^m)$. 
Note that, for conciseness of exposition, these graphs are assumed to be {\em directed}
but the original formulation in \cite{arroyo2019inference} is for undirected graphs. 
Our theoretical results nevertheless apply to both the undirected and directed settings,
see Remark~\ref{rm:undirect} and
Remark~\ref{rm:undirect2} for details. Also, as mentioned in 
Section~\ref{sec:COISIE}, multilayer SBMs are a special case of the COSIE model. 
More specifically the edge probabilities of multilayer SBMs are of the form $\mpp^{(i)} = \mz
\mb^{(i)} \mz^\top$, where $\mz \in \mathbb{R}^{n \times K}$ with
entries in $\{0,1\}$ and $\sum_{k=1}^K \mz_{sk} = 1$ for all $s \in
[n]$ represents the consensus community assignments (which do not
change across graphs), and $\{\mb^{(i)}\}_{i=1}^{m} \subset
\mathbb{R}^{K \times K}$ with entries in $[0,1]$ represent the varying
community-wise edge probabilities.  This is equivalent to 
setting $\muu = \mv = \mz (\mz^\top
\mz)^{-1/2}$ and $\mr^{(i)} = (\mz^\top \mz)^{1/2} \mb^{(i)} (\mz^\top
\mz)^{1/2}$ for the COSIE parameters; see Proposition~1 in \cite{arroyo2019inference} for more
details.

Given $\left(\mathbf{A}^{(1)}, \ldots, \mathbf{A}^{(m)}\right) \sim
\operatorname{COSIE}(\mathbf{U}, \mathbf{V},
\{\mathbf{R}^{(i)}\}_{i=1}^m)$, 
we can use a simplified variant of Algorithm~\ref{Alg_COISIE} 
to estimate $\muu,\mv$ and $\mr^{(i)}$; see Algorithm~\ref{Alg} in
Section~\ref{sec:formal COSIE} of the supplementary material for
more details. Expansions for the resulting estimates $\hat{\muu}$ and $\hat{\mv}$, their error bounds,
and row-wise normal approximations are then special cases of Theorem~\ref{thm:UhatW-U=EVR^{-1}+xxx_part2}, Proposition~\ref{thm:Vhat-VW_part2}, and Theorem~\ref{thm:WhatU_k-U_k->norm_part2}. See Assumption~\ref{ass:main}, Theorem~\ref{thm:UhatW-U=EVR^{-1}+xxx}, Proposition~\ref{thm:Vhat-VW}, and Theorem~\ref{thm:WhatU_k-U_k->norm} in Section~\ref{sec:formal COSIE} of the supplementary material for the formal statements. 
Our main focus in this subsection is the following result on the limiting distribution of
$\{\hat{\mr}^{(i)}\}_{i=1}^{m}$. 
\begin{theorem1}
  \label{thm:What_U Rhat What_v^T-R->norm}
   Consider $\left(\mathbf{A}^{(1)}, \ldots, \mathbf{A}^{(m)}\right) \sim \operatorname{COSIE}(\mathbf{U}, \mathbf{V}, \{\mathbf{R}^{(i)}\}_{i=1}^m)$ under the conditions in Assumption~\ref{ass:main} and furthermore assume that $\{\ma^{(i)}\}_{i=1}^{m}$ are mutually independent. Let $\hat{\mathbf{U}}$, $\hat{\mathbf{V}}$, and $\hat{\mathbf{R}}^{(i)}$ be the estimates of $\mathbf{U}$, $\mathbf{V}$, and $\mathbf{R}^{(i)}$ obtained by Algorithm~\ref{Alg}, and let $\mathbf{W}_{\mathbf{U}}$ and $\mathbf{W}_{\mathbf{V}}$ be the minimizers of $\|\hat{\mathbf{U}} \mathbf{O} - \mathbf{U}\|_{F}$ and $\|\hat{\mathbf{V}} \mathbf{O} - \mathbf{V}\|_{F}$ over all $d \times d$ orthogonal matrices $\mathbf{O}$, respectively.
  Define
  $\tilde{\mathbf{D}}^{(i)}$ and $\breve{\mathbf{D}}^{(i)}$ as the $n \times n$ diagonal matrices with
  entries
  $$\tilde{\mathbf{D}}^{(i)}_{kk} = \sum_{\ell=1}^n \mathbf{P}^{(i)}_{k \ell} (1 -
  \mathbf{P}^{(i)}_{k \ell}), \quad \breve{\mathbf{D}}^{(i)}_{kk} =
  \sum_{\ell=1}^n \mathbf{P}^{(i)}_{\ell k} (1 -
  \mathbf{P}^{(i)}_{\ell k}),$$
  and define $\mathbf{D}^{(i)}$ as the $n^2 \times n^2$ diagonal matrix
  with diagonal entries
  $$\mathbf{D}^{(i)}_{k_1 + (k_2 - 1)n, k_1 + (k_2 - 1)n} =
  \mathbf{P}_{k_1 k_2}^{(i)} ( 1 - \mathbf{P}_{k_1k_2}^{(i)})$$
  for any $k_1,k_2\in[n]$. 
  Now let $\bm{\mu}^{(i)} \in \mathbb{R}^{d^2}$ be given by
\begin{equation*}
\begin{aligned}
    \boldsymbol\mu^{(i)}
	&=\operatorname{vec}\Bigl(
    \frac{1}{m}\muu^\top\tilde\md^{(i)}\muu(\mr^{(i)\top})^{-1}
    -\frac{1}{2m^2}\sum_{j=1}^m\mr^{(i)}(\mr^{(j)})^{-1}\muu^\top\tilde\md^{(j)}\muu(\mr^{(j)\top})^{-1}\Bigr) \\
    &+ 
	\mathrm{vec}\Bigl(\frac{1}{m}(\mr^{(i)\top})^{-1}\mv^\top\breve\md^{(i)}\mv
	-\frac{1}{2m^2}\sum_{j=1}^m(\mr^{(j)\top})^{-1}\mv^\top\breve\md^{(j)}\mv(\mr^{(j)})^{-1}\mr^{(i)}\Bigr).
\end{aligned}
\end{equation*}
Note that 
$\|\boldsymbol\mu^{(i)}\|_{\max}\lesssim m^{-1}$.
Next define $\bm{\Sigma}^{(i)}$ as the $d^2 \times d^2$ symmetric matrix
$$\bm{\Sigma}^{(i)} = (\mathbf{V} \otimes \mathbf{U})^{\top}
\mathbf{D}^{(i)} (\mathbf{V} \otimes \mathbf{U}).$$
Note that $\|\bm{\Sigma}^{(i)}\| \lesssim \rho_n$. Suppose also that $\sigma_{\min}(\mSigma^{(i)}) \gtrsim \rho_n$.
Then for $n \rho_n = \omega(n^{1/2})$ we
have 
      \begin{equation*}
      \bigl(\bm{\Sigma}^{(i)}\bigr)^{-1/2} \bigl(\mathrm{vec}\bigl(\mathbf{W}_{\mathbf{U}}^{\top}
        \hat{\mathbf{R}}^{(i)} \mathbf{W}_{\mathbf{V}} -
        \mr^{(i)}\bigr) - \bm{\mu}^{(i)} \bigr) \rightsquigarrow
        \mathcal{N}\bigl(\bm{0}, \mathbf{I} \bigr)
      \end{equation*}
as $n \rightarrow \infty$. 
Furthermore, the
$\{ 
\mw_\muu^\top\hat\mr^{(i)}\mw_\mv 
\}_{i=1}^{m}$
are {\em asymptotically} mutually independent. Finally, if $n \rho_n =
O(n^{1/2})$ we have
$$
     \mathrm{vec}\bigl(\mathbf{W}_{\mathbf{U}}^{\top}
        \hat{\mathbf{R}}^{(i)} \mathbf{W}_{\mathbf{V}} -
        \mr^{(i)}\bigr) - \bm{\mu}^{(i)}
        \overset{\mathrm{p}}{\longrightarrow} \bm{0} $$
        as $n \rightarrow \infty$. 
\end{theorem1}

\begin{remark1}
  \label{rem:sqrt_n} 
  The normal approximation in
Theorem~\ref{thm:What_U Rhat What_v^T-R->norm} requires $n \rho_n =
\omega(n^{1/2})$, as opposed to the much weaker condition of $n \rho_n
= \Omega(\log n)$ in Theorem~\ref{thm:WhatU_k-U_k->norm}. The main
reason for this discrepancy is that
Theorem~\ref{thm:WhatU_k-U_k->norm} is a limit result for any given
row of $\hat{\muu}$ while Theorem~\ref{thm:What_U Rhat
What_v^T-R->norm} requires \emph{averaging} over all $n$ rows of
$\hat{\muu}$; indeed, $\hat{\mr}^{(i)} = \hat{\muu}^\top \ma^{(i)}
\hat{\mv}$ is a bilinear form in $\{\hat{\muu}, \hat{\mv}\}$.  The
main technical challenge for Theorem~\ref{thm:What_U Rhat
What_v^T-R->norm} lies in showing that $\hat{\mr}^{(i)}$ has
substantially smaller variability (compared to the variability in
any given row of $\hat{\muu}$) without incurring significant bias, and
currently we can only guarantee this for $n \rho_n \gg n^{1/2}$. While
this might seem, at first glance, disappointing, it is however
expected as the $n^{1/2}$ threshold also appears in many related limit
results that involve averaging over the rows of $\hat{\muu}$.  For
example, \cite{li2018two} considers testing whether the community
memberships of two graphs are the same, and their test statistic,
which is based on the $\sin$-$\Theta$ distance between the singular
subspaces of the two graphs, converges to a standard normal
distribution under the condition $n \rho_n \gtrsim n^{1/2 + \epsilon}$
for some $\epsilon > 0$; see Assumption~3 in \cite{li2018two}.  As
another example, \cite{han_fan} studies the asymptotic distributions
for the leading eigenvalues and eigenvectors of a symmetric matrix
$\mathbf{X}$ under the assumption that $\mathbf{X} = \mathbf{H} +
\mathbf{W}$, where $\mathbf{H}$ is an unobserved low-rank symmetric
matrix and $\mathbf{W}$ is an unobserved generalized Wigner matrix
(i.e., the upper triangular entries of $\mathbf{W}$ are independent
mean-zero random variables). Among the numerous conditions in their
paper, one \emph{sufficient} condition for several of their main
results is $ \min_{k \ell} (\mathrm{Var}[w_{k \ell}])^{1/2} \gg
\|\mathbb{E}[\mathbf{W}^2]\|^{1/2} \times
|\lambda_{r}(\mathbf{H})|^{-1}, $ for all $r \leq d$. Here, $w_{k
\ell}$ denotes the random variable for the $k\ell$th entry of
$\mathbf{W}$, and $\lambda_r(\mathbf{H})$ is the $r$th largest
eigenvalue (in modulus) of $\mathbf{H}$; see Eq.~(13) in
\cite{han_fan} for more details.  Suppose we fix an $i \in [m]$ and
let $\mathbf{X} = \ma^{(i)}$, $\mathbf{H} = \mathbf{P}^{(i)}$, and
$\mathbf{W} = \mathbf{E}^{(i)}$ (note that the
eigenvalues of $\mathbf{P}^{(i)}$ can be extracted from those of
$\hat{\mr}^{(i)}$). Then,
assuming the conditions in Assumption~\ref{ass:main}, we have $
\min_{k \ell} (\mathrm{Var}[w_{k \ell}])^{1/2} \lesssim \rho_n^{1/2}$,
$\|\mathbb{E}[\mathbf{W}^2]\|^{1/2} \asymp (n \rho_n)^{1/2}$,
$\lambda_r(\mathbf{H}) \asymp n \rho_n $, and thus the condition in \cite{han_fan} simplifies to $\rho_n^{1/2} \gg (n \rho_n)^{-1/2}$, or
equivalently, $n \rho_n \gg n^{1/2}$.
{\color{black}
We also note that increasing the number of layers $m$ does not relax this assumption required for the normal approximation of $\hat{\mr}^{(i)}$ in our Theorem~\ref{thm:What_U Rhat What_v^T-R->norm} as well as in other related existing results discussed above. Intuitively, this is because $\mr^{(i)}$ can differ across $i$, so that neither the estimation of each $\hat{\mr}^{(i)}$ (as discussed in Remark~\ref{rm:compare_single}) nor the assumption required for its normal approximation can benefit from averaging across layers. From a more technical standpoint, although larger $m$ improves the estimation of the common subspaces $\hat{\muu}$ and $\hat{\mv}$, the $n^{1/2}$ requirement for $\hat{\mr}^{(i)}$ arises from controlling the bias induced by the corresponding single-layer noise $\me^{(i)}$ (see Eq.~\eqref{eq:What_U Rhat What_v^T-R} in the supplementary material), and this bias cannot be reduced by averaging across layers.
}

In addition, Theorem~\ref{thm:What_U Rhat What_v^T-R->norm} assumes
that the minimum eigenvalue of $\mSigma^{(i)}$ grows at the rate
$\rho_n$. This condition is analogous to the condition for
$\bm{\Upsilon}^{(k)}_{\muu_c}$ and $\bm{\Upsilon}^{(i,k)}_{\muu_s}$ in
Theorem~\ref{thm:WhatU_k-U_k->norm_part2} and is satisfied whenever
the entries of $\mpp^{(i)}$ are homogeneous.  Furthermore, as we will
see in the two-sample testing problem below, both $\mSigma^{(i)}$ and
$(\mSigma^{(i)})^{-1}$ are generally unknown and need to be estimated,
and consistent estimation of $\mSigma^{(i)}$ does not necessarily
imply consistent estimation of $(\mSigma^{(i)})^{-1}$ (and vice versa)
unless we can control $\sigma_{\min}(\mSigma^{(i)})$.
\end{remark1}

\begin{remark1}
\label{rm:undirect2}
Theorem~\ref{thm:What_U Rhat What_v^T-R->norm} also holds under the undirected setting with $\mv=\muu$, and can be derived using the same arguments as those presented in the supplementary material with the main difference being that the covariance matrix $\bm{\Sigma}^{(i)}$ in Theorem~\ref{thm:What_U Rhat What_v^T-R->norm} now has to account for the symmetry in $\me^{(i)}$. 
More specifically, let $\mathrm{vech}$ denote the half-vectorization of a matrix, and let $\md^{(i)}$ denote the $\tbinom{n+1}{2} \times \tbinom{n+1}{2}$ diagonal matrix with diagonal entries 
$
\operatorname{diag}(\md^{(i)}) = \operatorname{vech}\big(\mpp^{(i)}_{k_1 k_2}(1-\mpp^{(i)}_{k_1 k_2})\big).
$
Denote by $\mathcal{D}_n$ the $n^2 \times \tbinom{n+1}{2}$ duplication matrix which, for any $n \times n$ symmetric matrix $\mm$, transforms $\operatorname{vech}(\mm)$ into $\operatorname{vec}(\mm)$. Define
$$
\mUpsilon^{(i)}=(\muu\otimes\muu)^\top\mathcal{D}_n\md^{(i)}\mathcal{D}_n^\top(\muu\otimes\muu),
$$
and Theorem~\ref{thm:What_U Rhat What_v^T-R->norm}, when stated for undirected graphs, becomes
$$
\bigl(\mathcal{L}_d \bm{\Sigma}^{(i)} \mathcal{L}_d^\top\bigr)^{-1/2} 
\Bigl(\mathrm{vech}\bigl(\mathbf{W}_{\mathbf{U}}^\top
        \hat{\mathbf{R}}^{(i)} \mathbf{W}_{\mathbf{U}} -
        \mr^{(i)}\bigr)-\mathcal{L}_d\bm{\mu}^{(i)} \Bigr)
        \rightsquigarrow \mathcal{N}\bigl(\mathbf{0}, \mi\bigr),
$$
as $n \to \infty$.      
Here, $\mathcal{L}_d$ denotes the $\tbinom{d+1}{2} \times d^2$ elimination matrix that, given any $d \times d$ symmetric matrix $\mathbf{M}$, transforms $\mathrm{vec}(\mathbf{M})$ into $\mathrm{vech}(\mathbf{M}).$
\end{remark1}

We now consider the problem of detecting similarities or differences between multiple graphs, which is of both practical and theoretical importance. One typical application is testing for similarity across brain networks; see, e.g., \cite{zalesky2010network, rubinov2010complex, he2008structural}.
A simple and natural formulation of two-sample hypothesis testing for graphs assumes that they are {\em edge-independent} random graphs on the same set of vertices, and given any two graphs, they are said to be from the same (resp. ``similar'') distribution if their edge probability matrices are the same (resp. ``close''); see, e.g., \cite{tang2017semiparametric, ginestet2017hypothesis, ghoshdastidar2020two, li2018two, levin2017central, durante2018bayesian} for several recent examples of this type of formulation.

However, many existing test statistics do not have known {\em non-degenerate} limiting distributions, especially when comparing only two graphs, and calibration of their rejection regions has to be performed either via bootstrapping (see, e.g., \cite{tang2017semiparametric}) or via non-asymptotic concentration inequalities (see, e.g., \cite{ghoshdastidar2020two}).
Both of these approaches can be sub-optimal: bootstrapping is computationally expensive and has inflated type-I error when the bootstrapped distribution exhibits larger variability compared to the true distribution while non-asymptotic concentration inequalities are overly conservative and thus result in a significant loss of power.

We now discuss two-sample testing in the context of the COSIE model.
More specifically, suppose we are given a collection of networks $\left(\mathbf{A}^{(1)}, \ldots, \mathbf{A}^{(m)}\right) \sim \operatorname{COSIE}\left(\muu, \mathbf{V}, \{\mathbf{R}^{(i)}\}_{i=1}^m\right)$ and are interested in testing the null hypothesis 
$
\mathbb{H}_0 \colon \mpp^{(i)} = \mpp^{(j)}
$
against the alternative hypothesis 
$
\mathbb{H}_A \colon \mpp^{(i)} \neq \mpp^{(j)}
$
for some indices $i \neq j$. Since $\mpp^{(i)} = \muu \mr^{(i)} \mv^\top$, this is equivalent to testing 
$
\mathbb{H}_0 \colon \mr^{(i)} = \mr^{(j)} 
$ against $ 
\mathbb{H}_A \colon \mr^{(i)} \neq \mr^{(j)}.
$
We emphasize that this reformulation transforms the problem from comparing $n \times n$ matrices to comparing $d \times d$ matrices.

Our test statistic is based on a
{\em Mahalanobis} distance between
$\operatorname{vec}(\hat{\mr}^{(i)})$ and
$\operatorname{vec}(\hat{\mr}^{(j)})$, i.e., by
Theorem~\ref{thm:What_U Rhat What_v^T-R->norm} we have
\begin{equation*}
	(\mathbf{\Sigma}^{(i)}+\mathbf{\Sigma}^{(j)})^{-1/2}
    \operatorname{vec}\big(\mathbf{W}_{\mathbf{U}}^{\top}
    \hat{\mathbf{R}}^{(i)}
    \mathbf{W}_{\mathbf{V}}-\mathbf{W}_{\mathbf{U}}^{\top}
    \hat{\mathbf{R}}^{(j)}
    \mathbf{W}_{\mathbf{V}}-\mathbf{R}^{(i)}+\mathbf{R}^{(j)} -
    \bm{\mu}^{(i)} + \bm{\mu}_j\big) 
	\rightsquigarrow \mathcal{N}\big(\bm{0}, \mathbf{I} \big)
  \end{equation*}
as $n \rightarrow \infty$.
Now suppose the null hypothesis $\mr^{(i)}=\mr^{(j)}$ is true. Then 
$\boldsymbol{\mu}^{(i)}=\boldsymbol{\mu}^{(j)}$ and, with $\mw_* =
\mw_{\mv} \otimes \mw_{\muu}$, we have 
\begin{equation}
\label{eq:HT:initial}
\begin{aligned}
  &
\mathrm{vec}(\hat{\mr}^{(i)} - \hat{\mr}^{(j)})^{\top}
    \mw_* (\mathbf{\Sigma}^{(i)}+\mathbf{\Sigma}^{(j)})^{-1} \mw_*^{\top}
\mathrm{vec}(\hat{\mr}^{(i)} - \hat{\mr}^{(j)})
\rightsquigarrow \chi_{d^2}^{2}
\end{aligned}
\end{equation}
as $n \rightarrow \infty$.
Our objective is to convert Eq.~\eqref{eq:HT:initial} into a
test statistic that depends only on estimates. Toward this
aim, we first define $\hat\mSigma^{(i)}$ as a $d^2\times d^2$ matrix
of the form 
\begin{equation}\label{eq:HT:hatsigma} 
\begin{aligned}
\hat\mSigma^{(i)}=(\hat\mv \otimes \hat\muu)^{\top} \hat\md^{(i)} (\hat\mv \otimes \hat\muu),
\end{aligned}
\end{equation}
where $\hat\md^{(i)}$ is a $n^2 \times n^2$ diagonal matrix whose
diagonal elements are
$$\hat\md^{(i)}_{k_1 + (k_2 - 1)n,
  k_1 + (k_2 - 1)n} = \hat\mpp^{(i)}_{k_1k_2}(1 -
\hat\mpp^{(i)}_{k_1k_2})$$
for any $k_1 \in [n], k_2\in[n]$; here we set
$\hat\mpp^{(i)}=\hat\muu\hat\mr^{(i)}\hat\mv^\top$. The
following lemma shows that
$\big(\hat\mSigma^{(i)}+\hat\mSigma^{(j)}\big)^{-1}$ is a consistent
estimate of
$(\mw_\mv\otimes\mw_\muu)(\mathbf{\Sigma}^{(i)}+\mathbf{\Sigma}^{(j)})^{-1}(\mw_\mv\otimes\mw_\muu)^\top$. 
\begin{lemma1}
\label{lemma:|sigma inverse-sigma|,order of sigma}
	Consider the setting in Theorem~\ref{thm:HT}. 
    We then have
	$$
    \begin{aligned}
        \rho_n\big\|(\mw_\mv\otimes\mw_\muu)(\mSigma^{(i)}+\mSigma^{(j)})^{-1}(\mw_\mv\otimes\mw_\muu)^\top
		-\big(\hat\mSigma^{(i)}+\hat\mSigma^{(j)}\big)^{-1}\big\|
		\lesssim d (n\rho_n)^{-1/2}(\log n)^{1/2}
    \end{aligned}
    $$
    with high probability.
\end{lemma1}
Given Lemma~\ref{lemma:|sigma inverse-sigma|,order of sigma}, the
following result provides a test statistic for $\mathbb{H}_0 \colon
\mr^{(i)} = \mr^{(j)}$ that converges to a central (resp. non-central)
$\chi^2$ under the null (resp. local alternative) hypothesis. 
\begin{theorem1}
\label{thm:HT}
Consider the setting in Theorem~\ref{thm:What_U Rhat What_v^T-R->norm}. Fix $i, j \in [m]$ with $i \neq j$, and let $\hat{\mr}^{(i)}$ and $\hat{\mr}^{(j)}$ be the estimates of $\mr^{(i)}$ and $\mr^{(j)}$ obtained from Algorithm~\ref{Alg}. Suppose $\sigma_{\min}(\mSigma^{(i)} + \mSigma^{(j)}) \asymp \rho_n$, and define the test statistic
\begin{equation*}
    T_{ij} = \mathrm{vec}^\top\big(\hat{\mathbf{R}}^{(i)} - \hat{\mathbf{R}}^{(j)}\big)
    \big(\hat{\mSigma}^{(i)} + \hat{\mSigma}^{(j)}\big)^{-1}
    \mathrm{vec}\big(\hat{\mathbf{R}}^{(i)} - \hat{\mathbf{R}}^{(j)}\big),
\end{equation*}
where $\hat{\mSigma}^{(i)}$ and $\hat{\mSigma}^{(j)}$ are given in Eq.~\eqref{eq:HT:hatsigma}.
Then under the null hypothesis $\mathbb{H}_0 \colon \mr^{(i)} = \mr^{(j)}$, we have 
$
T_{ij} \leadsto \chi^2_{d^2}$ as $n \to \infty.
$
Next, suppose that $\eta > 0$ is a finite constant and that $\mr^{(i)} \neq \mr^{(j)}$ satisfies a local alternative hypothesis such that
\begin{equation*}
    \mathrm{vec}^\top(\mathbf{R}^{(i)} - \mathbf{R}^{(j)})
    (\mathbf{\Sigma}^{(i)} + \mathbf{\Sigma}^{(j)})^{-1}
    \mathrm{vec}(\mathbf{R}^{(i)} - \mathbf{R}^{(j)}) 
    \to \eta.
\end{equation*}
We then have 
$
T_{ij} \leadsto \chi^2_{d^2}(\eta)$ as $n \to \infty,
$
where $\chi^2_{d^2}(\eta)$ is the noncentral chi-square distribution with $d^2$ degrees of freedom and noncentrality parameter $\eta$.
\end{theorem1}

\begin{remark1}\label{rm:Ttildeij}
Theorem~\ref{thm:HT} indicates that, for a chosen significance level $\alpha$, we reject $\mathbb{H}_0$ if $T_{ij} > c_{1-\alpha}$, where $c_{1-\alpha}$
is the $100 \times (1-\alpha)$ percentile of the $\chi^2$ distribution with $d^2$ degrees of freedom. 
Theorem~\ref{thm:HT} is derived based on the normal approximation of $\operatorname{vec}(\mw_\muu^\top \hat{\mr}^{(i)} \mw_\mv - \mr^{(i)})$ in Theorem~\ref{thm:What_U Rhat What_v^T-R->norm} and thus also has the assumption $n \rho_n = \omega(n^{1/2})$; see Remark~\ref{rem:sqrt_n} for further discussion on this $n^{1/2}$ threshold.
If the average degree grows at rate $O(n^{1/2})$, we still have 
$\operatorname{vec}(\mw_\muu^\top \hat{\mr}^{(i)} \mw_\mv - \mr^{(i)}) \to \bm{\mu}^{(i)}$, and thus 
$\operatorname{vec}(\hat{\mr}^{(i)} - \hat{\mr}^{(j)}) \to \mathbf{0}$ under $\mathbb{H}_0$. We can therefore use 
$\tilde{T}_{ij} = \|\hat{\mr}^{(i)} - \hat{\mr}^{(j)}\|_{F}$ as a test statistic and calibrate the rejection region for $\tilde{T}_{ij}$ via bootstrapping. 
We note that $\tilde{T}_{ij}$ is also used as a test statistic in \cite{arroyo2019inference}, but they only assume (and do not theoretically show) that 
$\|\hat{\mr}^{(i)} - \hat{\mr}^{(j)}\|_{F} \to 0$ under the null hypothesis.
\end{remark1}

Theorem~\ref{thm:HT} can also be extended to the multi-sample setting, i.e., testing $\mathbb{H}_0 \colon \mr^{(1)} = \mr^{(2)} = \dots = \mr^{(m)}$ against 
$\mathbb{H}_A \colon \mr^{(i)} \neq \mr^{(j)}$ for some (generally) unknown pair $(i, j)$. 
Our test statistic is then defined as the sum of the (empirical) Mahalanobis distances between $\hat{\mr}^{(i)}$ and 
$\bar{\hat{\mathbf{R}}} = m^{-1} \sum_{i=1}^m \hat{\mathbf{R}}^{(i)}$. More specifically, let
\begin{equation}\label{eq:T_multi}
    T = \sum_{i=1}^m \operatorname{vec}^\top\big(\hat{\mathbf{R}}^{(i)} - \bar{\hat{\mathbf{R}}}\big)
    \big(\bar{\hat{\mathbf{\Sigma}}}\big)^{-1}
    \operatorname{vec}\big(\hat{\mathbf{R}}^{(i)} - \bar{\hat{\mathbf{R}}}\big),
\end{equation}
where $\bar{\hat{\mathbf{\Sigma}}} = m^{-1} \sum_{i=1}^m \hat{\mathbf{\Sigma}}^{(i)}$. Let $\bar{{\mathbf{\Sigma}}} = m^{-1} \sum_{i=1}^m 
\mathbf{\Sigma}^{(i)}$ and suppose $\sigma_{\min}(\bar{\mSigma}) \asymp \rho_n$. Then, under $\mathbb{H}_0 \colon 
\mr^{(1)} = \cdots = \mr^{(m)}$, we have
$
T \leadsto \chi^2_{(m-1)d^2}$ as $n \to \infty.
$
Next, let $\eta > 0$ be a finite constant, and suppose that 
$\{\mr^{(i)}\}$ satisfies a local alternative hypothesis of the form
\begin{equation*}
    \sum_{i=1}^m \operatorname{vec}^\top(\mathbf{R}^{(i)} - \bar{\mathbf{R}})
    (\bar{\mathbf{\Sigma}})^{-1}
    \operatorname{vec}(\mathbf{R}^{(i)} - \bar{\mathbf{R}})
    \to \eta,
\end{equation*}
where $\bar{\mathbf{R}} = m^{-1} \sum_{i=1}^m \mathbf{R}^{(i)}$. 
Then, we also have $T \leadsto \chi^2_{(m-1)d^2}(\eta)$ as $n \to \infty$; see Section~\ref{sec:proof_thmHT} in the supplementary material for a proof sketch of these limiting results.

Thus, for a chosen significance level $\alpha$, we reject $\mathbb{H}_0 \colon \mr^{(1)} = \cdots = \mr^{(m)}$ 
if $T$ exceeds the $100 \times (1-\alpha)$ percentile of the $\chi^2$ distribution with $(m-1)d^2$ degrees of freedom. 
Furthermore, if we reject this $\mathbb{H}_0$, we can perform post-hoc analysis to identify pairs $(i,j)$ where $\mr^{(i)} \neq \mr^{(j)}$ 
by first computing the $p$-values of the test statistics $T_{ij}$ in Theorem~\ref{thm:HT} for all $i \neq j$, 
and then applying Bonferroni correction to these $\tbinom{m}{2}$ $p$-values.
The test statistic in Eq.~\eqref{eq:T_multi} also works for testing the hypothesis $\mathbb{H}_0 \colon \mr^{(i)} = 
\mr^{(i+1)}$ for all $1 \leq i \leq m-1$ against $\mathbb{H}_A \colon \mr^{(i)} \neq \mr^{(i+1)}$ for some possibly unknown 
$i$, which is useful in the context of change-point detection for time series of networks. Once again, if we reject this 
$\mathbb{H}_0$, we can identify the indices $i$ where $\mr^{(i)} \neq \mr^{(i+1)}$ by applying Bonferroni correction to the $p$-values 
of the $T_{i,i+1}$ in Theorem~\ref{thm:HT} for all $1 \leq i \leq m-1$.

\subsection{Related works}
\label{sec:related_works}

Some existing works on multiple networks assume common subspaces
across networks without individual subspaces, i.e., they can be
covered by the COSIE model
$\mathbf{P}^{(i)}=\mathbf{U}\mathbf{R}^{(i)}\mathbf{V}^\top$, but
their theoretical properties remain less complete than those presented here. For instance, when assuming $\mathbf{R}^{(i)}$ are diagonal and
considering undirected networks by setting $\mathbf{U}=\mathbf{V}$,
\cite{nielsen2018multiple,wang2019joint} estimate $\mathbf{U}$ via
alternating gradient descent but provide no error bounds for the
resulting estimates, except in the special case where
$\{\mathbf{R}^{(i)}\}$ are scalars. \cite{arroyo2019inference} 
proposes the COSIE model for undirected networks
and uses the same estimation procedure as Algorithm~\ref{Alg} 
but the theoretical results in \cite{arroyo2019inference} are much weaker than
those presented in the current paper. Indeed, for the estimation of
$\mathbf{U}$, \cite{arroyo2019inference} also provides a Frobenius
norm upper bound for $\hat{\mathbf{U}} \mathbf{W} - \mathbf{U}$ that
is slightly less precise than our Proposition~\ref{thm:Vhat-VW}, but
they do not provide more refined results such as those in
Section~\ref{sec:formal COSIE} (Theorem~\ref{thm:UhatW-U=EVR^{-1}+xxx}
and Theorem~\ref{thm:WhatU_k-U_k->norm}) for the $2 \to \infty$ norm
and row-wise fluctuations of $\hat{\mathbf{U}} \mathbf{W} -
\mathbf{U}$.  Meanwhile, for estimating $\mathbf{R}^{(i)}$,
\cite{arroyo2019inference} shows that
$\operatorname{vec}(\mathbf{W}\hat{\mathbf{R}}^{(i)}\mathbf{W}^\top-\mathbf{R}^{(i)}+\mathbf{H}^{(i)})$
converges to a multivariate normal distribution, but their result does
not yield a proper limiting distribution as it depends on a
non-vanishing and \emph{random} bias term $\mathbf{H}^{(i)}$ which
they can only bound by
$\mathbb{E}(\|\mathbf{H}^{(i)}\|_F)=O(dm^{-1/2})$. In contrast,
Theorem~\ref{thm:What_U Rhat What_v^T-R->norm} shows
$\mathrm{vec}(\mathbf{H}^{(i)}) = \bm{\mu}^{(i)} +
O_{p}((n\rho_n)^{-1/2})$, and thus $\mathbf{H}^{(i)}$ can be replaced
by the \emph{deterministic} term $\bm{\mu}^{(i)}$ in the limiting
distribution. This replacement is essential for subsequent inference;
for example, it allows us to derive the limiting distribution for
two-sample testing of the null hypothesis that two graphs have the
same edge probability matrices (see Section~\ref{sec:COSIE}). This is
also technically challenging as it requires detailed analysis of
$(\hat{\mathbf{U}} \mathbf{W}_{\mathbf{U}} -
\mathbf{U})^{\top}\mathbf{E}^{(i)} (\hat{\mathbf{V}}
\mathbf{W}_{\mathbf{V}}- \mathbf{V})$ using the expansions for
$\hat{\mathbf{U}} \mathbf{W}_{\mathbf{U}} - \mathbf{U}$ and
$\hat{\mathbf{V}}\mathbf{W}_{\mathbf{V}} - \mathbf{V}$ from
Theorem~\ref{thm:UhatW-U=EVR^{-1}+xxx} (see
Sections~\ref{proof:thm3_3} and \ref{sec:technical_lemmas2} for more
details).
{\color{black}
\cite{xie2024bias} also studies the COSIE model for undirected networks and proposes a bias-corrected joint spectral embedding that applies an iterative bias-correction scheme to $\sum_{i=1}^m (\ma^{(i)})^2$. \cite{xie2024bias} focuses on the estimation of the shared subspace $\mathbf{U}$, and establishes $2\to\infty$ norm bounds and entrywise normal approximations for $\hat{\mathbf{U}} \mathbf{W} - \mathbf{U}$. 
}


\cite{jones2020multilayer} considers multiple networks that share a common left subspace but can have possibly different right invariant subspaces, i.e., they assume $\mpp^{(i)}=\muu\mr^{(i)}\mv^{(i)\top}$ where $\muu$ is the common left subspace and $\mr^{(i)},\mv^{(i)}$ are possibly distinct across networks. The resulting model is then a special case of the COISIE model with $\muu^{(i)} = \muu_c$. 
Given a realization $\{\ma^{(i)}\}_{i=1}^{m}$ of these multiple GRDPGs, \cite{jones2020multilayer} defines $\hat{\muu}$ as the $n \times d$ matrix whose columns are the $d$ leading left singular vectors of the $n \times nm$ matrix $[\ma^{(1)} \mid \cdots \mid \ma^{(m)}]$ obtained by concatenating the columns of $\{\ma^{(i)}\}_{i=1}^{m}$, and also define $\hat{\my}$ as the $nm \times d$ matrix whose columns are the $d$ leading (right) singular vectors of $[\ma^{(1)} \mid \cdots \mid \ma^{(m)}]$; $\hat{\my}$ represents an estimate of the column space associated with $\{\mv^{(i)}\}$. They then derive $2 \to \infty$ norm bounds and normal approximations for the rows of $\hat{\muu}$ and $\hat{\my}$.
Their results, at least for estimation of $\hat{\muu}$, are qualitatively worse than ours. For example, Theorem~2 in \cite{jones2020multilayer} implies the bound
$$
\inf_{\mw\in \mathcal{O}_d} \|\hat{\muu} \mw - \muu\|_{2 \to \infty} \lesssim d^{1/2} (n\rho_n)^{-1} \log ^{1/2}{n},
$$
which is worse than the bound obtained from Proposition~\ref{thm:Vhat-VW_part2} by at least a factor of $\rho_n^{-1/2}$; recall that $\rho_n$ can converge to $0$ at rate $\rho_n \succsim n^{-1} \log n$. As another example, \cite{jones2020multilayer} assumes $m$ is fixed, and Theorem~3 in \cite{jones2020multilayer} yields a normal approximation for the rows of $\hat{\muu}$ that is identical to Theorem~\ref{thm:WhatU_k-U_k->norm_part2} of the current paper, but under the much more restrictive assumption $n \rho_n = \omega(n^{1/2})$ instead of $n \rho_n = \omega(\log n)$ in our paper. In addition, \cite{jones2020multilayer} does not discuss the estimation of $\{\mr^{(i)}\}$.

The MultiNeSS model \citep{macdonald2022latent} also assumes multiple networks are composed of the sum of common structure and individual structure, i.e.,
$\mpp^{(i)} = \mx_c \mi_{p_0,q_0} \mx_c^\top + \mx_s^{(i)} \mi_{p_i,q_i} \mx_s^{(i)\top}$ and provides upper bounds for $\min_{\mathbf{W}} \|\hat\mx_c \mw - {\mx}_c\|_{F}$ and $\min_{\mathbf{W}^{(i)}} \|\hat{\mx}_s^{(i)}\mw^{(i)}-\mx_s^{(i)}\|_{F}$; see Remark~\ref{rk:generalized model} for details and a comparison between the results in \cite{macdonald2022latent} and our paper. 
{\color{black}
\cite{tian2026efficient} also considers multiple heterogeneous networks with shared and distinct latent spaces with the model $g^{-1}(\mpp^{(i)}) = \mx_c \mx_c^\top + \mx_s^{(i)} \mx_s^{(i)\top}$ for a link function $g$, and proposes a two-stage procedure consisting of a spectral initialization followed by a likelihood-based refinement to estimate the latent embeddings. \cite{tian2026efficient} establishes rate-optimal error bounds for estimating both the shared and individual latent vectors.
}

There are also some existing works on multilayer SBMs. Recall that
multilayer SBMs assume
$\mathbf{P}^{(i)}=\mathbf{Z}\mathbf{B}^{(i)}\mathbf{Z}^\top$ where
$\mathbf{Z}$ represents community assignments for vertices and
$\mathbf{B}^{(i)}$ are block probability matrices for individual
networks, and this is a special case of the undirected COSIE model
with $\mathbf{U} = (\mathbf{Z}^{\top}\mathbf{Z})^{-1/2} \mathbf{Z}$.
We emphasize that the estimation of $\muu$ in both our
paper and \cite{arroyo2019inference} is based on an 
"estimate-then-aggregate" approach, i.e., we first obtain individual
estimates $\hat{\mathbf{U}}^{(i)}$ of $\mathbf{U}$ from each
$\mathbf{A}^{(i)}$, then aggregate all $\hat{\mathbf{U}}^{(i)}$ to
obtain $\hat{\mathbf{U}}$. In contrast, existing works on multilayer
SBMs (e.g., \cite{paul2020spectral,jing2021community,lei2020bias})
primarily use "aggregate-then-estimate" approaches, i.e., they
aggregate all $\mathbf{A}^{(i)}$ first and then obtain
$\hat{\mathbf{U}}$. For example, \cite{lei2020bias} uses the leading
eigenvectors of the debiased $\sum_{i=1}^m (\mathbf{A}^{(i)})^2$ to
obtain $\hat{\mathbf{U}}$. In general, these two types of methods have
their respective advantages and are complementary to each other.  The
advantage of "aggregate-then-estimate" approaches is that they can
have weaker requirements on the sparsity $\rho_n$ when the number of networks $m$ increases. 
For example, \cite{paul2020spectral} requires $m n
\rho_n = \omega(\log n)$, and \cite{jing2021community} requires $m n
\rho_n = \omega(\log^{4}n)$. In contrast, our "estimate-then-aggregate"
approach needs to guarantee that each individual
$\hat{\mathbf{U}}^{(i)}$ is a consistent estimate of $\mathbf{U}$ and
thus requires $n \rho_n = \Omega(\log n)$.  
If $m$ is bounded then our conditions 
are comparable to those of the "aggregate-then-estimate" approaches. Note that the setting of bounded $m$ is practically relevant as, for
many real-world applications, we only have a small number of graphs even when the number of vertices
in these graphs can be quite large. 

One important advantage of the "estimate-then-aggregate" is that
it is a distributed method and is thus applicable even when, due to certain constraints, the "aggregate-then-estimate" approaches
are infeasible. For instance, when each network is large and stored in
different locations, aggregation of the raw data may be
impractical due to high communication costs, privacy constraints, or
storage limitations at the aggregation site.  
Another important advantage is that both "aggregate-then-estimate"
and "estimate-then-aggregate" approaches can achieve accurate estimation for models with only common subspaces and no individual subspaces (such as multilayer SBMs), but the "aggregate-then-estimate" approaches will fail when individual subspaces are present, while the "estimate-then-aggregate" approach remains effective.
For instance, in the COISIE model, which assumes
$\mathbf{U}^{(i)} = [\mathbf{U}_c \,|\, \mathbf{U}_s^{(i)}]$ to
include possibly distinct individual subspaces $\mathbf{U}_s^{(i)}$,
using the leading eigenvectors of $\sum_{i=1}^m (\mathbf{A}^{(i)})^2$
fails to provide an estimate of $\mathbf{U}_c$; see the
simulation results in Section~\ref{sec:simu_comp} for compelling evidence supporting this claim.
{\color{black}
\cite{ma2026optimal} studies when such ``aggregate-then-estimate'' approaches fail in the presence of individual subspaces.
A similar situation arises for degree-corrected multiple networks with layer-specific degree heterogeneity, where the ``estimate-then-aggregate'' approach can accommodate the degree heterogeneity within each individual estimate while ``aggregate-then-estimate'' approaches again fail.
For further discussion of these points, together with a comparison of the theoretical results in this paper with those of existing works on multilayer SBMs, see Section~\ref{sec:related_works_technical} in the supplementary material.
}




\section{Distributed PCA}
\label{sec:distributed_PCA}

Principal component analysis (PCA) \citep{hotelling1933analysis} is
the most classical and widely applied dimension reduction technique for
high-dimensional data. Standard uses of PCA
involve computing the leading singular vectors of a matrix and thus
generally assume that the data can be stored in memory and/or allowed
for random access. However, massive datasets are now quite prevalent
and these data are often stored
across multiple machines in possibly distant geographic
locations. The communication cost for applying traditional
PCA on these datasets can be rather prohibitive if all the data are
sent to a central location, not to mention that (1) the central location may not have the capability to store and process such
large datasets or (2) due to privacy constraints the raw data cannot
be shared between machines. To meet these challenges, significant efforts have been spent on
designing and analyzing algorithms for PCA in either distributed or
streaming environments; see \cite{garber2017communication,
charisopoulos2021communication, chen2021distributed,
fan2019distributed,pmlr-streaming-PCA} for several recent developments in
this area. 

A succinct description of distributed PCA is as follows. Let
$\{X_j\}_{j=1}^N$ be $N$ iid random vectors in $\mathbb{R}^{D}$ with $X_j \sim \mathcal{N}(\mathbf{0},
\bm{\Sigma})$, and suppose $\{X_j\}$ are scattered across $m$
computing nodes with each node $i$ storing $n_i$ samples.
We denote by $\mathbf{X}^{(i)}$ the $D\times n_i$
matrix formed by the samples stored on the $i$th
node. A natural distributed procedure (see e.g., \cite{fan2019distributed})
for estimating the $d$ leading principal
components $\mathbf{U}$ of $\bm{\Sigma}$ is: (1) each node computes the $D \times d$ matrix $\hat{\mathbf{U}}^{(i)}$ whose columns are the leading eigenvectors of the
sample covariance matrix $\hat{\bm{\Sigma}}^{(i)}= n_i^{-1} \mathbf{X}^{(i)}\mathbf{X}^{(i)\top}$;
(2) $\{\hat{\mathbf{U}}^{(i)}\}_{i=1}^{m}$ are sent to a central node; (3) the central node computes the $D \times d$ matrix
$\hat{\mathbf{U}}$ whose columns are the leading $d$ left
singular vectors of the $D \times dm$ matrix
$\bigl[\hat{\mathbf{U}}^{(1)} \mid \dots \mid \hat{\mathbf{U}}^{(m)}\big]$.

The distributed PCA described above considers the same covariance matrix $\bm{\Sigma}$ across all $m$ computing nodes. We extend this to allow possible heterogeneity across different nodes by assuming that the covariance matrix $\bm{\Sigma}^{(i)}$ for node $i$ shares a common $d_0$-dimensional subspace $\mathbf{U}_c$, but may have possibly distinct $(d_i-d_0)$-dimensional individual subspaces $\mathbf{U}_s^{(i)}$. 
More specifically, we investigate the theoretical properties of distributed PCA assuming 
a spiked covariance structure for $\bm{\Sigma}^{(i)}$, i.e., \begin{equation}
  \label{eq:spiked}
  \bm{\Sigma}^{(i)}=\mathbf{U}^{(i)}\bm{\Lambda}^{(i)}\mathbf{U}^{(i)\top}+\sigma_i^2 (\mathbf{I} - \mathbf{U}^{(i)} \mathbf{U}^{(i)\top}),
\end{equation} 
where $\mathbf{U}^{(i)}=[\mathbf{U}_c|\mathbf{U}_s^{(i)}]\in\mathcal{O}_{D\times d_i}$ and $\bm{\Lambda}^{(i)}$ is a diagonal matrix with diagonal entries $\lambda^{(i)}_1,\ldots,\lambda^{(i)}_{d_i}$ satisfying $\min_{k\in[d_i]}\lambda^{(i)}_k>\sigma_i^2>0$. 
The corresponding distributed PCA estimator 
is presented in Algorithm~\ref{Alg_disPCA}.
\begin{algorithm}[htbp!]
\caption{Distributed PCA}	
\label{Alg_disPCA}
\begin{algorithmic}
\REQUIRE $D\times n_i$ data matrix $\mathbf{X}^{(i)}$ formed by the samples stored on the $i$th node, subspace dimensions $d_1, \dots, d_m$, and common subspace dimension $d_0$. 
\begin{enumerate}
	\item Each node $i \in [m]$ computes the $D\times d_i$ matrix $\hat{\mathbf{U}}^{(i)}$ whose columns are the $d_i$ leading eigenvectors of the sample covariance matrix $\hat{\bm{\Sigma}}^{(i)}= n_i^{-1} \mathbf{X}^{(i)}\mathbf{X}^{(i)\top}$, and sends $\hat{\mathbf{U}}^{(i)}$ to a central node.
	
	\item The central node computes $\hat{\mathbf{U}}_c$ as the $D \times d_0$ matrix whose columns are the leading left singular vectors of $[\hat{\mathbf{U}}^{(1)} \mid \cdots \mid \hat{\mathbf{U}}^{(m)}]$, and sends $\hat{\mathbf{U}}_c$ to all nodes.

    \item Each node $i \in [m]$ computes $\hat{\mathbf{U}}_s^{(i)}$ as the $D \times (d_i - d_{0})$ matrix whose columns are the leading left singular vectors of $(\mathbf{I} - \hat{\mathbf{U}}_c \hat{\mathbf{U}}_c^\top) \hat{\mathbf{U}}^{(i)}$.
	
\end{enumerate} 
\ENSURE $\hat{\mathbf{U}}_c, \{\hat{\mathbf{U}}_s^{(i)}
\}_{i=1}^m$.
\end{algorithmic}
\end{algorithm}

Covariance matrices of the form in Eq.~\eqref{eq:spiked} are studied extensively in the high-dimensional statistics literature; see e.g., \cite{Spiked,sparse_AOS,sparse_AOS1,sparse,Tony_sparse,bai_sample_cov} and the references therein. A common assumption for $\mathbf{U}^{(i)}$ is that it is sparse, e.g., the $\ell_{q}$ quasi-norms, for some $q \in [0,1]$, of the columns of $\mathbf{U}^{(i)}$ are bounded. Note that sparsity of $\mathbf{U}^{(i)}$ also implies sparsity of $\bm{\Sigma}^{(i)}$. In this paper we do not impose sparsity constraints on $\mathbf{U}^{(i)}$ but instead assume that $\mathbf{U}^{(i)}$ has bounded coherence, i.e., $\|\mathbf{U}^{(i)}\|_{2 \to \infty} \lesssim D^{-1/2}$. The resulting $\bm{\Sigma}^{(i)}$ will no longer be sparse. Bounded coherence is also a natural condition in the context of covariance matrix estimation; see e.g., \cite{cape2019two, yan2021inference, chen2021distributed, xie2018bayesian}, as it allows for the spiked eigenvalues $\bm{\Lambda}^{(i)}$ to grow with $D$ while also guaranteeing that the entries of the covariance matrix $\bm{\Sigma}^{(i)}$ remain bounded, i.e., there is a large gap between the spiked eigenvalues and the remaining eigenvalues. In contrast, if $\mathbf{U}^{(i)}$ is sparse then the spiked eigenvalues $\bm{\Lambda}^{(i)}$ grow with $D$ if and only if the variances and covariances in $\bm{\Sigma}^{(i)}$ also grow with $D$, and this can be unrealistic in many settings as increasing the dimension of the $X_j$ (e.g., by adding more features) should not change the magnitude of the existing features.

We now state the analogues of Theorem~\ref{thm:UhatW-U=EVR^{-1}+xxx_part2}, Theorem~\ref{thm:WhatU_k-U_k->norm_part2} and Proposition~\ref{thm:Vhat-VW_part2} in the setting of distributed PCA. For simplicity (and with minimal loss of generality), we assume $n_i \equiv n=\lfloor N/m \rfloor$. We emphasize that these results should be interpreted in the regime where both $n$ and $D$ are arbitrarily large or $n, D \rightarrow \infty$.
\begin{theorem1}
  \label{thm:PCA_UhatW-U=(I-UU^T)EVR^{-1}+xxx_previous}
Suppose we have $m$ computing nodes and each node $i$ stores $n$ iid mean zero $D$-dimensional multivariate Gaussian random samples with covariance matrix $\bm{\Sigma}^{(i)}$ of the form in Eq.~\eqref{eq:spiked} with common subspace $\mathbf{U}_c$ and individual subspace $\mathbf{U}_s^{(i)}$. Let $\hat{\mathbf{U}}_c$ be the estimate of $\mathbf{U}_c$ obtained by Algorithm~\ref{Alg_disPCA}. Suppose $\sigma_i^2\lesssim 1$, $\|\mathbf{U}^{(i)}\|_{2\to\infty}\lesssim \sqrt{d_i/D}$, $\lambda^{(i)}_1\asymp \lambda^{(i)}_{d_i}\asymp D^{\gamma}$ for some $\gamma\in(0,1]$, and suppose there exists a constant $c_s>0$ such that $\|m^{-1}\sum_{i=1}^m\muu_s^{(i)}\muu_s^{(i)\top}\|\leq 1-c_s$.
Let $r_i = \mathrm{tr}(\bm{\Sigma}^{(i)})/\lambda^{(i)}_1$ be the effective rank of $\bm{\Sigma}^{(i)}$ and $r=\max_{i\in[m]}{r_i}$. Define $ \varphi = \Bigl({\max\{r, \log D\}/n}\Bigr)^{1/2}.$ Let $\mathbf{W}_{\muu_c}$ minimize $\|\hat{\mathbf{U}}_c \mathbf{O} - \mathbf{U}_c\|_{F}$ over all $d_0 \times d_0$ orthogonal matrix $\mathbf{O}$. Then when $n=\omega (\max\{D^{1 - \gamma},\log D\})$ we have
\begin{equation}
  \label{eq:expansion_PCA1}
  \hat\muu_c \mw_{\muu_c}-\muu_c
=\frac{1}{m}\sum_{i=1}^m(\mi-\muu^{(i)}\muu^{(i)\top})(\hat{\mSigma}^{(i)} -\mSigma^{(i)}) \muu_c(\mLambda_c^{(i)})^{-1}+\mq_{\muu_c},
\end{equation}
 where $\mLambda^{(i)}_c$ is the principal submatrix of $\mLambda^{(i)}$ containing only the eigenvalues corresponding to the common subspace $\muu_c$, and $\mq_{\muu_c}$ is a random matrix satisfying
$$
\begin{aligned}
	&\|\mq_{\muu_c}\|
	\lesssim  D^{-\gamma/2} \tilde{\varphi}+D^{-\gamma} \varphi + \varphi^2
\end{aligned}
$$
with high probability.  
Furthermore, when $n=\omega(D^{2-2\gamma}\log D)$, we have
\begin{equation}
  \label{eq:mq_2toinf_pca}
	\|\mq_{\muu_c}\|_{2\to\infty}
\lesssim d_{\max}^{1/2}\tilde{\varphi}
[D^{-(\gamma+1)/2} +D^{-3\gamma/2} (1 + D \tilde{\varphi})]
\end{equation}
with high probability, where $d_{\max}=\max_{i\in[m]}d_i$ and $\tilde\varphi = n^{-1/2} \log^{1/2}{D}$. 

For each $i\in[m]$, let $\hat\muu_s^{(i)}$ be the estimation of $\muu_s^{(i)}$ obtained by Algorithm~\ref{Alg_disPCA}, and let $\mw_{\muu_s}^{(i)}$ be the minimizer of $\|\hat\muu_s^{(i)}\mo-\muu_s^{(i)}\|_{F}$ over all $(d_i-d_0)\times (d_i-d_0)$ orthogonal matrices $\mo$. Then
$$
\begin{aligned}
	\hat\muu_s^{(i)}\mw_{\muu_s}^{(i)}-\muu_s^{(i)}
	=(\mi-\muu^{(i)}\muu^{(i)\top})(\hat{\mSigma}^{(i)} -\mSigma^{(i)}) \muu_s^{(i)}(\mLambda_s^{(i)})^{-1}
+\mq_{\muu_s}^{(i)},
\end{aligned}
$$
where $\mLambda^{(i)}_s$ is the principal submatrix of $\mLambda^{(i)}$ containing only the eigenvalues corresponding to the common subspace $\muu_s^{(i)}$, and the random matrix
$\mq_{\muu_s}^{(i)}$ satisfies the same upper bounds as those for $\mq_{\muu_c}$.

Note that if all nodes share the common subspace $\muu_c$ without any individual subspaces $\{\muu_s^{(i)}\}$, then the bounds become
	$\|\mq_{\muu_c}\|
	\lesssim  D^{-\gamma}\varphi+\varphi^2
	$,
$\|\mq_{\muu_c}\|_{2\to\infty}
	\lesssim d_{\max}^{1/2}D^{-3\gamma/2}\tilde\varphi(1 + D \tilde{\varphi})$
with high probability.
\end{theorem1}
\begin{remark1}
Theorem~\ref{thm:PCA_UhatW-U=(I-UU^T)EVR^{-1}+xxx_previous} assumes that the $d_i$ leading (spiked) eigenvalues of $\bm{\Sigma}^{(i)}$ grow with $D$ at rate $D^{\gamma}$ for some $\gamma \in (0,1]$ while the remaining (non-spiked) eigenvalues remain bounded. Under this condition the effective rank of $\mSigma^{(i)}$ satisfies $r_i = \mathrm{tr}(\mSigma^{(i)})/\lambda_1^{(i)} \asymp D^{1 - \gamma}$ and thus $\gamma < 1$ and $\gamma \geq 1$ correspond to the cases where $r_i$ is growing with $D$ and remains bounded, respectively. The effective rank $r_i$ serves as a measure of the complexity of $\bm{\Sigma}^{(i)}$; see e.g., \cite{vershynin_hdp,tropp2,luo_xiao_bunea}. The condition $n = \omega(\max\{D^{1-\gamma}, \log D\})$ assumed for Eq.~\eqref{eq:expansion_PCA1} is thus very mild as we are only requiring the sample size in each node to grow slightly faster than the effective ranks $\{r_i\}$. Similarly the slightly more restrictive condition $n = \omega(D^{2 - 2 \gamma} \log D)$ for Eq.~\eqref{eq:mq_2toinf_pca} is also quite mild as it leads to much stronger (uniform) row-wise concentration for $\mq$. If $\gamma = 1$ then the above two conditions both simplify to $n = \omega(\log D)$ and thus allow for the dimension $D$ to grow exponentially with $n$. Finally, Theorem~\ref{thm:PCA_UhatW-U=(I-UU^T)EVR^{-1}+xxx_previous} also holds for $\gamma > 1$, with the only minor change being that the sample size requirement for Eq.~\eqref{eq:mq_2toinf_pca} continues to be $n = \omega(\log D)$ for $\gamma > 1$.
\end{remark1}

\begin{remark1}\label{rmk:T0T,condition}
The proof of Theorem~\ref{thm:PCA_UhatW-U=(I-UU^T)EVR^{-1}+xxx_previous} (see Section~\ref{Appendix:A_thm:PCA_UhatW-U=(I-UU^T)EVR^{-1}+xxx_previous}) is almost identical to that of Theorem~\ref{thm:UhatW-U=EVR^{-1}+xxx_part2} for the COISIE model (see Section~\ref{Appendix:thm:UhatW-U=EVR^{-1}+xxx_part2}). More specifically, after deriving an expansion for $\hat{\muu}^{(i)} \mw^{(i)} - \muu^{(i)}$ for each $i \in [m]$ (see Lemma~\ref{lemma:PCA_UhatW-U_previous} in Section~\ref{Appendix:A_thm:PCA_UhatW-U=(I-UU^T)EVR^{-1}+xxx_previous}), we apply Theorem~\ref{thm3_general} to obtain expansions for $\hat{\muu}_c$ and $\hat\muu_s^{(i)}$ based on these individual expansions for $\{\hat{\muu}^{(i)}\}$. 
{\color{black} 
Concretely, the individual expansions required to apply Theorem~\ref{thm3_general} take the form $\hat{\muu}^{(i)} \mw^{(i)} - \muu^{(i)} = \mt_0^{(i)} + \mt^{(i)}$, where the leading term $\mt_0^{(i)}$ equals $\me^{(i)} \mv^{(i)} (\mr^{(i)})^{-1}$ with $\me^{(i)} = \ma^{(i)} - \mpp^{(i)}$ for the COISIE model, and $(\mi-\muu^{(i)}\muu^{(i)\top})(\hat{\bm{\Sigma}}^{(i)} - \bm{\Sigma}^{(i)}) \muu^{(i)}(\bm{\Lambda}^{(i)})^{-1}$ for distributed PCA, while $\mt^{(i)}$ collects the higher-order remainders. In both settings, a simple sufficient condition for Eq.~\eqref{eq:cond_pis} is $\| m^{-1} \sum_{i=1}^m \muu_s^{(i)} \muu_{s}^{(i)\top}\| \leq 1 - c_s$ for some constant $c_s>0$, which is guaranteed by the last condition in Assumption~\ref{ass:main_part2} for the COISIE model and by the assumptions of Theorem~\ref{thm:PCA_UhatW-U=(I-UU^T)EVR^{-1}+xxx_previous} for distributed PCA; see Remark~\ref{rem:assumptions} for an intuitive explanation of this condition.
}
 We also note that the main difference between the leading terms in Theorem~\ref{thm:UhatW-U=EVR^{-1}+xxx_part2} and Theorem~\ref{thm:PCA_UhatW-U=(I-UU^T)EVR^{-1}+xxx_previous} is the appearance of the projection matrix $(\mi-\muu^{(i)}\muu^{(i)\top})$ (note that $\muu_c(\mLambda_c^{(i)})^{-1}=\muu(\mLambda^{(i)})^{-1}\muu^\top\muu_c$ and $\muu_s^{(i)}(\mLambda_s^{(i)})^{-1}=\muu(\mLambda^{(i)})^{-1}\muu^\top\muu_s^{(i)}$). This difference arises from the individual expansions for $\hat{\muu}^{(i)}$, and this is because for the COISIE model, $\mpp^{(i)} = \muu^{(i)} \mr^{(i)} \mv^{(i)\top}$ are low-rank matrices while for distributed PCA the matrices $\mSigma^{(i)} = \muu^{(i)} \bm{\Lambda}^{(i)} \muu^{(i)\top} + \sigma_i^2\muu^{(i)}_{\perp}\muu^{(i)\top}_{\perp}$ are not necessarily low-rank.

\end{remark1}

\begin{proposition1} \label{prop:PCA_||hatUW-U||_previous}
	Consider the setting and assumptions ($n=\omega(\max\{D^{1-\gamma},\log D\})$) in Theorem~\ref{thm:PCA_UhatW-U=(I-UU^T)EVR^{-1}+xxx_previous}. We then have
	$$
	\begin{aligned}
		&\|\hat\muu_c\mw_{\muu_c}-\muu_c\|_F
\lesssim \sqrt{\frac{d_0\max\{r,\log D\}}{mn}}
+\sqrt{\frac{d_0\log D}{D^\gamma n}}
	+ \sqrt{\frac{d_0\max\{r,\log D\}}{D^{2\gamma}n}}
	+\frac{d_0^{1/2}\max\{r,\log D\}}{n},\\
	&\|\hat\muu_s^{(i)}\mw_{\muu_s}^{(i)}-\muu_s^{(i)}\|_F
\lesssim \sqrt{\frac{(d_i-d_0)\max\{r,\log D\}}{n}}
	\end{aligned}
$$
	with high probability.
	Furthermore, if $m=O(D^{\gamma}\max\{r,\log D\}/\log D)$, $m=O(D^{2\gamma})$, and $m=O(n/\max\{r,\log D\})$ we have
	$$
	\|\hat\muu_c\mw_{\muu_c}-\muu_c\|_F
	\lesssim \sqrt{\frac{d_0\max\{r,\log D\}}{N}}
	$$
	with high probability.
	
	Note that if the nodes have no individual subspaces, then the bound for $\|\hat\muu_c\mw_{\muu_c}-\muu_c\|_F$ becomes
	$
	\sqrt{\frac{d_0\max\{r,\log D\}}{mn}}
	+ \sqrt{\frac{d_0\max\{r,\log D\}}{D^{2\gamma}n}}
	+\frac{d_0^{1/2}\max\{r,\log D\}}{n}
	$.
\end{proposition1}

\begin{remark1}
For the case where the covariance matrix is shared across all nodes, i.e., $\mSigma^{(i)}\equiv\mSigma$, we have $\muu^{(i)}\equiv \muu_c$, and Proposition~\ref{prop:PCA_||hatUW-U||_previous} becomes almost identical to Theorem~4 in \cite{fan2019distributed}, except that \cite{fan2019distributed} presented their results in terms of the $\psi_1$ Orlicz norm for $\|\hat\muu\mw-\muu\|_F$.
We note that for fixed $D$ and $\gamma$, the error bound in Proposition~\ref{prop:PCA_||hatUW-U||_previous} converges to zero at rate $N^{-1/2}$, where $N$ is the total number of samples, and is thus reminiscent of the error rate for traditional PCA (where $m = 1$) in the low-dimensional setting; see also the asymptotic covariance matrix $\bm{\Upsilon}_{\muu_c}$ in the following Theorem~\ref{thm:PCA_(UhatW-U)_k->normal_previous} (specifically, when $\mSigma^{(i)}\equiv\mSigma$, we have $\bm{\Upsilon}_{\muu_c}=\frac{1}{N}\sigma^2\mLambda^{-1}$).
\end{remark1}

\begin{theorem1}
	\label{thm:PCA_(UhatW-U)_k->normal_previous}
	Consider the setting in
    Theorem~\ref{thm:PCA_UhatW-U=(I-UU^T)EVR^{-1}+xxx_previous}. 
    Define $\bm{\Upsilon}_{\muu_c}$ as the $d_0\times d_0$ symmetric matrix $\bm{\Upsilon}_{\muu_c}
    =\frac{1}{Nm}\sum_{i=1}^m \sigma_i^2(\mLambda^{(i)}_c)^{-1}$.
    Then for the $k$th row $\hat u_{c,k}$ and $u_{c,k}$ of $\hat\muu_c$ and $\muu_c$, when $m=o(D^{2\gamma}/\log D)$ and $m=o(n/(D^{2-2\gamma}\log^2 D))$, we have 
	$$
\bm{\Upsilon}_{\muu_c}^{-1/2}\big(\mw_{\muu_c}^\top\hat u_{c,k}-u_{c,k}\big)
\leadsto \mathcal{N}\big(\mathbf{0},\mathbf{I}\big)
$$
as $n,D \rightarrow \infty$.

For each $i\in[m]$, define $\bm{\Upsilon}_{\muu_s}^{(i)}$ as the $(d_i-d_0)\times (d_i-d_0)$ symmetric matrix 
$\bm{\Upsilon}_{\muu_s}^{(i)}
    =\frac{1}{n} \sigma_i^2(\mLambda^{(i)}_s)^{-1}$.
    Then for the $k$th row $\hat u_{s,k}^{(i)}$ and $u_{s,k}^{(i)}$ of $\hat\muu_s^{(i)}$ and $\muu_s^{(i)}$, when $n=\omega(D^{2-2\gamma}\log^2 D)$ we have 
	$$
(\bm{\Upsilon}_{\muu_s}^{(i)})^{-1/2}\big(\mw_{\muu_s}^{(i)\top}\hat u_{s,k}^{(i)}-u_{s,k}^{(i)}\big)
\leadsto \mathcal{N}\big(\mathbf{0},\mathbf{I}\big)
$$
as $n,D \rightarrow \infty$.

\end{theorem1}


\begin{remark1}
The condition that the number of distributed machines $m$ cannot be too large also appears in other distributed estimation settings, including $M$-estimation and PCA. 
More specifically, suppose a dataset is split across $m$ nodes with each node having $n$ observations. Theorems~4.1 and 4.2 in \cite{huo2019aggregated} present error bounds for distributed $M$-estimation, and the optimal rate $N^{-1/2}$ is achieved and the central limit theorem holds when $m = O(n)$.  
Similarly, Eq.~(4.6) and Eq.~(4.7) of \cite{fan2019distributed} show that distributed PCA where all nodes share the common covariance matrix achieves the same estimation error rate as that of traditional PCA when $m = O(n)$.

In addition, the condition $m=o(n/(D^{2-2\gamma}\log^2 D))$ stated in Theorem~\ref{thm:PCA_(UhatW-U)_k->normal_previous} is imposed purely for ease of exposition, as the normal approximation for $\mw_{\muu_c}^\top\hat u_{c,k}-u_{c,k}$ when $m=o(n/(D^{2-2\gamma}\log D))$ requires more tedious book-keeping of $\|q_k\|$. See Remark~\ref{rem:condition m} for Theorem~\ref{thm:UhatW-U=EVR^{-1}+xxx_part2} for similar discussions.
\end{remark1}



\begin{remark1}
  \label{rem:demean}
  For ease of exposition, the previous results are stated under the assumption that $\mathbb{E}[X_j^{(i)}]$ is known and thus, without loss of generality, we can assume $\mathbb{E}[X_j^{(i)}] = \mathbf{0}$. If $\mathbb{E}[X_j^{(i)}]$ is unknown, we have to demean the data before performing PCA.
More specifically, let 
$\tilde\mOmega^{(i)}=\frac{1}{n}\sum_{j=1}^n(X_j^{(i)}-\bar{X}^{(i)})(X_j^{(i)}-\bar{X}^{(i)})^{\top}$
be the sample covariance matrix for the $i$th server, where
$\bar{X}^{(i)}=\frac{1}{n}\sum_{j=1}^nX_j^{(i)}$.
Then, with
$\hat\mOmega^{(i)}=\frac{1}{n}\sum_{j=1}^n(X_j-\bm{\mu}^{(i)})(X_j-\bm{\mu}^{(i)})^{\top}$, we have
$$
\underbrace{\tilde\mOmega^{(i)}-\mOmega^{(i)}}_{\me^{(i)}}
=\underbrace{\hat\mOmega^{(i)}-\mOmega^{(i)}}_{\me_1^{(i)}}
-\underbrace{(\bar{X}^{(i)}-\bm{\mu}^{(i)})(\bar{X}^{(i)}-\bm{\mu}^{(i)})^\top}_{\me_2^{(i)}}.
$$
Bounds for $\me_1^{(i)}$ are provided in the proof of Theorem~\ref{thm:PCA_UhatW-U=(I-UU^T)EVR^{-1}+xxx_previous}. Since $\bar X^{(i)}\sim \mathcal{N}(\bm\mu^{(i)},\mSigma^{(i)}/n)$, we have
$$
\begin{aligned}
	&\|\me_2^{(i)}\|\lesssim n^{-1/2}  D^\gamma \varphi,
	\quad \|\me_2^{(i)}\|_{\infty}
	\lesssim n^{-1/2} D^\gamma \tilde\varphi
\end{aligned}
$$
with high probability. We thus obtain, from 
Eq.~(8.16)
in \cite{chen2021classification}, that
$$
\begin{aligned}
	\|\me_2^{(i)}\muu^{(i)}\|_{2\to\infty}
	\lesssim d_i^{1/2} \Big(\frac{D^\gamma}{n}\sqrt{\frac{d_i}{D}}
	+\frac{\max\{r,\log D\}}{n}D^{\gamma/2}\Big)
\end{aligned}
$$
with high probability. Therefore, $\|\me_2^{(i)}\|$, $\|\me_2^{(i)}\|_{\infty}$, and $\|\me_2^{(i)}\muu^{(i)}\|_{2\to\infty}$ are all of smaller order than the corresponding terms for $\me_1^{(i)}$. Consequently, we can ignore all terms depending on $\me_2^{(i)}$ in the proofs of Theorem~\ref{thm:PCA_UhatW-U=(I-UU^T)EVR^{-1}+xxx_previous}, Theorem~\ref{thm:PCA_(UhatW-U)_k->normal_previous}, and Proposition~\ref{prop:PCA_||hatUW-U||_previous}; that is, these results continue to hold even when $\mathbb{E}[X_j^{(i)}]$ is unknown.
\end{remark1}

\begin{remark1}
  \label{rem:general_covariance}
The theoretical results in this section can be easily extended to the case where 
$\mSigma^{(i)}=\muu^{(i)}\mLambda^{(i)}\muu^{(i)\top} + \mathbf{U}^{(i)}_{\perp}
\bm{\Lambda}_{\perp}^{(i)} \mathbf{U}_{\perp}^{(i)\top}$
with $\muu^{(i)}=[\muu_c|\muu_s^{(i)}]$, $\lambda_1^{(i)}\asymp \lambda_d^{(i)}\asymp D^\gamma$ for all $i$, and $\max_{i} \|\mLambda_{\perp}^{(i)}\| \leq M$ for some finite constant $M > 0$ that does not depend on $m$, $n$, and $D$. 
Under this setting, the expansions in Theorem~\ref{thm:PCA_UhatW-U=(I-UU^T)EVR^{-1}+xxx_previous} still hold, while the limit result in Theorem~\ref{thm:PCA_(UhatW-U)_k->normal_previous} holds with covariance matrices
$\bm{\Upsilon}_{\muu_c}
= \frac{1}{Nm}\sum_{i=1}^m \zeta_{kk}^{(i)} (\mLambda_c^{(i)})^{-1},
\bm{\Upsilon}_{\muu_s}^{(i)}
= \frac{1}{n}\zeta_{kk}^{(i)} (\mLambda_s^{(i)})^{-1},
$
where $\zeta_{kk}^{(i)}$ is the $k$-th diagonal element of $\mathbf{U}^{(i)}_{\perp}
\bm{\Lambda}_{\perp}^{(i)} \mathbf{U}_{\perp}^{(i)\top}$.

Finally, all results in this section can also be generalized to the case where the $X$ are only sub-Gaussian. Indeed, the same bounds (up to constant factors) for $\hat{\mSigma}^{(i)} - \mSigma^{(i)}$ as those presented in the current paper are also available in the sub-Gaussian setting; see, e.g., \cite{Lounici2017,chen2021classification,chen2021spectral}. Thus, the arguments presented in the supplementary material still carry through. 
The only minor change is in the expressions of covariance matrices in Theorem~\ref{thm:PCA_(UhatW-U)_k->normal_previous}. 
Specifically, if $X^{(i)}$ has mean $\bm{0}$ and is sub-Gaussian, then
$$
\begin{aligned}
	\bm{\Upsilon}_{\muu_c} &= \frac{1}{Nm}\sum_{i=1}^m[\zeta_{i,k}
\otimes \muu_c (\mLambda_c^{(i)})^{-1}]^\top 
\bm{\Xi}^{(i)} [\zeta_{i,k}
\otimes \muu_c (\mLambda_c^{(i)})^{-1}],
	\\
	\bm{\Upsilon}_{\muu_s}^{(i)} &= \frac{1}{n}[\zeta_{i,k}
\otimes \muu_s^{(i)} (\mLambda_s^{(i)})^{-1}]^\top 
\bm{\Xi}^{(i)} [\zeta_{i,k}
\otimes \muu_s^{(i)} (\mLambda_s^{(i)})^{-1}],
\end{aligned}
$$
where $\zeta_{i,k}$ is the $k$-th row of $\mi - \muu^{(i)} \muu^{(i)\top}$ and $\bm{\Xi}^{(i)} = \mathrm{Var}[\mathrm{vec}(X^{(i)}X^{(i)\top})]$ contains the fourth-order (mixed) moments of $X^{(i)}$ and thus need not depend only on $\mSigma^{(i)}$. 
In the special case when $X^{(i)} \sim \mathcal{N}(\bm{0}, \bm{\Sigma}^{(i)})$, we have $\mathrm{Var}[\mathrm{vec}(X^{(i)}X^{(i)\top})] = (\mSigma^{(i)} \otimes \mSigma^{(i)})(\mathbf{I}_{D^2} + \mathcal{K}_{D})$, where $\mathcal{K}_{D}$ is the $D^2 \times D^2$ commutation matrix, and this implies the expressions for $\bm{\Upsilon}_{\muu_c}$ and $\bm{\Upsilon}_{\muu_s}^{(i)}$ in Theorem~\ref{thm:PCA_(UhatW-U)_k->normal_previous} (see Eq.~\eqref{eq:var_pca_1}).
\end{remark1}

\subsection{Related works}
\label{sec:related_works_pca}

We begin by comparing our results for distributed PCA in the setting where $\mSigma^{(i)} \equiv \mSigma = \muu\mLambda\muu^\top + \sigma^2\mi$ with $\muu \in \mathbb{R}^{D \times d}$, against the minimax bound for traditional PCA (where all $N = nm$ observations are centralized on a single node) provided in \cite{cai2013sparse}. For ease of exposition, we state these comparisons in terms of the $\sin$-$\Theta$ distance between subspaces, as these are equivalent to the corresponding Procrustes distances. 
Let $\Theta$ be the family of spiked covariance matrices of the form
\begin{equation*}
\begin{aligned}
	\muu\mLambda\muu^\top+\sigma^2 \mi \colon 
	&C_2 D^{\gamma} \leq \lambda_d\leq\dots\leq\lambda_1\leq
   C_1D^{\gamma}, \muu\in\mathbb{R}^{D\times d},\muu^\top\muu=\mi_d,
\end{aligned}
\end{equation*}
where $C_1$, $C_2$, $\sigma^2$, and $\gamma \in (0,1]$ are fixed constants. 
Then for any $\mSigma \in \bm{\Theta}$, we have from Proposition~\ref{prop:PCA_||hatUW-U||_previous} that
\begin{equation}
	\label{Eq:PCA_rate1_previous}
\|\sin \Theta(\hat{\muu}, \muu)\|_{F}^2 \lesssim \frac{\sigma^2 d \max\{D^{1-\gamma},\log D\}}{N}
\end{equation}
with high probability, provided that $\muu$ has \emph{bounded coherence}. 
Meanwhile, by Theorem~1 in \cite{cai2013sparse}, the \emph{minimax} error rate for the class $\Theta$ is 
\begin{equation}\label{Eq:PCA_rate2_previous}
\begin{aligned}
\inf_{\tilde{\muu}} \sup_{\mSigma\in\Theta}
  \mathbb{E} \|\sin \Theta(\tilde{\muu}, \muu)\|_{F}^2
\asymp \frac{\sigma^2 dD^{1-\gamma}}{N},
\end{aligned}
\end{equation}
where the infimum is taken over all estimators $\tilde{\muu}$ of $\muu$. 
If $\gamma<1$, then the error rate in Eq.~\eqref{Eq:PCA_rate1_previous} for distributed PCA is the same as that in Eq.~\eqref{Eq:PCA_rate2_previous} for traditional PCA, while if $\gamma=1$, then there is a (multiplicative) gap of order at most $\log D$ between the two error rates. Note, however, that Eq.~\eqref{Eq:PCA_rate1_previous} provides a high-probability bound for $\|\hat\muu\mw-\muu\|_F^2$, which is a slightly stronger guarantee than the expected value in Eq.~\eqref{Eq:PCA_rate2_previous}.

We now compare our results with existing results for distributed PCA 
in \cite{garber2017communication, charisopoulos2021communication, chen2021distributed}. 
Note that the existing literature on distributed PCA assumes that all nodes share a common covariance matrix, therefore we compare the results under the setting $\mSigma^{(i)} \equiv \mSigma = \muu\mLambda\muu^\top + \sigma^2 (\mi - \muu\muu^\top)$ where $\muu \in \mathbb{R}^{D \times d}$, and thus $\muu^{(i)} \equiv \muu_c = \muu$ and $d_i \equiv d_0$. 
We remark at the outset that our $\| \cdot \|_{2 \to \infty}$ norm 
bound for $\hat{\muu}$ in Theorem~\ref{thm:PCA_UhatW-U=(I-UU^T)EVR^{-1}+xxx_previous} and the
row-wise normal approximations of $\hat{u}_k$ in
Theorem~\ref{thm:PCA_(UhatW-U)_k->normal_previous} are, to the best of
our knowledge, novel. Previous theoretical analyses for distributed
PCA have focused exclusively on the coarser Frobenius norm error of $\hat{\muu}$ and
$\muu$. 
\cite{garber2017communication} proposes a
procedure for estimating the 
leading eigenvector of $\muu$ by aligning all local estimates
(using sign-flips) to a reference solution and then averaging the
aligned local estimates. \cite{charisopoulos2021communication} 
extends this procedure to handle multiple eigenvectors by employing orthogonal
Procrustes transformations to align the local estimates. Let $\hat{\muu}^{(P)}$ denote the resulting
estimate of $\muu$. Theorem~4 
in \cite{charisopoulos2021communication} gives
\begin{equation}
  \label{eq:charisopoulos2021communication}
  \|\sin \Theta(\hat{\muu}^{(P)}, \muu)\| \lesssim \sqrt{\frac{d(r+\log n)}{N}} + \frac{\sqrt{d}(r+\log m)}{n},
\end{equation}
with high probability.
 The error rates for $\hat{\muu}$ and $\hat{\muu}^{(P)}$ are therefore 
almost identical; cf. Eq.~\eqref{Eq:PCA_rate1_previous}. 
\cite{chen2021distributed} considers distributed estimation of $\muu$ by aggregating the eigenvectors $\{\hat{\muu}^{(i)}\}_{i=1}^{m}$ associated with subspaces of $\{\hat{\mSigma}^{(i)}\}_{i=1}^{m}$ whose dimensions are slightly larger than that of $\muu$. While the aggregation scheme in \cite{chen2021distributed} is considerably more complicated than that studied in \cite{fan2019distributed} and the current paper, it also requires possibly weaker eigengap conditions, and thus a detailed comparison between the two sets of results is perhaps not meaningful. Nevertheless, if we assume the above setting, then Theorem~3.3 in \cite{chen2021distributed} yields an error bound for $\sin \Theta(\hat{\muu}, \muu)$ equivalent to Eq.~\eqref{Eq:PCA_rate1_previous}.

In this paper, we assume that $D$ grows with $n$, as the case where $D$ is fixed has been addressed in several classic works. For example, Theorem~13.5.1 in \cite{anderson1962introduction} states that $\mathrm{vec}(\hat{\mathbf{U}} - \mathbf{U})$ converges to a multivariate normal distribution in $\mathbb{R}^{D^2}$, provided that the eigenvalues of $\bm{\Sigma}$ are distinct. This result is subsequently extended to the case where the $\{X_j\}_{j=1}^{N}$ are from an elliptical distribution with possibly non-distinct eigenvalues (see Sections~3.1.6 and 3.1.8 of \cite{kollo2005advanced}) or when they only have finite fourth-order moments \citep{davis1977asymptotic}. These cited results are for the joint distribution of \emph{all} rows of $\hat{\mathbf{U}}$ and are thus slightly stronger than the row-wise results presented in this paper, which currently only imply that the joint distribution for any \emph{finite} collection of rows of $\hat{\muu}$ converges to multivariate normal.
{\color{black}
We also discuss in Section~\ref{sec:related_works_technical_pca} of the supplementary material how our row-wise normal approximations for distributed PCA relate to the classical eigenvector and eigenvalue central limit theorems for spiked covariance models in the regime where the dimension $D$ grows proportionally with the sample size $n$ \citep{paul2007,wang_fan_2017,fan_li_xia_zheng}.}
  

Finally, we present another variant of Theorem~\ref{thm:PCA_UhatW-U=(I-UU^T)EVR^{-1}+xxx_previous} and Theorem~\ref{thm:PCA_(UhatW-U)_k->normal_previous}, but with different assumptions on $n$ and $D$.
More specifically, rather than basing our
analysis on the sample covariances $\hat{\mSigma}^{(i)} =
\frac{1}{n} \sum_{j} X_j^{(i)} X_j^{(i)\top}$, we instead view each
$X_j^{(i)}$ as $Y_j^{(i)} + Z_j^{(i)}$ where $Y_j^{(i)} \stackrel{\mathrm{iid}}{\sim}
\mathcal{N}(\mathbf{0},\muu^{(i)}(\mLambda^{(i)}-\sigma_i^2\mi)\muu^{(i)\top})$ and
$Z_j^{(i)} \stackrel{\mathrm{iid}}{\sim}
\mathcal{N}(\mathbf{0},\sigma_i^2\mi)$ represent the ``signal'' and
``noise'' components, respectively. 
Let $\my^{(i)}=(Y^{(i)}_1,\dots,Y^{(i)}_n)$, 
$\mz^{(i)}=(Z^{(i)}_1,\dots,Z^{(i)}_n)$ and note that 
$$
\my^{(i)}=\muu^{(i)}(\mLambda^{(i)}-\sigma_i^2\mi)^{1/2}\mathbf{F}^{(i)},
\text{ where }\mathbf{F}^{(i)}=(F^{(i)}_1,\dots,F^{(i)}_n) \in
\mathbb{R}^{d \times n} \,\, \text{with } 
F^{(i)}_k\stackrel{\mathrm{iid}}{\sim} \mathcal{N}(\mathbf{0},\mi_d).
$$
The column space of $\my^{(i)}$ is, almost surely, the same as that
spanned by $\muu^{(i)}$. Furthermore the leading eigenvectors $\hat{\muu}^{(i)}$ 
of $\hat{\mSigma}^{(i)}$ are also the leading left singular vectors of
$\mathbf{X}^{(i)}$ and thus they can be considered as a noisy perturbation of
the leading left singular vectors of $\my^{(i)}$ (see Section~3 in
\cite{yan2021inference} for more details). Note that $\mz^{(i)}$ has 
  mutually independent entries; in contrast, the entries of
  $\hat{\mSigma}^{(i)} - \mSigma^{(i)}$ are
{\em dependent}. 
We then have the following results. 
\begin{theorem1}
\label{thm:PCA_UhatW-U=(I-UU^T)EVR^{-1}+xxx}
Consider the same setting as that in Theorem~\ref{thm:PCA_UhatW-U=(I-UU^T)EVR^{-1}+xxx_previous}.
Then when $\frac{\log^3(n+D)}{\min\{n,D\}}\lesssim 1,\phi :=
  \frac{(n+D)\log(n+D)}{nD^\gamma}\ll 1.$ we have
  \begin{equation*}
  \begin{aligned}
  	 &\hat\muu_c \mw_{\muu_c}-\muu_c
=\frac{1}{m}\sum_{i=1}^m\mz^{(i)}(\my^{(i)})^{\dagger} \muu_c +\mq_{\muu_c}  \end{aligned}
 \end{equation*}
where $(\cdot)^{\dagger}$ denotes the Moore-Penrose pseudo-inverse 
and the residual matrix $\mq_{\muu_c}$ satisfies
$$
\begin{aligned}
	\|\mq_{\muu_c}\|_{2\to\infty}
	&\lesssim  \frac{d_{\max} \phi}{(n+D)^{1/2}}
    	+\frac{d_{\max} \phi}{D^{1/2}\log(n+D)}
    	+\frac{ d_{\max}  \phi^{3/2} D^{1/2} \log^{1/2}(n+D)}{(n+D)}
    	+\frac{ d_{\max} \phi^{1/2} \log^{1/2}(n+D)}{(n+D)^{1/2}D^{1/2}}
\end{aligned}
$$
with probability at least $1-O((n+D)^{-10})$. 
  
 For each $i\in[m]$, we have
 \begin{equation*}
  \begin{aligned}
  	 &\hat\muu_s^{(i)}\mw_{\muu_s}^{(i)}-\muu_s^{(i)}
=\mz^{(i)}(\my^{(i)})^{\dagger} \muu_s^{(i)} +\mq_{\muu_s}^{(i)}  \end{aligned}
 \end{equation*}
 where the random matrix
$\mq_{\muu_s}^{(i)}$ satisfies the same upper bound as that for $\mq_{\muu_c}$. 
\end{theorem1}
\begin{theorem1}
	\label{thm:PCA_(UhatW-U)_k->normal}
	Consider the setting in
    Theorem~\ref{thm:PCA_UhatW-U=(I-UU^T)EVR^{-1}+xxx}. Then when
    $m=o\Big(\frac{nD^\gamma}{(n+D)\log^2(n+D)}\big)$ and $
m=o\Big(D^{1+\gamma}/n\Big),$
    we have
$$
\bm{\Upsilon}_{\muu_c}^{-1/2}\big(\mw_{\muu_c}^\top\hat u_{c,k}-u_{c,k}\big)
\leadsto \mathcal{N}\big(\mathbf{0},\mathbf{I}\big)
$$
as $n,D \rightarrow \infty$. 

And for each $i\in[m]$, when
    $\frac{(n+D)\log^2(n+D)}{nD^\gamma}=o(1)$ and $
n/D^{1+\gamma}=o(1),$
we have 
	$$
(\bm{\Upsilon}_{\muu_s}^{(i)})^{-1/2}\big(\mw_{\muu_s}^{(i)\top}\hat u_{s,k}^{(i)}-u_{s,k}^{(i)}\big)
\leadsto \mathcal{N}\big(\mathbf{0},\mathbf{I}\big)
$$
as $n,D \rightarrow \infty$.
Here $\bm{\Upsilon}_{\muu_c}$ and $\bm{\Upsilon}_{\muu_s}^{(i)}$ are defined in Theorem~\ref{thm:PCA_(UhatW-U)_k->normal_previous}.
\end{theorem1}
  As we mentioned above, the conclusions of
  Theorems~\ref{thm:PCA_UhatW-U=(I-UU^T)EVR^{-1}+xxx_previous} and
  \ref{thm:PCA_(UhatW-U)_k->normal_previous} are the same as those in 
 Theorem~\ref{thm:PCA_UhatW-U=(I-UU^T)EVR^{-1}+xxx}
  and Theorem~\ref{thm:PCA_(UhatW-U)_k->normal}. In particular, for the estimate error for $\muu_c$, the leading
  term $m^{-1} \sum_{i=1}^{m} (\mi-\muu^{(i)}\muu^{(i)\top})(\hat{\mSigma}^{(i)} -\mSigma^{(i)}) \muu_c^{(i)}(\mLambda_c^{(i)})^{-1}$ in
Theorem~\ref{thm:PCA_UhatW-U=(I-UU^T)EVR^{-1}+xxx_previous} is
equivalent to the leading term $m^{-1} \sum_{i=1}^{m} \mz^{(i)}
  (\my^{(i)})^{\dagger} \muu_c$ in
  Theorem~\ref{thm:PCA_UhatW-U=(I-UU^T)EVR^{-1}+xxx}; see the
  derivations in Section~\ref{sec:equivalence_thm7_thm9} for more
  details. And the covariance matrix $\bm{\Upsilon}_{\muu_c}$ in
Theorem~\ref{thm:PCA_(UhatW-U)_k->normal_previous} is identical to
that in Theorem~\ref{thm:PCA_(UhatW-U)_k->normal}. Thus, the only
difference between these sets of results is in the conditions assumed for $n$,$D$ and $m$ (see Table~\ref{tb:PCA}). 
More specifically, Theorem~\ref{thm:PCA_UhatW-U=(I-UU^T)EVR^{-1}+xxx_previous} and 
Theorem~\ref{thm:PCA_(UhatW-U)_k->normal_previous} only
requires $n$ to be large compared to $D$, e.g. in Theorem~\ref{thm:PCA_(UhatW-U)_k->normal_previous} $n =
\omega(mD^{2-2\gamma}\log^2 D)$, while
Theorem~\ref{thm:PCA_UhatW-U=(I-UU^T)EVR^{-1}+xxx}  and
Theorem~\ref{thm:PCA_(UhatW-U)_k->normal} require $n$ to be large
but not too large compared to $D$, e.g. in Theorem~\ref{thm:PCA_(UhatW-U)_k->normal}
$mD^{1 - \gamma} \log^{2} D \ll n \ll D^{1+ \gamma}/m$.
The main reason behind these discrepancies is in the noise structure in
$\mz^{(i)}$ (independent entries) compared to $\me^{(i)} = \hat{\mSigma}^{(i)} -
\mSigma^{(i)}$ (dependent entries). For example, if $D$ is fixed then $\|\me^{(i)}\|
\rightarrow 0$ in probability and $\|\mSigma^{(i)}\| \asymp D^{\gamma}$. 
In contrast, for a fixed $D$ we have $n^{-1/2} \|\mz^{(i)}\|
\rightarrow \sigma_i^2$ as $n \rightarrow \infty$ but $n^{-1/2} \|\my^{(i)}\|
\asymp D^{\gamma/2}$ with high probability. The signal to noise
ratio ($\|\mSigma^{(i)}\|/\|\me^{(i)}\|$) in 
Theorem~\ref{thm:PCA_UhatW-U=(I-UU^T)EVR^{-1}+xxx_previous} thus behaves
quite differently from the signal to noise ratio ($\|\my^{(i)}\|/\|\mz^{(i)}\|$) 
in Theorem~\ref{thm:PCA_UhatW-U=(I-UU^T)EVR^{-1}+xxx} as $n$
increases. Finally, for fixed $m$ if $\gamma > 1/3$ then $D^{1 + \gamma} \gg D^{2 -
  2\gamma}$ and, by combining
Theorems~\ref{thm:PCA_(UhatW-U)_k->normal_previous} and
\ref{thm:PCA_(UhatW-U)_k->normal}, 
$\bm{\Upsilon}^{-1/2}_{\muu_c}(\mathbf{W}^{\top}_{\muu_c} \hat{u}_{c,k} - u_{c,k})
\rightsquigarrow \mathcal{N}(0, \mathbf{I})$ under the very mild condition of $n \gg D^{1 - \gamma}$. 
Similar remarks also hold for the estimation error of $\muu_s^{(i)}$.
\begin{table}[htbp]
\centering
	\begin{tabular}{cc}
	\toprule
	Result &Conditions\\
	\midrule
	Theorem~\ref{thm:PCA_UhatW-U=(I-UU^T)EVR^{-1}+xxx_previous} &\tabincell{c}{$\frac{D^{2-2\gamma}\log D}{n}=o(1)$}\\
	Theorem~\ref{thm:PCA_UhatW-U=(I-UU^T)EVR^{-1}+xxx} &\tabincell{c}{
                                               $\frac{D^{1-\gamma}\log D}{n}=o(1)$, $\frac{\log^3 D}{n}=O(1)$,
                                               $\frac{\log n}{D^\gamma}=o(1)$ and $\frac{\log^3 n}{D}=O(1)$
                                                }
      \\\midrule
Theorem~\ref{thm:PCA_(UhatW-U)_k->normal_previous}&\tabincell{c}{$m=o\Big(\frac{n}{D^{2-2\gamma}\log^2 D}\Big)$ and $m=o\Big(\frac{D^{2\gamma}}{\log D}\Big)$}\\
	Theorem~\ref{thm:PCA_(UhatW-U)_k->normal}&\tabincell{c}{
                                               $m=o\Big(\frac{n}{D^{1-\gamma}\log^2 D}\Big)$,
                                               $m=o\Big(\frac{D^{\gamma}}{\log^2 n}\Big)$,
                                               $m=o\Big(\frac{D^{1+\gamma}}{n}\Big)$ and $\frac{\log^3D}{n}
                                               = O(1)$ }
      \\

	\bottomrule
\end{tabular}
\caption{Relationship between $n,D$ and $m$ assumed for $2\to\infty$
  norm bounds and asymptotic normality of $\mw^{\top}_{\muu_c} \hat{u}_{c,k} - u_{c,k}$.}\label{tb:PCA}
\end{table}




\section{Simulation Results and Real Data Experiments}
\label{sec:simulation and real data}


\subsection{COISIE model}\label{sec:simu_COISIE}
We now present simulations to validate our theoretical results for the COISIE model.
We consider the COISIE model with $n = 1000$, $m = 3$, $d_i \equiv 4$, and $d_{0,\muu} = d_{0,\mv} = 2$, resulting in $\muu_c$, $\mv_c$, $\muu^{(i)}_s$, and $\mv^{(i)}_s$ each being $1000 \times 2$ matrices. 
The orthonormal matrices $\muu_c$ and $\mv_c$ are randomly generated. For each $i$, orthonormal matrices $\muu^{(i)}_s$ and $\mv^{(i)}_s$, which are orthogonal to $\muu_c$ and $\mv_c$, respectively, are also randomly generated, and $\mr^{(i)}$ is constructed, with its entries independently drawn from the uniform distribution $U(0,n)$. We then obtain the underlying matrices $\mpp^{(i)}=[\muu_c|\muu^{(i)}_s]\mr^{(i)}[\mv_c|\mv^{(i)}_s]^\top$.
As mentioned in Remark~\ref{rmk:general error}, the COISIE model can be generalized to settings with bounded error or sub-Gaussian error, and our theoretical results remain valid in these cases. For each Monte Carlo replicate, we generate error matrices $\me^{(i)}$ whose entries are independently drawn from the Gaussian distribution with mean zero and variance $0.5^2$. The observed matrices are then given by $\ma^{(i)}=\mpp^{(i)}+\me^{(i)}$.
We then apply Algorithm~\ref{Alg_COISIE} to obtain estimated common subspaces and individual subspaces.

\begin{figure}[htbp!]
\centering
\subfigure
{\includegraphics[height=4.2cm]{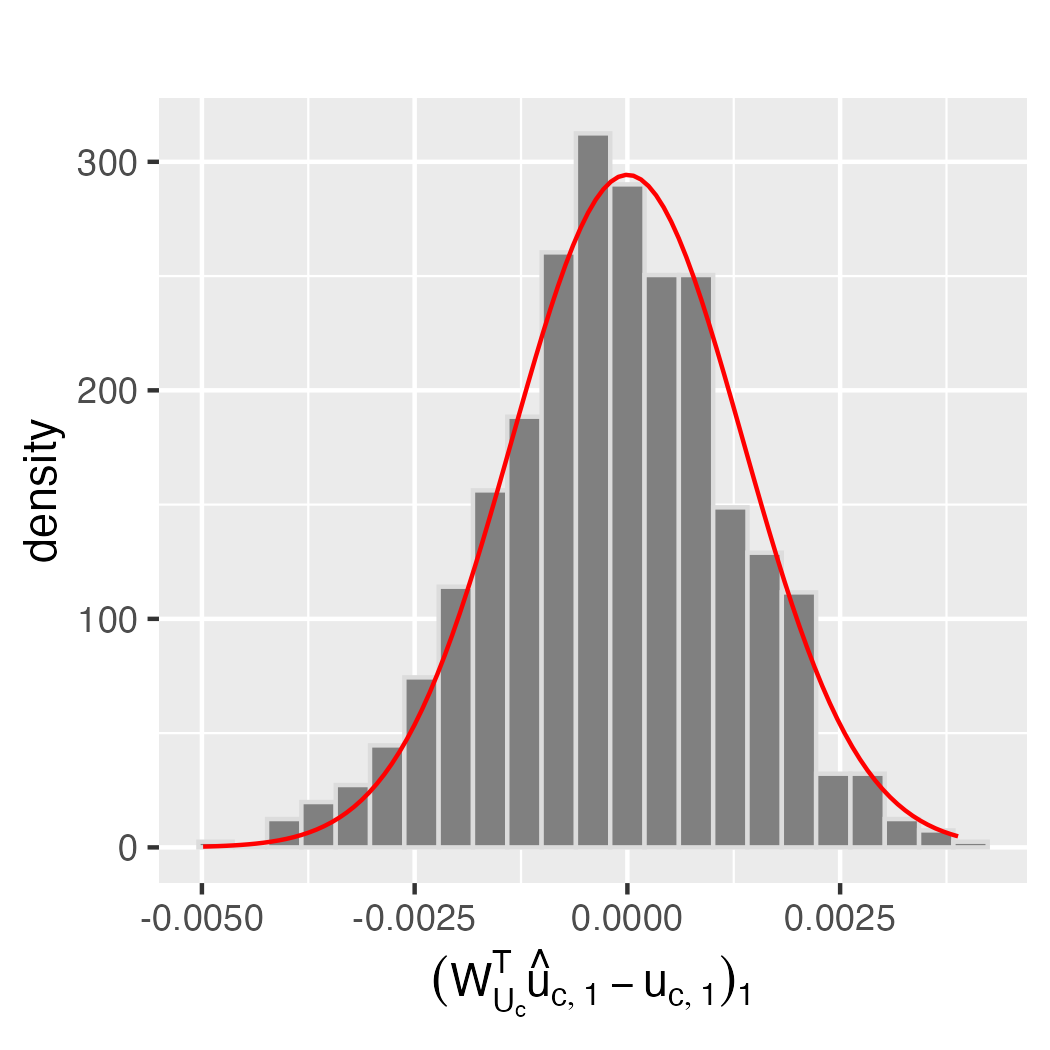}}
\subfigure
{\includegraphics[height=4.2cm]{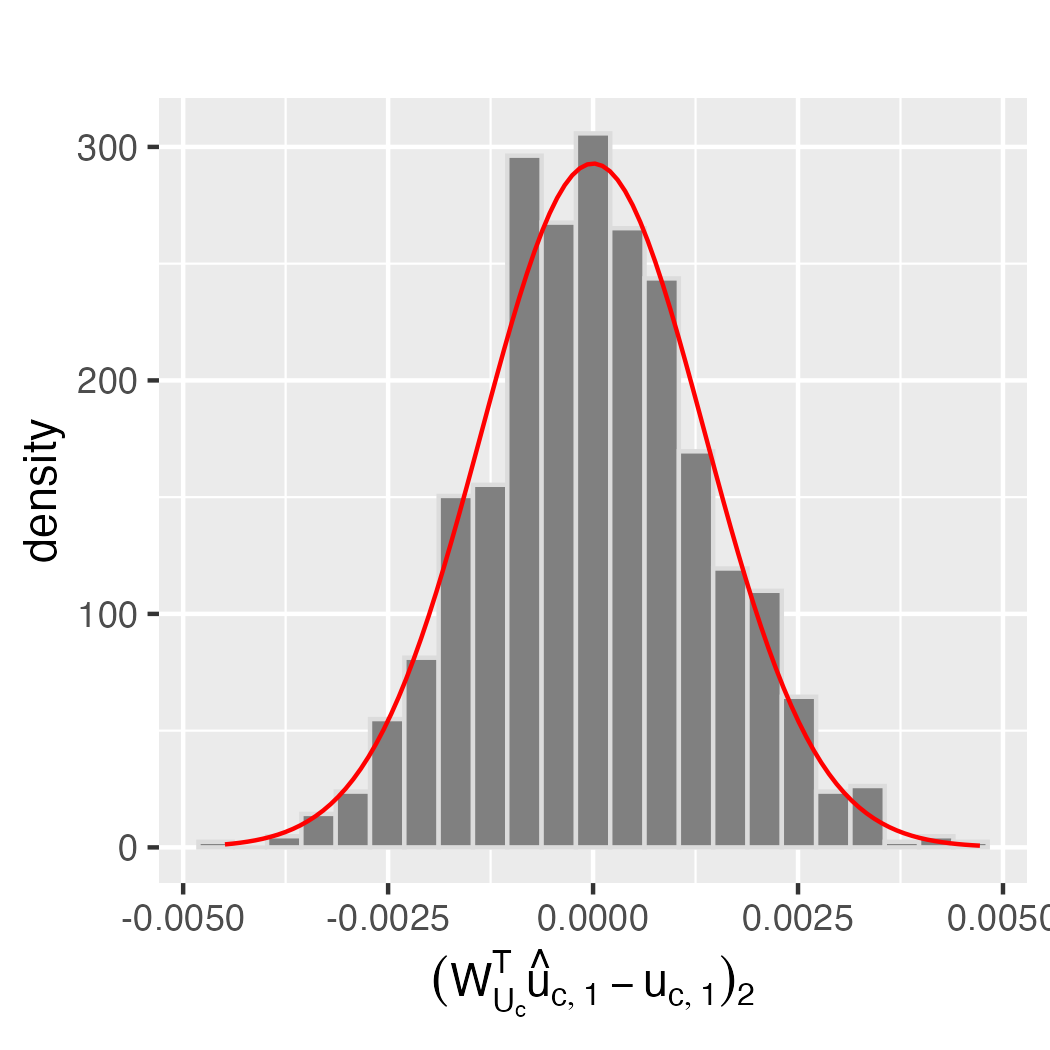}}
\subfigure
{\includegraphics[height=4.2cm]{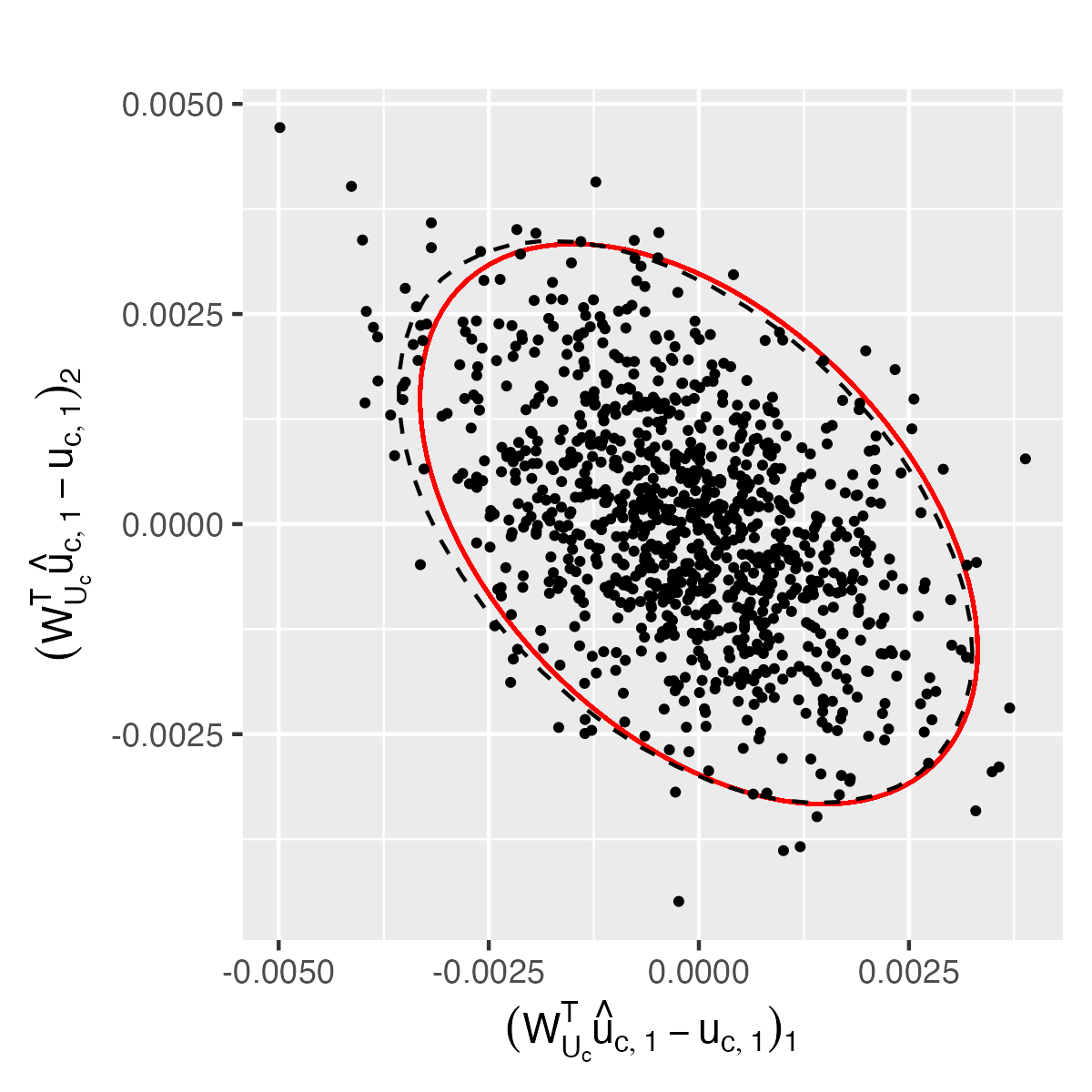}}

\caption{\footnotesize 
The left two panels are histograms of the empirical distributions of the entries of the estimation error $\mw^\top_{\muu_c}\hat{u}_{c,k} - {u}_{c,k}$ for $k = 1$. These histograms are based on $1000$ independent Monte Carlo replicates of the COISIE model with $n = 1000$, $m = 3$, $d_i \equiv 4$, $d_{0,\muu} = d_{0,\mv} = 2$. The red lines represent the probability density functions of the normal distributions with parameters specified in Theorem~\ref{thm:WhatU_k-U_k->norm_part2}. 
The right panel displays a bivariate plot of the empirical distributions of the entries. The dashed black ellipses represent 95\% level curves for the empirical distributions, while the solid red ellipses represent 95\% level curves for the theoretical distributions as specified in Theorem~\ref{thm:WhatU_k-U_k->norm_part2}.
}
\label{fig:simulation_CLT1}
\end{figure}

\begin{figure}[htbp!]
\centering
\subfigure
{\includegraphics[height=4.2cm]{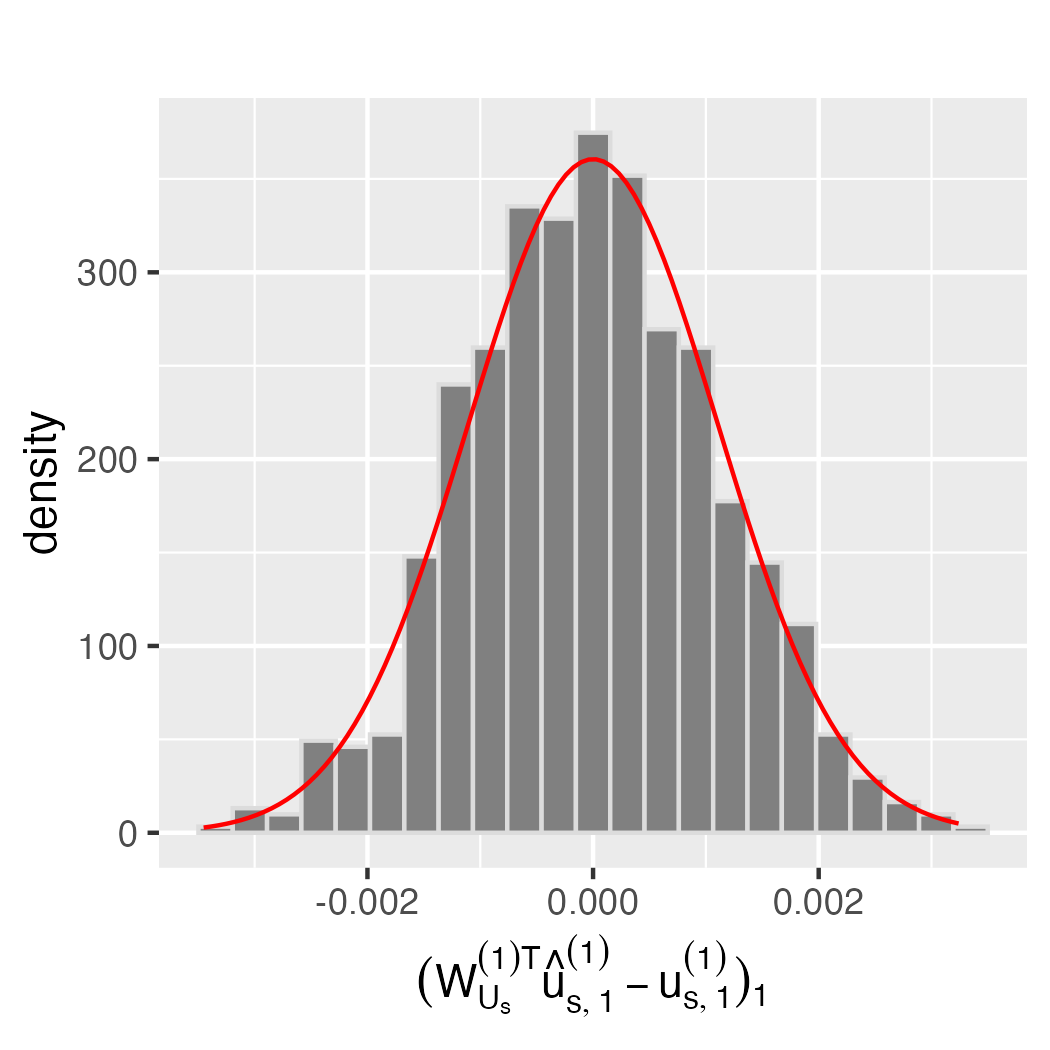}}
\subfigure
{\includegraphics[height=4.2cm]{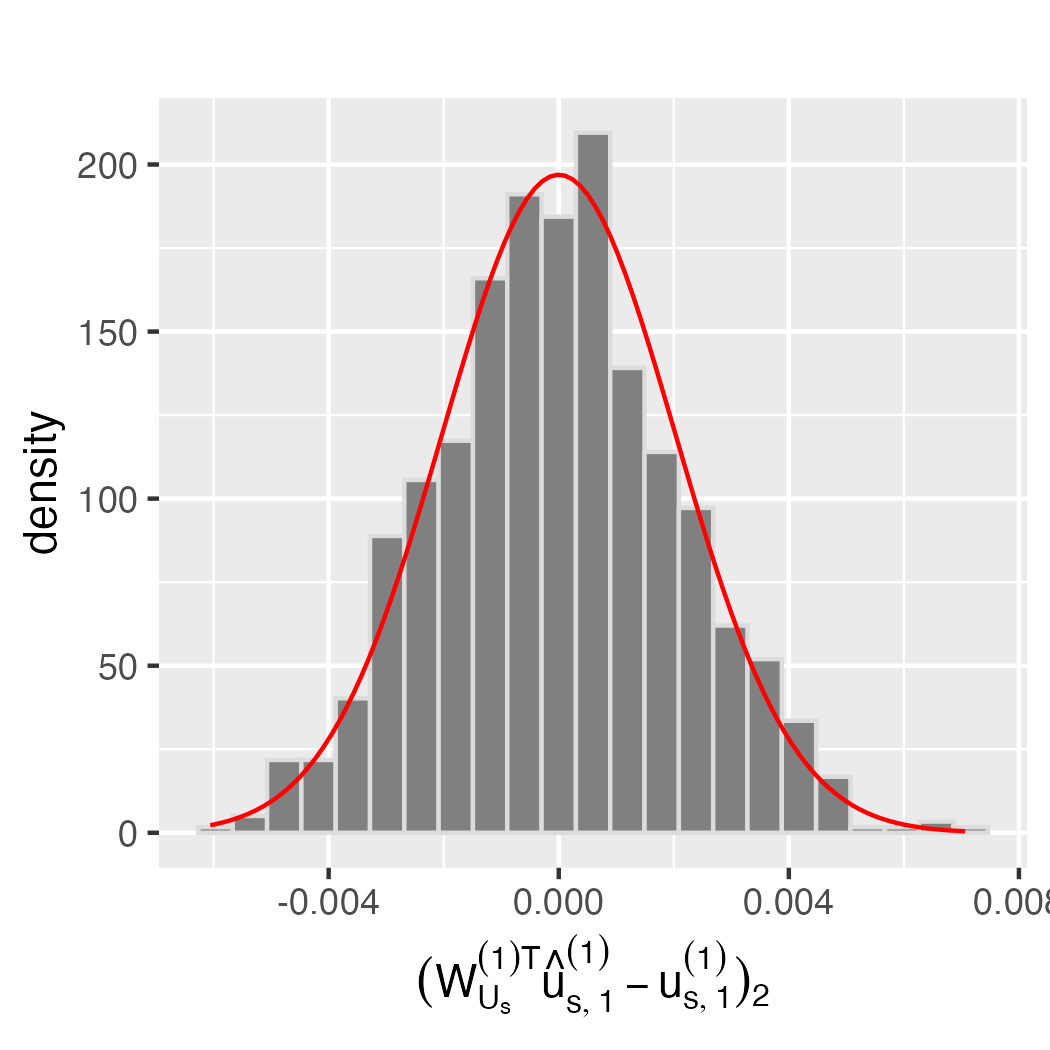}}
\subfigure
{\includegraphics[height=4.2cm]{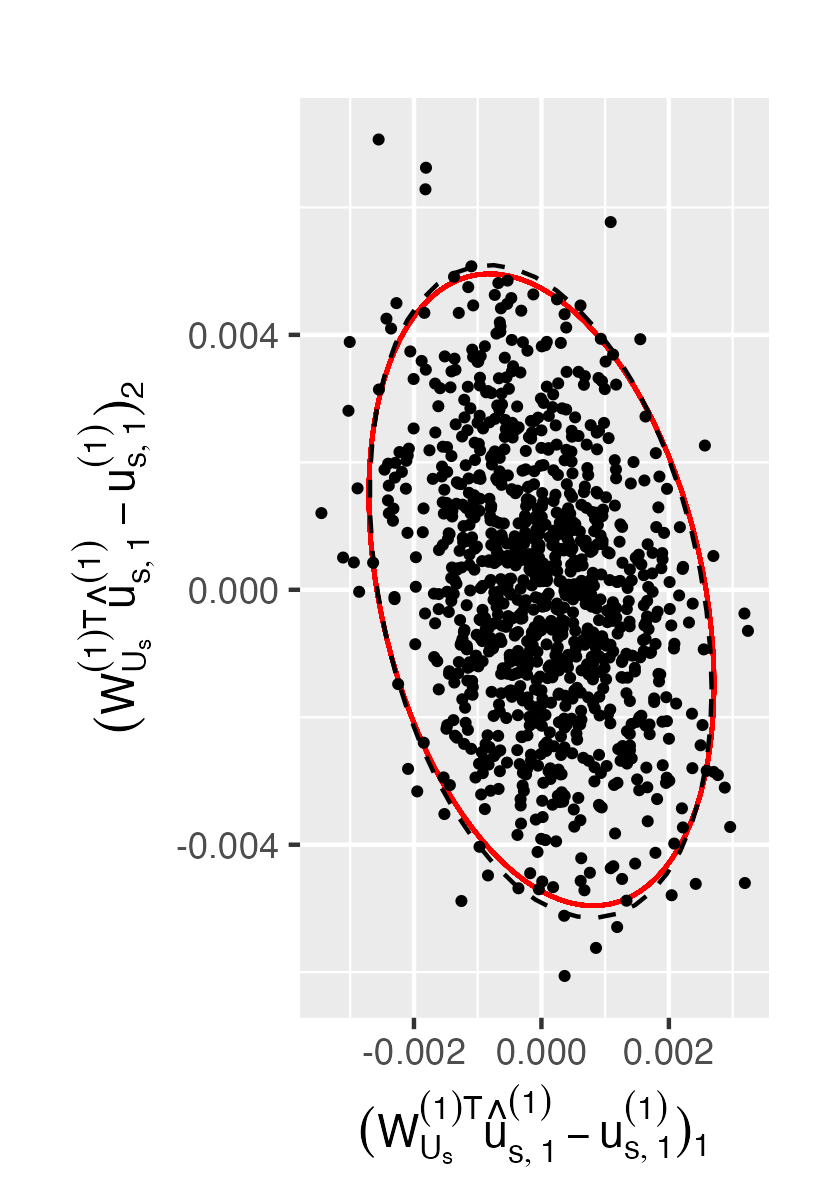}}

\caption{\footnotesize 
Histograms and a bivariate plot of the empirical distributions of the entries of the estimation error $\mw^{(i)\top}_{\muu_s}\hat{u}^{(i)}_{s,k} - {u}^{(i)}_{s,k}$ for $i = 1$ and $k = 1$ are presented. Refer to Figure~\ref{fig:simulation_CLT1} for more details.
}
\label{fig:simulation_CLT2}
\end{figure}

We conduct $1000$ independent Monte Carlo replicates to obtain empirical distributions of the estimation errors $\mw^\top_{\muu_c}\hat{u}_{c,k}-{u}_{c,k}$ for $k=1$ and $\mw^{(i)\top}_{\muu_s}\hat{u}^{(i)}_{s,k}-{u}^{(i)}_{s,k}$ for $i=1,k=1$, which we then compare against the limiting distribution given in Theorem~\ref{thm:WhatU_k-U_k->norm_part2}. The results are summarized in Figure~\ref{fig:simulation_CLT1} and Figure~\ref{fig:simulation_CLT2}. Henze-Zirkler's normality test indicates that the empirical distributions of $\mw^\top_{\muu_c}\hat{u}_{c,k}-{u}_{c,k}$ and $\mw^{(i)\top}_{\muu_s}\hat{u}^{(i)}_{s,k}-{u}^{(i)}_{s,k}$ are well-approximated by multivariate normal distributions, and the figures furthermore show that the empirical covariance matrices are very close to the theoretical covariance matrices.

\subsection{COSIE model and the two-sample hypothesis testing}
We next demonstrate the theoretical results for the COSIE model. Specifically, we consider the setting of \emph{directed} multilayer SBMs on $n = 2000$ vertices, with $m = 3$ graphs and $K = 3$ blocks. For each vertex $v$, we randomly generate the \emph{outgoing} and \emph{incoming} community assignments $\tau(v)$ and $\phi(v)$, i.e., the $\tau(v)$ are iid random variables with $\mathbb{P}[\tau(v) = k] = 1/3$ for $k \in \{1,2,3\}$, and similarly for $\phi(v)$.  
Next let $\mathbf{Z}_{\tau}$ be the $n \times 3$ matrix where $(\mathbf{Z}_{\tau})_{vk} = 1$ if $\tau(v) = k$ and $(\mathbf{Z}_{\tau})_{vk} = 0$ otherwise, 
and define $\mathbf{Z}_{\phi}$ analogously. 
Then for each $i$, we randomly generate the $3 \times 3$ block probability matrix $\mathbf{B}^{(i)}$, with entries independently drawn from $U(0,1)$, and set
$\mpp^{(i)} = \mathbf{Z}_{\tau} \mathbf{B}^{(i)} \mathbf{Z}_{\phi}^\top$.
For each Monte Carlo replicate, we randomly generate observed adjacency matrices $\mathbf{A}^{(i)}$ according to $\mpp^{(i)}$, and estimate $\muu = \mathbf{Z}_{\tau} (\mathbf{Z}_{\tau}^\top \mathbf{Z}_{\tau})^{-1/2},  
\mv = \mathbf{Z}_{\phi} (\mathbf{Z}_{\phi}^\top \mathbf{Z}_{\phi})^{-1/2},  
\mr^{(i)} = (\mathbf{Z}_{\tau}^\top \mathbf{Z}_{\tau})^{1/2} \mathbf{B}^{(i)} (\mathbf{Z}_{\phi}^\top \mathbf{Z}_{\phi})^{1/2}$ using Algorithm~\ref{Alg}.

\begin{figure}[htbp!]
\centering
\subfigure
{\includegraphics[width=3.5cm]{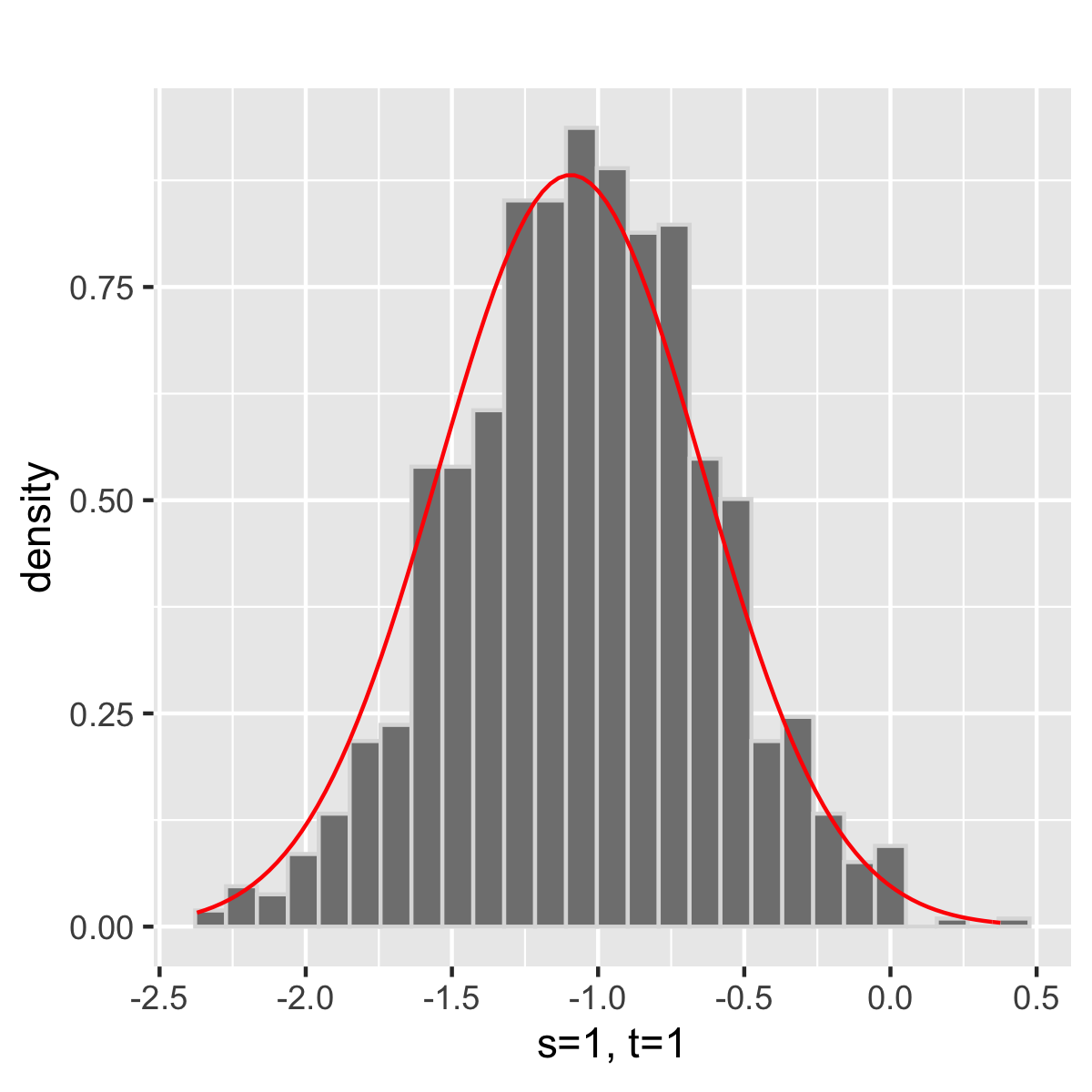}}
\subfigure
{\includegraphics[width=3.5cm]{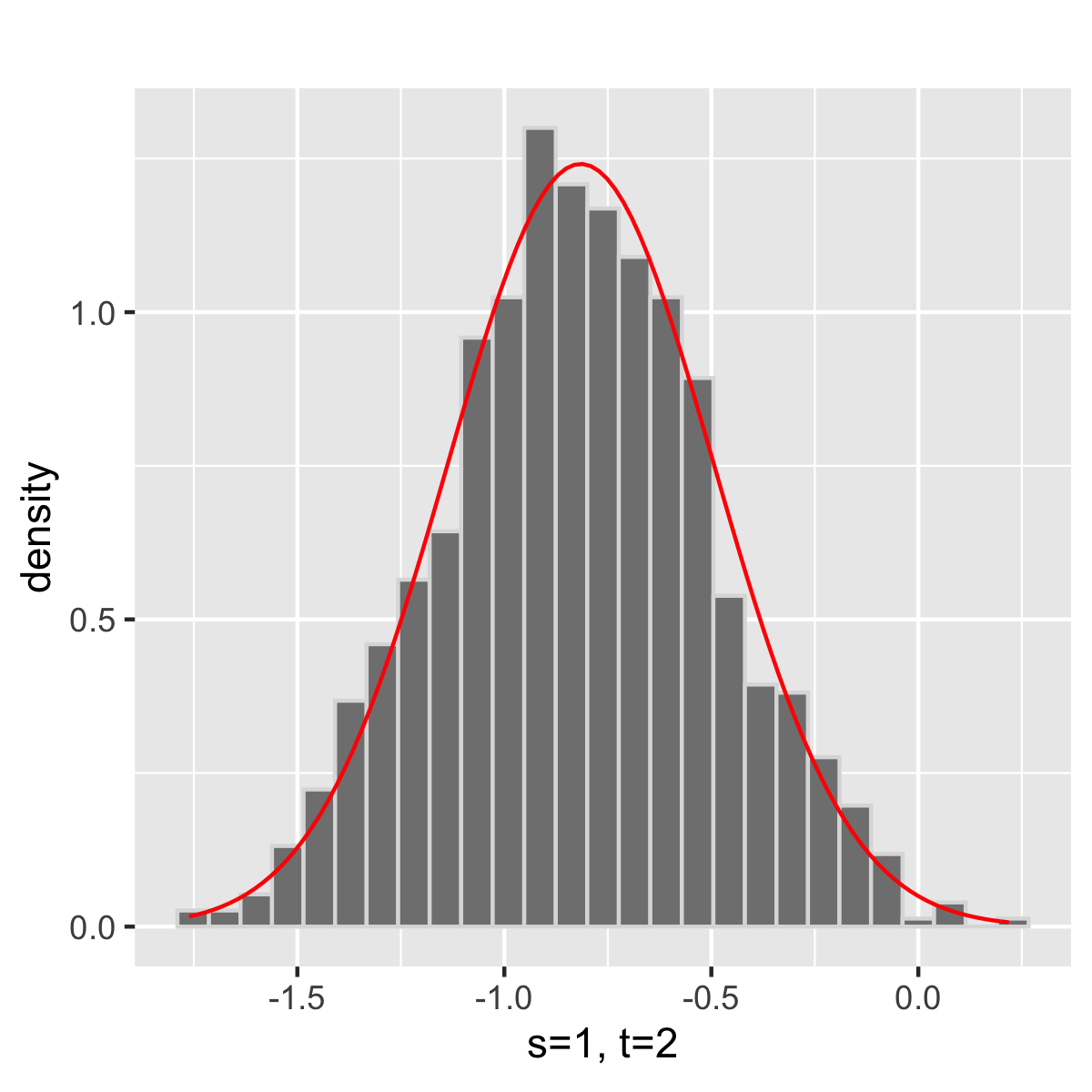}}
\subfigure
{\includegraphics[width=3.5cm]{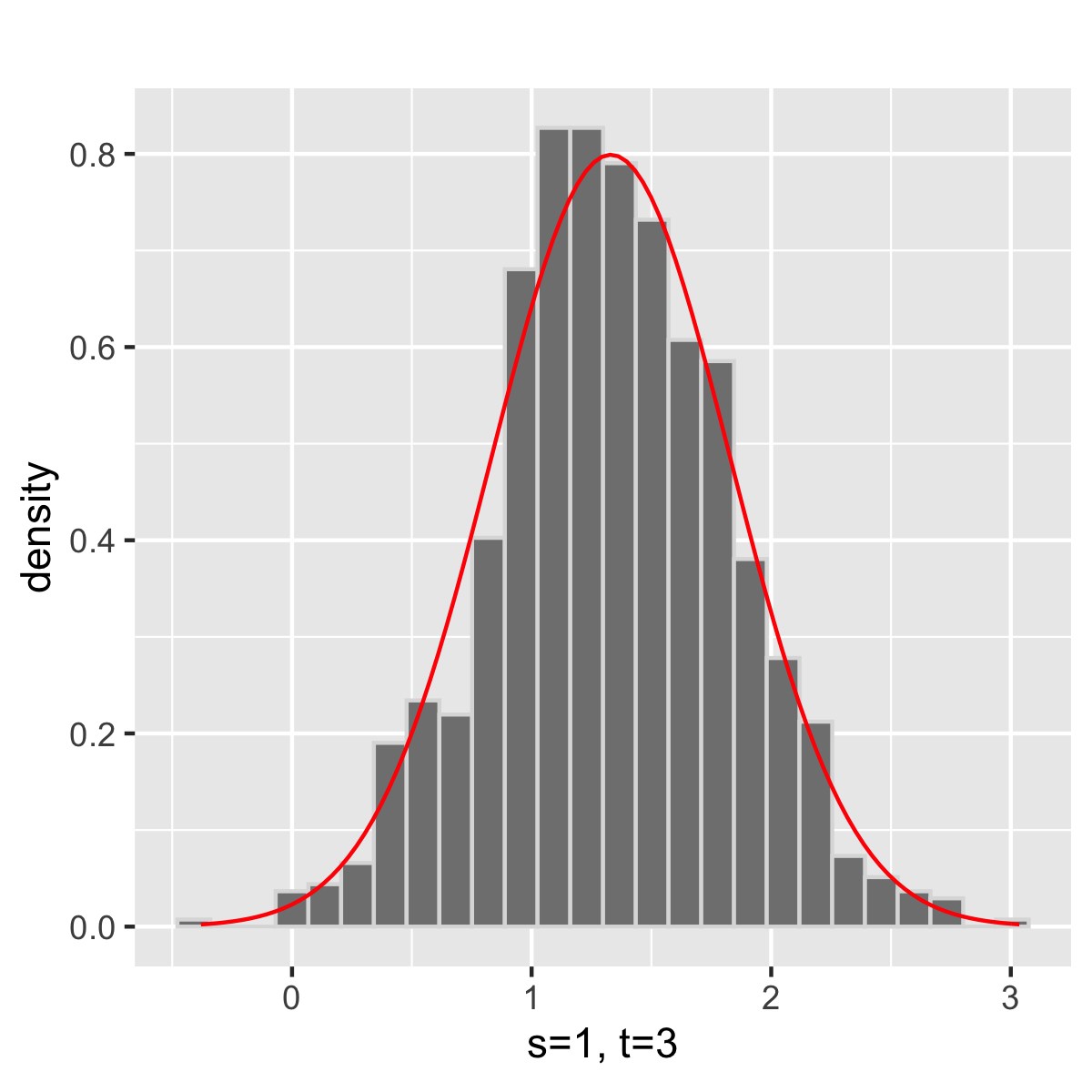}}

\subfigure
{\includegraphics[width=3.5cm]{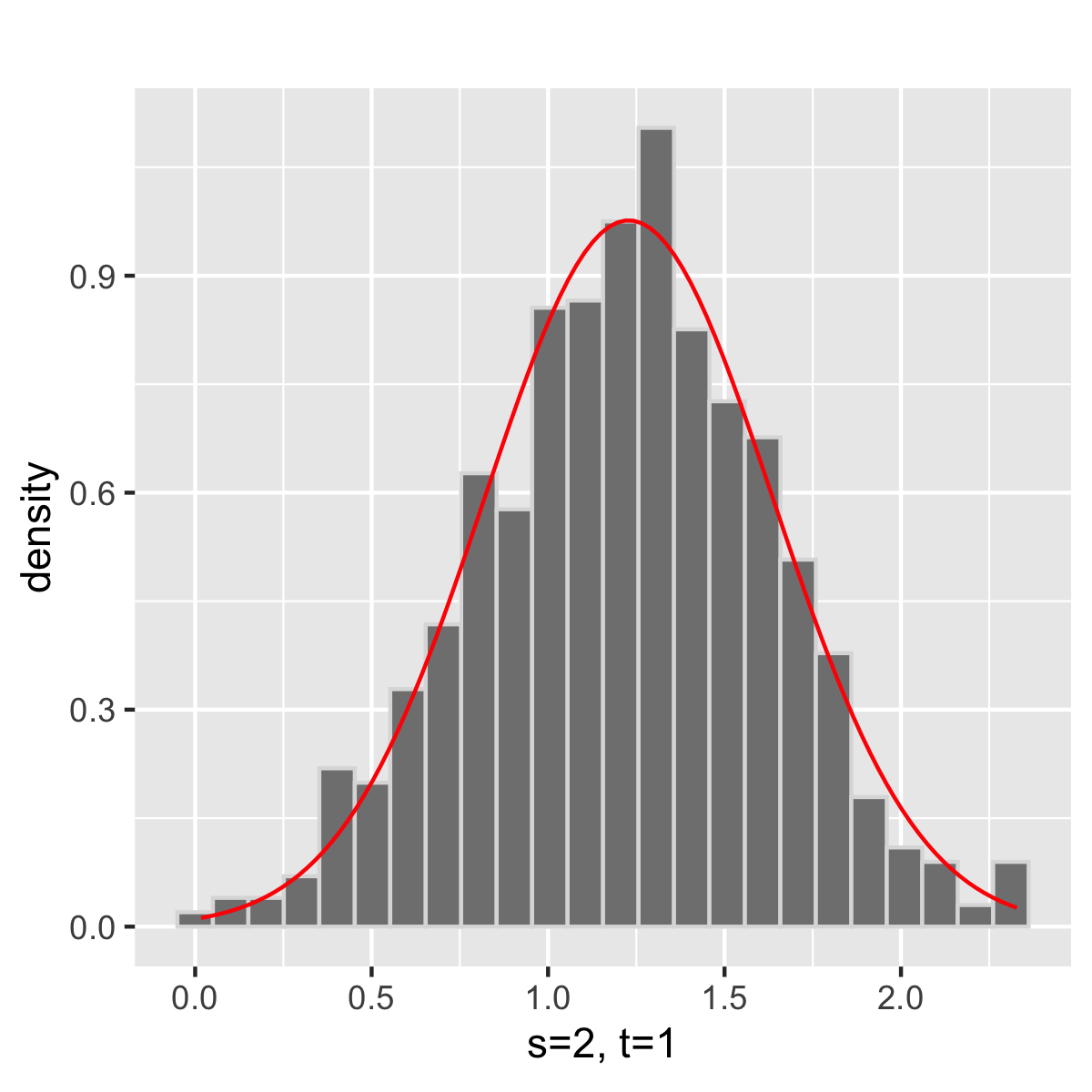}}
\subfigure
{\includegraphics[width=3.5cm]{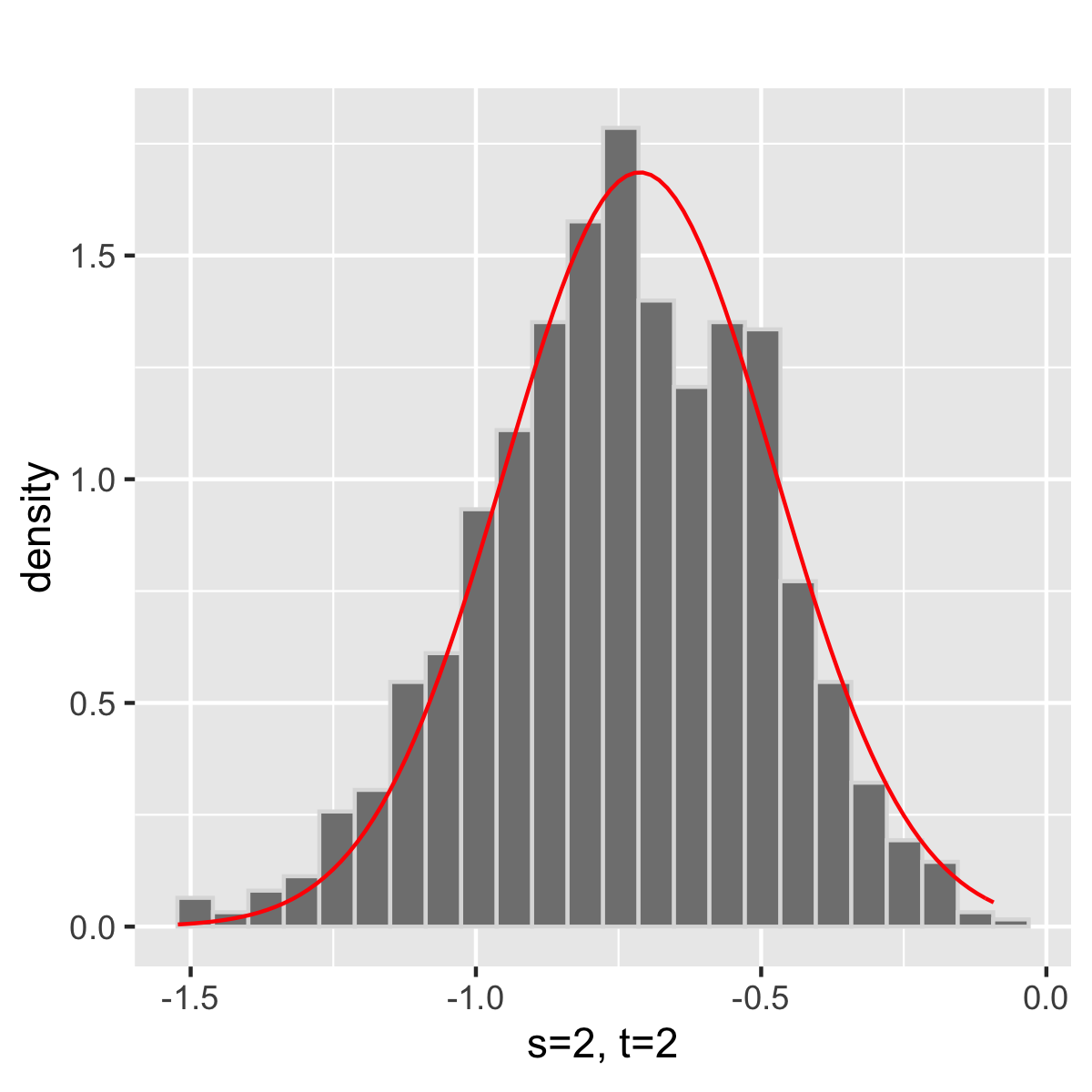}}
\subfigure
{\includegraphics[width=3.5cm]{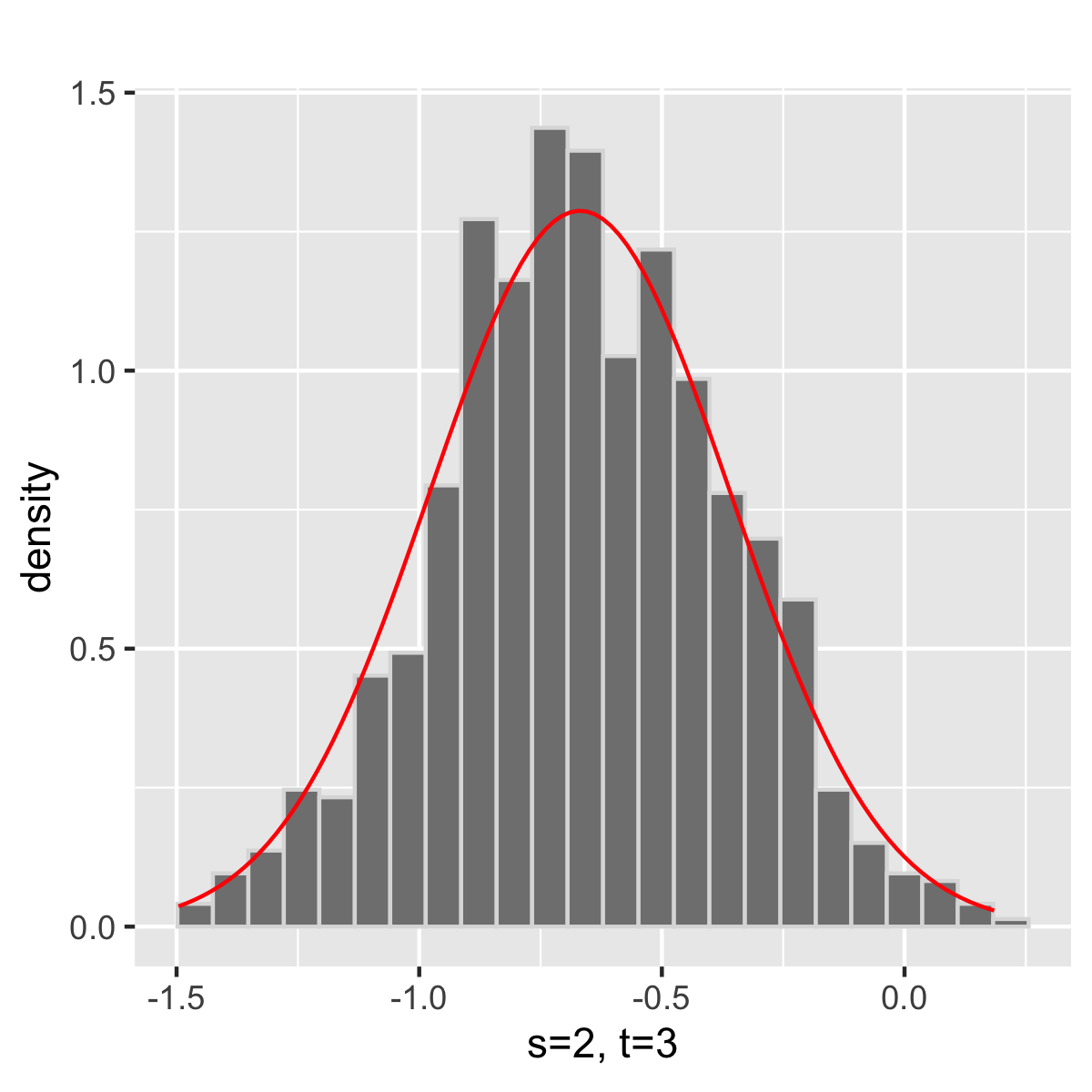}}

\subfigure
{\includegraphics[width=3.5cm]{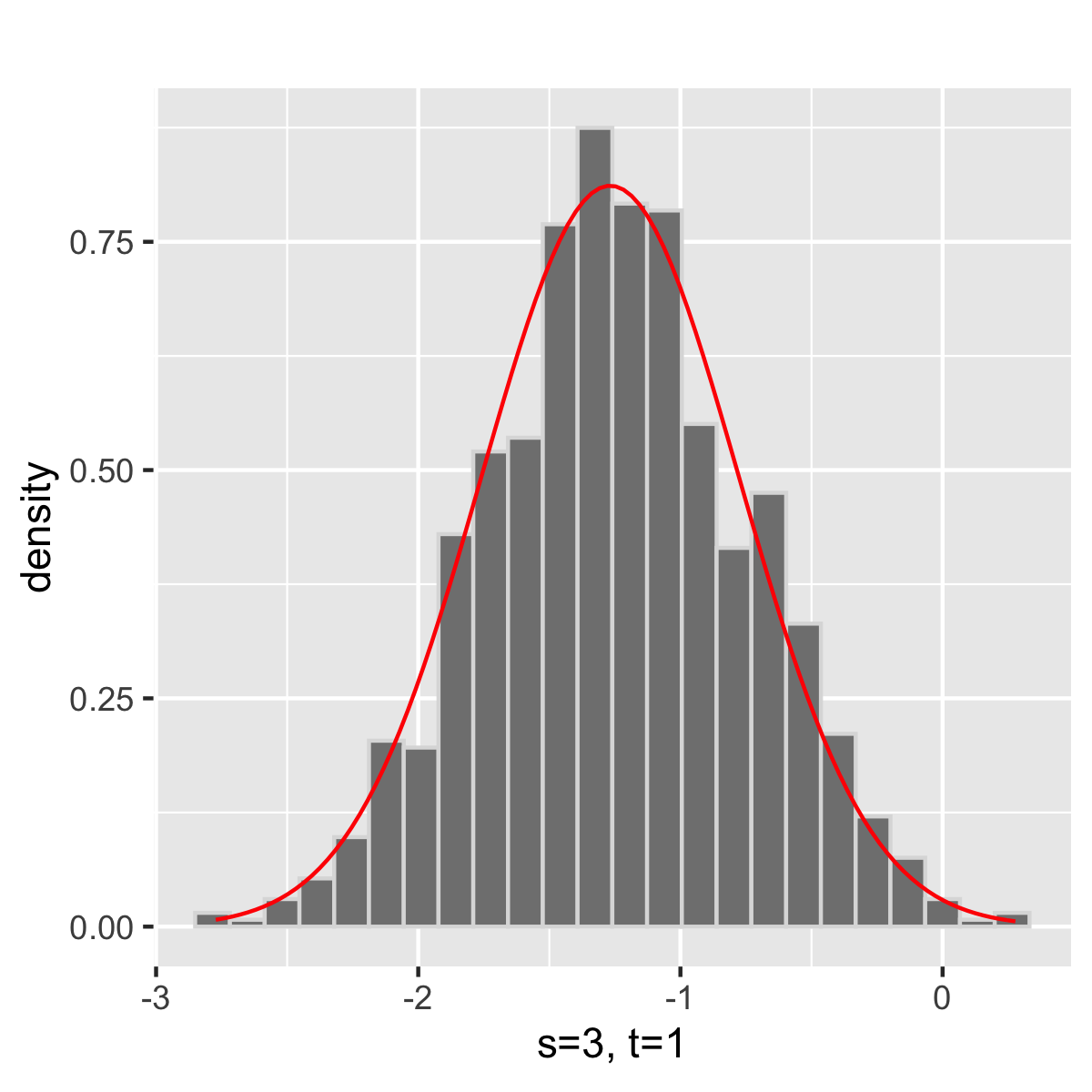}}
\subfigure
{\includegraphics[width=3.5cm]{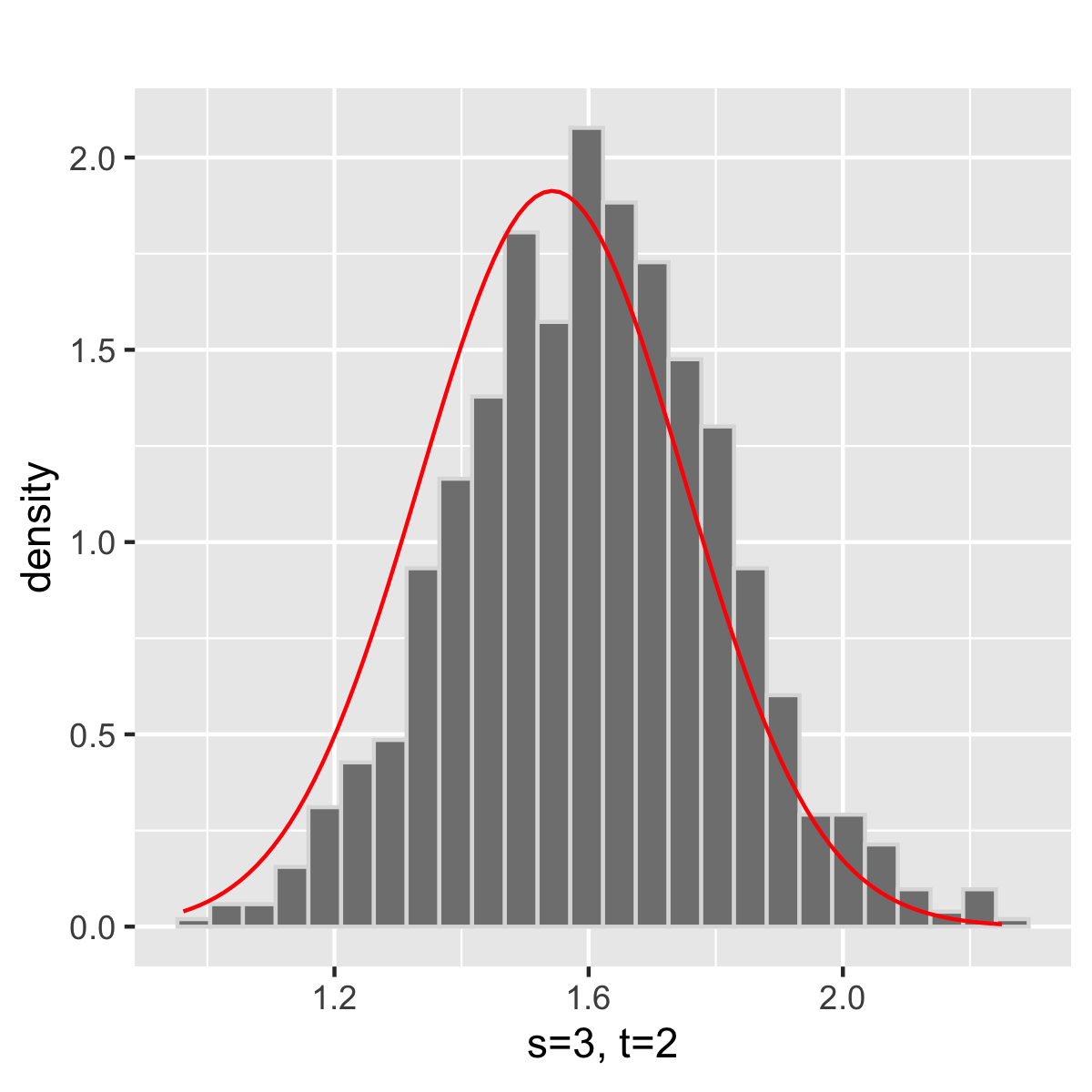}}
\subfigure
{\includegraphics[width=3.5cm]{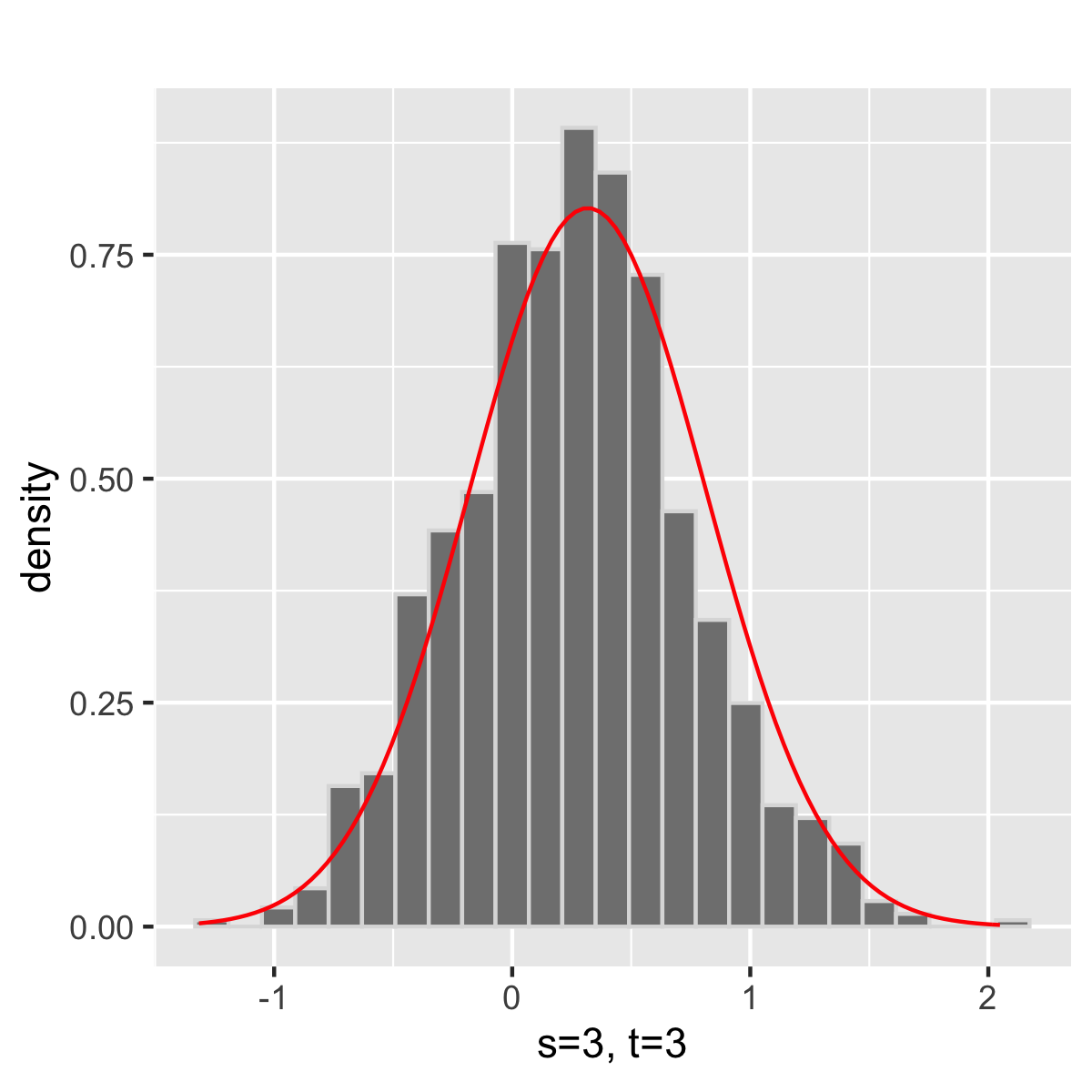}}

\caption{\footnotesize 
Histograms for the empirical distributions of the entries $(\mw_\muu^\top\hat\mr^{(1)}\mw_\mv-\mr^{(1)})_{st}$. The histograms are based on $1000$ samples of multilayer SBM graphs on $n=2000$ vertices with $m = 3$ layers and $K=3$ blocks. 
The red lines represent the probability density functions of the normal distributions with parameters specified in Theorem~\ref{thm:What_U Rhat What_v^T-R->norm}.
}
\label{fig:normality of R}
\end{figure}

We conduct $1000$ independent Monte Carlo replicates to obtain an empirical distribution of $\operatorname{vec}(\mw_\muu^\top\hat\mr^{(i)}\mw_\mv-\mr^{(i)})$, which we then compare against the limiting distribution given in Theorem~\ref{thm:What_U Rhat What_v^T-R->norm}. The results are summarized in Figure~\ref{fig:normality of R}. 
The Henze-Zirkler normality test indicates that the empirical distribution for $\operatorname{vec}(\mw_\muu^\top\hat\mr^{(i)}\mw_\mv-\mr^{(i)})$ is well-approximated by a multivariate normal distribution, and furthermore the empirical covariances for $\operatorname{vec}(\mw_\muu^\top\hat\mr^{(i)}\mw_\mv-\mr^{(i)})$ are very close to the theoretical covariances.

\begin{figure}[htbp!]
\centering
\subfigure[null hypothesis]
{\includegraphics[width=5cm]{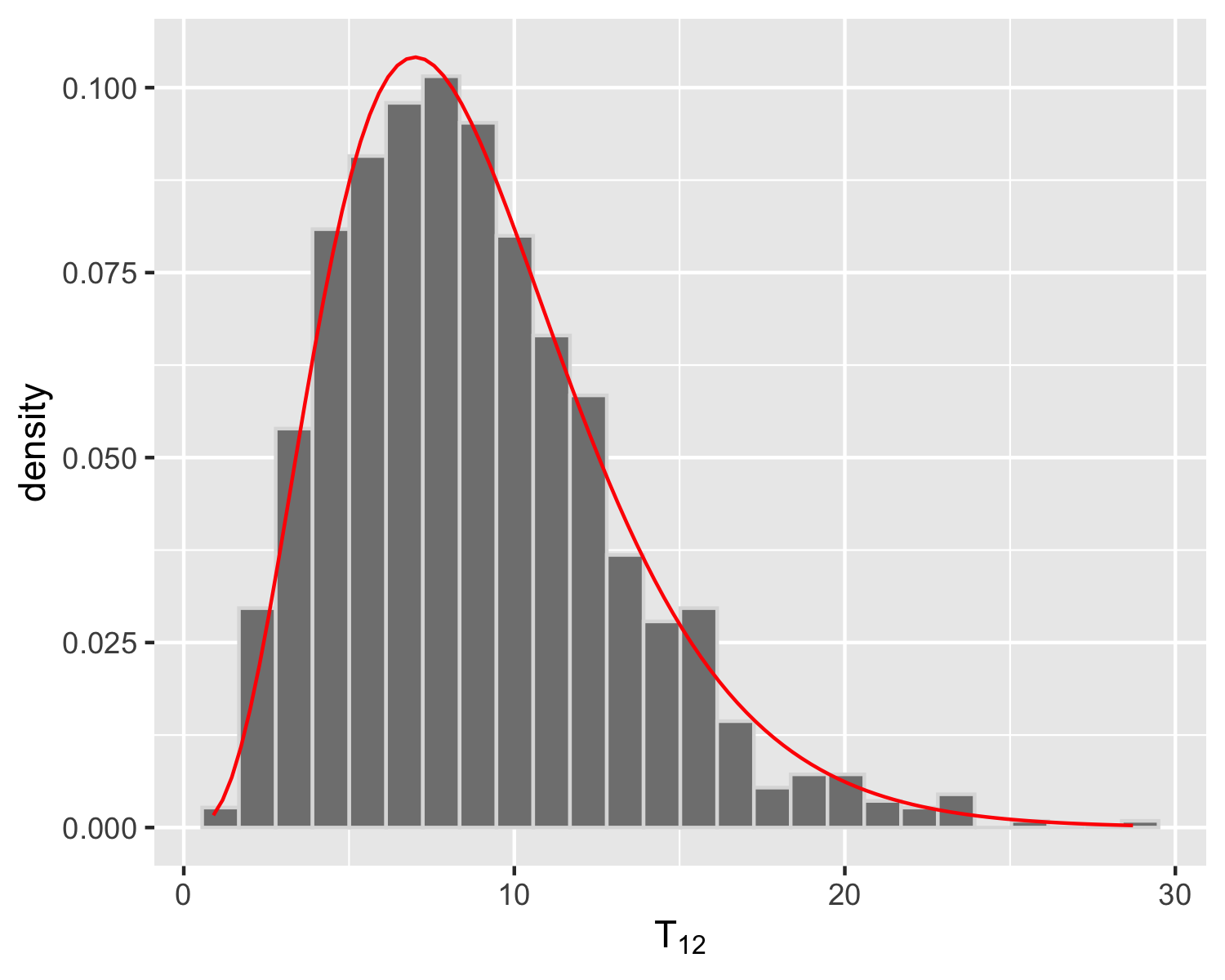}}
\quad
\subfigure[alternative hypothesis]
{\includegraphics[width=5cm]{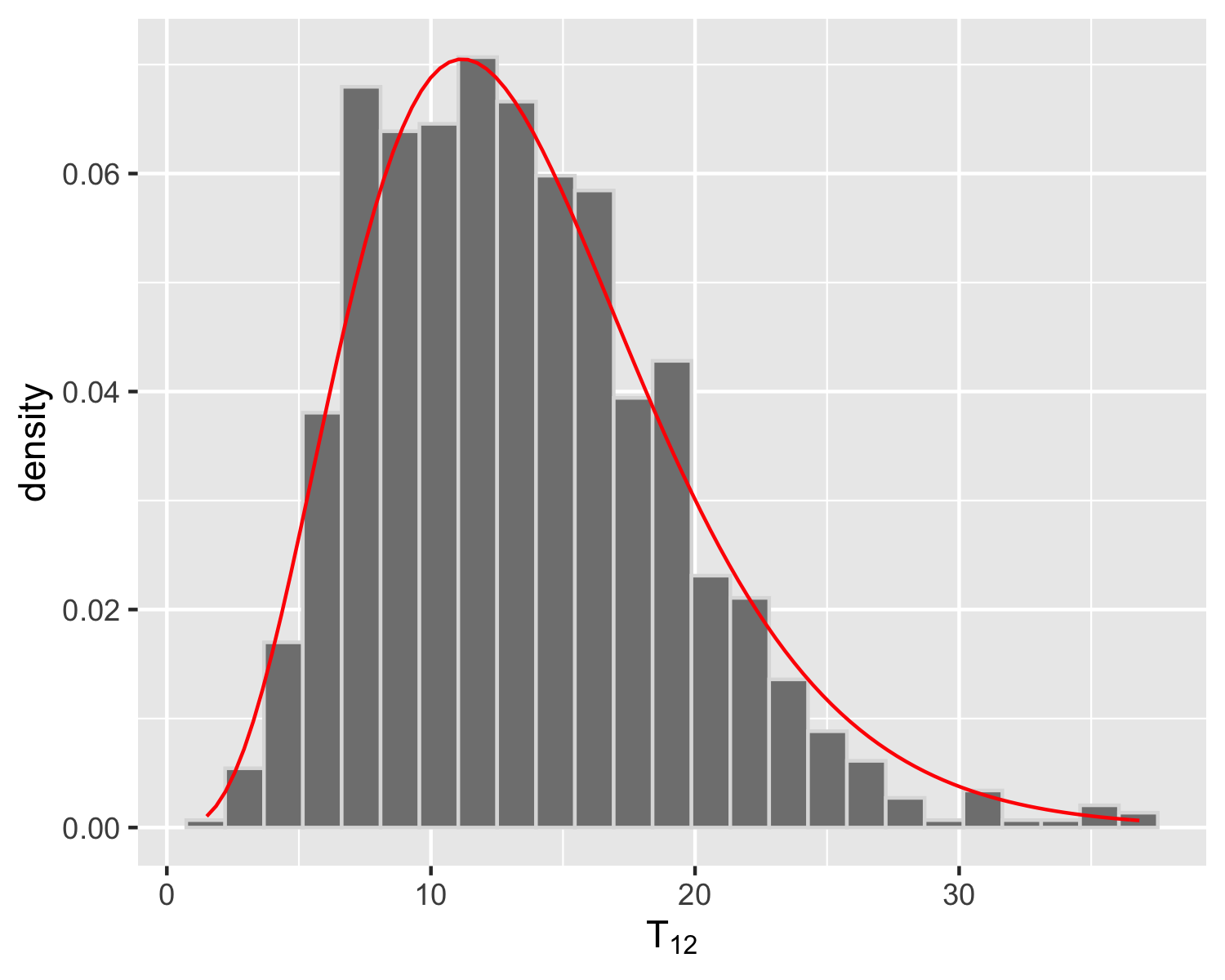}}

\caption{\footnotesize 
Histograms for the empirical distributions of $T_{12}$ under either the null or local alternative hypothesis. 
Refer to Figure~\ref{fig:normality of R} for more details.
The red lines represent the probability density functions for the central and non-central chi-square distributions with degrees of freedom and non-centrality parameters specified in Theorem~\ref{thm:HT}.
}
\label{fig:HT}
\end{figure}

We next consider the problem of determining whether or not two graphs $\mathbf{A}^{(i)}$ and $\mathbf{A}^{(j)}$ have the same distribution, i.e., we wish to test $\mathbb{H}_0 \colon \mr^{(i)} = \mr^{(j)}$ against $\mathbb{H}_A \colon \mr^{(i)} \neq \mr^{(j)}$.
We once again generate $1000$ Monte Carlo replicates where, for each replicate, we generate a \emph{directed} multilayer SBM with $m = 3$ graphs, $K = 3$ blocks using a similar setting to that described above, except now we set either $\mathbf{B}^{(2)} = \mathbf{B}^{(1)}$ or $\mb^{(2)}=\mb^{(1)}+\frac{1}{n}\mathbf{1}\mathbf{1}^\top$. These two choices for $\mathbf{B}^{(2)}$ correspond to the null and \emph{local} alternative, respectively. For each Monte Carlo replicate we compute the test statistic in Theorem~\ref{thm:HT}. 
We compare its empirical distributions under the null and alternative hypotheses against the central and non-central $\chi^2$ distributions with degrees of freedom $9=3^2$ and non-centrality parameters specified in Theorem~\ref{thm:HT} in Figure~\ref{fig:HT}.

\subsection{MultiNeSS model}
\label{sec:MultiNeSS}

We now evaluate the accuracy of Algorithm~\ref{Alg_COISIE} for
recovering the common and individual structures in a collection of
matrices generated from the
MultiNeSS model \citep{macdonald2022latent} with Gaussian errors. More specifically, for any $i\in[m]$, let $\mpp^{(i)}$
be a $n\times n$ matrix of the
form $$\mpp^{(i)}=\mx_c\mx_c^\top+\mx_s^{(i)}\mx_s^{(i)\top},$$ where
$\mx_c\in\mathbb{R}^{n\times d_1},\mx_s^{(i)}\in\mathbb{R}^{n\times
d_2}$.  Let $\mathbf{F}:=\mx_c\mx_c^\top$ be the common structure
across all $\{\mpp^{(i)}\}$, and let
$\mathbf{G}^{(i)}:=\mx_s^{(i)}\mx_s^{(i)\top}$ be the individual
structure for $\mpp^{(i)}$. We then generate $\ma^{(i)} = \mpp^{(i)} +
\me^{(i)}$ where $\me^{(i)}$ is a symmetric
random matrix whose upper triangular entries are iid $N(0, \sigma^2)$ random
variables. See Remark~\ref{rk:generalized model} for further discussion of the MultiNeSS model and
its relevance to the current paper. 

Given $\{\ma^{(i)}\}_{i=1}^{m}$ we first compute
$\hat\muu^{(i)}$ as the $n\times
(d_1+d_2)$ matrix whose columns are the $d_1+d_2$ leading eigenvectors
of $\ma^{(i)}$ for each $i \in [m]$. Next we let $\hat\muu_c$ be the
$n\times d_1$ matrix whose columns are the $d_1$ leading left singular vectors of
$[\hat\muu^{(1)}|\cdots |\hat\muu^{(m)}]$ and 
$\hat\muu_s^{(i)}$ be the $n\times d_2$ matrix containing the $d_2$
leading left singular vectors of
$(\mi-\hat\muu_c\hat\muu_c^\top)\hat\muu^{(i)}$ for all $i \in
[m]$. Finally we compute the
estimates of $\mathbf{F},\mathbf{G}^{(i)},\mathbf{P}^{(i)}$ via
$$
\begin{aligned}
	&\hat{\mathbf{F}}=\hat\muu_c\hat\muu_c^\top\bar{\ma}\hat\muu_c\hat\muu_c^\top, \quad
	\hat{\mathbf{G}}^{(i)}=\hat\muu_s^{(i)}\hat\muu_s^{(i)\top}\ma^{(i)}\hat\muu_s^{(i)}\hat\muu_s^{(i)\top}, \quad
	\hat\mpp^{(i)}=\hat\muu_{c,s}^{(i)}\hat\muu_{c,s}^{(i)\top}\ma^{(i)}\hat\muu_{c,s}^{(i)}\hat\muu_{c,s}^{(i)\top},
\end{aligned}
$$
where $\bar{\ma}=m^{-1}\sum_{i=1}^m \ma^{(i)}$ and $\hat\muu_{c,s}^{(i)}=[\hat\muu_c|\hat\muu_s^{(i)}]$.

\begin{figure}[htbp!]
\centering
\subfigure
{\includegraphics[width=4.2cm]{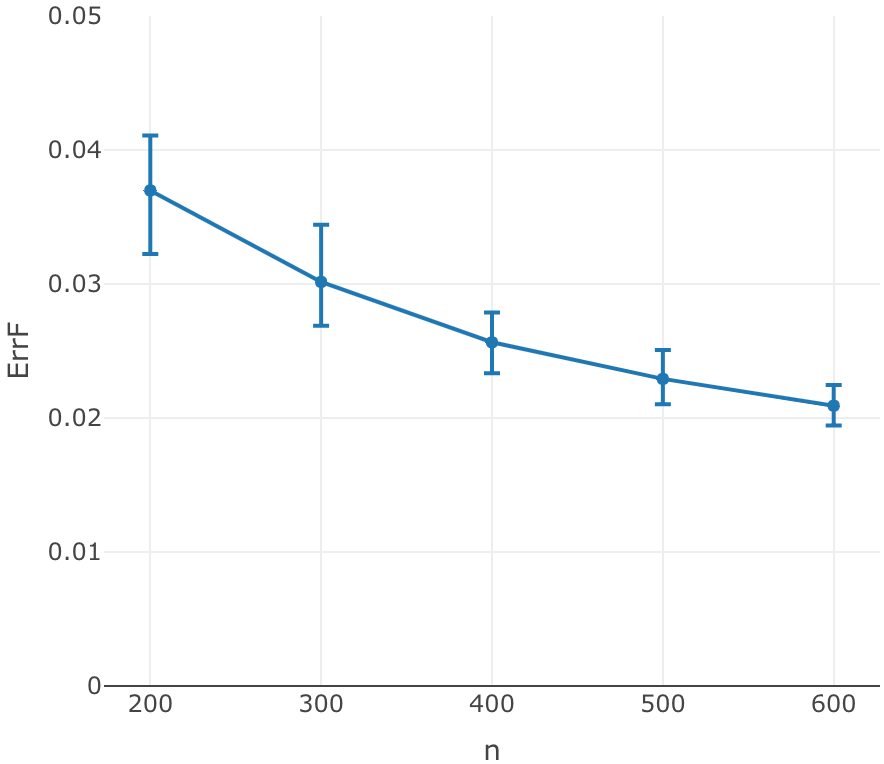}}
\subfigure
{\includegraphics[width=4.2cm]{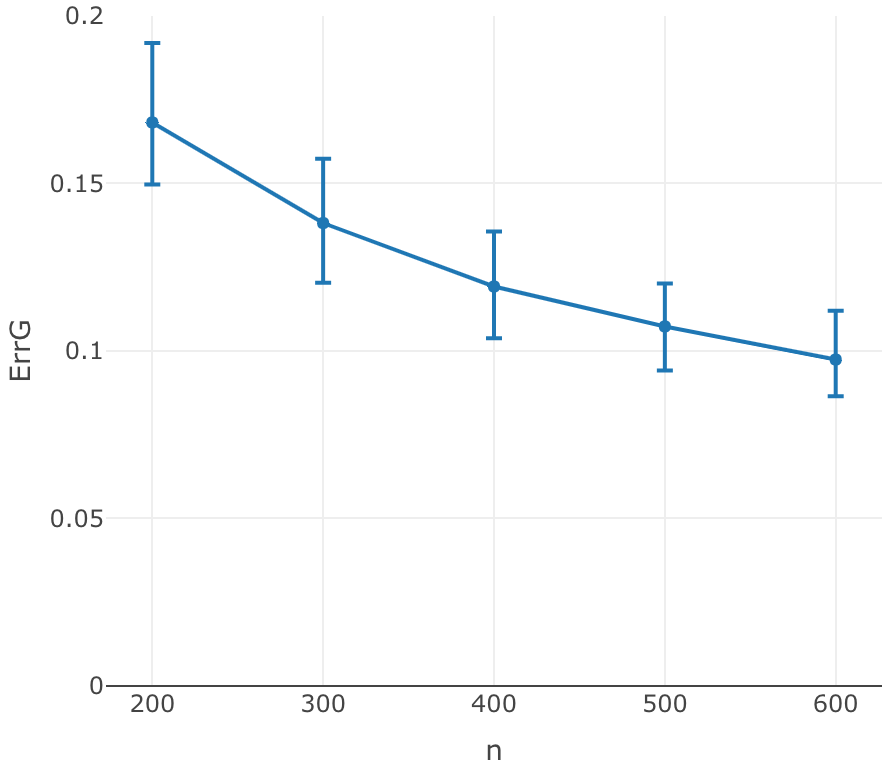}}
\subfigure
{\includegraphics[width=4.2cm]{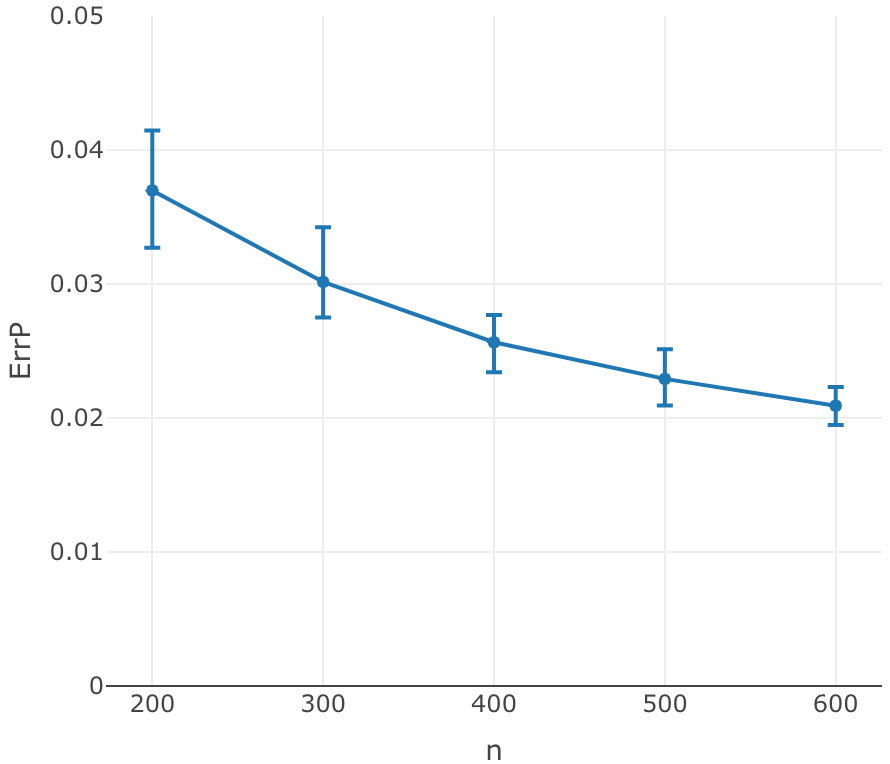}}

\caption{\footnotesize 
Relative Frobenius norm errors for the common structure (left panel), the individual structure (middle panel), and the overall expectation of the matrix (right panel) with $d_1=d_2=2,\sigma=1,m=8$ and $n\in\{200,300,400,500,600\}$. The figures display the mean, $0.05$ and $0.95$ quantile points, over $100$ independent Monte Carlo replications.
}
\label{fig:MN_n}
\end{figure}

\begin{figure}[htbp!]
\centering
\subfigure
{\includegraphics[width=4.2cm]{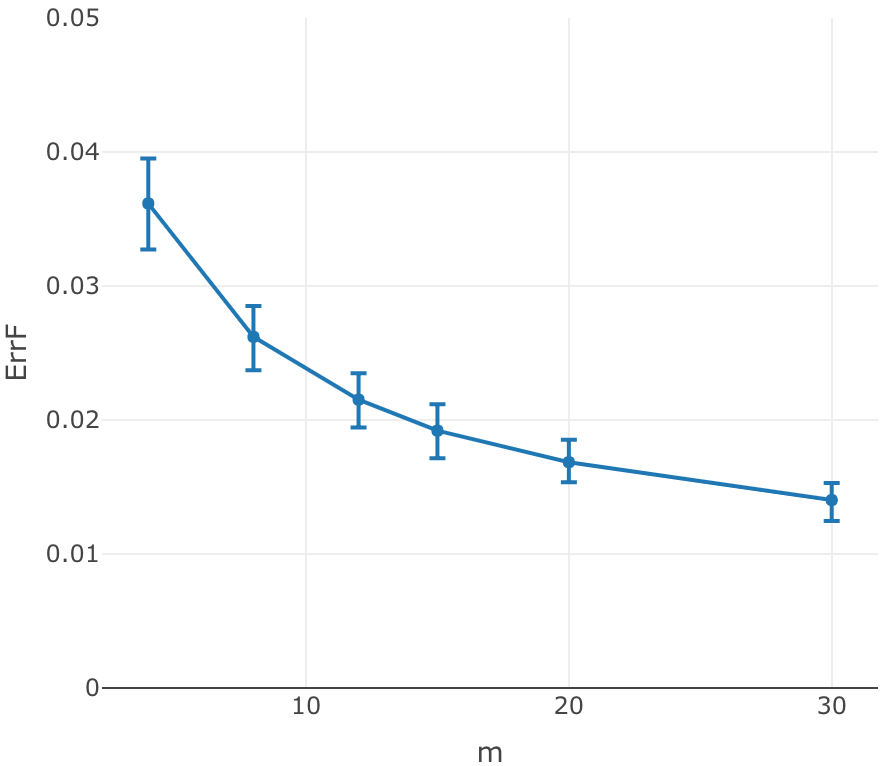}}
\subfigure
{\includegraphics[width=4.2cm]{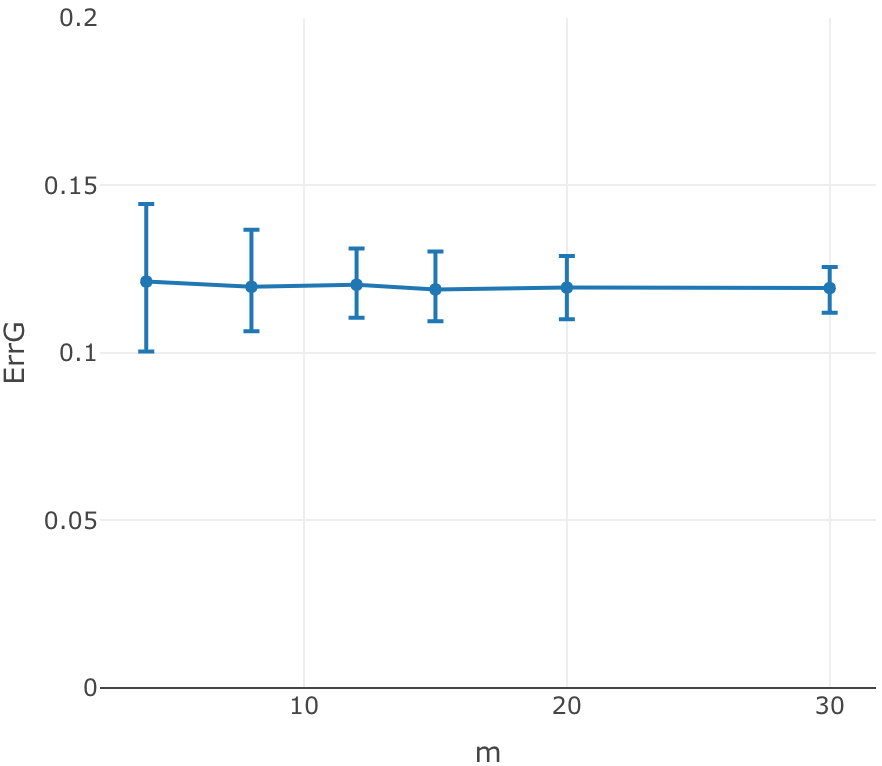}}
\subfigure
{\includegraphics[width=4.2cm]{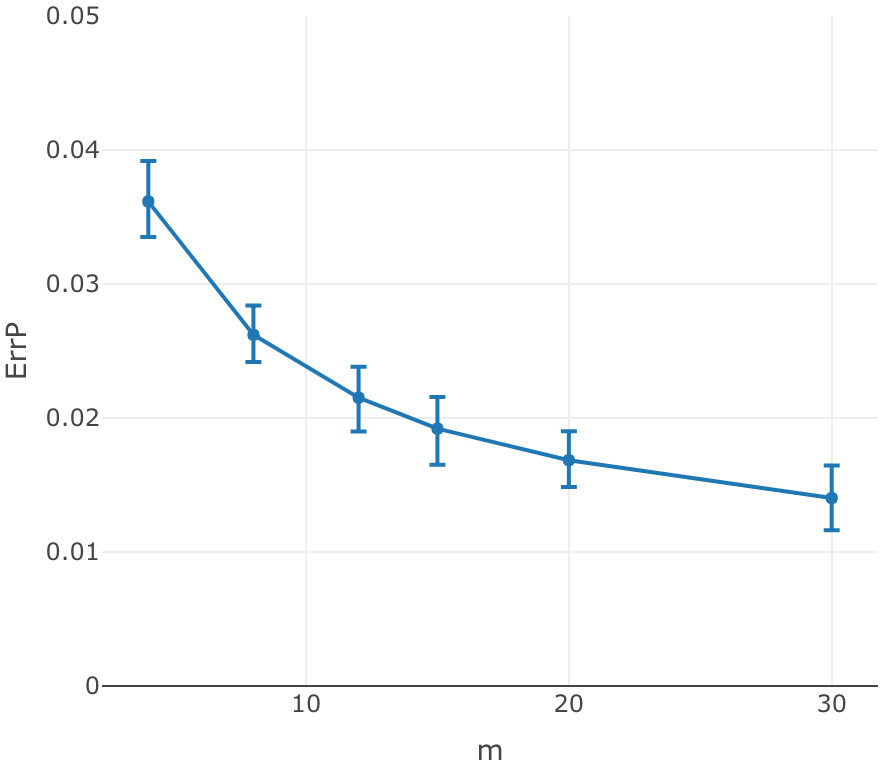}}

\caption{\footnotesize 
Relative Frobenius norm errors for the common structure (left panel), the individual structure (middle panel), and the overall expectation of the matrix (right panel) with $d_1=d_2=2,\sigma=1,n=400$ and $m\in\{4,8,12,15,20,30\}$. The figures display the mean, $0.05$ and $0.95$ quantile points, over $100$ independent Monte Carlo replications.
}
\label{fig:MN_m}
\end{figure}

We use the same setting as that in Section~5.2 in
\cite{macdonald2022latent}. More specifically we fix
$d_1=d_2=2,\sigma=1$, and either fix $m = 8$ and vary $n\in\{200,300,400,500,600\}$
or fix $n = 400$ and vary $m\in\{4,8,12,15,20,30\}$. The estimation error for $\hat{\mathbf{F}},
\{\hat{\mathbf{G}}^{(i)}\}$ and $\{\hat{\mathbf{P}}^{(i)}\}$ are also
evaluated using the same metric as that in \cite{macdonald2022latent}, i.e., we compute
$$
\text{ErrF}=\frac{\|\hat{\mathbf{F}}-\mathbf{F}\|_{\tilde{F}}}{\|\mathbf{F}\|_{\tilde{F}}},\quad
\text{ErrG}=\frac{1}{m}\sum_{i=1}^m\frac{\|\hat{\mathbf{G}}^{(i)}-\mathbf{G}^{(i)}\|_{\tilde{F}}}{\|\mathbf{G}^{(i)}\|_{\tilde{F}}},\quad
\text{ErrP}=\frac{1}{m}\sum_{i=1}^m\frac{\|\hat{\mathbf{P}}^{(i)}-\mathbf{P}^{(i)}\|_{\tilde{F}}}{\|\mathbf{P}^{(i)}\|_{\tilde{F}}},
$$
where $\|\cdot\|_{\tilde{F}}$ denote the Frobenius norm of a matrix
after setting its diagonal entries to $0$.
The results are summarized in Figure~\ref{fig:MN_n} and
Figure~\ref{fig:MN_m}. Comparing
the relative Frobenius norm errors in Figure~\ref{fig:MN_n} and
Figure~\ref{fig:MN_m} with those in Figure~2 of
\cite{macdonald2022latent}, we see that the two set of estimators have
comparable performance. Nevertheless, our algorithm is slightly better
for recovering the
common structure (smaller ErrF) while
the algorithm in \cite{macdonald2022latent} is slightly better for
recovering individual structure  (smaller ErrG).  
Finally for recovering the overall edge probabilities $\{\mpp^{(i)}\}$,
our ErrPs are always smaller than theirs. Indeed, as $n$ varies from
$200$ to $600$, the mean of our ErrP varies from about $0.076$ to
$0.044$ while the mean in \cite{macdonald2022latent} varies from about
$0.08$ to $0.05$.
Simialrly, as $m$ varies from $4$ to $30$, the mean of our ErrP varies from about
$0.056$ to $0.051$ while the mean in \cite{macdonald2022latent} varies from about $0.07$ to
$0.06$. In summary, while the two algorithms yield estimates with
comparable performance, our algorithm has some computational advantage
as (1) it is not an interative procedure and (2) it does not require any tuning parameters
(note that the embedding dimensions $d_1$ and $d_2$ are generally not tuning parameters but rather chosen via
some dimension selection procedure). 

\subsection{Comparison of estimation methods}\label{sec:simu_comp}
\begin{figure}[htbp!]
\centering
\subfigure
{\includegraphics[width=6cm]{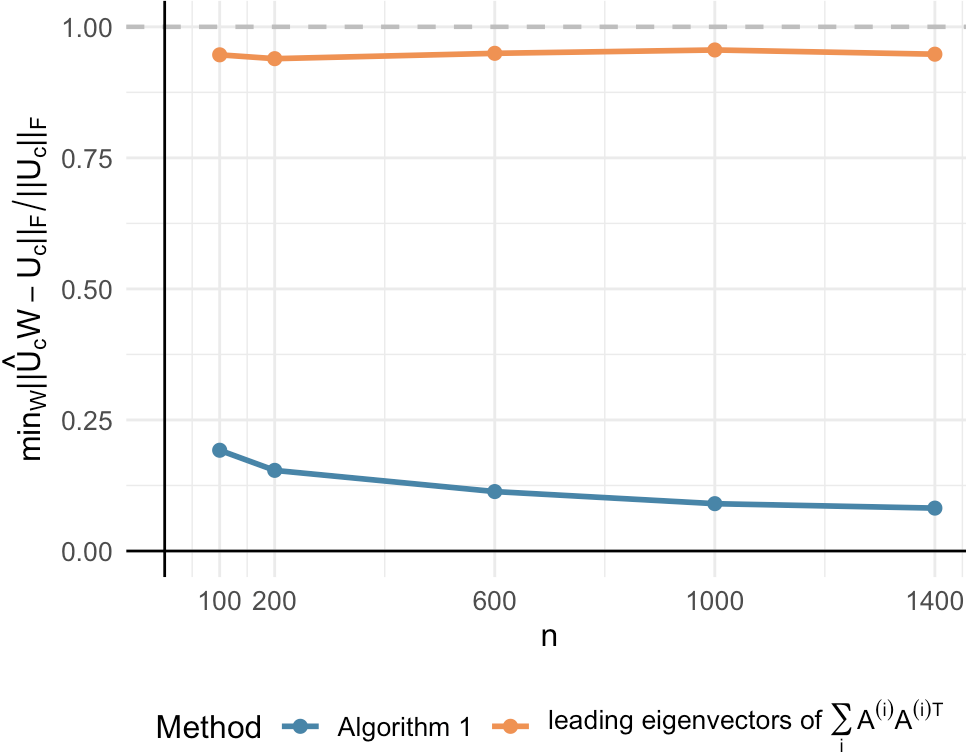}}

\caption{\footnotesize 
Empirical relative Frobenius norm $\min_\mw\|\hat\muu_c\mw-\muu_c\|_F/\|\muu_c\|_F$ for Algorithm~\ref{Alg_COISIE} and the "aggregate-then-estimate" approach that uses the leading eigenvectors of $\sum_{i}\ma^{(i)}\ma^{(i)\top}$ as $\hat\muu_c$ for the COISIE model when varying $n \in \{100,200,600,1000,1400\}$ while fixing $m = 3$, $d_i \equiv 4$ and $d_{0,\muu} = d_{0,\mv} = 2$. Additional details of the settings are provided in Section~\ref{sec:simu_comp}. The lines represent the means of $100$ independent Monte Carlo replicates.
}
\label{fig:compare}
\end{figure}

In Section~\ref{sec:related_works}, we mention that although
"aggregate-then-estimate" approaches allow for milder conditions if $m$ goes to infinity
and there are only common subspaces with no individual subspaces,
they can fail to be consistent when individual subspaces are present. We now provide
some simulation results to support this claim. Consider the setting in Section~\ref{sec:simu_COISIE} and suppose
that the $\mpp^{(i)}$ are also randomly generated in each Monte Carlo replicate. We compare $\hat\muu_c$ obtained by Algorithm~\ref{Alg_COISIE} with the "aggregate-then-estimate" approach that uses the leading eigenvectors of $\sum_{i}\ma^{(i)}\ma^{(i)\top}$ as $\hat\muu_c$. We measure estimation accuracy using the relative Frobenius norm $\min_\mw\|\hat\muu_c\mw-\muu_c\|_F/\|\muu_c\|_F$. As shown in Figure~\ref{fig:compare}, this "aggregate-then-estimate" approach fails to provide accurate subspace estimation, while Algorithm~\ref{Alg_COISIE} is effective.

\subsection{Distributed PCA} 
\begin{figure}[htbp!]
\centering
\subfigure
{\includegraphics[height=4.2cm]{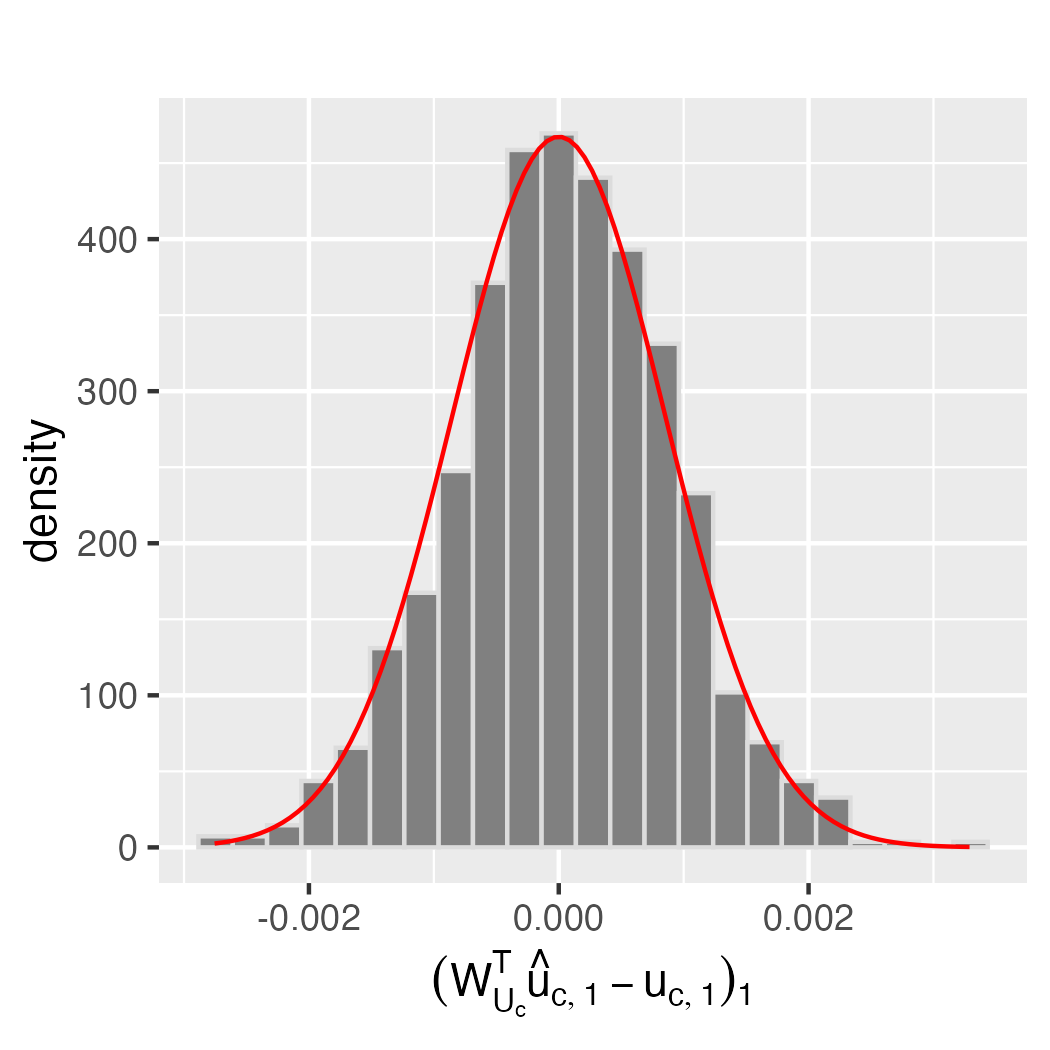}}
\subfigure
{\includegraphics[height=4.2cm]{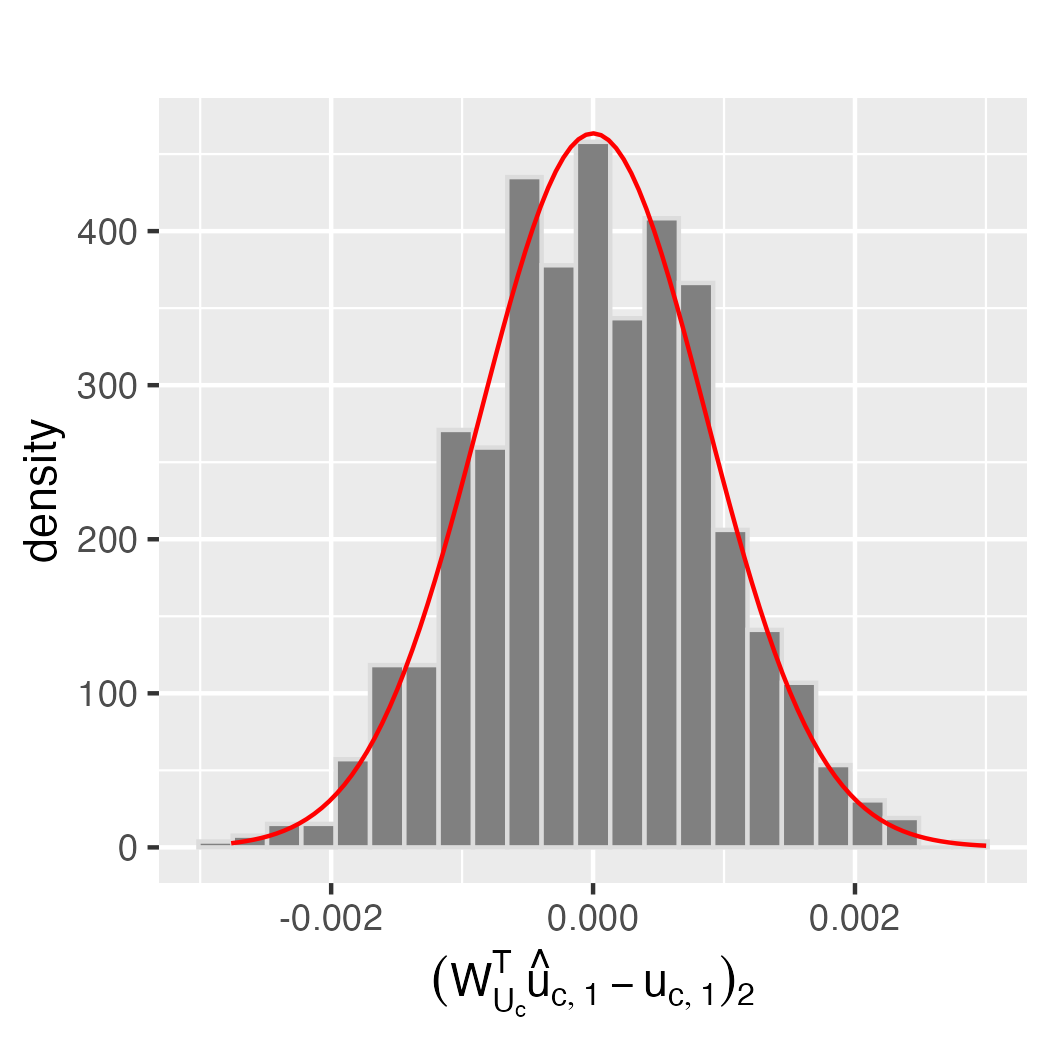}}
\subfigure
{\includegraphics[height=4.2cm]{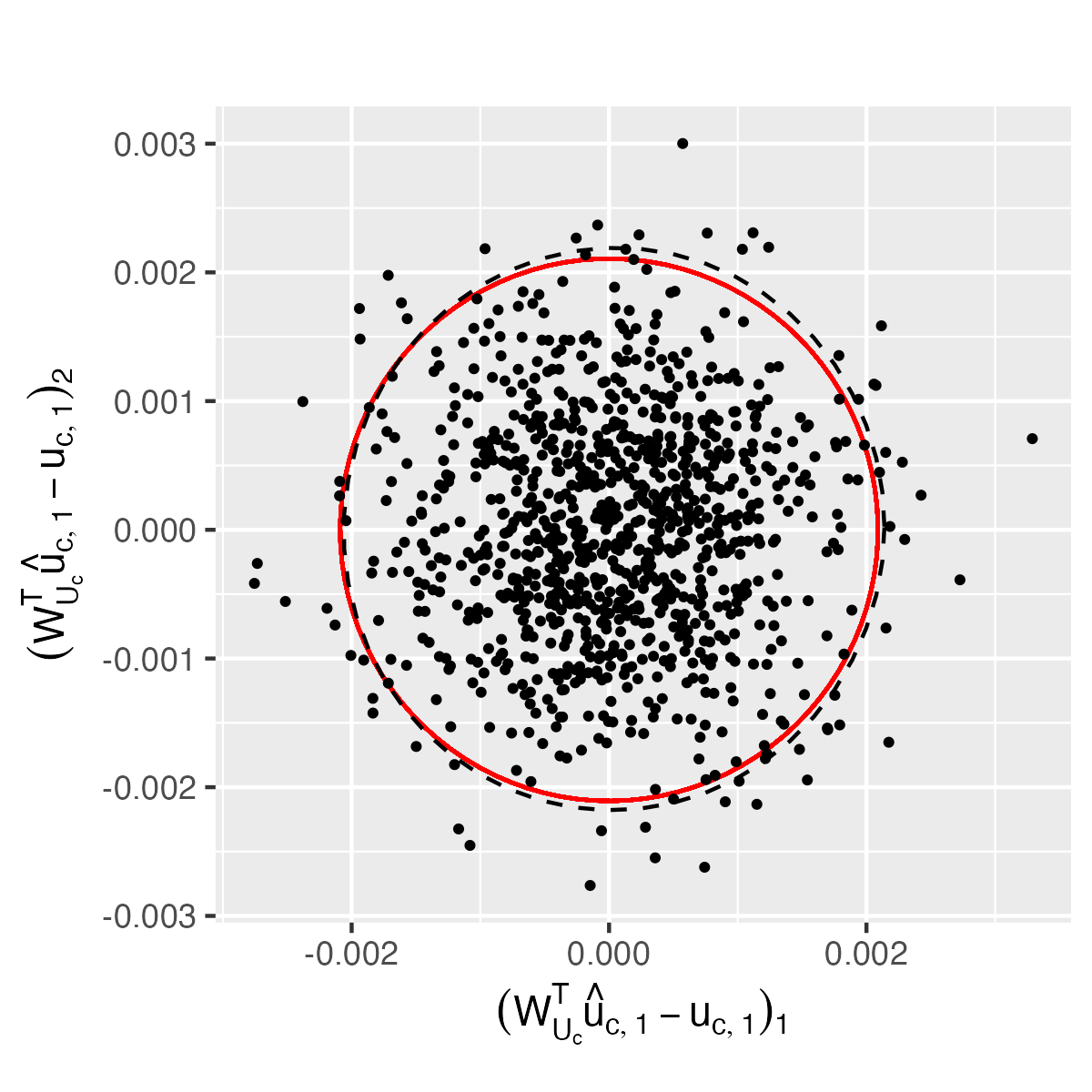}}

\caption{\footnotesize 
The left two panels are histograms of the empirical distributions of the entries of the estimation error $\mw^\top_{\muu_c}\hat{u}_{c,k} - {u}_{c,k}$ for $k = 1$. These histograms are based on $1000$ independent Monte Carlo replicates of the distributed PCA with $n_i \equiv n= 4000$, $m = 10$, $D=1000$, $d_i \equiv 4$, $d_{0} = 2$, $\sigma_i\equiv 1$, $\max\lambda^{(i)}_\ell=50$ and $\min\lambda^{(i)}_\ell=20$. The red lines represent the probability density functions of the normal distributions with parameters specified in Theorem~\ref{thm:PCA_(UhatW-U)_k->normal_previous}. 
The right panel displays a bivariate plot of the empirical distributions of the entries. The dashed black ellipses represent 95\% level curves for the empirical distributions, while the solid red ellipses represent 95\% level curves for the theoretical distributions as specified in Theorem~\ref{thm:PCA_(UhatW-U)_k->normal_previous}.
}
\label{fig:simulation_CLT1_pca}
\end{figure}

\begin{figure}[htbp!]
\centering
\subfigure
{\includegraphics[height=4.2cm]{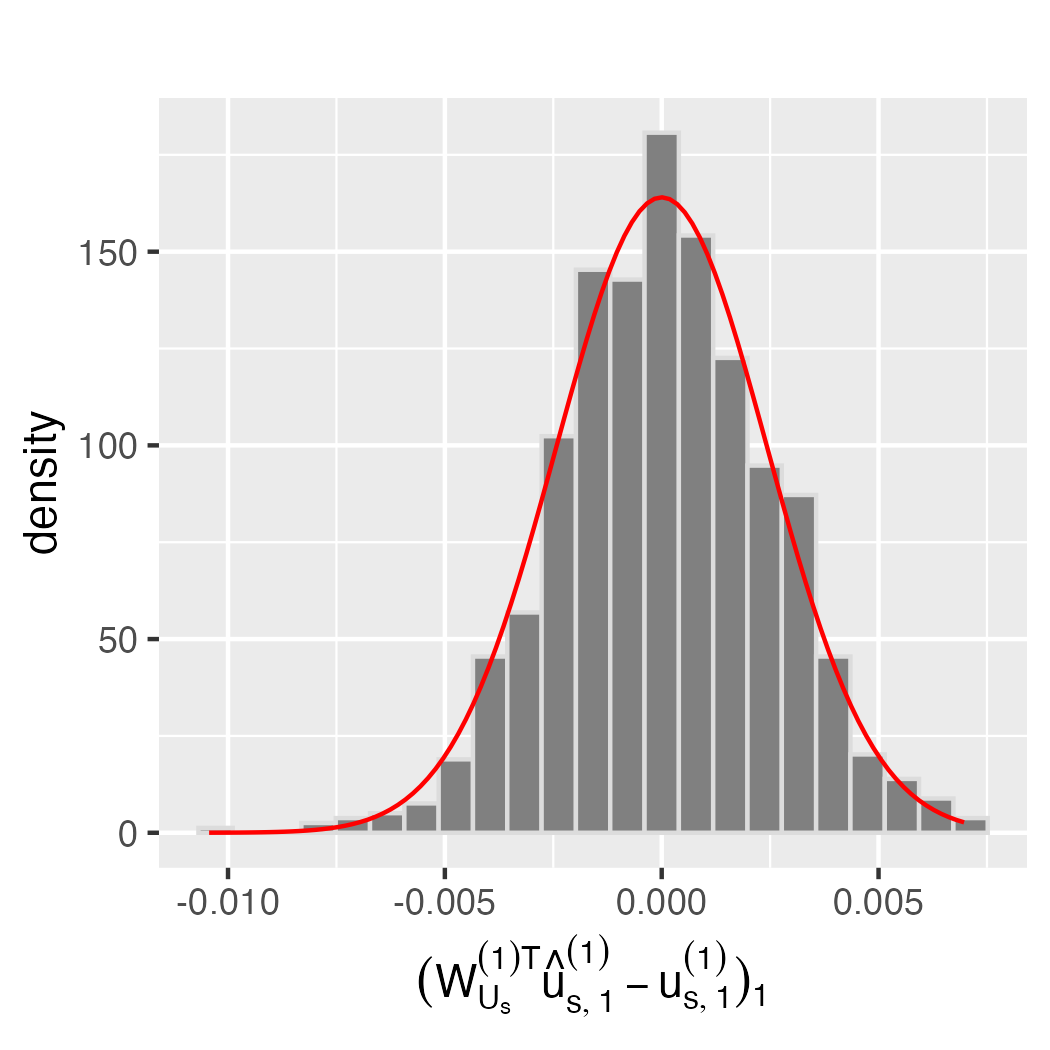}}
\subfigure
{\includegraphics[height=4.2cm]{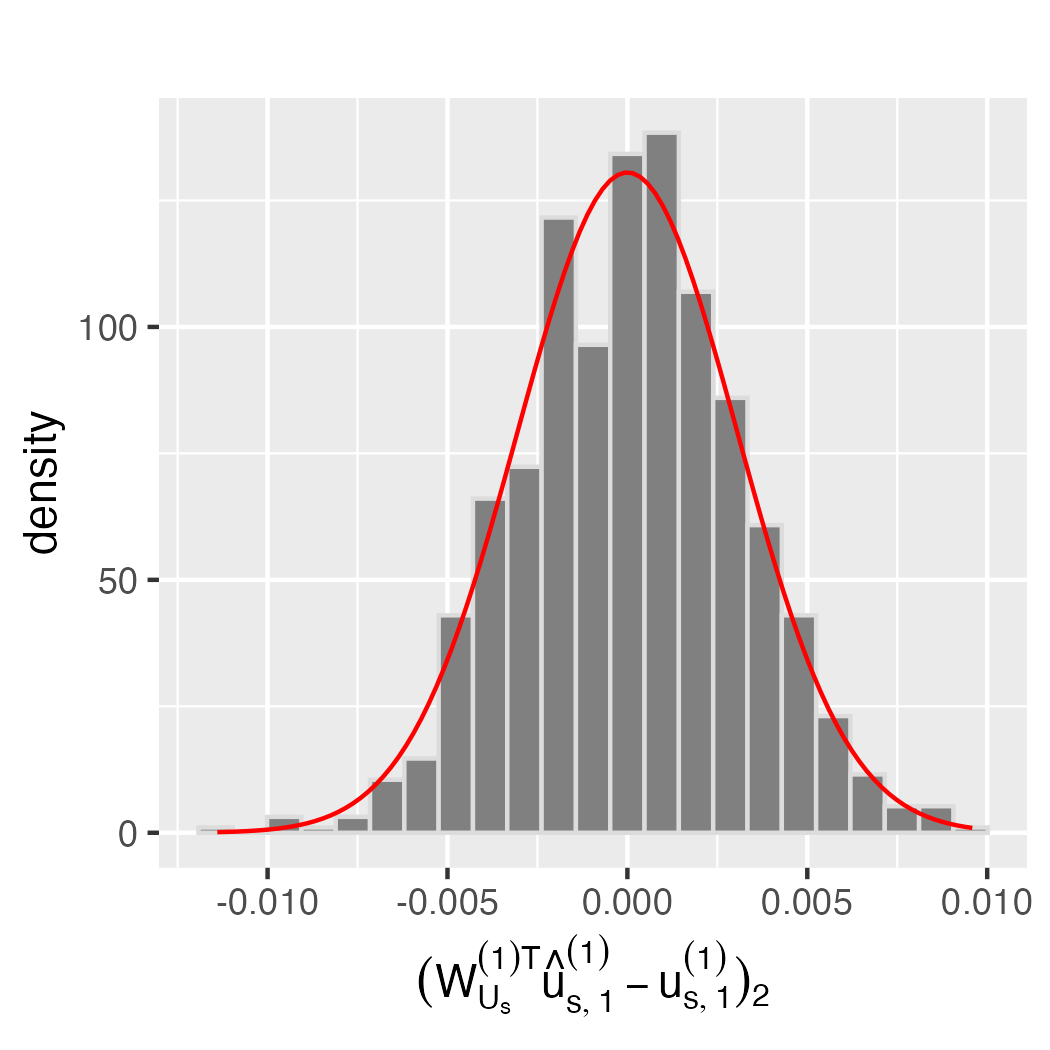}}
\subfigure
{\includegraphics[height=4.2cm]{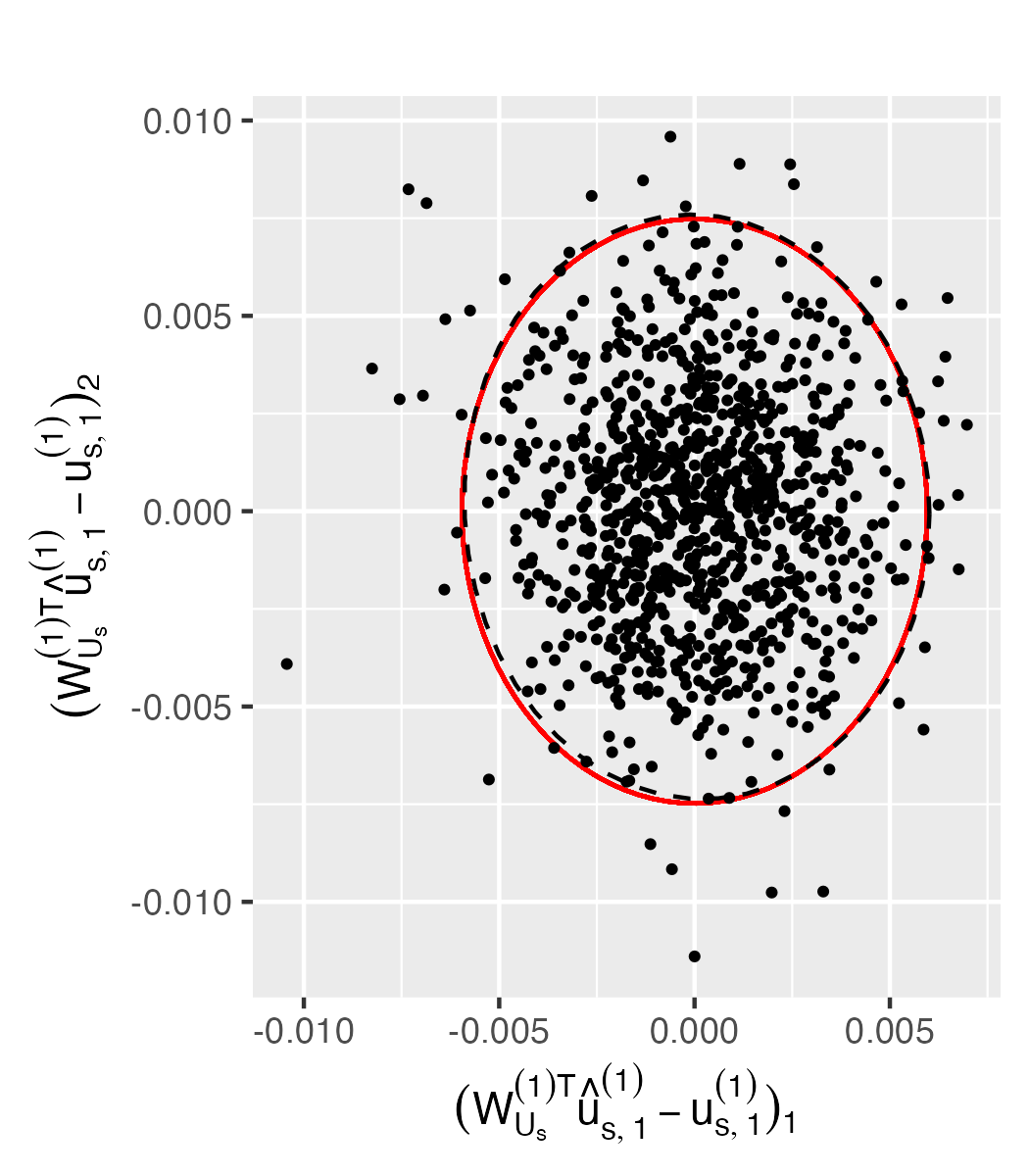}}

\caption{\footnotesize 
Histograms and a bivariate plot of the empirical distributions of the entries of the estimation error $\mw^{(i)\top}_{\muu_s}\hat{u}^{(i)}_{s,k} - {u}^{(i)}_{s,k}$ for $i = 1$ and $k = 1$ are presented. Refer to Figure~\ref{fig:simulation_CLT1_pca} for more details.
}
\label{fig:simulation_CLT2_pca}
\end{figure}

We now present simulations to validate
our theoretical results for distributed PCA. We consider the setting
with $m = 10$, $D=1000$, $d_i \equiv 4$, and $d_{0} = 2$, resulting in
$\muu_c$ and $\muu^{(i)}_s$ being $1000 \times 2$ matrices, and
$\mLambda^{(i)}$ being $2\times 2$ matrices. The orthonormal matrix
$\muu_c$ is randomly generated. For each $i$, orthonormal matrix $\muu^{(i)}_s$, which is orthogonal to $\muu_c$, is also randomly generated. We
generate the diagonal entries of $\bm{\Lambda}^{(i)}$ as iid random variables from the uniform
distribution $U(20,50)$. We then set $\muu^{(i)}=[\muu_c|\muu^{(i)}_s]$, 
$\sigma_i\equiv 1$, and 
$\mSigma^{(i)}=\muu^{(i)}\mLambda^{(i)}\muu^{(i)\top}+(\mi-\muu^{(i)}\muu^{(i)\top})$.
With $n_i \equiv n= 4000$, for each Monte Carlo replicate, we generate
$1000\times 4000$ data matrices $\mx^{(i)}$ whose columns are
independently drawn from the multivariate Gaussian distribution with
mean $\mathbf{0}$ and covariance matrix $\mSigma^{(i)}$. We then apply
Algorithm~\ref{Alg_disPCA} to obtain estimated common subspaces and
individual subspaces. Comparison of the resulting empirical distributions, 
based on $1000$ independent Monte Carlo replicates, against the
limiting distribution given in
Theorem~\ref{thm:PCA_(UhatW-U)_k->normal_previous} is summarized in
Figures~\ref{fig:simulation_CLT1_pca} and
\ref{fig:simulation_CLT2_pca}.

\subsection{Connectivity of brain networks}
\label{sec:brain_data}
In this section, we use the test statistic $T_{ij}$ in Section~\ref{sec:COSIE} to measure similarities between different connectomes constructed from the HNU1 study \citep{zuo2014open}. The data consists of diffusion magnetic resonance imaging (dMRI) records for $30$ healthy adult subjects, where each subject received $10$ dMRI scans over the span of one month. The resulting $m = 300$ dMRIs are then converted into undirected and unweighted graphs on $n = 200$ vertices by registering the brain regions for these images to the CC200 atlas of \cite{craddock2012whole}. 

Taking the $m = 300$ graphs as one realization from an undirected COSIE model, we first apply Algorithm~\ref{Alg} to extract the parameter estimates $\hat{\muu}$, $\hat{\mv}$, $\{\hat{\mr}^{(i)}\}_{i=1}^{300}$ associated with these graphs. 
The initial embedding dimensions $\{d_i\}_{i=1}^{300}$, which range from $5$ to $18$, and the final embedding dimension $d = 11$ are all selected using the (automatic) dimensionality selection procedure described in \cite{zhu2006automatic}.
Given the quantities $\hat{\muu}$, $\hat{\mv}$, and $\{\hat{\mr}^{(i)}\}$, we compute $\hat\mpp^{(i)}=\hat\muu\hat\mr^{(i)}\hat\muu^\top$ for each graph $i$ (and truncate the entries of the resulting $\hat{\mpp}^{(i)}$ to lie in $[0,1]$) before computing $\{\hat\mSigma^{(i)}\}_{i=1}^{300}$ using the formula in Remark~\ref{rm:undirect2}. 
Finally, we compute the test statistic $T_{ij}$ for all pairs $i,j \in [m]$, $i \neq j$, as defined in Theorem~\ref{thm:HT}. 

\begin{figure}[htbp!]
\centering
\subfigure[$T_{ij}$]
{\includegraphics[width=7.5cm]{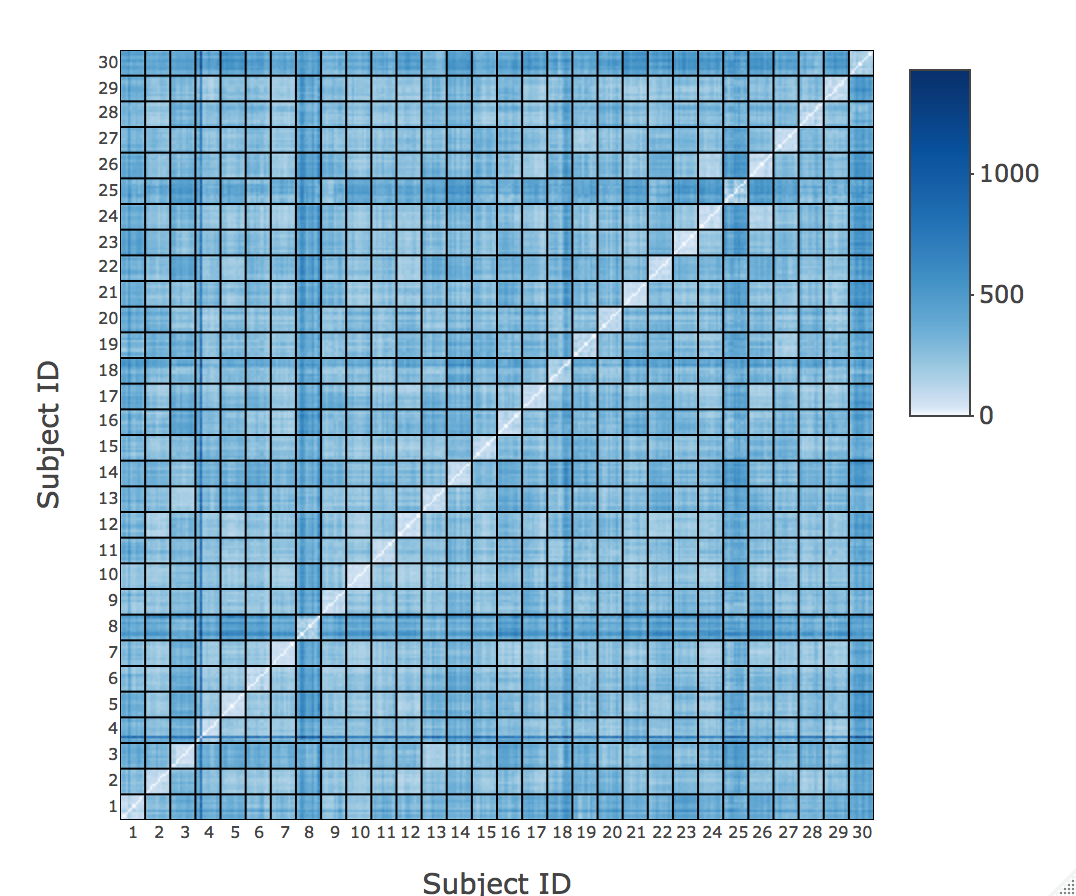}}
\subfigure[$p$-values]
{\includegraphics[width=7.5cm]{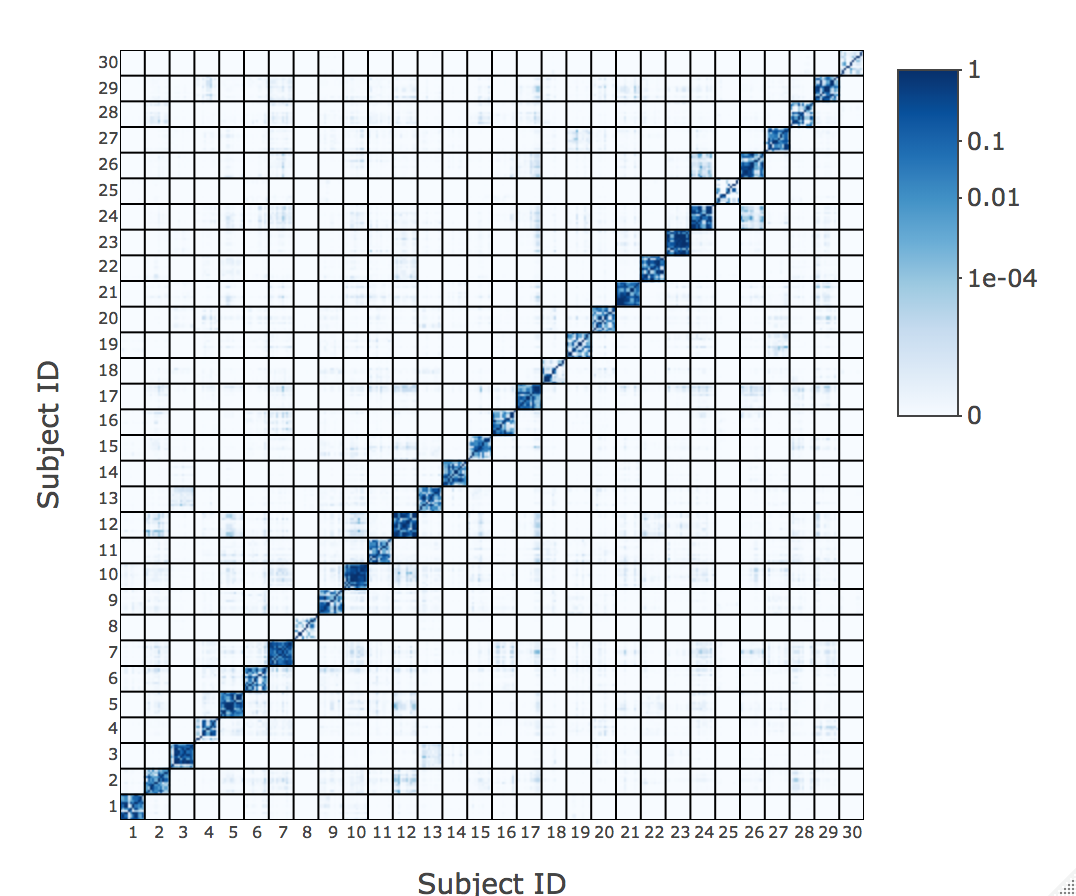}}

\caption{\footnotesize 
Left panel: Test statistic $T_{ij}$ for each pair of brain connectivity networks. Right panel: $p$-values for $T_{ij}$ computed using the $\chi^2$ distribution with $66$ degrees of freedom.
}
\label{fig:brain_Tij}
\end{figure}

\begin{figure}[htbp!]
\centering
\subfigure
{\includegraphics[width=5cm]{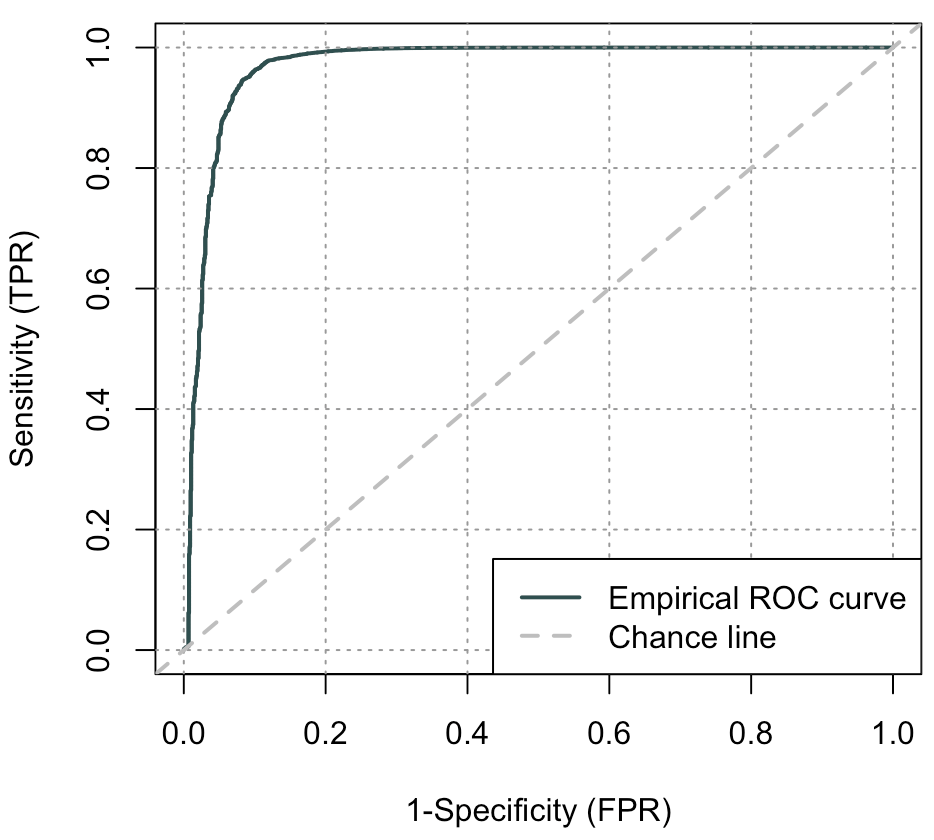}}

\caption{\footnotesize 
ROC curve for classifying whether a pair of graphs represent scans from the same subject (specificity) or from different subjects (sensitivity) as determined by thresholding the values of $T_{ij}$. The corresponding AUC is $0.970$.
}
\label{fig:brain_pval}
\end{figure}

The left panel of Figure~\ref{fig:brain_Tij} shows the matrix of $T_{ij}$ values for all pairs $(i,j) \in [m] \times [m]$ with $i \neq j$, while the right panel presents the $p$-values associated with these $T_{ij}$ (as computed using the $\chi^2$ distribution with $\binom{d}{2} = 66$ degrees of freedom). Note that for ease of presentation, we have rearranged the $m = 300$ graphs so that graphs for the same subject are grouped together, and furthermore we only include on the $x$ and $y$ axes the labels for the subjects but not the individual scans within each subject. We see that our test statistic $T_{ij}$ can discern between scans from the same subject (where $T_{ij}$ are generally small) and scans from different subjects (where $T_{ij}$ are quite large). Indeed, given any two scans $i$ and $j$ from different subjects, the $p$-value for $T_{ij}$ (under the null hypothesis that $\mr^{(i)} = \mr^{(j)}$) is always smaller than $0.01$.
Figure~\ref{fig:brain_pval} shows the ROC curve when we use $T_{ij}$ to classify whether a pair of graphs represents scans from the same subject (specificity) or from different subjects (sensitivity). The corresponding AUC is $0.970$ and is thus close to optimal.

The HNU1 data have also been analyzed in \cite{arroyo2019inference}. In particular, \cite{arroyo2019inference} proposes $\|\hat\mr^{(i)}-\hat\mr^{(j)}\|_F^2$ as a test statistic, and instead of computing $p$-values from some limiting distribution directly, \cite{arroyo2019inference} calculates empirical $p$-values using: 1) a parametric bootstrap approach; 2) the asymptotic null distribution of $\|\hat\mr^{(i)}-\hat\mr^{(j)}\|_F^2$. By neglecting the effect of the bias term $\mh^{(i)}$, \cite{arroyo2019inference} approximates the null distribution of $\|\hat\mr^{(i)}-\hat\mr^{(j)}\|_F^2$ as a generalized $\chi^2$ distribution and estimate it by Monte Carlo simulations of a mixture of normal distributions with the estimates $\hat\mSigma^{(i)}$ and $\hat\mSigma^{(j)}$.

Comparing the $p$-values of our test in Figure~\ref{fig:brain_Tij} with the results obtained by their two methods in Figure~15, we see that for different methods, the ratios of the $p$-values for pairs from the same subject to those for pairs from different subjects are very similar. Thus, both test statistics can detect whether pairs of graphs are from the same subject well. Our test statistic, however, has the benefit that its $p$-value is computed using a large-sample $\chi^2$ approximation and is thus much less computationally intensive compared to test procedures that use bootstrapping and other Monte Carlo simulations.

\subsection{Worldwide food trade networks}
\label{sec:food_trading}
For the next example, we use the trade networks between countries for different food and agriculture products during the year $2018$. The data is collected by the Food and Agriculture Organization of the United Nations and is available at \url{https://www.fao.org/faostat/en/#data/TM}. We construct a collection of networks, one for each product, where vertices represent trade entities (countries or regions) and the edges in each network represent trade relationships between trade entities; the resulting adjacency matrices $\{\mathbf{A}^{(i)}\}$ are directed but unweighted as we (1) set $\ma^{(i)}_{rs}=1$ if trade entity $r$ exports product $i$ to trade entity $s$, and (2) ignore any links between trade entities $r$ and $s$ in $\mathbf{A}^{(i)}$ if their total trade amount for the $i$-th product is less than two hundred thousand US dollars. Finally, we extract the \emph{intersection} of the \emph{largest connected components} of $\{\mathbf{A}^{(i)}\}$ and obtain $56$ networks on a set of $75$ shared vertices.



\begin{figure}[htbp!]
\centering
\subfigure
{\includegraphics[width=17cm]{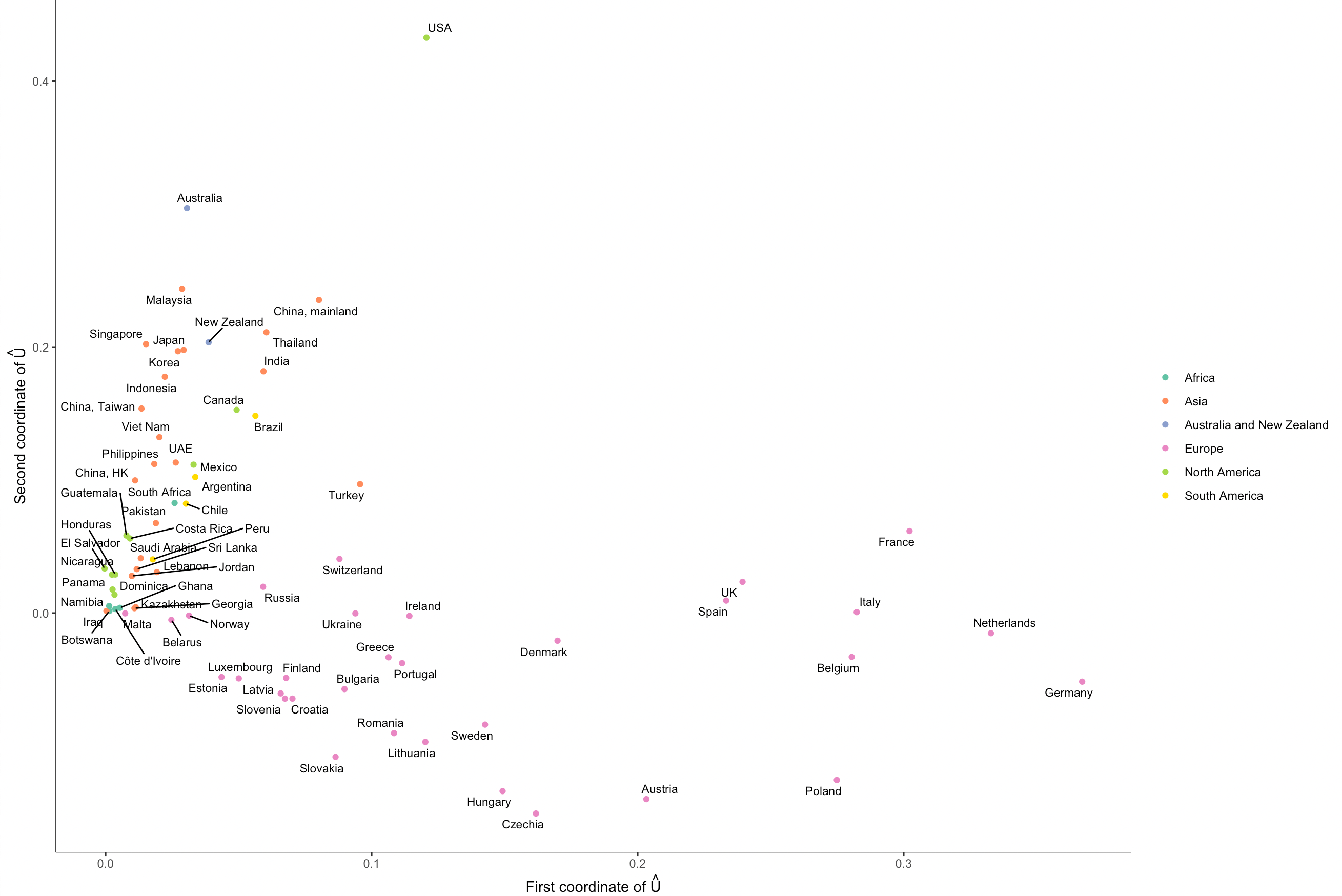}}

\caption{\footnotesize 
Latent positions of trade entities as exporter
}
\label{fig:food_countries_U}
\end{figure}

\begin{figure}[htbp!]
\centering
\subfigure
{\includegraphics[width=17cm]{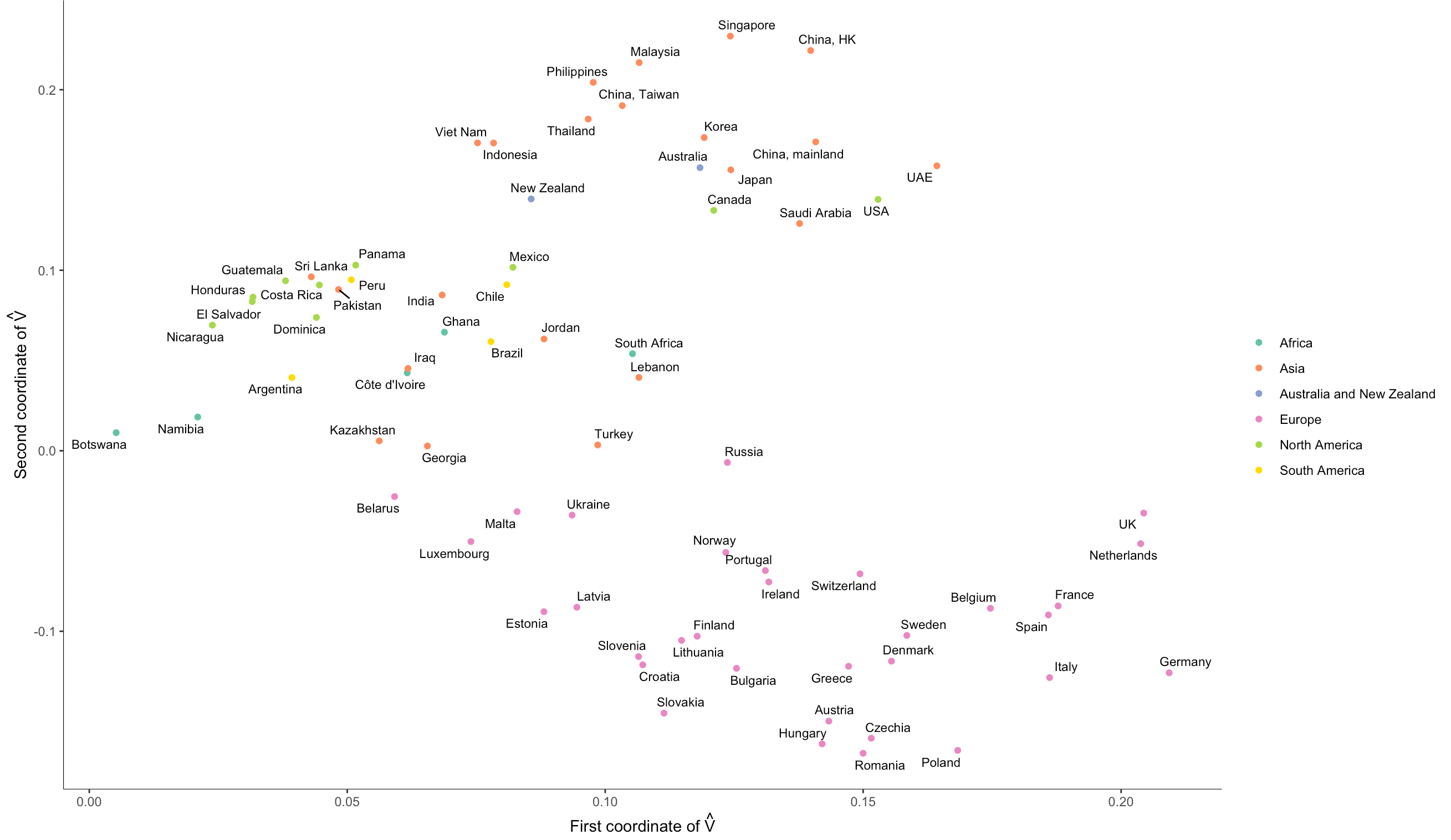}}

\caption{\footnotesize 
Latent positions of trade entities as importer
}
\label{fig:food_countries_V}
\end{figure}

Taking the $m = 56$ networks as one realization from a directed COSIE model, we apply Algorithm~\ref{Alg} to compute the parameter estimates $\hat{\muu}$, $\hat{\mv}$, $\{\hat{\mr}^{(i)}\}_{i=1}^{56}$ associated with these graphs with initial embedding dimensions $\{d_i\}_{i=1}^{56}$ as well as the final embedding dimension $d$ all chosen to be $2$. Figure~\ref{fig:food_countries_U} and Figure~\ref{fig:food_countries_V} present scatter plots for the rows of $\hat{\mathbf{U}}$ and $\hat{\mathbf{V}}$, respectively; we interpret the $r$th row of $\hat\muu$ (resp. $\hat\mv$) as representing the estimated latent position for this country as an exporter (resp. importer). We see that there is a high degree of correlation between these estimated latent positions and the true underlying geographic proximities, e.g., countries in the same continent are generally placed close together in Figure~\ref{fig:food_countries_U} and Figure~\ref{fig:food_countries_V}.

\begin{figure}[htbp!]
\centering
\subfigure
{\includegraphics[width=11cm]{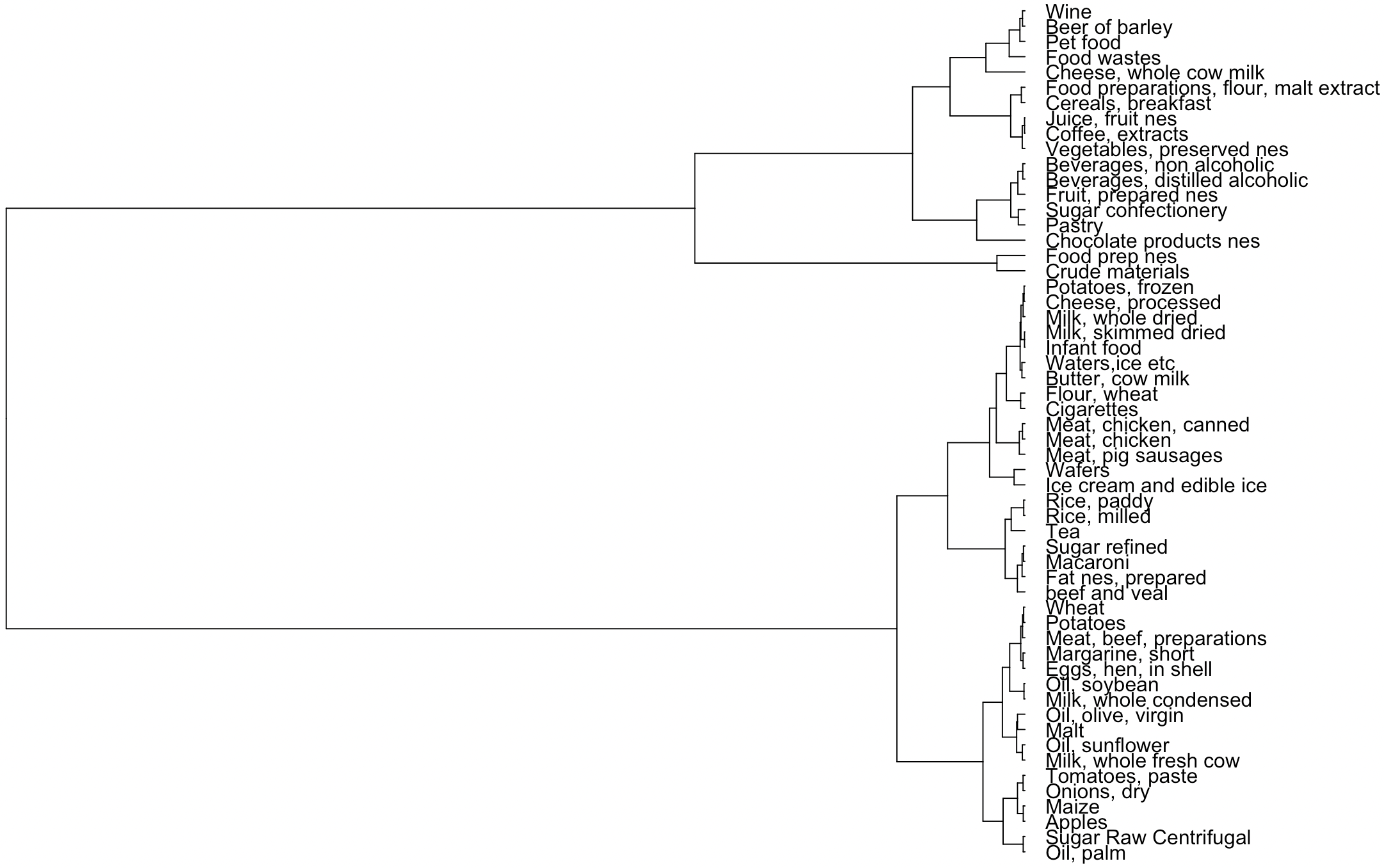}}
\caption{\footnotesize 
Hierarchical clustering of food products
}
\label{fig:food_items_hcluster}
\end{figure}

Next, we compute the statistic $T_{ij}$ in Theorem~\ref{thm:HT} to measure the differences between $\hat{\mathbf{R}}^{(i)}$ and $\hat{\mathbf{R}}^{(j)}$ for all pairs of products $\{i,j\}$.
Viewing $(T_{ij})$ as a distance matrix, we organize the food products using hierarchical clustering \citep{johnson1967hierarchical}; see the dendrogram in Figure~\ref{fig:food_items_hcluster}. There appear to be two main clusters formed by raw/unprocessed products (bottom cluster) and processed products (top cluster), which suggest discernible differences in the trade patterns for these types of products.

The trade dataset (but for 2010) has also been analyzed in \cite{jing2021community}. In particular, \cite{jing2021community} studies the mixture multilayer SBM and propose a tensor-based algorithm to reveal memberships of vertices and memberships of layers. For the food trading networks, \cite{jing2021community} first groups the layers, i.e., the food products, into two clusters, and then obtains the embeddings and the clustering result of the trade entities for each food cluster. Our results are similar to theirs. In particular, their clustering of the food products also shows a difference in the trade patterns for unprocessed and processed foods, while their clustering of the trade entities is also related to geographical location. However, as we also compute the test statistic $T_{ij}$ for each pair of products, we obtain a more detailed analysis of the product relationships. In addition, as we keep the orientation of the edges (and thus our graphs are directed), we can also analyze the trade entities in terms of both their export and import behavior, and Figures~\ref{fig:food_countries_U} and \ref{fig:food_countries_V} show that there is indeed some difference between these behaviors, e.g., the USA and Australia are outliers as exporters but are clustered with other trade entities as importers.

\subsection{Distributed PCA and MNIST}
\begin{figure}[htbp!]
\centering
\subfigure
{\includegraphics[width=5.5cm]{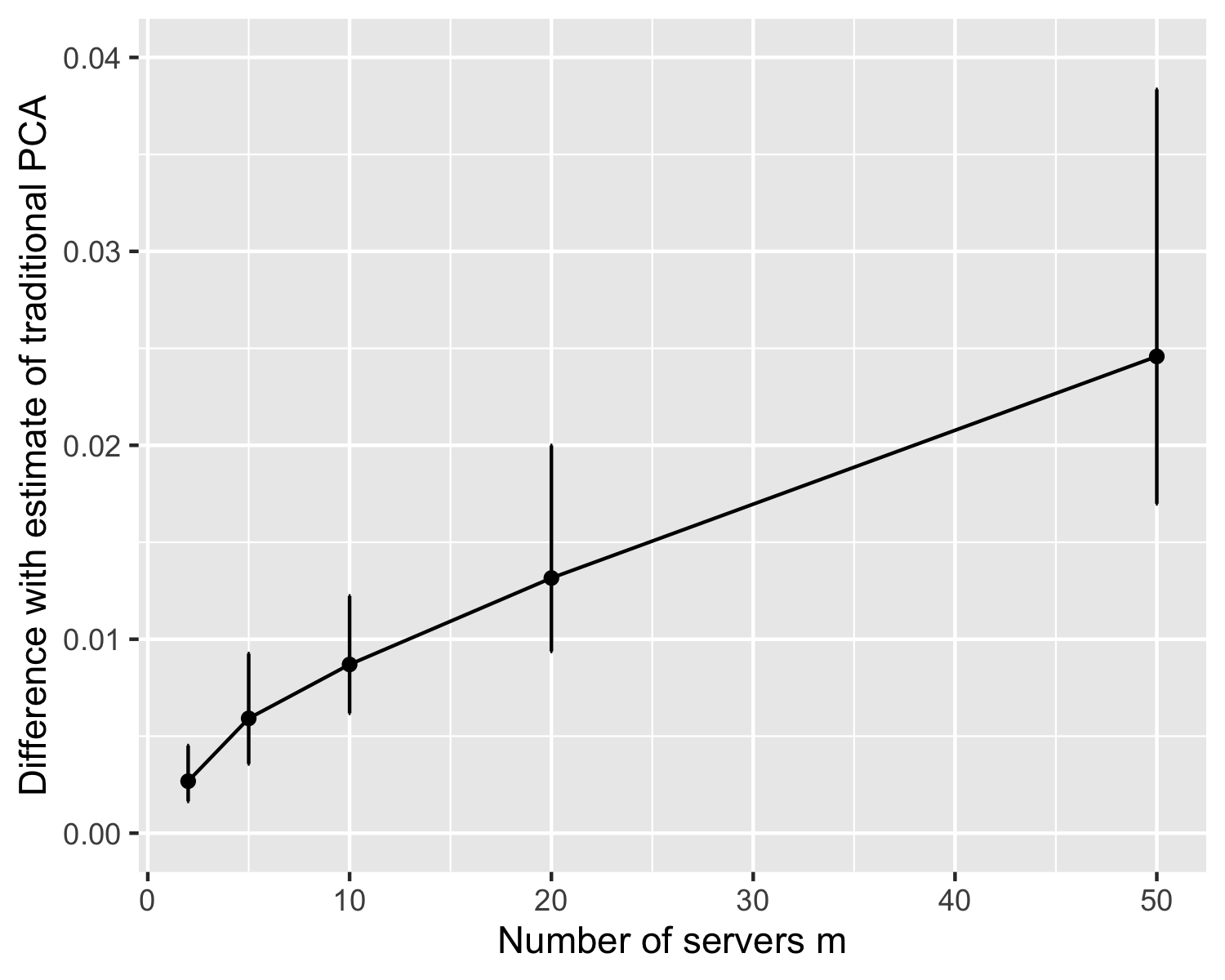}}
\caption{\footnotesize 
Empirical estimates for the difference between the $d = 9$ leading principal components of the MNIST data as computed by traditional PCA and by distributed PCA with $m \in \{2,5,10,20,50\}$. The difference is quantified by $\min_{\mw\in\mathcal{O}_d}\|\hat\muu^{(m)}\mw-\hat\muu\|_F$. The estimates (together with the $95\%$ confidence intervals) are based on $100$ independent Monte Carlo replicates.
}
\label{fig:disPCA_realdata_change_m}
\end{figure}

We now perform dimension reduction on the MNIST dataset using distributed PCA for the case where the covariance matrix is shared across $m \geq 2$ nodes and compare the result against traditional PCA ($m=1$) on the full dataset. The MNIST data consists of $60{,}000$ grayscale images of handwritten digits of the numbers $0$ through $9$. Each image is of size $28 \times 28$ pixels and can be viewed as a vector in $\mathbb{R}^{784}$ with entries in $[0,255]$.
Letting $\mathbf{X}$ be the $60{,}000 \times 784$ matrix whose rows represent the images, we first extract the matrix $\hat{\muu}$ whose columns are the $d = 9$ leading principal components of $\mathbf{X}$. The choice $d = 9$ is arbitrary and is chosen purely for illustrative purposes. Next, we approximate $\hat\muu$ using distributed PCA by randomly splitting $\mathbf{X}$ into $m \in \{2,5,10,20,50\}$ subsamples. Letting $\hat{\muu}^{(m)}$ be the resulting approximation, we compute $\min_{\mw\in\mathcal{O}_d}\|\hat\muu^{(m)}\mw-\hat\muu\|_{F}$. We repeat these steps for $100$ independent Monte Carlo replicates and summarize the results in Figure~\ref{fig:disPCA_realdata_change_m}, which shows that the errors between $\hat\muu^{(m)}$ and $\hat{\muu}$ are always substantially smaller than $\|\hat{\muu}\|_F = \|\muu\|_F=3$.
We emphasize that while the errors in Figure~\ref{fig:disPCA_realdata_change_m} do increase with $m$, this is mainly an artifact of the experimental setup as there is no underlying ground truth and we are only using $\hat{\muu}$ as a surrogate for some unknown (or possibly non-existent) $\muu$. In other words, $\hat{\muu}$ is noise-free in this setting while $\hat{\muu}^{(m)}$ is inherently noisy, and thus it is reasonable for the noise level in $\hat{\muu}^{(m)}$ to increase with $m$. Finally, we note that for this experiment, we have assumed that the rows of $\mathbf{X}$ are iid samples from a \emph{mixture} of $10$ multivariate Gaussians with each component corresponding to a number in $\{0,1,\dots,9\}$. As a mixture of multivariate Gaussians is sub-Gaussian, the results in Section~\ref{sec:distributed_PCA} remain relevant in this setting; see Remark~\ref{rem:general_covariance}.

\section{Discussion}
\label{sec:discussion}

In this paper, we present a general framework for deriving limit results for distributed estimation of the leading singular vectors for a collection of matrices with shared invariant subspaces and possibly distinct individual subspaces, and apply this framework to multiple heterogeneous network inference and distributed PCA. 

We now mention several potential related directions for future research on multiple network inference and distributed PCA.
First, the COISIE model has low-rank edge probability matrices $\{\mpp^{(i)}\}_{i=1}^{m}$, while for distributed PCA, the intrinsic rank of $\bm{\Sigma}$ grows at order $D^{1 - \gamma}$ for some $\gamma \in (0,1]$ and can thus be arbitrarily close to ``full'' rank. This suggests that we can extend our results to general edge-independent random graphs where the ranks of each $\mpp^{(i)}$ grow with $n$. The main challenge is then in formulating a sufficiently general and meaningful yet still tractable model under these constraints. 
Second, the results for distributed PCA in this paper assume that for each $i \in [m]$, the estimate $\hat{\muu}^{(i)}$ is given by the leading eigenvectors of the sample covariance matrix $\hat{\mSigma}^{(i)}$. If the eigenvectors in $\muu^{(i)}$ are known to be sparse, then it might be more desirable to let each $\hat{\muu}^{(i)}$ be computed from $\hat{\mSigma}^{(i)}$ using some sparse PCA algorithm (see, e.g., \cite{amini,vu_lei_rohe,d_apresmont}) and then aggregate these estimates to yield a final $\hat{\muu}$. Recently, \cite{agterberg2022entrywise} derives $\ell_{2\to\infty}$ bounds for sparse PCA given a single sample covariance $\hat{\mSigma}$ under a general high-dimensional subgaussian design and thus, by combining their analysis with ours, it may be possible to also obtain limit results for $\hat{\muu}$ in distributed \emph{sparse} PCA.
{\color{black}
Third, we are interested in extending Theorem~\ref{thm:What_U Rhat What_v^T-R->norm} and Theorem~\ref{thm:HT} to the $o(n^{1/2})$ regime but, as we discussed in Remark~\ref{rem:sqrt_n}, this appears to be highly challenging as related existing results all require $\omega(n^{1/2})$. Nevertheless, while the asymptotic bias for $\mw_{\muu}^{\top} \hat{\mr}^{(i)} \mw_{\mv} - \mr^{(i)}$ is important for the normal approximation of a single $\hat{\mr}^{(i)}$, we surmise that part of this asymptotic bias is not essential for two-sample testing. Indeed, the test statistic in Theorem~\ref{thm:HT} depends on the difference $\hat{\mr}^{(i)} - \hat{\mr}^{(j)}$, and under the null hypothesis $\mathbb{H}_0 \colon \mr^{(i)} = \mr^{(j)}$ a more refined analysis of this difference may cancel some of the asymptotic bias terms, thereby allowing for a weaker condition on $n\rho_n$. We therefore expect that Theorem~\ref{thm:HT} will continue to hold even in the $o(n^{1/2})$ regime, although a more careful and technically involved analysis is required and is left for future work.
}

{\color{black}
}

Finally, as we alluded to in the introduction, our framework can also be applied to other matrix estimation problems,
such as the joint and individual variation explained (JIVE) model for integrative data analysis \citep{jive,feng2018angle} and population value decomposition for the analysis of image populations \citep{population_svd}. Taking JIVE as a specific example, recall that the JIVE model assumes that there are $m$ data matrices $\{\mathbf{X}^{(i)}\}_{i=1}^{m}$ where each $\mathbf{X}^{(i)}$ is of dimension $d_i \times n$; the columns of $\mathbf{X}^{(i)}$ correspond to experimental subjects while the rows correspond to features. 
Furthermore, $\mathbf{X}^{(i)}$ are modeled as $\mathbf{X}^{(i)} = \mathbf{J}^{(i)} + \mathbf{I}^{(i)} + \mathbf{N}^{(i)}$ where $\{\mathbf{J}^{(i)}\}_{i=1}^m$ share a common row space (denoted as $\mathbf{J}_*$), $\mathbf{I}^{(i)}$ represent individual structures, and $\mathbf{N}^{(i)}$ denote additive noise perturbations. The estimation of $\mathbf{J}_*$ and $\{\mathbf{I}^{(i)}\}_{i=1}^{m}$ can be done using the aJIVE procedure \citep{feng2018angle}
that is very similar to Algorithm~\ref{Alg_COISIE} in our paper. 
While \cite{feng2018angle} presents criteria for choosing the dimensions for $\mathbf{J}_*$ and $\{\mathbf{I}^{(i)}\}$, it does not provide theoretical guarantees for the estimation of $\mathbf{J}_*$ and $\{\mathbf{I}^{(i)}\}$; this is partly because they did not consider any noise model for $\mathbf{N}^{(i)}$. We surmise that if the entries of each  $\mathbf{N}^{(i)}$ are independent
mean-zero sub-Gaussian variables then $2 \to \infty$ norm error bounds for estimating
$\mathbf{J}_*$ and $\{\mathbf{I}^{(i)}\}$ can be obtained following the same analysis as that done for the COISIE model.

\end{sloppypar}

\newpage

\bibliographystyle{unsrtnat}
\bibliography{sample}

@book{bhatia2013matrix,
  title={Matrix analysis},
  author={R.~Bhatia},
  year={2013},
  publisher={Springer}
}

@article{distributedPCA_review,
author = {S.~Wu and H.-T~Wai and L.~Lin and A.~Scaglione},
title = {A review of distributed algorithms for principal component analysis},
year = {2018},
journal = {Proceedings of the {IEEE}},
volume = {106},
pages = {1321--1340}
}

@article{bertrand_moonen,
author = {A.~Bertrand and M.~Moonen},
title = {Distributed adaptive estimation of covariance matrix eigenvectors in wireless sensor networks with applications to distributed {PCA}},
journal = {Signal Processing},
volume = 104,
pages = {120--135},
year = 2014
}

@article{lin_scaglione2,
author = {L.~Lin and A.~Scaglione and J.~H.~Manton},
title = {Distributed Principal Subspace Estimation in Wireless Sensor Networks},
journal = {IEEE Journal of Selected Topics in Signal Processing},
year = {2011},
pages = {725--738},
volume = {5},
}

@article{jiao_gu,
author = {Y.~Jiao and Y.~Gu},
title = {Communication-Efficient Decentralized Subspace Estimation},
journal = {IEEE Journal of Selected Topics in Signal Processing},
year = {2022},
pages = {516--531},
volume = {16},
}

@article{roy_vetterli,
title = {Dimensionality Reduction for Distributed Estimation in the Infinite Dimensional Regime},
author = {O.~Roy and M.~Vetterli},
journal = {IEEE Transactions on Information Theory},
volume = {54},
pages = {1655--1669},
year = {2008},
}

@article{jive,
author = {E.~F.~Lock and K.~A.~Hoadley and J.~S.~Marron and A.~B.~Nobel},
title = {Joint and individual variation explained ({JIVE}) for integrated analysis of multiple data types},
journal = {Annals of Applied Statistics},
volume = {7},
year = {2013},
pages = {523--542}
}

@article{hector,
author = {E.~Hector and P.~X.~K.~Song},
title = { A distributed and integrated method of moments for high-dimensional correlated data analysis},
journal = {Journal of the American Statistical Association},
volume = {116},
pages = {805--818},
year = {2021}}

@article{amini,
author = {A.~A.~Amini and M.~J.~Wainwright},
title = {High-dimensional analysis of semidefinite relaxations for sparse principal components},
journal = {Annals of Statistics},
year = {2009},
volume = {37},
pages = {2877--2921}
}

@article{vu_lei_rohe,
author = {V.~Q.~Vu and J.~Cho and J.~Lei and K.~Rohe},
title = {Fantope projection and selection: a near-optimal convex relaxation of sparse {PCA}},
journal = {Advances in Neural Informations Processing Systems},
year = {2013},
pages = {2670--2678},
volume = {26}
}

@article{d_apresmont,
author = {A. d’Aspremont and L. El Ghaoui and M. I. Jordan and G. R. G. Lanckriet},
title = {A Direct Formulation for Sparse PCA Using Semidefinite Programming},
journal = {SIAM Review},
volume = {49},
pages = {434–448},
year = {2007}
}

@article{Spiked,
author = {Johnstone, Iain},
year = {2001},
pages = {295--327},
title = {On the Distribution of the Largest Eigenvalue in Principal Components Analysis},
volume = {29},
journal = {Annals of Statistics},
doi = {10.1214/aos/1009210544}
}

@InProceedings{pmlr-streaming-PCA,
  title = 	 {Streaming Principal Component Analysis in Noisy Setting},
  author =       {Marinov, Teodor Vanislavov and Mianjy, Poorya and Arora, Raman},
  booktitle = 	 {Proceedings of the 35th International Conference on Machine Learning},
  pages = 	 {3413--3422},
  year = 	 {2018},
  volume = 	 {80},
}

@article{zalesky2010network,
  title={Network-based statistic: identifying differences in brain networks},
  author={Zalesky, Andrew and Fornito, Alex and Bullmore, Edward T},
  journal={Neuroimage},
  volume={53},
  pages={1197--1207},
  year={2010},
  publisher={Elsevier}
}

@article{rubinov2010complex,
  title={Complex network measures of brain connectivity: uses and interpretations},
  author={Rubinov, Mikail and Sporns, Olaf},
  journal={Neuroimage},
  volume=52,
  pages={1059--1069},
  year=2010,
  publisher={Elsevier}
}

@article{he2008structural,
  title={Structural insights into aberrant topological patterns of large-scale cortical networks in Alzheimer's disease},
  author={He, Yong and Chen, Zhang and Evans, Alan},
  journal={Journal of Neuroscience},
  volume={28},
  pages={4756--4766},
  year={2008},
  publisher={Soc Neuroscience}
}

@article{Tony_sparse,
author = {Cai, T. and Ma, Zongming and Wu, Yihong},
year = {2013},
pages = {781--815},
title = {Optimal Estimation and Rank Detection for Sparse Spiked Covariance Matrices},
volume = {161},
journal = {Probability Theory and Related Fields},
doi = {10.1007/s00440-014-0562-z}
}

@article{sparse_AOS,
author = {Birnbaum, Aharon and Johnstone, Iain and Nadler, Boaz and Paul, Debashis},
year = {2012},
pages = {1055--1084},
title = {Minimax bounds for sparse {PCA} with noisy high-dimensional data},
volume = {41},
journal = {Annals of Statistics},
doi = {10.1214/12-AOS1014}
}

@article{sparse_AOS1,
author = {Berthet, Quentin and Rigollet, Philippe},
year = {2012},
pages = {1780--1815},
title = {Optimal detection of sparse principal components in high dimension},
volume = {41},
journal = {Annals of Statistics},
doi = {10.1214/13-AOS1127}
}

@book{bai_sample_cov,
author = {J.~Yao and S.~Zheng and Z.~Bai},
title = {Large Sample Covariance Matrices and High-Dimensional Data Analysis},
year = {2015},
publisher = {Cambridge University Press}
}

@inproceedings{sparse,
author = {Vu, Vincent and Lei, Jing},
year = {2012},
pages = {1276--1286},
title = {Minimax Rates of Estimation for Sparse PCA in High Dimensions},
booktitle = {Proceedings of the fifteenth international conference on artificial intelligence and statistics}
}

@article{yu2015useful,
  title={A useful variant of the Davis--Kahan theorem for statisticians},
  author={Yu, Yi and Wang, Tengyao and Samworth, Richard J},
  journal={Biometrika},
  volume={102},
  number={2},
  pages={315--323},
  year={2015},
  publisher={Oxford University Press}
}

@article{Lounici2017,
  title={Concentration inequalities and moment bounds for sample covariance operators},
  author={Koltchinskii, Vladimir and Lounici, Karim},
  journal={Bernoulli},
  volume={23},
  number={1},
  pages={110--133},
  year={2017},
  publisher={Bernoulli Society for Mathematical Statistics and Probability}
}

@article{arroyo2019inference,
  title={Inference for multiple heterogeneous networks with a common invariant subspace},
  author={Arroyo, Jes{\'u}s and Athreya, Avanti and Cape, Joshua and Chen, Guodong and Priebe, Carey E and Vogelstein, Joshua T},
  journal={Journal of Machine Learning Research},
  volume={22},
  number={142},
  pages={1--49},
  year={2021}
}

@article{rinaldo_2013,
  author =     {J.~Lei and A. Rinaldo},
  title =    {Consistency of spectral clustering in stochastic blockmodels},
  journal =     {Annals of Statistics},
  year =     {2015},
  volume = {43},
  pages = {215--237}
}

@article{mao_sarkar,
  title={Estimating mixed memberships with sharp eigenvector deviations},
  author={Mao, Xueyu and Sarkar, Purnamrita and Chakrabarti, Deepayan},
  journal={Journal of the American Statistical Association},
  volume={116},
  number={536},
  pages={1928--1940},
  year={2021},
  publisher={Taylor \& Francis}
}

@book{stewart_sun,
author = {G.~W.~Stewart and J.-G.~Sun},
title = {Matrix perturbation theory},
publisher = {Academic Press},
year = {1990}
}

@unpublished{oliveira2009concentration,
  title =    {Concentration of the adjacency matrix and of the
                  {L}aplacian in random graphs with independent edges},
  author =     {R.~I.~Oliveira},
  note =     {\url{http://arxiv.org/abs/0911.0600}},
  year =     {2009}
}

@article{cai2021subspace,
  title={Subspace estimation from unbalanced and incomplete data matrices: $\ell_{2,\infty}$ statistical guarantees},
  author={Cai, Changxiao and Li, Gen and Chi, Yuejie and Poor, H Vincent and Chen, Yuxin},
  journal={Annals of Statistics},
  volume={49},
  pages={944--967},
  year={2021},
  publisher={Institute of Mathematical Statistics}
}

@article{davis70,
  title={The rotation of eigenvectors by a perturbation. III},
  author={Davis, Chandler and Kahan, William Morton},
  journal={SIAM Journal on Numerical Analysis},
  volume={7},
  number={1},
  pages={1--46},
  year={1970},
  publisher={SIAM}
}

@article{cape2019two,
  title={The two-to-infinity norm and singular subspace geometry with applications to high-dimensional statistics},
  author={Cape, Joshua and Tang, Minh and Priebe, Carey E},
  journal={Annals of Statistics},
  volume={47},
  number={5},
  pages={2405--2439},
  year={2019},
  publisher={Institute of Mathematical Statistics}
}

@unpublished{carlsson2018perturbation,
  title={Perturbation theory for the matrix square root and matrix modulus},
  author={Carlsson, Marcus},
  note={arXiv preprint at \url{https://arxiv.org/abs/1810.01464}},
  year={2018}
}

@article{bandeira2016sharp,
  title={Sharp nonasymptotic bounds on the norm of random matrices with independent entries},
  author={Bandeira, Afonso S and Van Handel, Ramon},
  journal={Annals of Probability},
  volume={44},
  pages={2479--2506},
  year={2016},
  publisher={Institute of Mathematical Statistics}
}

@book{van2000asymptotic,
  title={Asymptotic statistics},
  author={Van der Vaart, Aad W},
  year={2000},
  publisher={Cambridge University Press}
}

@article{cape2019signal,
  title={Signal-plus-noise matrix models: eigenvector deviations and fluctuations},
  author={Cape, Joshua and Tang, Minh and Priebe, Carey E},
  journal={Biometrika},
  volume={106},
  number={1},
  pages={243--250},
  year={2019},
  publisher={Oxford University Press}
}

@article{tang2018eigenvalues,
  title={Eigenvalues of stochastic blockmodel graphs and random graphs with low-rank edge probability matrices},
  author={A.~Athreya and J.~Cape and M.~Tang},
  journal={Sankhya A},
  volume = {84},
  pages = {36--63},
  year={2022}
}

@article{zuo2014open,
  title={An open science resource for establishing reliability and reproducibility in functional connectomics},
  author={X.~Zuo and others},
  journal={Scientific Data},
  volume={1},
  pages={1--13},
  year={2014},
  publisher={Nature Publishing Group}
}

@article{craddock2012whole,
  title={A whole brain fMRI atlas generated via spatially constrained spectral clustering},
  author={R.~C.~Craddock and G.~A.~James and P.~E.~Holtzheimer III and X.~P.~Hu and H.~S.~Mayberg},
  journal={Human Brain Mapping},
  volume={33},
  number={8},
  pages={1914--1928},
  year={2012},
  publisher={Wiley Online Library}
}

@article{zhu2006automatic,
  title={Automatic dimensionality selection from the scree plot via the use of profile likelihood},
  author={M.~Zhu and A.~Ghodsi},
  journal={Computational Statistics \& Data Analysis},
  volume={51},
  number={2},
  pages={918--930},
  year={2006},
  publisher={Elsevier}
}

@article{magnus1979commutation,
  title={The commutation matrix: some properties and applications},
  author={J.~R.~Magnus and H.~Neudecker},
  journal={Annals of Statistics},
  volume={7},
  number={2},
  pages={381--394},
  year={1979},
  publisher={Institute of Mathematical Statistics}
}

@article{neudecker1986symmetry,
  title={Symmetry, 0-1 matrices and {J}acobians},
  author={Neudecker, H},
  journal={Econometric Theory},
  volume=2,
  pages={157--190},
  year=1986,
  publisher={Citeseer}
}

@article{fan2019distributed,
  title={Distributed estimation of principal eigenspaces},
  author={J.~Fan and D.~Wang and K.~Wang and Z.~Zhu},
  journal={Annals of Statistics},
  volume={47},
  pages={3009--3031},
  year={2019},
}

@article{van2008visualizing,
  title   = {Visualizing Data using t-SNE},
  author  = {L.~Van der Maaten and G.~Hinton},
  journal = {Journal of Machine Learning Research},
  year    = {2008},
  volume  = {9},
  number  = {86},
  pages   = {2579-2605}
}

@article{johnson1967hierarchical,
  title={Hierarchical clustering schemes},
  author={S.~C.~Johnson},
  journal={Psychometrika},
  volume={32},
  number={3},
  pages={241--254},
  year={1967},
  publisher={Springer}
}

@article{chen2021spectral,
  title={Spectral methods for data science: a statistical perspective},
  author={Chen, Yuxin and Chi, Yuejie and Fan, Jianqing and Ma, Cong},
  journal={Foundations and Trends{\textregistered} in Machine Learning},
  volume={14},
  number={5},
  pages={566--806},
  year={2021},
  publisher={Now Publishers, Inc.}
}

@article{abbe2020entrywise,
  title={Entrywise eigenvector analysis of random matrices with low expected rank},
  author={Abbe, Emmanuel and Fan, Jianqing and Wang, Kaizheng and Zhong, Yiqiao},
  journal={Annals of Statistics},
  volume={48},
  number={3},
  pages={1452--1474},
  year={2020},
  publisher={NIH Public Access}
}

@unpublished{lei2019unified,
  title={Unified $\ell_{2\to\infty}$ eigenspace perturbation theory for symmetric random matrices},
  author={Lei, Lihua},
  note={arXiv preprint at \url{https://arxiv.org/abs/1909.04798}},
  year={2019}
}

@article{damle,
  author = {A.~Damle and Y.~Sun},
  title = {Uniform bounds for invariant subspace perturbations},
  journal = {SIAM Journal on Matrix Analysis and Applications},
  volume={41},
  number={3},
  pages={1208--1236},
  year={2020},
  publisher={SIAM}
}

@article{fan2018eigenvector,
  title={An $\ell_{\infty}$ eigenvector perturbation bound and its application to robust covariance estimation},
  author={Fan, Jianqing and Wang, Weichen and Zhong, Yiqiao},
  journal={Journal of Machine Learning Research},
  volume={18},
  number={207},
  pages={1--42},
  year={2018},
  publisher={Microtome Publishing}
}

@unpublished{chen2021classification,
  title={Classification of high-dimensional data with spiked covariance matrix structure},
  author={Y.~Chen and M.~Tang},
  note={arXiv preprint at \url{https://arxiv.org/abs/2110.01950}},
journal={Transactions on Machine Learning Research},
  year={2026}
}

@article{zhong_boumal,
  title = {Near-optimal bounds for phase synchronization},
  author = {Y.~Zhong and N.~Boumal},
  journal = {SIAM Journal on Optimization},
  volume={28},
  number={2},
  pages={989--1016},
  year={2018},
  publisher={SIAM}
}

@article{javanmard_montanari,
  title = {Debiasing the lasso: optimal sample size for Gaussian designs},
  author = {A.~Javanmard and A.~Montanari},
  journal = {Annals of Statistics},
  volume={46},
  number={6A},
  pages={2593--2622},
  year={2018},
  publisher={Institute of Mathematical Statistics}
}

@article{cai2013sparse,
  title={Sparse PCA: optimal rates and adaptive estimation},
  author={Cai, T Tony and Ma, Zongming and Wu, Yihong},
  journal={Annals of Statistics},
  volume={41},
  number={6},
  pages={3074--3110},
  year={2013},
  publisher={Institute of Mathematical Statistics}
}

@article{jing2021community,
  title={Community detection on mixture multilayer networks via regularized tensor decomposition},
  author={Jing, Bing-Yi and Li, Ting and Lyu, Zhongyuan and Xia, Dong},
  journal={Annals of Statistics},
  volume={49},
  number={6},
  pages={3181--3205},
  year={2021},
  publisher={Institute of Mathematical Statistics}
}

@article{lei2020bias,
  title={Bias-adjusted spectral clustering in multi-layer stochastic block models},
  author={Lei, Jing and Lin, Kevin Z},
  journal={Journal of the American Statistical Association},
  year={2022+},
  publisher={Taylor \& Francis}
}

@article{grdpg1,
  author = {P.~Rubin-Delanchy and J.~Cape and M.~Tang and C.~E.~Priebe},
  title = {A statistical interpretation of spectral embedding: the generalised random dot product graph},
  journal = {Journal of the Royal Statistical Society, Series B},
  volume = {84},
  year = {2022},
  pages = {1446--1473}
}

@article{wang2019joint,
  title={Joint embedding of graphs},
  author={Wang, Shangsi and Arroyo, Jes{\'u}s and Vogelstein, Joshua T and Priebe, Carey E},
  journal={IEEE Transactions on Pattern Analysis and Machine Intelligence},
  year={2021},
  volume = {43},
  pages = {1324--1336},
  publisher={IEEE}
}

@unpublished{nielsen2018multiple,
  title={The multiple random dot product graph model},
  author={Nielsen, Agnes Martine and Witten, Daniela},
  note={arXiv preprint at \url{https://arxiv.org/abs/1811.12172}},
  year={2018}
}

@unpublished{jones2020multilayer,
  title={The multilayer random dot product graph},
  author={Jones, Andrew and Rubin-Delanchy, Patrick},
  note={arXiv preprint at \url{https://arxiv.org/abs/2007.10455}},
  year={2020}
}

@unpublished{draves2020bias,
  title={Bias-variance tradeoffs in joint spectral embeddings},
  author={Draves, Benjamin and Sussman, Daniel L},
  note={arXiv preprint at \url{https://arxiv.org/abs/2005.02511}},
  year={2020}
}

@article{paul2020spectral,
  title={Spectral and matrix factorization methods for consistent community detection in multi-layer networks},
  author={Paul, Subhadeep and Chen, Yuguo},
  journal={Annals of Statistics},
  volume={48},
  number={1},
  pages={230--250},
  year={2020},
  publisher={Institute of Mathematical Statistics}
}

@inproceedings{tang2009clustering,
  title={Clustering with multiple graphs},
  author={Tang, Wei and Lu, Zhengdong and Dhillon, Inderjit S},
  booktitle={2009 Ninth IEEE International Conference on Data Mining},
  pages={1016--1021},
  year={2009}
}

@inproceedings{kumar2011co,
  title={Co-regularized multi-view spectral clustering},
  author={Kumar, Abhishek and Rai, Piyush and Daume, Hal},
  booktitle={Advances in Neural Information Processing Systems},
  pages={1413--1421},
  year={2011}
}

@inproceedings{han2015consistent,
  title={Consistent estimation of dynamic and multi-layer block models},
  author={Han, Qiuyi and Xu, Kevin and Airoldi, Edoardo},
  booktitle={International Conference on Machine Learning},
  pages={1511--1520},
  year={2015}
}

@article{holland1983stochastic,
  title={Stochastic blockmodels: first steps},
  author={Holland, Paul W and Laskey, Kathryn Blackmond and Leinhardt, Samuel},
  journal={Social Networks},
  volume={5},
  number={2},
  pages={109--137},
  year={1983},
  publisher={Elsevier}
}

@article{paul2007,
title = { Asymptotics of sample eigenstructure for a large dimensional spiked covariance model},
author = {D.~Paul},
year = {2007},
journal = {Statistical Sinica},
volume = {17},
pages = {1617--1642}
}

@article{wang_fan_2017,
title = {Asymptotics of empirical eigenstructure for high-dimensional spiked covariance},
year = {2017},
journal = {Annals of Statistics},
author = {W.~Wang and J.~Fan},
volume = {45},
pages = {1342--1374},
}

@article{fan_li_xia_zheng,
author = {J.~Fan and Y.~Li and N.~Xia and X.~Zheng},
title = {Tests for principal eigenvalues and eigenvectors},
year = {2026},
journal = {Journal of the American Statistical Association},
}

@article{paul2016consistent,
  title={Consistent community detection in multi-relational data through restricted multi-layer stochastic blockmodel},
  author={Paul, Subhadeep and Chen, Yuguo},
  journal={Electronic Journal of Statistics},
  volume={10},
  number={2},
  pages={3807--3870},
  year={2016},
  publisher={Institute of Mathematical Statistics and Bernoulli Society}
}

@article{chen2021distributed,
  title={Distributed estimation for principal component analysis: An enlarged eigenspace analysis},
  author={Chen, Xi and Lee, Jason D and Li, He and Yang, Yun},
  journal={Journal of the American Statistical Association},
  year={2022},
  volume = {117},
  pages = {1775--1786},
}

@article{charisopoulos2021communication,
  title={Communication-efficient distributed eigenspace estimation},
  author={Charisopoulos, Vasileios and Benson, Austin R and Damle, Anil},
  journal={SIAM Journal on Mathematics of Data Science},
  volume={3},
  number={4},
  pages={1067--1092},
  year={2021},
  publisher={SIAM}
}

@article{xie2021entrywise,
  title={Entrywise limit theorems for eigenvectors of signal-plus-noise matrix models with weak signals},
  author={Xie, Fangzheng},
  journal={Bernoulli},
  year={2023+}
  
}

@article{han_fan,
author = {J.~Fan and Y.~Fan and X.~Han and J.~Lv},
title = {Asymptotic theory of eigenvectors for random matrices with diverging spikes},
journal = {Journal of the American Statistical Association},
year = {2022},
volume = {117},
pages = {996--1009},
}

@unpublished{yan2021inference,
  title={Inference for heteroskedastic PCA with missing data},
  author={Yan, Yuling and Chen, Yuxin and Fan, Jianqing},
  note={arXiv preprint at \url{https://arxiv.org/abs/2107.12365}},
  year={2021}
}

@inproceedings{garber2017communication,
  title={Communication-efficient algorithms for distributed stochastic principal component analysis},
  author={Garber, Dan and Shamir, Ohad and Srebro, Nathan},
  booktitle={International Conference on Machine Learning},
  pages={1203--1212},
  year={2017}
}

@article{battiston2017multilayer,
  title={Multilayer motif analysis of brain networks},
  author={Battiston, Federico and Nicosia, Vincenzo and Chavez, Mario and Latora, Vito},
  journal={Chaos: An Interdisciplinary Journal of Nonlinear Science},
  volume={27},
  number={4},
  pages={047404},
  year={2017},
  publisher={AIP Publishing LLC}
}

@article{de2017multilayer,
  title={Multilayer modeling and analysis of human brain networks},
  author={De Domenico, Manlio},
  journal={GigaScience},
  volume={6},
  number={5},
  pages={gix004},
  year={2017},
  publisher={Oxford University Press}
}

@unpublished{kong2021multiplex,
  title={Multiplex graph networks for multimodal brain network analysis},
  author={Kong, Zhaoming and Sun, Lichao and Peng, Hao and Zhan, Liang and Chen, Yong and He, Lifang},
  note={arXiv preprint at \url{https://arxiv.org/abs/2108.00158}},
  year={2021}
}

@article{bullmore2009complex,
  title={Complex brain networks: graph theoretical analysis of structural and functional systems},
  author={Bullmore, Ed and Sporns, Olaf},
  journal={Nature Reviews Neuroscience},
  volume={10},
  number={3},
  pages={186--198},
  year={2009},
  publisher={Nature Publishing Group}
}

@inproceedings{papalexakis2013more,
  title={Do more views of a graph help? community detection and clustering in multi-graphs},
  author={Papalexakis, Evangelos E and Akoglu, Leman and Ience, Dino},
  booktitle={Proceedings of the 16th International Conference on Information Fusion},
  pages={899--905},
  year={2013}
}

@inproceedings{greene2013producing,
  title={Producing a unified graph representation from multiple social network views},
  author={Greene, Derek and Cunningham, P{\'a}draig},
  booktitle={Proceedings of the 5th Annual ACM Web Science Conference},
  pages={118--121},
  year={2013}
}

@article{schweitzer2009economic,
  title={Economic networks: the new challenges},
  author={Schweitzer, Frank and Fagiolo, Giorgio and Sornette, Didier and Vega-Redondo, Fernando and Vespignani, Alessandro and White, Douglas R},
  journal={Science},
  volume={325},
  number={5939},
  pages={422--425},
  year={2009},
  publisher={American Association for the Advancement of Science}
}

@article{lee2016strength,
  title={Strength of weak layers in cascading failures on multiplex networks: case of the international trade network},
  author={Lee, Kyu-Min and Goh, K-I},
  journal={Scientific Reports},
  volume={6},
  number={1},
  pages={1--9},
  year={2016},
  publisher={Nature Publishing Group}
}

@article{hotelling1933analysis,
  title={Analysis of a complex of statistical variables into principal components},
  author={Hotelling, Harold},
  journal={Journal of Education Psychology},
  volume={24},
  pages={417--441},
  year={1933},
}

@article{xie2018bayesian,
  title={Bayesian estimation of sparse spiked covariance matrices in high dimensions},
  author={Xie, Fangzheng and Xu, Yanxun and Priebe, Carey E and Cape, Joshua},
  year={2022},
  journal = {Bayesian Analysis},
  pages = {1193--1217},
}

@article{candes_recht,
author = {E.~J.~Candes and B.~Recht},
title = {Exact matrix completion via convex optimization},
journal = {Foundations of Computational Mathematics},
year = {2009},
volume = {9},
pages = {717--772}
}

@article{tang2017semiparametric,
  title={A semiparametric two-sample hypothesis testing problem for random graphs},
  author={Tang, Minh and Athreya, Avanti and Sussman, Daniel L and Lyzinski, Vince and Park, Youngser and Priebe, Carey E},
  journal={Journal of Computational and Graphical Statistics},
  volume={26},
  number={2},
  pages={344--354},
  year={2017},
  publisher={Taylor \& Francis}
}

@inproceedings{levin2017central,
  title={A central limit theorem for an omnibus embedding of multiple random dot product graphs},
  author={Levin, Keith and Athreya, Avanti and Tang, Minh and Lyzinski, Vince and Priebe, Carey E},
  booktitle={2017 IEEE International Conference on Data Mining Workshops},
  pages={964--967},
  year={2017}
}

@article{ginestet2017hypothesis,
  title={Hypothesis testing for network data in functional neuroimaging},
  author={Ginestet, Cedric E and Li, Jun and Balachandran, Prakash and Rosenberg, Steven and Kolaczyk, Eric D},
  journal={Annals of Applied Statistics},
  volume={11},
  number={2},
  pages={725--750},
  year={2017},
  publisher={JSTOR}
}

@unpublished{li2018two,
  title={Two-sample test of community memberships of weighted stochastic block models},
  author={Li, Yezheng and Li, Hongzhe},
  note={arXiv preprint at \url{https://arxiv.org/abs/1811.12593}},
  year={2018}
}

@article{ghoshdastidar2020two,
  title={Two-sample hypothesis testing for inhomogeneous random graphs},
  author={Ghoshdastidar, Debarghya and Gutzeit, Maurilio and Carpentier, Alexandra and Von Luxburg, Ulrike},
  journal={Annals of Statistics},
  volume={48},
  number={4},
  pages={2208--2229},
  year={2020},
  publisher={Institute of Mathematical Statistics}
}

@unpublished{short_note_concentration,
author = {C.~Jin and P.~Netrapalli and R.~Ge and S.~M.~Kakade and M.~I.~Jordan},
title = {A short note on concentration inequalities for random vectors with sub-Gaussian norms},
note = {arXiv preprint at \url{https://arxiv.org/abs/1902.03736}},
year = {2019}
}

@article{macdonald2022latent,
  title={Latent space models for multiplex networks with shared structure},
  author={MacDonald, Peter W and Levina, Elizaveta and Zhu, Ji},
  journal={Biometrika},
  volume={109},
  number={3},
  pages={683--706},
  year={2022},
  publisher={Oxford University Press}
}

@inBook{vershynin_hdp,
author = {R.~Vershynin},
chapter = {Introduction to the non-asymptotic analysis of random matrices},
pages = {210--268},
title = {Compressed Sensing: Theory and Applications},
publisher = {Cambridge University Press},
year = {2012}
}

@article{luo_xiao_bunea,
author = {F.~Bunea and L.~Xiao},
title = {On the sample covariance matrix estimation of reduced effective rank population matrices, with applications to fPCA},
journal = {Bernoulli},
pages = {1200--1230},
year = {2015},
volume = {21}
}

@article{tropp2,
author = {J.~A.~Tropp},
title = {An Introduction to Matrix Concentration Inequalities},
journal = {Foundations and Trends in Machine Learning},
volume = {8},
pages = {1--230},
year = {2015}
}

@article{lyzinski2014perfect,
  title={Perfect clustering for stochastic blockmodel graphs via adjacency spectral embedding},
  author={Lyzinski, Vince and Sussman, Daniel L and Tang, Minh and Athreya, Avanti and Priebe, Carey E},
  journal={Electronic Journal of Statistics},
  volume={8},
  number={2},
  pages={2905--2922},
  year={2014},
  publisher={Institute of Mathematical Statistics and Bernoulli Society}
}

@article{durante2018bayesian,
  title={Bayesian inference and testing of group differences in brain networks},
  author={Durante, Daniele and Dunson, David B},
  journal={Bayesian Analysis},
  volume={13},
  number={1},
  pages={29--58},
  year={2018},
  publisher={International Society for Bayesian Analysis}
}

@article{feng2018angle,
  title={Angle-based joint and individual variation explained},
  author={Feng, Qing and Jiang, Meilei and Hannig, Jan and Marron, JS},
  journal={Journal of multivariate analysis},
  volume={166},
  pages={241--265},
  year={2018},
  publisher={Elsevier}
}

@unpublished{ponzi2021rajive,
  title={RaJIVE: robust angle based JIVE for integrating noisy multi-source data},
  author={Ponzi, Erica and Thoresen, Magne and Ghosh, Abhik},
  note={arXiv preprint at \url{https://arxiv.org/abs/2101.09110}},
  year={2021}
}

@unpublished{zhang2022perturbation,
  title={Perturbation Analysis of Randomized SVD and its Applications to High-dimensional Statistics},
  author={Zhang, Yichi and Tang, Minh},
  note={arXiv preprint at \url{https://arXiv:2203.10262}},
  year={2022}
}

@article{agterberg_dcsbm,
title = {Joint spectral clustering in multilayer degree-corrected stochastic blockmodels},
author = {J.~Agterberg and Z.~Lubberts and J.~Arroyo},
journal = {Journal of the American Statistical Asociation},
volume = {120},
pages = {1607--1620},
year = {2025}
}

@inproceedings{agterberg2022entrywise,
  title={Entrywise recovery guarantees for sparse PCA via sparsistent Algorithms},
  author={Agterberg, Joshua and Sulam, Jeremias},
  booktitle={International Conference on Artificial Intelligence and Statistics},
  pages={6591--6629},
  year={2022}
}

@article{population_svd,
author = {C.~M.~Crainiceanu and B.~S.~Caffo and S.~Luo and N.~M.~Punjabi},
title = {Population value decomposition, a framework for the analysis of image populations},
journal = {Journal of the American Statistical Association},
year = {2011},
volume = {106},
pages = {775--790}
}

@book{vershynin2018high,
  title={High-dimensional probability: An introduction with applications in data science},
  author={Vershynin, Roman},
  year={2018},
  publisher={Cambridge University Press}
}

@article{tropp2012user,
  title={User-friendly tail bounds for sums of random matrices},
  author={Tropp, Joel A},
  journal={Foundations of computational mathematics},
  volume={12},
  number={4},
  pages={389--434},
  year={2012},
  publisher={Springer}
}

@article{abbe2015exact,
  title={Exact recovery in the stochastic block model},
  author={Abbe, Emmanuel and Bandeira, Afonso S and Hall, Georgina},
  journal={IEEE Transactions on information theory},
  volume={62},
  number={1},
  pages={471--487},
  year={2015},
  publisher={IEEE}
}

@inproceedings{mossel2015consistency,
  title={Consistency thresholds for the planted bisection model},
  author={Mossel, Elchanan and Neeman, Joe and Sly, Allan},
  booktitle={Proceedings of the forty-seventh annual ACM Symposium on Theory of Computing},
  pages={69--75},
  year={2015}
}

@article{federated_survey,
author = {C.~Zhang and Y.~Xie and H.~Bai and B.~Yu and W.~Liand Y.~Gao},
title = {A survey on federated learning},
journal = {Knowledge-based systems},
volume = 216,
year = 2021,
pages = 106775
}

@article{ridge_regression_dobriban,
author = {E.~Dobriban and Y.~Sheng},
title = {WONDER: Weighted one-shot distributed ridge regression in high dimensions},
journal = {Journal of Machine Learning Research},
volume = {21},
year = {2020}
}

@article{aggregated_inference_huo,
author = {X.~Huo and S.~Cao},
title = {Aggregated inference},
journal = {WIREs Computational Statistics},
volume = {11},
year = {2019}
}

@book{anderson1962introduction,
  title={An introduction to multivariate statistical analysis},
  author={Anderson, Theodore Wilbur},
  year={2003},
edition = {3rd},
  publisher={Wiley New York}
}

@book{kollo2005advanced,
  title={Advanced multivariate statistics with matrices},
  author={Kollo, Tonu},
  year={2005},
  publisher={Springer}
}

@article{matrix_quadratic,
author = {M.~Singull and T.~Koski},
title = {On the distribution of matrix quadratic forms},
journal = {Communications in Statistics: Theory and Methods},
year = {2012},
pages = {3403--3415},
volume = {41}
}

@article{davis1977asymptotic,
  title={Asymptotic theory for principal component analysis: non-normal case},
  author={A.~W.~Davis},
  journal={Australian Journal of Statistics},
  volume={19},
  number={3},
  pages={206--212},
  year={1977},
  publisher={Wiley Online Library}
}

@article{zheng2022vertex,
  title={Vertex nomination between graphs via spectral embedding and quadratic programming},
  author={Zheng, Runbing and Lyzinski, Vince and Priebe, Carey E and Tang, Minh},
  journal={Journal of Computational and Graphical Statistics},
  volume={31},
  number={4},
  pages={1254--1268},
  year={2022},
  publisher={Taylor \& Francis}
}

@inproceedings{liang2014improved,
  title={Improved distributed principal component analysis},
  author={Liang, Yingyu and Balcan, Maria-Florina F and Kanchanapally, Vandana and Woodruff, David},
  booktitle={Proceedings of the 27th International Conference on Neural Information Processing Systems},
  year={2014},
  pages = {3113--3121},
}

@article{tang2021integrated,
  title={Integrated principal components analysis},
  author={Tang, Tiffany M and Allen, Genevera I},
  journal={Journal of Machine Learning Research},
  volume={22},
  number={1},
  pages={8953--9023},
  year={2021},
  publisher={JMLRORG}
}

@inproceedings{sagonas2017robust,
  title={Robust joint and individual variance explained},
  author={Sagonas, Christos and Panagakis, Yannis and Leidinger, Alina and Zafeiriou, Stefanos},
  booktitle={Proceedings of the IEEE Conference on Computer Vision and Pattern Recognition},
  pages={5267--5276},
  year={2017}
}

@article{cai2022non,
  title={On the non-asymptotic concentration of heteroskedastic Wishart-type matrix},
  author={Cai, T Tony and Han, Rungang and Zhang, Anru R},
  journal={Electronic Journal of Probability},
  volume={27},
  pages={1--40},
  year={2022},
  publisher={The Institute of Mathematical Statistics and the Bernoulli Society}
}

@article{huo2019aggregated,
  title={Aggregated inference},
  author={Huo, Xiaoming and Cao, Shanshan},
  journal={Wiley Interdisciplinary Reviews: Computational Statistics},
  volume={11},
  number={1},
  pages={e1451},
  year={2019},
  publisher={Wiley Online Library}
}

@article{lei2023bias,
  title={Bias-adjusted spectral clustering in multi-layer stochastic block models},
  author={Lei, Jing and Lin, Kevin Z},
  journal={Journal of the American Statistical Association},
  volume={118},
  number={544},
  pages={2433--2445},
  year={2023},
  publisher={Taylor \& Francis}
}

@article{lei2024computational,
  title={Computational and statistical thresholds in multi-layer stochastic block models},
  author={Lei, Jing and Zhang, Anru R and Zhu, Zihan},
  journal={The Annals of Statistics},
  volume={52},
  number={5},
  pages={2431--2455},
  year={2024},
  publisher={Institute of Mathematical Statistics}
}

@article{weylandt2022multivariate,
  title={Multivariate Analysis for Multiple Network Data via Semi-Symmetric Tensor PCA},
  author={Weylandt, Michael and Michailidis, George},
  journal={arXiv preprint arXiv:2202.04719},
  year={2022}
}

@article{mheich2020brain,
  title={Brain network similarity: methods and applications},
  author={Mheich, Ahmad and Wendling, Fabrice and Hassan, Mahmoud},
  journal={Network Neuroscience},
  volume={4},
  number={3},
  pages={507--527},
  year={2020},
  publisher={MIT Press One Rogers Street, Cambridge, MA 02142-1209, USA journals-info~…}
}

@article{fan2015similarity,
  title={Similarity between community structures of different online social networks and its impact on underlying community detection},
  author={Fan, Wei and Yeung, Kai-Hau},
  journal={Communications in Nonlinear Science and Numerical Simulation},
  volume={20},
  number={3},
  pages={1015--1025},
  year={2015},
  publisher={Elsevier}
}

@article{zalesky2012connectivity,
  title={Connectivity differences in brain networks},
  author={Zalesky, Andrew and Cocchi, Luca and Fornito, Alex and Murray, Micah M and Bullmore, ED},
  journal={Neuroimage},
  volume={60},
  number={2},
  pages={1055--1062},
  year={2012},
  publisher={Elsevier}
}

@article{xie2024bias,
  title={Bias-corrected joint spectral embedding for multilayer networks with invariant subspace: entrywise eigenvector perturbation and inference},
  author={Xie, Fangzheng},
  journal={IEEE Transactions on Information Theory},
  volume={70},
  number={12},
  pages={9036--9083},
  year={2024},
  publisher={IEEE}
}

@article{ma2026optimal,
  title={Optimal estimation of shared singular subspaces across multiple noisy matrices},
  author={Ma, Zhengchi and Ma, Rong},
  journal={IEEE Transactions on Information Theory},
  year={2026},
  publisher={IEEE}
}

@article{tian2026efficient,
  title={Efficient analysis of latent spaces in heterogeneous networks},
  author={Tian, Yuang and Sun, Jiajin and He, Yinqiu},
  journal={Journal of the American Statistical Association},
  pages={1--13},
  year={2026},
  publisher={Taylor \& Francis}
}

\newpage

\appendix

\begin{center}%
    {\Large Supplementary Material for ``Limit results for distributed estimation of invariant subspaces in multiple networks inference and PCA"\par}%
  \end{center}

\counterwithin{figure}{section}
\counterwithin{theorem1}{section}
\counterwithin{lemma1}{section}
\counterwithin{definition1}{section}
\counterwithin{proposition1}{section}
\counterwithin{assumption}{section}

\begin{sloppypar}

\section{Proofs of Main Results}
\label{Appendix:A}
\subsection{Proof of Theorem~\ref{thm3_general}}
\label{sec:thm3_general}
From the assumption on $\hat{\muu}^{(i)}$ we have
 \begin{equation}
  \label{eq:UU^T2_separete_2}
\begin{split}
\frac{1}{m}\sum_{i=1}^m\hat\muu^{(i)}(\hat\muu^{(i)})^\top
  	&=\frac{1}{m}\sum_{i=1}^m\big(\hat\muu^{(i)}\mw_\muu^{(i)}\big)\big(\hat\muu^{(i)}\mw^{(i)}_\muu\big)^\top
  	\\&=\frac{1}{m}\sum_{i=1}^m\big(\muu^{(i)}+\mt^{(i)}_0+\mt^{(i)}\big)\big(\muu^{(i)}+\mt^{(i)}_0+\mt^{(i)}\big)^\top
   \\ &= \frac{1}{m} \sum_{i=1}^m \muu^{(i)} \muu^{(i)\top} + \tilde{\me}
  	=\muu_c \muu_c^{\top} + \frac{1}{m} \sum_{i=1}^m \muu_s^{(i)} \muu_{s}^{(i)\top} + \tilde\me,
\end{split}
\end{equation}
where the matrix $\tilde\me$ is defined as
\begin{equation}
  \label{eq:def_e_l_thm3}
\begin{aligned}
    \tilde\me
    =&\frac{1}{m}\sum_{i=1}^m\Big[\mt_0^{(i)} \muu^{(i)\top}+
       \muu^{(i)}\mt_{0}^{(i)\top}] + \ml, \\
	\ml =&\frac{1}{m}\sum_{i=1}^m \bigl[
           \mt^{(i)} \muu^{(i)\top} +\muu^{(i)} \mt^{(i)\top} + (\mt_{0}^{(i)}+\mt^{(i)})(\mt_{0}^{(i)}+\mt^{(i)})^\top \bigr].
\end{aligned}
\end{equation}

Write the eigendecomposition for
$\tfrac{1}{m}\sum_{i=1}^m\hat\muu^{(i)}\hat\muu^{(i)\top}$ as
\begin{equation}
  \label{eq:UU^T_separate}
\begin{split}
\frac{1}{m}\sum_{i=1}^m\hat\muu^{(i)}\hat\muu^{(i)\top}
=\hat\muu_c\hat\mLambda \hat\muu_c^\top+\hat\muu_{c\perp}\hat\mLambda_{\perp} \hat\muu_{c\perp}^\top
=\muu_c\muu_c^{\top}+ \bm{\Pi}_{s} + \tilde\me,
\end{split}
\end{equation}
where $\bm{\Pi}_{s} = m^{-1} \sum_{i=1}^m \muu_s^{(i)} \muu_{s}^{(i)\top}$.
Here, $\hat{\mLambda}$ is the diagonal matrix containing the
$d_{0}$ largest eigenvalues, and $\hat{\muu}_c$ is the matrix
whose columns are the corresponding eigenvectors. The final equality follows from Eq.~\eqref{eq:UU^T2_separete_2}.
Now, as each $\muu^{(i)}$ has orthonormal columns, we have $\muu_c^{\top}
\muu^{(i)}_s = \mathbf{0}$ for all $i$ and hence $\muu_c^{\top}
\bm{\Pi}_{s} = \mathbf{0}$. In summary $\muu_c\muu_c^{\top} +
\bm{\Pi}_{s}$ has $d_{0,\muu}$ eigenvalues equal to $1$ and the
remaining eigenvalues are at most $\|\bm{\Pi}_s\|$. 
By Weyl's inequality, we have
\begin{equation}
  \begin{split}
  \label{eq:hatLambda_separate_2}
\max_{i \leq d_{0,\muu}} |\hat{\mLambda}_{ii} - 1| &\leq \|\tilde{\me}\|
 \leq  \frac{2}{m} \sum_{j=1}^{m} \|\mt_0^{(j)} + \mt^{(j)}\| +
 \frac{1}{m} \sum_{i=1}^{m} \|\mt_0^{(j)} +
 \mt^{(j)}\|^2,
    \end{split}
\end{equation}
where the last equality is from the definition of $\tilde\me$ in Eq.~\eqref{eq:epsilon_e_thm3}.
Eq.~\eqref{eq:UU^T_separate} also implies
$\hat\muu_c\hat\mLambda-(\tilde\me +
\bm{\Pi}_{s})\hat\muu_c=\muu_c\muu_c^{\top}\hat\muu_c$. 
And hence, under the conditions in Eq.~\eqref{eq:cond_pis}, the
eigenvalues of $\hat{\bm{\Lambda}}$ are disjoint from the eigevalues
of $\tilde{\me} + \bm{\Pi}_s$. Therefore $\hat{\muu}_c$ has a
von Neumann series expansion \cite{bhatia2013matrix} as
\begin{equation}
  \label{eq:von_Neumann1_seperate_2}
\hat{\muu}_c
	=\sum_{k=0}^{\infty} (\tilde\me + \bm{\Pi}_{s})^{k} \muu_c \muu_c^{\top} \hat{\muu}_c \hat{\mLambda}^{-(k+1)}.
\end{equation}
Now for any $d_{0} \times d_{0}$ orthogonal matrix $\mathbf{W}$, we define the matrices
\begin{equation}
\label{eq:Q1-Q5}
  \begin{aligned}
\mathbf{Q}_{\muu_c,1} &= \muu_c \muu_c^\top \hat{\muu}_c \hat{\mLambda}^{-1}\mw-\muu_c,
\\
\mathbf{Q}_{\muu_c,2} &=
\frac{1}{m}\sum_{i=1}^m\mt_0^{(i)} \muu^{(i)\top} \muu_c \left(\muu_c^{\top}
  \hat{\muu}_c \hat{\mLambda}^{-2}-\mw^\top\right)\mw, \\
\mathbf{Q}_{\muu_c,3} &= \frac{1}{m}\sum_{i=1}^m\muu^{(i)}                        \mt_0^{(i)\top} \muu_c\muu_c^\top
\hat{\muu}_c \hat{\mLambda}^{-2}\mw, \\
\mathbf{Q}_{\muu_c,4} &= \ml \muu_c\muu_c^\top \hat{\muu}_c
      \hat{\mLambda}^{-2}\mw, \\
\mathbf{Q}_{\muu_c,5} &= 	\sum_{k=2}^{\infty} (\tilde\me + \bm{\Pi}_{s})^{k}\muu_c\muu_c^\top \hat{\muu}_c \hat{\mLambda}^{-(k+1)}\mw.
\end{aligned}
\end{equation}
Notice $\bm{\Pi}_{s} \mathbf{U}_c = \bm{0}$, and recall the definition of $\tilde{\me}$ and $\ml$ in Eq.~\eqref{eq:def_e_l_thm3}. Then by the expansion of $\hat\muu_c$ in Eq.~\eqref{eq:von_Neumann1_seperate_2}
 we have
\begin{equation}
  \label{eq:series of hatU-UW_seperate_2}
\begin{split}
	\hat{\muu}_c\mw-\muu_c
	&=\mathbf{Q}_{\muu_c,1} +\tilde\me \muu_c \muu_c^\top \hat{\muu}_c \hat{\mLambda}^{-2}\mw
	+\sum_{k=2}^{\infty} (\tilde\me + \bm{\Pi}_{s})^{k} \muu_c \muu_c^\top \hat{\muu}_c \hat{\mLambda}^{-(k+1)}\mw\\
&=\frac{1}{m}\sum_{i=1}^m \mt_{0}^{(i)} \muu^{(i)\top} \muu_c + \mq_{\muu_c},
\end{split}
\end{equation}
where we let $\mq_{\muu_c} = \mq_{\muu_c,1} + \mq_{\muu_c,2} + \dots + \mq_{\muu_c,5}$. 

Let $\mw_{\muu_c}$ denote the minimizer of $\|\hat{\muu}_c^{\top} \mo - \muu_c\|_{F}$ over all $d_{0} \times d_{0}$
orthogonal matrices $\mo$. 
We now bound $\mq_{\muu_c,1}$ through $\mq_{\muu_c,5}$ for this choice of $\mw =
\mw_{\muu_c}$. 
We first define the quantities associated with
$\tilde{\me}$ and $\ml$
\begin{equation*}
	\begin{aligned}
		\epsilon_{\ml}=\|\ml\|,
		\quad \zeta_{\ml}=\|\ml\|_{2\to\infty},
		\quad \epsilon_{\tilde\me}=\|\tilde\me\|,
		\quad \zeta_{\tilde\me}=\|\tilde\me\|_{2\to\infty}.
	\end{aligned}
\end{equation*}
Under the condition in Eq.~\eqref{eq:cond_pis} we have $\zeta_{\mt_0}\leq\epsilon_{\mt_0}<1,\zeta_{\mt}\leq\epsilon_{\mt}<1$. Then we have
\begin{equation}
\begin{aligned}
	\label{eq:epsilon_e_thm3}
  &\epsilon_{\ml} 
  \leq 2 \epsilon_{\mt} + (\epsilon_{\mt_0}+\epsilon_{\mt})^2
  \lesssim \epsilon_{\mt_0}^2+\epsilon_{\mt},\\
  &
  \epsilon_{\tilde{\me}} 
  \leq 2\epsilon_{\mt_0}+ \epsilon_{\ml}
  \lesssim \epsilon_{\mt_0}+\epsilon_{\mt}, 
  \\& \zeta_{\ml} 
  \leq \zeta_{\mt} +\zeta_{\muu} \epsilon_{\mt}  +
  (\zeta_{\mt_0} + \zeta_{\mt})  (\epsilon_{\mt}+\epsilon_{\mt})
  \lesssim \zeta_{\muu} \epsilon_{\mt}+\zeta_{\mt_0} (\epsilon_{\mt_0}+\epsilon_{\mt}) + \zeta_{\mt}
  ,\\&
  \zeta_{\tilde{\me}} 
  \leq \zeta_{\mt_0}+\zeta_{\muu}\epsilon_{\mt_0}  + \zeta_{\ml}
  \lesssim \zeta_{\muu} (\epsilon_{\mt_0}+\epsilon_{\mt}) + \zeta_{\mt_0}+ \zeta_{\mt}.
\end{aligned}
\end{equation}

{\bf Bounding $\mq_{\muu_c,1}$:}
Let $\delta_{c} = \min_{i > d_{0,\muu}} |1 -
\hat{\lambda}_{i}|$ where $\hat{\lambda}_{i}$ for $i > d_{0,\muu}$
are the eigenvalues in $\hat{\bm{\Lambda}}_{\perp}$ in
Eq.~\eqref{eq:UU^T_separate}. By similar reasoning to that for
Eq.~\eqref{eq:hatLambda_separate_2}, we have
$\delta_c \geq 1 - \|\bm{\Pi}_s\| - \epsilon_{\tilde\me}$. 
Now, by the general form of the Davis-Kahan Theorem (see Theorem~VII.3 in \cite{bhatia2013matrix}) we have
\begin{equation*}
  \begin{split}
    \bigl\|\sin \Theta(\hat{\muu}_c, \muu_c) \bigr\|  &= \|(\mi - \hat{\muu}_c \hat{\muu}_c^{\top}) \muu_c \muu_c^{\top}\| 
    \\&\leq\frac{\|(\mi - \hat{\muu}_c \hat{\muu}_c^{\top})(\tilde\me + \bm{\Pi}_{s}) \muu_c \muu_c^{\top}\|}{\mathrm{\delta}_c}
     \leq \frac{\|\tilde{\me}\|}{\delta_c} \leq
     \frac{\epsilon_{\tilde{\me}}}{1 - \epsilon_{\tilde{\me}} - \|\bm{\Pi}_s\|}.
   \end{split}
\end{equation*}
And hence
\begin{equation}
	\label{eq:|(I-UcUc^)hat Uc|}
	\begin{aligned}
		\|(\mi-\muu_c\muu_c^\top)\hat\muu_c\|
		\leq \sqrt{2}\bigl\|\sin \Theta(\hat{\muu}_c, \muu_c) \bigr\|
		\leq \frac{2^{1/2}\epsilon_{\tilde\me}}{1 - \|\bm{\Pi}_s\| - \epsilon_{\tilde\me}}.
	\end{aligned}
\end{equation}
As $\mw_{\muu_c}$ is the solution of orthogonal Procrustes problem, we have
\begin{equation}
  \label{eq:|U^T hatU-W|_F_seperate}
\begin{split}
\begin{aligned}
	\|\muu_c^\top\hat\muu_c-\mw_{\muu_c}^\top\|
	&\leq 1-\sigma^2_{\min}(\muu_c^\top\hat\muu_c)
	 \leq \bigl\|\sin \Theta(\hat{\muu}_c, \muu_c) \bigr\|^2
	\leq \frac{\epsilon_{\tilde\me}^2}{(1 -
       \|\bm{\Pi}_s\| - \epsilon_{\tilde\me})^2}.
\end{aligned}
\end{split}
\end{equation}
Rewrite $\mq_{\muu,1}$ as
\begin{equation}
  \label{eq:q_u1_decomp_thm3}
\begin{aligned}
	\mq_{\muu_c,1}
	&=\muu_c\big(\muu_c^{\top} \hat{\muu}_c \hat{\mLambda}^{-1}-\mw_{\muu_c}^\top\big)\mw_{\muu_c}\\
	&=-\muu_c\big(\muu_c^\top \hat{\muu}_c^\top\hat{\mLambda}-\muu_c^\top \hat{\muu}_c^\top\big)\hat{\mLambda}^{-1}\mw_{\muu_c}
	+\muu_c\big(\muu_c^\top \hat{\muu}_c-\mw_{\muu_c}^\top\big)\mw_{\muu_c}
	\\ &=-\muu_c \muu_c^{\top} (\tilde{\me} + \bm{\Pi}_s) \hat{\muu}_c \hat{\mLambda}^{-1}\mw_{\muu_c}
	+\muu_c\big(\muu_c^\top \hat{\muu}_c-\mw_{\muu_c}^\top\big)\mw_{\muu_c} \\
 &= - \muu_c \muu_c^{\top} \tilde{\me} \hat{\muu}_c \hat{\mLambda}^{-1}\mw_{\muu_c}
	+\muu_c\big(\muu_c^\top
   \hat{\muu}_c-\mw_{\muu_c}^\top\big)\mw_{\muu_c} \\
  &= - \muu_c \muu_c^{\top} \bigl(\tilde{\me} \muu_c \muu_c^{\top} \hat{\muu}_c + 
    \tilde{\me} (\mi - \muu_c \muu_c^{\top}) \hat{\muu}_c \bigr) \hat{\mLambda}^{-1}\mw_{\muu_c} + \muu_c\big(\muu_c^\top
   \hat{\muu}_c-\mw_{\muu_c}^\top\big)\mw_{\muu_c}.
\end{aligned}
\end{equation}
Recalling the expression for $\tilde{\me}$ in Eq.~\eqref{eq:def_e_l_thm3}, we have
\begin{equation}
  \label{eq:5}
 \begin{aligned}
	\muu_c^{\top} \tilde\me \muu_c	&=\frac{1}{m}\sum_{i=1}^m \muu_c^\top [ \muu^{(i)}
      \mt_0^{(i)\top} + \mt_0^{(i)} \muu^{(i)\top}] \muu_c +
      \muu_c^{\top} \ml \muu_c  
\end{aligned}
\end{equation}
and hence
\begin{equation}
  \label{eq:|UcTtildeEUc|}
 \begin{aligned}
\|\muu_c^{\top} \tilde{\me} \muu_c\| \leq \frac{2}{m}
\sum_{i=1}^{m} \|\muu_c^{\top} \mt_0^{(i)}\| + \|\ml\| \leq 2
\epsilon_{\star} + \epsilon_{\ml}.
\end{aligned}
\end{equation}
Plugging Eq.~\eqref{eq:hatLambda_separate_2}, Eq.~\eqref{eq:|(I-UcUc^)hat Uc|}, Eq.~\eqref{eq:|U^T hatU-W|_F_seperate} and Eq.\eqref{eq:|UcTtildeEUc|}
into Eq.~\eqref{eq:q_u1_decomp_thm3} yields
\begin{gather*}
\begin{aligned}
	\|\mq_{\muu_c,1}\| 
	 &\leq (\|\muu_c^{\top} \tilde{\me} \muu_c\| + \|\tilde{\me}\|
       \cdot \|(\mi - \muu_c \muu_c^{\top}) \hat{\muu}_c\|)\cdot
\|\hat{\bm{\Lambda}}^{-1}\| + \|\muu_c^\top
\hat{\muu}_c-\mw_{\muu_c}^\top\| 
  \\ &\leq \frac{2 \epsilon_{\star} + \epsilon_{\ml}}{1 -
       \epsilon_{\tilde{\me}}} 
       +\frac{2^{1/2}\epsilon_{\tilde\me}^2}{(1 - \|\bm{\Pi}_s\| - \epsilon_{\tilde\me})(1 -
       \epsilon_{\tilde{\me}})}
       + 
       \frac{\epsilon_{\tilde{\me}}^2}{(1 - \|\bm{\Pi}_s\| - \epsilon_{\tilde{\me}})^2},
       \\
	\|\mq_{\muu_c,1}\|_{2 \to \infty}  &\leq \|\muu_c\|_{2\to \infty}\big[(\|\muu_c^{\top} \tilde{\me} \muu_c\| + \|\tilde{\me}\|
       \cdot \|(\mi - \bm{\Pi}_c) \hat{\muu}_c\|)\cdot
\|\hat{\bm{\Lambda}}^{-1}\| + \|\muu_c^\top
\hat{\muu}_c-\mw_{\muu_c}^\top\|\big]  
	 \\ & \leq \zeta_{\muu} \Bigl( \frac{2 \epsilon_{\star} + \epsilon_{\ml}}{1 -
       \epsilon_{\tilde{\me}}} 
       +\frac{2^{1/2}\epsilon_{\tilde\me}^2}{(1 - \|\bm{\Pi}_s\| - \epsilon_{\tilde\me})(1 -
       \epsilon_{\tilde{\me}})}
       + 
       \frac{\epsilon_{\tilde{\me}}^2}{(1 - \|\bm{\Pi}_s\| - \epsilon_{\tilde{\me}})^2}\Bigr).
\end{aligned}
\end{gather*}

{\bf Bounding $\mq_{\muu_c,2}$:}
We first have
\begin{equation*}
\begin{aligned}
	\muu^{\top}_c \hat{\muu}_c \hat{\bm{\Lambda}}^{-2} -
\mw^\top_{\muu_c} 
&= (\muu_c^{\top} \hat{\muu}_c-\muu^{\top}_c \hat{\muu}_c
\hat{\bm{\Lambda}}^{2} )\hat{\bm{\Lambda}}^{-2} + (\muu_c^{\top}
\hat{\muu}_c - \mw^\top_{\muu_c})\\  
 & =\big[\muu^{\top}_c \hat{\muu}_c-\muu^{\top}_c (\muu_c \muu_c^{\top} + \tilde{\me} + \bm{\Pi}_{s})^2 \hat{\muu}_c \big]
    \hat{\bm{\Lambda}}^{-2}
    + (\muu_c^{\top}\hat{\muu}_c - \mw^\top_{\muu_c})\\
& =-\muu_c^{\top}\big(\tilde\me+ \tilde\me\muu_c\muu_c^\top+\tilde\me^2 + \tilde{\me} \bm{\Pi}_{s}\big)\hat{\muu}_c
    \hat{\bm{\Lambda}}^{-2}
    + (\muu_c^{\top}\hat{\muu}_c - \mw^\top_{\muu_c}),
\end{aligned}
\end{equation*}
where the final equality follows from the fact that $\bm{\Pi}_{s}
\muu_c = \mathbf{0}$. 
By
Eq.~\eqref{eq:hatLambda_separate_2} and
Eq.~\eqref{eq:|U^T hatU-W|_F_seperate},
we have
\begin{equation}
  \label{eq:Q2_part1_seperate}
\begin{aligned}
	\|\muu^{\top}_c \hat{\muu}_c \hat{\bm{\Lambda}}^{-2} - \mw^\top_{\muu_c}\|
  & \leq 
    \frac{4\epsilon_{\tilde{\me}}}{(1-\epsilon_{\tilde{\me}})^2}
     +
    \frac{\epsilon_{\tilde{\me}}^2}{(1 - \|\bm{\Pi}_s\| - \epsilon_{\tilde{\me}})^2}.
\end{aligned}
\end{equation}
Eq.~\eqref{eq:Q2_part1_seperate} then implies
$$
\begin{aligned}
	&\|\mq_{\muu_c,2}\| \leq \frac{1}{m}\sum_{i=1}^m
      \|\mt_{0}^{(i)}\| \cdot \|\muu_c^{\top}\hat{\muu}_c
    \hat{\bm{\Lambda}}^{-2} - \mw^\top_{\muu_c}\| 
      \leq  
      \epsilon_{\mt_0}\Bigl(\frac{4\epsilon_{\tilde{\me}}}{(1-\epsilon_{\tilde{\me}})^2}
     +
    \frac{\epsilon_{\tilde{\me}}^2}{(1 - \|\bm{\Pi}_s\| - \epsilon_{\tilde{\me}})^2}\Bigr)\\
    &\|\mq_{\muu_c,2}\|_{2 \to \infty} \leq
  \frac{1}{m} \sum_{i=1}^m \|\mt_0^{(i)}\|_{2 \to \infty} \cdot
 \|\muu_c^{\top}\hat{\muu}_c
      \hat{\bm{\Lambda}}^{-2} - \mw_{\muu_c}^\top\| 
      \leq \zeta_{\mt_0} \Bigl(\frac{4\epsilon_{\tilde{\me}}}{(1-\epsilon_{\tilde{\me}})^2}
     +
    \frac{\epsilon_{\tilde{\me}}^2}{(1 - \|\bm{\Pi}_s\| - \epsilon_{\tilde{\me}})^2}\Bigr).
    \end{aligned}
$$

{\bf Bounding $\mq_{\muu_c,3}$ and $\mq_{\muu_c,4}$:}  
By Eq.~\eqref{eq:hatLambda_separate_2}, these terms can be
controlled using
$$
\begin{aligned}
	&\|\mq_{\muu_c,3}\|
	\leq \frac{1}{m}\sum_{i=1}^m \|\mt_0^{(i) \top} \muu_c\| \cdot
      \|\hat{\bm{\Lambda}}^{-1}\|^2 
      \leq \frac{\epsilon_{\star}}{(1 - \epsilon_{\tilde{\me}})^2}, \\
    &\|\mq_{\muu_c,3}\|_{2\to\infty} 
    \leq
      \frac{1}{m}\sum_{i=1}^m\|\muu^{(i)}\|_{2\to\infty}\cdot
      \|\mt_0^{(i)\top} \muu_c\| \cdot\| \hat{\mLambda}^{-1}\|^2
	\leq \frac{\zeta_{\muu}\epsilon_{\star}}{(1 - \epsilon_{\tilde{\me}})^2}, \\
 &\|\mq_{\muu_c,4}\| \leq \|\ml\|\cdot\|\hat{\bm{\Lambda}}^{-1}\|^2
   \leq \frac{\epsilon_{\ml}}{(1 - \epsilon_{\tilde{\me}})^2}, \\
    &\|\mq_{\muu_c,4}\|_{2\to\infty} \leq \|\ml\|_{2\to\infty}\cdot
\|\hat{\bm{\Lambda}}^{-1}\|^2 \leq \frac{\zeta_{\ml}}{(1 - \epsilon_{\tilde{\me}})^2}.
\end{aligned}
$$

{\bf Bounding $\mq_{\muu_c,5}$:} 
First note that, as $\bm{\Pi}_{s} \muu_c =
\bm{0}$, we have $$(\tilde\me + \bm{\Pi}_{s})^{2} \muu_c 
= (\tilde\me + \bm{\Pi}_{s})\tilde\me \muu_c
=\tilde\me^2 \muu_c+ \bm{\Pi}_{s}\tilde\me \muu_c,$$
and thus $$\|(\tilde\me + \bm{\Pi}_{s})^{2} \muu_c \|\leq \epsilon_{\tilde{\me}}^2+\|\bm{\Pi}_{s}\tilde\me \muu_c\|.$$
Then for any $k \geq 2$ we have $$\|(\tilde\me + \bm{\Pi}_{s})^{k} \muu_c \|
\leq\|(\tilde\me + \bm{\Pi}_{s})^{k-2}\|\cdot\| (\tilde\me + \bm{\Pi}_{s})^2\muu_c \|
\leq (\epsilon_{\tilde{\me}}+\| \bm{\Pi}_{s}\|)^{k-2}[\epsilon_{\tilde{\me}}^{2}+\|\bm{\Pi}_{s}\tilde\me \muu_c\|].$$

Let $\hat{\lambda}^{-1} =
\|\hat{\bm{\Lambda}}^{-1}\|$. We then have
\begin{equation}
  \label{eq:thm3_mq_5_part1_2}
    \begin{split}
    \|\mq_{\muu_c,5}\| 
    &\leq \sum_{k=2}^{\infty} 
                       \hat{\lambda}^{-(k+1)} (\epsilon_{\tilde{\me}}+\| \bm{\Pi}_{s}\|)^{k-2}[\epsilon_{\tilde{\me}}^{2}+\|\bm{\Pi}_{s}\tilde\me \muu_c\|] \\
    &\leq [\epsilon_{\tilde{\me}}^{2}+\|\bm{\Pi}_{s}\tilde\me \muu_c\|]\hat{\lambda}^{-3}\sum_{\ell=0}^{\infty} 
                       \hat{\lambda}^{-\ell} (\epsilon_{\tilde{\me}}+\| \bm{\Pi}_{s}\|)^{\ell} \\
    &\leq [\epsilon_{\tilde{\me}}^{2}+\|\bm{\Pi}_{s}\tilde\me \muu_c\|]\hat{\lambda}^{-3}
    \frac{1}{1-\hat{\lambda}^{-1}(\epsilon_{\tilde{\me}}+\| \bm{\Pi}_{s}\|)}.
    \end{split}
\end{equation}
Notice that under the conditions in Eq.~\eqref{eq:cond_pis} we have $\hat{\lambda}^{-1}(\epsilon_{\tilde{\me}}+\| \bm{\Pi}_{s}\|)<c'$ for some constant $c'<1$.
Recalling the definition of
$\bm{\Pi}_{s} = m^{-1} \sum_{i=1}^m \muu^{(i)}_s
\muu^{(i)\top}_s$, and following the similar argument for Eq.~\eqref{eq:5}, we have
\begin{equation}
  \label{eq:pi_stildee_uc_2}
\begin{split}
\|\bm{\Pi}_{s} \tilde{\me} \muu_c\| 
                                     \leq \frac{1}{m^2}
  \sum_{i=1}^{m} \sum_{j=1}^{m} \|\muu_s^{(i)\top} \mt_0^{(j)}\| +
  \frac{1}{m}  \sum_{i=1}^{m}
  \|\muu_c^{\top} \mt_0^{(i)}\| + \|\ml\| 
  \leq 2\epsilon_{\star} + \epsilon_{\ml}.
\end{split}
\end{equation}
Substituting Eq.~\eqref{eq:pi_stildee_uc_2} into Eq.~\eqref{eq:thm3_mq_5_part1_2},
and then using Eq.~\eqref{eq:hatLambda_separate_2} to bound $\hat{\lambda}^{-1}$
we obtain
$$
\|\mq_{\muu_c,5}\| 
\leq \frac{\epsilon_{\tilde{\me}}^{2}+2\epsilon_{\star} + \epsilon_{\ml}}{(1 - \epsilon_{\tilde{\me}})^3[1-\hat{\lambda}^{-1}(\epsilon_{\tilde{\me}}+\| \bm{\Pi}_{s}\|)]}.
$$

For $\|\mq_{\muu_c,5}\|_{2\to\infty}$,
we note that
$$
\begin{aligned}
	\|(\tilde\me + \bm{\Pi}_{s})^{2} \muu_c\|_{2\to\infty} 
=\|\tilde\me^2 \muu_c+ \bm{\Pi}_{s}\tilde\me \muu_c\|_{2\to\infty}
\leq \zeta_{\tilde\me}\epsilon_{\tilde\me}+\zeta_{\muu}\frac{1}{m}
                                           \sum_{i=1}^{m}
                                           \|\muu_{s}^{(i)\top}
                                      \tilde{\me} \muu_c\|
\leq \zeta_{\tilde\me}\epsilon_{\tilde\me}+\zeta_{\muu}(2\epsilon_{\star} + \epsilon_{\ml}),
\end{aligned}
$$
and for any $k > 2$ we have 
$$
\begin{aligned}
	\|(\tilde\me + \bm{\Pi}_{s})^{k} \muu_c \|_{2\to\infty}
&\leq\|\tilde\me + \bm{\Pi}_{s}\|_{2\to\infty}\cdot\|(\tilde\me + \bm{\Pi}_{s})^{k-3}\|\cdot\| (\tilde\me + \bm{\Pi}_{s})^2\muu_c \|\\
&\leq (\zeta_{\tilde\me}+\zeta_{\muu})(\epsilon_{\tilde{\me}}+\| \bm{\Pi}_{s}\|)^{k-3}(\epsilon_{\tilde{\me}}^{2}+2\epsilon_{\star} + \epsilon_{\ml}).
\end{aligned}
$$
Then using the same reasoning as
that for Eq.~\eqref{eq:thm3_mq_5_part1_2}, we have
\begin{equation}
  \label{eq:thm3_mq_5_part2_2}
    \begin{split}
    \|\mq_{\muu_c,5}\|_{2 \to \infty}
    &\leq 
    \hat{\lambda}^{-3}\|(\tilde\me + \bm{\Pi}_{s})^{2} \muu_c\|_{2\to\infty} 
    +
    \sum_{k=3}^{\infty} 
                       \hat{\lambda}^{-(k+1)} \|(\tilde\me + \bm{\Pi}_{s})^{k} \muu_c\|_{2\to\infty} \\
    &\leq 
    \hat{\lambda}^{-3}[\zeta_{\tilde\me}\epsilon_{\tilde\me}+\zeta_{\muu}(2\epsilon_{\star} + \epsilon_{\ml})] 
    +
    \sum_{k=3}^{\infty} 
                       \hat{\lambda}^{-(k+1)} (\zeta_{\tilde\me}+\zeta_{\muu})(\epsilon_{\tilde{\me}}+\| \bm{\Pi}_{s}\|)^{k-3}(\epsilon_{\tilde{\me}}^{2}+2\epsilon_{\star} + \epsilon_{\ml}) \\
  &\leq \hat{\lambda}^{-3}[\zeta_{\tilde\me}\epsilon_{\tilde\me}+\zeta_{\muu}(2\epsilon_{\star} + \epsilon_{\ml})]
    +\hat{\lambda}^{-4}(\zeta_{\tilde\me}+\zeta_{\muu})
    (\epsilon_{\tilde{\me}}^{2}+2\epsilon_{\star} + \epsilon_{\ml})
    \sum_{\ell=0}^{\infty} 
                       \hat{\lambda}^{-\ell} (\epsilon_{\tilde{\me}}+\| \bm{\Pi}_{s}\|)^{\ell}
                       \\
    &\leq  \frac{\zeta_{\tilde\me}\epsilon_{\tilde\me}+\zeta_{\muu}(2\epsilon_{\star} + \epsilon_{\ml})}{(1-\epsilon_{\tilde\me})^3}
   + \frac{(\zeta_{\tilde\me}+\zeta_{\muu})
    (\epsilon_{\tilde{\me}}^{2}+2\epsilon_{\star} + \epsilon_{\ml})}{(1-\epsilon_{\tilde\me})^4[1-\hat{\lambda}^{-1}(\epsilon_{\tilde{\me}}+\| \bm{\Pi}_{s}\|)]}.
    \end{split}
    \end{equation}

We now combine the above bounds for $\mq_{\muu_c,1}$ through $\mq_{\muu_c,5}$. 
Notice under the conditions in Eq.~\eqref{eq:cond_pis} we have $(1-\epsilon_{\tilde\me})\gtrsim 1, (1 - \|\bm{\Pi}_s\| - \epsilon_{\tilde\me})\gtrsim 1,
[1-\hat{\lambda}^{-1}(\epsilon_{\tilde{\me}}+\| \bm{\Pi}_{s}\|)]\gtrsim 1,\epsilon_{\tilde{\me}}\lesssim 1, \epsilon_{\mt_0}\lesssim 1,\epsilon_{\mt}\lesssim 1$.
And recall the bounds in Eq.~\eqref{eq:epsilon_e_thm3}. We then have
\begin{equation}
	\label{eq:mq_muuc_bounds_2}
\begin{aligned}
	\|\mq_{\muu_c}\| 
	&\leq  \sum_{k=1}^{5} \|\mq_{\muu_c,k}\|
       \lesssim \epsilon_{\star}+\epsilon_\ml+\epsilon_{\tilde\me}(\epsilon_{\tilde\me}+\epsilon_{\mt_0})
       \lesssim \epsilon_{\star}+\epsilon_{\mt_0}^2+\epsilon_{\mt},\\
       \|\mq_{\muu_c}\|_{2\to\infty}
       &\leq \sum_{k=1}^{5} \|\mq_{\muu_c,k}\|_{2\to\infty}
       \lesssim \zeta_{\muu}(\epsilon_{\star}+\epsilon_\ml+\epsilon_{\tilde\me}^2)
       +\zeta_{\mt_0}\epsilon_{\tilde\me}
       +\zeta_\ml
       +\zeta_{\tilde\me}(\epsilon_{\star}+\epsilon_\ml+\epsilon_{\tilde\me})\\
       &\lesssim \zeta_\muu(\epsilon_{\star}+\epsilon_{\mt_0}^2+\epsilon_{\mt})
       +\zeta_{\mt_0}(\epsilon_{\star}+\epsilon_{\mt_0}+\epsilon_{\mt})
       +\zeta_{\mt}.
       \end{aligned}
\end{equation}
The expansion for $\hat\muu_c$ in Theorem~\ref{thm3_general} follows directly.

{\bf Estimation of $\{\muu_s^{(i)}\}$:}

We estimate $\muu_s^{(i)}$ using the $d_i - d_{0}$ leading left singular vectors of
$(\mi - \hat{\muu}_c \hat{\muu}_c^{\top}) \hat{\muu}^{(i)}$. Let $\bm{\Pi}_c = \muu_c \muu_c^{\top}$ and $\widehat{\bm{\Pi}}_c =
\hat{\muu}_c \hat{\muu}_c^{\top}$. 
 Let $\mathbf{M}^{(i)} = (\widehat{\bm{\Pi}}_c - \bm{\Pi}_c)
 \hat{\muu}^{(i)}\mw_{\muu}^{(i)}$. From the previous expansion for
 $\hat{\muu}_c \mw_c - \muu_c$, 
 we have
 \begin{equation*}
   \begin{split}
  \widehat{\bm{\Pi}}_c - \bm{\Pi}_c
  &= \hat{\muu}_c \mw_{\muu_c}
 \mw_{\muu_c}^{\top} \hat{\muu}_c^\top - \muu_c \muu_c^{\top} \\
   &= (\hat{\muu}_c \mw_{\muu_c}-\muu_c)
 (\hat{\muu}_c \mw_{\muu_c}-\muu_c)^\top
 +(\hat{\muu}_c \mw_{\muu_c}-\muu_c )\muu_c^\top
 +\muu_c
 (\hat{\muu}_c \mw_{\muu_c}-\muu_c)^\top\\
  &= \frac{1}{m}
 \sum_{j=1}^m \mt_0^{(j)} \muu^{(j)\top} \muu_c 
 \muu_c^{\top} 
+ \frac{1}{m} \sum_{j=1}^{m} \muu_c \muu_c^{\top}
 \muu^{(j)} \mt_0^{(j)\top} +
 \tilde{\mq}_{\muu},
   \end{split}
 \end{equation*}
 where $\tilde{\mq}_{\muu}$ satisfies the same bound as $\mq_{\muu_c}$
 given in Eq.~\eqref{eq:mq_muuc_bounds_2}. Now define
 $$
 {\mathcal{L}} = \frac{1}{m} \sum_{j=1}^m \mt_0^{(j)}
 \muu^{(j)\top} \muu_c \muu_c^{\top}.$$
 The above expansion for $\widehat{\bm{\Pi}}_c - \bm{\Pi}_c$ can then
 be written as
 $$\widehat{\bm{\Pi}}_c - \bm{\Pi}_c = {\mathcal{L}} +
 {\mathcal{L}}^{\top} + \tilde{\mq}_{\muu},$$
and we have
 \begin{equation}
   \begin{split}
 \label{eq:M_0}
 	\mm^{(i)}&={\mathcal{L}}\hat\muu^{(i)}\mw_{\muu}^{(i)} +
 {\mathcal{L}}^{\top}\hat\muu^{(i)}\mw_{\muu}^{(i)} + \tilde{\mq}_{\muu}\hat\muu^{(i)}\mw_{\muu}^{(i)}.
     \end{split}
 \end{equation}
 We now analyze the terms on the right hand of Eq.~\eqref{eq:M_0}.
 For ${\mathcal{L}}\hat\muu^{(i)}\mw_{\muu}^{(i)}$, recalling the expansion for $\hat{\muu}^{(i)}$, by the assumption about $\hat{\muu}^{(i)}$ we have
 \begin{equation*}
   \begin{split}
   {\mathcal{L}}\hat{\muu}^{(i)}\mw_{\muu}^{(i)} &=
  {\mathcal{L}}\bigl(\muu^{(i)} + \mt_0^{(i)}  + \mt^{(i)}\bigr)  
    = {\mathcal{L}} \muu^{(i)} + \mathbf{T}_{\mathcal{L},1}^{(i)},
   \end{split}
 \end{equation*}
 where $\mathbf{T}_{\mathcal{L},1}^{(i)}$ is a $n_l \times d_i$ residual
 matrix satisfying
 \begin{equation*}
   \|\mathbf{T}_{\mathcal{L},1}^{(i)}\| \lesssim \epsilon_{\mt_0}(\epsilon_{\mt_0} + \epsilon_{\mt}), \quad
   \|\mathbf{T}_{\mathcal{L},1}^{(i)}\|_{2 \to \infty} \lesssim
   \zeta_{\mt_0}(\epsilon_{\mt_0} + \epsilon_{\mt}).
 \end{equation*}
 Similarly for ${\mathcal{L}}^{\top}\hat\muu^{(i)}\mw_{\muu}^{(i)}$, we have
 \begin{equation*}
   \begin{split}
  {\mathcal{L}}^{\top} \hat{\muu}^{(i)}\mw_{\muu}^{(i)} &=
 {\mathcal{L}}^{\top}  \muu^{(i)} +
 \mathbf{T}_{\mathcal{L},2}^{(i)},
   \end{split}
 \end{equation*}
 where $\mathbf{T}_{\mathcal{L},2}^{(i)}$ also satisfies 
 $$
 \|\mathbf{T}_{\mathcal{L},2}^{(i)}\| \lesssim \epsilon_{\mt_0}
 (\epsilon_{\mt_0} + \epsilon_{\mt}) \quad
   \|\mathbf{T}_{\mathcal{L},2}^{(i)}\|_{2 \to \infty} \lesssim
   \zeta_{\muu} \epsilon_{\mt_0} (\epsilon_{\mt_0} + \epsilon_{\mt}),$$
and for ${\mathcal{L}}^{\top}  \muu^{(i)}$, because 
 $${\mathcal{L}}^{\top}  \muu^{(i)}=\frac{1}{m} \sum_{j=1}^{m} \muu_c \muu_c^{\top} \muu^{(j)}
 (\mt_0^{(j)\top} \muu^{(i)}) ,
 $$
 we have
 \begin{gather*}
  \bigl\|{\mathcal{L}}^{\top}  \muu^{(i)}\bigr\|
 \lesssim \epsilon_{\star}, \quad 
   \bigl\|{\mathcal{L}}^{\top}  \muu^{(i)}\bigr\|_{2 \to
   \infty} \lesssim \zeta_{\muu} \epsilon_{\star}.
 \end{gather*}
 For $\tilde{\mq}_{\muu} \hat{\muu}^{(i)}\mw_{\muu}^{(i)}$, its bounds are the same as
 those for $\tilde{\mq}_{\muu}$.
 Combining the above results for terms in Eq.~\eqref{eq:M_0},
  we conclude that
 \begin{equation}
   \label{eq:m0form}
 \mm^{(i)}  =
{\mathcal{L}}  \muu^{(i)}  +
 \mt_{\mathcal{L}}^{(i)},
 \end{equation}
 where $\mt_{\mathcal{L}}^{(i)} = \mt_{\mathcal{L},1}^{(i)} + {\mathcal{L}}^{\top}  \muu^{(i)}+\mt_{\mathcal{L},2}^{(i)}
 + \tilde{\mq}_{\muu} \hat{\muu}^{(i)}\mw_{\muu}^{(i)}$ satisfies the same bounds as
 that for $\mq_{\muu_c}$.
 
From the expansion for $\hat{\muu}^{(i)}$ we have
\begin{equation*}
  \begin{split} (\mi - \hat{\bm{\Pi}}_c) \hat{\muu}^{(i)}\mw_{\muu}^{(i)}  
   &= (\mi - \bm{\Pi}_c) \big(\muu^{(i)} + \mt_0^{(i)} + \mt^{(i)}\big) + \mathbf{M}^{(i)}
\\ &= \bigl[\begin{matrix} \bm{0} \mid \muu_s^{(i)} \end{matrix}
    \bigr] + (\mi - \bm{\Pi}_c) \mt_0^{(i)} + (\mi - \bm{\Pi}_c) \mt^{(i)} + \mathbf{M}^{(i)}\\
     &=\bigl[\begin{matrix} \bm{0} \mid \muu_s^{(i)} \end{matrix}
    \bigr]+\mt_{0,s}^{(i)}+\mt_{s}^{(i)}+\mm^{(i)},
  \end{split}
\end{equation*}
where we define $\mt_{0,s}^{(i)} =  (\mi - \bm{\Pi}_c) \mt_0^{(i)}$ and $\mt_s^{(i)} = (\mi - \bm{\Pi}_c) 
\mt^{(i)}$. 
Now $\hat{\muu}_s^{(i)}$ is the leading left singular vectors of $(\mi
- \hat{\bm{\Pi}}_c) \hat{\muu}^{(i)}$ and is thus the
leading eigenvectors of
\begin{equation*}
  \begin{split}
(\mi - \hat{\bm{\Pi}}_c) \hat{\muu}^{(i)}\hat{\muu}^{(i)\top}
\bigl(\mi - \hat{\bm{\Pi}}_c)
& = [(\mi - \hat{\bm{\Pi}}_c) \hat{\muu}^{(i)}
\mw_{\muu}^{(i)}]
 [(\mi - \hat{\bm{\Pi}}_c) \hat{\muu}^{(i)}
\mw_{\muu}^{(i)}]^\top
\\ &= \muu_s^{(i)} \muu_s^{(i)\top} + (\mt_{0,s}^{(i)} + \mt_s^{(i)} +
\mathbf{M}^{(i)}) [\bm{0} \mid \muu_s^{(i)}]^{\top} 
+ 
 [\bm{0} \mid \muu_s^{(i)}] (\mt_{0,s}^{(i)} + \mt_s^{(i)} + \mathbf{M}^{(i)})^{\top}  
\\ &+ (\mt_{0,s}^{(i)} + \mt_s^{(i)} + \mathbf{M}^{(i)})(\mt_{0,s}^{(i)} + 
\mt_s^{(i)} + \mathbf{M}^{(i)})^{\top}.
\end{split}
\end{equation*}
From Eq.~\eqref{eq:m0form} we have
\begin{equation*}
  \begin{split}
 \mathbf{M}^{(i)} [\bm{0} \mid \muu_s^{(i)}]^{\top} &=
{\mathcal{L}} \muu^{(i)} [\bm{0} \mid \muu_s^{(i)}]^{\top} +
 \mathbf{T}_{\mathcal{L}}  [\bm{0} \mid
 \muu_s^{(i)}]^{\top} \\
 &=  \frac{1}{m} \sum_{j=1}^m \mt_0^{(j)} \muu_c
 \muu_c^{\top} \muu^{(i)} [\bm{0} \mid \muu_s^{(i)}]^{\top} +
\mathbf{T}_{\mathcal{L}} [\bm{0} \mid
\muu_s^{(i)}]^{\top}  \\
&= \mathbf{T}_{\mathcal{L}}  [\bm{0} \mid
\muu_s^{(i)}]^{\top},
  \end{split}
\end{equation*}
where the final equality is because
$$\muu_c \muu_c^{\top} \muu^{(i)} = \muu_c \muu_c^{\top} [\muu_c \mid
\muu_s^{(i)}] = [\muu_c \mid \bm{0}]$$
and $\muu_c^{\top} \muu_s^{(i)} = 0$ for all $i\in[m]$. 
We therefore have
$$
(\mi - \hat{\bm{\Pi}}_c)
\hat{\muu}^{(i)}
\hat{\muu}^{(i)\top}  (\mi - \hat{\bm{\Pi}}_c)
=\muu_s^{(i)} \muu_s^{(i)\top}+\tilde{\me}_s^{(i)},
$$
where we define
\begin{gather*}
  \tilde{\me}_s^{(i)} = [\bm{0} \mid \muu_s^{(i)}] \mt_{0,s}^{(i)\top}+\mt_{0,s}^{(i)} [\bm{0} \mid \muu_s^{(i)}]^{\top}+\tilde{\ml}_s^{(i)}
  ,\\
 {\ml}_s^{(i)} = [\bm{0} \mid \muu_s^{(i)}] 
(\mt_s^{(i)} + \mathbf{T}_{\mathcal{L}}^{(i)}
 )^{\top}+(\mt_s^{(i)} + \mathbf{T}_{\mathcal{L}}^{(i)}
 )
  [\bm{0} \mid \muu_s^{(i)}]^{\top}
   + (\mt_{0,s}^{(i)} + \mt_s^{(i)} +
  \mathbf{M}^{(i)})(\mt_{0,s}^{(i)} + \mt_s^{(i)} + \mathbf{M}^{(i)})^{\top}.
\end{gather*}
Note that, following similar derivations for $\tilde{\me}$ and $\ml$, we know $\tilde{\me}_s^{(i)}$ and ${\ml}_s^{(i)}$ have the same bounds with $\tilde{\me}$ and $\ml$ in Eq.~\eqref{eq:epsilon_e_thm3}.

Now write the eigen-decomposition of $(\mi - \hat{\bm{\Pi}}_c)
\hat{\muu}^{(i)}
\hat{\muu}^{(i)\top}  (\mi - \hat{\bm{\Pi}}_c)$ as 
$$\hat{\muu}_s^{(i)} \hat{\bm{\Lambda}}_s^{(i)} \hat{\muu}_s^{(i)\top}
+ \hat{\muu}_{s,\perp}^{(i)} \hat{\bm{\Lambda}}_{s,\perp}^{(i)}
\hat{\muu}_{s,\perp}^{(i)\top} = (\mi - \hat{\bm{\Pi}}_c)
\hat{\muu}^{(i)}
\hat{\muu}^{(i)\top}  (\mi - \hat{\bm{\Pi}}_c)=
\muu_s^{(i)} \muu_s^{(i)\top} + \tilde{\me}^{(i)}_s.$$
Once again $\hat{\muu}_s^{(i)}$ has a von-Neumann series expansion as
$$\hat{\muu}_s^{(i)} = \sum_{k=0}^{\infty} \tilde{\me}_s^{(i)}
\muu_s^{(i)} \muu_s^{(i)\top} \hat{\muu}_s^{(i)} (\hat{\bm{\Lambda}}_s^{(i)})^{-(k+1)}.$$
We can finally follow the exact same argument as that in the previous derivations for
$\hat{\muu}_c$, with $\tilde{\me}, 
\ml$ and $\bm{\Pi}_s$ there replaced by $\tilde{\me}_s^{(i)}, 
\tilde{\ml}_s^{(i)}$ and $\bm{0}$, respectively. We omit the
straightforward but tedious technical
details. 
In summary we obtain
\begin{equation*}
  \begin{split}
  \hat{\muu}_s^{(i)} \mw_{\muu_s}^{(i)} - \muu_s^{(i)} &= (\mi -
\bm{\Pi}_c) \mt_0^{(i)} [\bm{0} \mid
\muu_s^{(i)}]^{\top} \muu_s^{(i)} + {\mq}_{\muu_s}^{(i)} \\ &= (\mi -
\bm{\Pi}_c) \mt_0^{(i)} \muu^{(i)\top}
\muu_s^{(i)}  + \mq_{\muu_s}^{(i)} \\
&= \mt_0^{(i)} \muu^{(i)\top} \muu_s^{(i)} +
\mq_{\muu_s}^{(i)} + \muu_c \muu_c^{\top} \mt_0^{(i)} \muu^{(i)\top}
\muu_s^{(i)},
  \end{split}
\end{equation*}
where $\mq_{\muu_s}^{(i)}$ satisfies the same upper bound as that for $\mq_{\muu_c}$, and the term $\muu_c \muu_c^{\top} \mt_0^{(i)} \muu^{(i)\top}
\muu_s^{(i)}$
satisfies the same upper bound as
$\mq_{\muu_s}^{(i)}$ and can thus be subsumed by $\mq_{\muu_s}^{(i)}$. The expansion for $\hat\muu_s^{(i)}$ in Theorem~\ref{thm3_general} follows 
directly.

For the expansion for $\muu^{(i)}=[\muu_c \mid
\muu_s^{(i)}]$, combining the expansion for $\hat{\muu}_c$ and $\hat{\muu}_s^{(i)}$, we
conclude that there exists a block orthogonal matrix $\mw_{\muu}$ such
that 
\begin{equation*}
  \label{eq:final_decomp_thm3}
 [\hat{\muu}_c \mid \hat{\muu}_s^{(i)}] \mw_{\muu}^{(i)} - [\muu_c \mid
\muu_s^{(i)}] =
\Bigl[\frac{1}{m} \sum_{j=1}^{m} \mt_0^{(j)} 
\muu^{(j)\top} \muu_c \mid \mt_0^{(i)} 
\muu^{(i)\top} \muu_s^{(i)}\Bigr] + \mq_{\muu}^{(i)},
\end{equation*}
where $\mq_{\muu}^{(i)}$ satisfies the same upper bound as that for
$\mq_{\muu_c}$ and $\mq_{\muu_s}^{(i)}$
\hspace*{\fill} \qedsymbol

\subsection{Proof of Theorem~\ref{thm:UhatW-U=EVR^{-1}+xxx_part2}}
\label{Appendix:thm:UhatW-U=EVR^{-1}+xxx_part2}

Theorems~\ref{thm:UhatW-U=EVR^{-1}+xxx_part2} and~\ref{thm:WhatU_k-U_k->norm_part2} remain valid under a more general noise model for $\me^{(i)}$ as described in Assumption~\ref{ass:general}. Our proofs of these theorems (along with the corresponding technical lemmas) are based on this generalized noise model. See also~\cite{xie2021entrywise} for similar assumptions.
\begin{assumption}
  \label{ass:general}
  For each $i \in [m]$, $\me^{(i)}$ is an $n \times n$ matrix that can be decomposed as $\me^{(i)} = \me^{(i,1)} + \me^{(i,2)}$, with finite constants $C_1$, $C_2$, and $C_3$ independent of $m$ and $n$, such that
  \begin{enumerate}
  \item The entries $\{\me^{(i,1)}_{rs}\}_{r,s}$ are independent random variables with mean $0$, satisfying
    \begin{itemize}
      \item $\max_{i\in[m],r,s\in[n]} |\me^{(i,1)}_{rs}| \leq C_{1}$ almost
        surely.
      \item $\max_{i\in[m],r,s\in[n]}
        \mathbb{E}\bigl[(\me^{(i,1)}_{rs})^2\bigr] \leq C_{2}
        \rho_n$.
    \end{itemize}
    \item The entries $\{\me^{(i,2)}_{rs}\}_{r,s}$ are independent sub-Gaussian random variables with mean $0$, satisfying
      $$\max_{i\in[m],r,s\in[n]} \|\me^{(i,2)}_{rs}\|_{\psi_2} \leq C_3
      \rho_n^{1/2},$$
      where $\|\cdot\|_{\psi_2}$ represents the Orlicz-$2$ norm. 
    \item The matrices $\me^{(i,1)}$ and $\me^{(i,2)}$ are independent.
    \end{enumerate}
\end{assumption}

We begin by stating several fundamental bounds that will be consistently used in the proofs of Theorems~\ref{thm:UhatW-U=EVR^{-1}+xxx_part2} through~\ref{thm:HT}. 
Note that the proofs for Theorems~\ref{thm:UhatW-U=EVR^{-1}+xxx_part2} through~\ref{thm:HT} are primarily written for directed graphs; however, the same arguments apply to undirected graphs, where we assume $\muu_c = \mv_c$, $\muu_s^{(i)} = \mv_s^{(i)}$, and matrices $\ma^{(i)}$, $\mr^{(i)}$, $\mpp^{(i)}$, $\me^{(i)}$ are symmetric. 
The only step requiring additional attention arises in the proof of Lemma~\ref{lemma:U^T EE^T U R^-1->}, as the dependency among the entries of $\{\me^{(i)}\}$ leads to slightly more involved book-keeping.

\begin{lemma1}
  \label{lemma:|E|_2,|UEV|F}
  Consider the setting in Theorem~\ref{thm:UhatW-U=EVR^{-1}+xxx_part2}, where, for each $i \in [m]$, the noise matrix 
$\mathbf{E}^{(i)} = \mathbf{A}^{(i)} - \mathbf{P}^{(i)}$ is of the form described in Assumption~\ref{ass:general}. Then for any $i,j\in[m]$, we have
$$
\begin{aligned}
&\|\mathbf{E}^{(i)}\|\lesssim (n \rho_{n})^{1 / 2},\quad
\|\muu^{(i)\top} \mathbf{E}^{(j)}\mv^{(j)}\| \lesssim d_i^{1/2}( \rho_{n} \log n)^{1/2},\\
&\|\me^{(i)}\mv^{(i)} \|_{2\to\infty}\lesssim d_i^{1/2}(\rho_n \log n)^{1/2},\quad
\|\me^{(i)\top}\muu^{(i)}  \|_{2\to\infty}\lesssim d_i^{1/2}(\rho_n \log n)^{1/2}
\end{aligned}
$$
with high probability.
If we further assume $\{\me^{(i)}\}_i$ are independent then
$$\Bigl\|\frac{1}{m}\sum_{i=1}^m\me^{(i)}\mv^{(i)} (\mr^{(i)})^{-1}\Bigr\|_{2\to\infty}\lesssim d_i^{1/2}(mn)^{-1/2}(n\rho_n)^{-1/2}\log^{1/2} n$$
with high probability.
\end{lemma1}
We now present an essential technical lemma for bounding the error of $\hat{\muu}^{(i)}$ as an estimate of the true $\muu^{(i)}$, for each $i \in [m]$.
\begin{lemma1}
  \label{lemma:Uhat-UW}
  Consider the setting in Theorem~\ref{thm:UhatW-U=EVR^{-1}+xxx}, where, for each $i \in [m]$, the noise matrix 
$\mathbf{E}^{(i)} = \mathbf{A}^{(i)} - \mathbf{P}^{(i)}$ is of the form described in Assumption~\ref{ass:general}. Fix an $i \in [m]$ and write the singular value decomposition of $\ma^{(i)}$ as 
$
\ma^{(i)} = \hat{\muu}^{(i)} \hat{\mSigma}^{(i)} \hat{\mv}^{(i)\top} + \hat{\muu}_\perp^{(i)} \hat{\mSigma}_\perp^{(i)} \hat{\mv}_\perp^{(i)\top}.
$
Next define $\mw_{\muu}^{(i)}$ as the minimizer of $\|\hat{\muu}^{(i)} \mo - \muu^{(i)}\|_F$ over all $d_i \times d_i$ orthogonal matrices $\mo$, and define $\mw_{\mv}^{(i)}$ analogously.
We then have
  $$
  \hat\muu^{(i)}\mw_{\muu}^{(i)}-\muu^{(i)}=\me^{(i)}\mv^{(i)}(\mr^{(i)})^{-1}+\mt^{(i)},
  $$
  where $\mt^{(i)}$ is a $n \times d_i$ matrix satisfying
  \begin{equation*}
    \begin{aligned}
  &\|\mt^{(i)}\| \lesssim (n\rho_n)^{-1}\max\{1,d_i^{1/2}(\rho_n \log
      n)^{1/2}\},  \\
  &\|\mt^{(i)}\|_{2\to\infty} \lesssim d_i^{1/2}n^{-1/2}(n\rho_n)^{-1}\log n
    \end{aligned}
  \end{equation*}
  with high probability. 
  An analogous result holds for 
$\hat{\mv}^{(i)} \mw_{\mv}^{(i)} - \mv$, where $\me^{(i)}$, $\mr^{(i)}$, and $\mv^{(i)}$ are replaced by 
$\me^{(i)\top}$, $\mr^{(i)\top}$, and $\muu^{(i)}$, respectively.
\end{lemma1}

The proofs of Lemma~\ref{lemma:|E|_2,|UEV|F} and
Lemma~\ref{lemma:Uhat-UW} are presented in Section~\ref{Appendix:proof
  of lemma:|E|_2,|UEV|F} and Section~\ref{Appendix:proof of lemma:Uhat-UW},
respectively. 

We now apply Theorem~\ref{thm3_general} to derive the expansions for the estimations of the invariant
subspace $\muu_c$ as well as the possibly distinct subspaces $\{\muu_s^{(i)}\}$.
The expansions for $\mv_c$ and $\{\mv_s^{(i)}\}$ follow almost identical arguments and are therefore omitted.

For each $i\in[m]$, by Lemma~\ref{lemma:Uhat-UW} we have the expansion
\begin{equation*}
  \hat\muu^{(i)}\mw^{(i)}_\muu-\muu^{(i)}=\mt^{(i)}_0+\mt^{(i)}
\end{equation*}
for some orthogonal matrix $\mw_\muu^{(i)}$, 
where $\mt^{(i)}_0=\me^{(i)}\mv^{(i)}(\mr^{(i)})^{-1}$ and
$\mt^{(i)}$ satisfies
$$
\begin{aligned}
	&\|\mt^{(i)}\|\lesssim (n\rho_n)^{-1}\max\{1,d^{1/2}(\rho_n \log n)^{1/2}\}
,\\
&\|\mt^{(i)}\|_{2\to\infty}\lesssim d_i^{1/2}n^{-1/2}(n\rho_n)^{-1}\log n
\end{aligned}
$$
with high probability,
so by Lemma~\ref{lemma:|E|_2,|UEV|F} we have
$$
\begin{aligned}
	&\|\mt^{(i)}_0\|
\leq\|\me^{(i)}\|\cdot \|(\mr^{(i)})^{-1}\|
\lesssim (n\rho_n)^{1/2}\cdot (n\rho_n)^{-1}
\lesssim (n\rho_n)^{-1/2},\\
&\|\mt^{(i)}_0\|_{2\to\infty}
\leq\|\me^{(i)}\mv^{(i)}\|_{2\to\infty}\cdot \|(\mr^{(i)})^{-1}\|
\lesssim d_i^{1/2}n^{-1/2}(n\rho_n)^{-1/2}\log^{1/2}n
\end{aligned}
$$
with high probability. Thus we have
$$
\begin{aligned}
	\max_{i\in[m]}\big(2\|\mt^{(i)}_0\|
	+2\|\mt^{(i)}\|
	+\|\mt^{(i)}_0+\mt^{(i)}\|^2\big)
	\lesssim (n\rho_n)^{-1/2}
\end{aligned}
$$
with high probability.
Then with the assumption $n\rho_n=\Omega(\log n)$, we have $\max_{i\in[m]}\big(2\|\mt^{(i)}_0\|
	+2\|\mt^{(i)}\|
	+\|\mt^{(i)}_0+\mt^{(i)}\|^2\big)=o_p(1)$.
	Let $\bm{\Pi}_{s} = m^{-1} \sum_{i=1}^m \muu_s^{(i)} \muu_{s}^{(i)\top}$.
Under the assumption that $\|\mathbf\Pi_s\|=\| m^{-1} \sum_{i=1}^m \muu_s^{(i)} \muu_{s}^{(i)\top}\|\leq 1-c_s$ for some constant $0<c_s\leq 1$, we have 
$\frac{1}{2}(1-\|\mathbf\Pi_s\|)\geq \frac{c_s}{2}$. Then for large enough $n$
 we have
 $$
\begin{aligned}
	\max_{i\in[m]}\big(2\|\mt^{(i)}_0\|
	+2\|\mt^{(i)}\|
	+\|\mt^{(i)}_0+\mt^{(i)}\|^2\big)
	\leq c(1-\|\mathbf\Pi_s\|)
	< \frac{1}{2}(1-\|\mathbf\Pi_s\|)
\end{aligned}
$$
with high probability for any constant $c<\frac{1}{2}$.
Let $\vartheta_n=\max\{1,d_{\max}^{1/2}(\rho_n \log n)^{1/2}\}$. Then we have
\begin{equation}
\label{eq:ep_t,ze_t}
	\begin{aligned}
		&\epsilon_{\mt_0}
		=\max_{i\in[m]}\|\mt^{(i)}_0\|
		\lesssim (n\rho_n)^{-1/2},\\
		&\zeta_{\mt_0}
		=\max_{i\in[m]}\|\mt^{(i)}_0\|_{2\to\infty}
		\lesssim d_{\max}^{1/2}n^{-1/2}(n\rho_n)^{-1/2}\log^{1/2}n,\\
		&\epsilon_{\mt}
		=\max_{i\in[m]}\|\mt^{(i)}\|
		\lesssim (n\rho_n)^{-1}\vartheta_n,\\
		&\zeta_{\mt}
		=\max_{i\in[m]}\|\mt^{(i)}\|_{2\to\infty}
		\lesssim d_{\max}^{1/2}n^{-1/2}(n\rho_n)^{-1}\log n
	\end{aligned}
\end{equation}
with high probability.
By the assumption about $\muu^{(i)}$, we have
\begin{equation}
\label{eq:ze_u}
	\begin{aligned}
		\zeta_\muu=\max_{i\in[m]}\|\muu^{(i)}\|_{2\to\infty}
		\lesssim d_{\max}^{1/2}n^{-1/2}.
	\end{aligned}
\end{equation}
By Lemma~\ref{lemma:|E|_2,|UEV|F} we have
\begin{equation}
\label{eq:ep_cs}
	\begin{aligned}
		\epsilon_{\star}
		&=\max_{i,j\in[m]}\|\muu^{(i)\top}\mt^{(j)}_0\|
          \leq\max_{i\in[m]}\|\muu^{(i)\top}\me^{(j)}\mv^{(j)}\|\cdot\|(\mr^{(j)})^{-1}\|
          \\ & \lesssim d_{\max}^{1/2}(\rho_n\log n)^{1/2}\cdot (n\rho_n)^{-1}
\lesssim d_{\max}^{1/2}n^{-1/2}(n\rho_n)^{-1/2}\log^{1/2} n
	\end{aligned}
\end{equation}
with high probability.

Therefore by Theorem~\ref{thm3_general}, for the estimation of the invariant subspace $\muu_c$ we have
$$
\begin{aligned}
	\hat\muu_c\mw_{\muu_c}-\muu_c
	=\frac{1}{m}\sum_{i=1}^m\mt_0^{(i)}\muu^{(i)\top}\muu_c+\mq_{\muu_c}
	=\frac{1}{m}\sum_{i=1}^m\me^{(i)}\mv^{(i)}(\mr^{(i)})^{-1}\muu^{(i)\top}\muu_c+\mq_{\muu_c},
\end{aligned}
$$
where $\mw_{\muu_c}$ is a minimizer of $\|\hat{\muu}_c \mo -
\muu_c\|_{F}$ over all orthogonal matrix $\mo$,
and by Eq.~\eqref{eq:ep_t,ze_t}, Eq.~\eqref{eq:ze_u} and Eq.~\eqref{eq:ep_cs}, $\mq_{\muu_c}$ satisfies
$$
\begin{aligned}
	\|\mq_{\muu_c}\|
	&\lesssim \epsilon_{\star} + \epsilon_{\mt_0}^2+ \epsilon_{\mt}\\
	&\lesssim d_{\max}^{1/2}n^{-1/2}(n\rho_n)^{-1/2}\log^{1/2} n
	+[(n\rho_n)^{-1/2}]^2
	+(n\rho_n)^{-1}\vartheta_n\\
	&\lesssim (n\rho_n)^{-1}\vartheta_n,\\
	\|\mq_{\muu_c}\|_{2\to\infty}
	&\lesssim \zeta_\muu(\epsilon_{\star}+\epsilon_{\mt_0}^2+\epsilon_{\mt})
       +\zeta_{\mt_0}(\epsilon_{\star}+\epsilon_{\mt_0}+\epsilon_{\mt})
       +\zeta_{\mt}\\
    &\lesssim d_{\max}^{1/2}n^{-1/2}\cdot(n\rho_n)^{-1}\vartheta_n
    +d_{\max}^{1/2}n^{-1/2}(n\rho_n)^{-1/2}\log^{1/2}n\cdot (n\rho_n)^{-1/2}\\
    &+d_{\max}^{1/2}n^{-1/2}(n\rho_n)^{-1}\log n\\
    &\lesssim d_{\max}^{1/2}n^{-1/2}(n\rho_n)^{-1}\log n
\end{aligned}
$$
with high probability.
And for each $i\in[m]$, the estimation for the possibly distinct subspace $\muu_s^{(i)}$ has the expansion
$$
\begin{aligned}
	\hat\muu_s^{(i)}\mw_{\muu_s}^{(i)}-\muu_s^{(i)}
	=\mt_0^{(i)}\muu^{(i)\top}\muu_{s}^{(i)}+\mq_{\muu_s}^{(i)}
	=\me^{(i)}\mv^{(i)}(\mr^{(i)})^{-1}\muu^{(i)\top}\muu_{s}^{(i)}+\mq_{\muu_s}^{(i)},
\end{aligned}
$$
where $\mw^{(i)}_{\muu_s}$ is a minimizer of $\|\hat{\muu}_s^{(i)} \mo -
\muu_s^{(i)}\|_{F}$ over all orthogonal matrix $\mo$, and $\mq_{\muu_s}^{(i)}$ satisfies the same upper bounds as that for $\mq_{\muu_c}$. 

Finally, for any fixed $k \in [n]$, the bound $\|q_{\muu_c,k}\| \lesssim d_{\max}^{1/2} n^{-1/2} (n \rho_n)^{-1} t$, which holds with probability at least $1 - n^{-c} - O(me^{-t})$ for any $c > 0$, can be derived as follows. 
First, we can replace the upper bound in Lemma~\ref{lemma:||e_h^top E (hat V[h]hat V[h]^topV-V)||}
with the bound $d_i^{1/2} n^{-1/2} t$ which holds with probability $1 - O(e^{-t})$.
Similarly, the $2 \to \infty$ norm bounds in
Lemmas~\ref{lemma:||E(hatVtildeW-V)||2toinfty} and \ref{lemma:T4}, which hold uniformly for all $n$
rows with high probability,
can be replaced by bounds for a single row of the form $d_i^{1/2} n^{-1/2} (n \rho_n)^{-1} t$ which holds with probability at least $1 - n^{-c}- O(e^{-t})$ for any $c>0$. Combining these modified bounds we can show that a single row of
$\mt^{(i)}$ in Lemma~\ref{lemma:Uhat-UW} is upper bounded by $d_i^{1/2} n^{-1/2} (n \rho_n)^{-1} t$
with probability at least $1 - n^{-c}- O(e^{-t})$ for any $c>0$ under the condition $m = O(n^{c'})$ for some 
finite constant $c' > 0$. Finally, by careful book-keeping we can show that
$\max_{1 \leq r \leq 5} \|q_{\muu_c,k,r}\|$ is also upper bounded by $d_{\max}^{1/2} n^{-1/2} (n \rho_n)^{-1} t$ with probability at least $1 - n^{-c}- O(me^{-t})$ for any $c>0$ under the assumption $m=O(n^{c'})$ for some $c'>0$; here $q_{\muu_c,k,r}$ is the $k$th row of the matrix $\mq_{\muu_c,r}$ defined
in the proof of Theorem~\ref{thm3_general}. The analysis for the bound on $\|q_{\muu_s,k}^{(i)}\|$ follows similar arguments. We omit the details as they are mostly technical and tedious.
\hspace*{\fill} \qedsymbol
	
\subsection{Proof of Proposition~\ref{thm:Vhat-VW_part2}}

Eq.~\eqref{eq:expansion_hatU_2inf2_part2} follows directly from
Eq.~\eqref{eq:expansion_hatU_seperate_3} and Lemma~\ref{lemma:|E|_2,|UEV|F}.

For Eq.~\eqref{eq:expansion_hatU_F_part2}, by Theorem~\ref{thm:UhatW-U=EVR^{-1}+xxx_part2} we have
$$
\begin{aligned}
	\hat{\mathbf{U}}_c\mathbf{W}_{\mathbf{U}_c} -\mathbf{U}_c =\frac{1}{m}\sum_{i=1}^m\me^{(i)}\mv^{(i)}(\mr^{(i)})^{-1}\muu^{(i)\top}\muu_c+\mq_{\muu_c},
\end{aligned}
$$
where $\|\mq_{\muu_c}\|_F\leq d_{0,\muu}^{1/2}\|\mq_{\muu_c}\|\lesssim d_{0,\muu}^{1/2}(n\rho_n)^{-1}\max\{1,d_{\max}^{1/2}\rho_n^{1/2}\log^{1/2} n\}$ with high probability.
Furthermore we have
$$
\begin{aligned}
	&\Bigl\|\frac{1}{m}\sum_{i=1}^m\me^{(i)}\mv^{(i)}(\mr^{(i)})^{-1}\muu^{(i)\top}\muu_c\Bigr\|_F^{2}\\
	&= \frac{1}{m^{2}}\mathrm{tr} \Bigl[\sum_{i=1}^{m} \sum_{j=1}^{m} \muu_c^{\top}\muu^{(i)}\bigl(\mathbf{R}^{(i) \top}\bigr)^{-1} \mathbf{V}^{(i)\top} \mathbf{E}^{(i) \top} \mathbf{E}^{(j)} \mathbf{V}^{(i)}\bigl(\mathbf{R}^{(j)}\bigr)^{-1}\muu^{(i)\top}\muu_c\Bigr]\\
	&=\frac{1}{m^2}\sum_{i=1}^m\bigl\| \mathbf{E}^{(i)} \mathbf{V}^{(i)}\bigl(\mathbf{R}^{(i)}\bigr)^{-1}\muu^{(i)\top}\muu_c\bigr\|_F^2
	+\frac{1}{m^{2}} \sum_{i\neq j}\mathrm{tr}\bigl[\muu_c^{\top}\muu^{(i)}\bigl(\mathbf{R}^{(i) \top}\bigr)^{-1} \mathbf{V}^{(i)\top} \mathbf{E}^{(i) \top} \mathbf{E}^{(j)} \mathbf{V}^{(i)}\bigl(\mathbf{R}^{(j)}\bigr)^{-1}\muu^{(i)\top}\muu_c\bigl]\\
	&\lesssim m^{-1}\cdot d_{\max}(n\rho_n)^{-1}
	+m^{-1}\cdot d_{0,\muu}\cdot d_{\max}^3n^{-1/2}(n\rho_n)^{-1} 
	\lesssim d_{\max}m^{-1}(n\rho_n)^{-1}
\end{aligned}
$$
with high probability. Indeed, for any $i\in[m]$ we have
$$
\begin{aligned}
	\bigl\|\me^{(i)}\mv^{(i)}(\mr^{(i)})^{-1}\muu^{(i)\top}\muu_c\bigr\|_F
	&\leq d_i^{1/2} \|\me^{(i)}\|\cdot \|(\mr^{(i)})^{-1}\|
	\cdot \|\muu^{(i)\top}\muu_c\|
	\lesssim d_i^{1/2}(n\rho_n)^{-1/2}
\end{aligned}
$$
with high probability, 
and with the similar analysis as the proof of Lemma \ref{lemma:U^T EE^T U R^-1->} we have, for any $i\neq
j$ and $s \in [d_{0,\muu}]$
$$
\begin{aligned}
	 \Bigl[\muu_c^{\top}\muu^{(i)}\bigl(\mathbf{R}^{(i) \top}\bigr)^{-1} \mathbf{V}^{(i)\top} \mathbf{E}^{(i) \top} \mathbf{E}^{(j)} \mathbf{V}^{(i)}\bigl(\mathbf{R}^{(j)}\bigr)^{-1}\muu^{(i)\top}\muu_c\Bigr]_{ss}
	 \lesssim d_{\max}^3n^{-1/2}(n\rho_n)^{-1}
\end{aligned}
$$
with high probability. 
In summary, we have
$$
 \Bigl\|\frac{1}{m}\sum_{i=1}^m\me^{(i)}\mv(\mr^{(i)})^{-1}\muu^{(i)\top}\muu_c\Bigr\|_F
 \lesssim d_{\max}^{1/2}m^{-1/2}(n\rho_n)^{-1/2}
$$
with high probability, and the desired result of $\|\hat{\mathbf{U}}_c\mathbf{W}_{\mathbf{U}_c} -\mathbf{U}_c \|_F$ is obtained.
The analysis for the bound of $\|\hat{\mathbf{U}}_s^{(i)}\mathbf{W}_{\mathbf{U}_s}^{(i)} -\mathbf{U}_s^{(i)} \|_F$ follows similar arguments.
\hspace*{\fill} \qedsymbol

\subsection{Proof of Theorem~\ref{thm:WhatU_k-U_k->norm_part2}}
\label{sec:proof_normality1_part2}
We emphasize once again that the following proof is written for the
more general noise model described in Assumption~\ref{ass:general}. 

We now derive Eq.~\eqref{eq:clt_COSIE_part2} for $\hat{u}_{c,k}$. The result for $\hat{u}_{s,k}^{(i)}$ follows from similar arguments.  
According to Eq.~\eqref{eq:expansion_hatU_seperate_3}, we have 
\begin{equation}\label{eq:thm2_proof1}
	\begin{aligned}
		\mw_{\muu_c}^\top \hat u_{c,k}-u_{c,k}=\sum_{i=1}^m\sum_{\ell=1}^n\my_{i,\ell}^{(k)}+q_{\muu_c,k},
	\end{aligned}
\end{equation}
where we define
$$
\my_{i,\ell}^{(k)}:= \sum_{\ell=1}^n \me^{(i)}_{k\ell}\muu_c^\top\muu^{(i)}(\mr^{(i)\top})^{-1}v_\ell.
$$
Note that $\{\my_{i,\ell}^{(k)}\}_{i\in[m],\ell\in[n]}$ are independent mean $\bm{0}$ random vectors.
For  any $i\in[m],\ell\in[n]$, the variance of $\my_{i,\ell}^{(k)}$ is
\begin{equation*}
	\begin{aligned}
\mathrm{Var}\big[\my_{i,\ell}^{(k)}\big]
=m^{-2}\operatorname{Var}[\me^{(i)}_{k\ell}]\cdot \muu_c^\top\muu^{(i)}
(\mr^{(i)\top})^{-1}v_\ell v_\ell^\top(\mr^{(i)})^{-1}\muu^{(i)\top}\muu_c.
	\end{aligned}
\end{equation*}
and hence 
\begin{equation*}
  \begin{split}
\sum_{i=1}^{m}\sum_{\ell=1}^{n}
\mathrm{Var}\big[\my_{i,\ell}^{(k)}\big]
=\sum_{i=1}^{m}m^{-2}\muu_c^\top\muu^{(i)}
(\mr^{(i)\top})^{-1}\mv^\top\bm\Xi^{(i,k)}\mv(\mr^{(i)})^{-1}\muu^{(i)\top}\muu_c
\end{split}
\end{equation*}
where, for each $(k,i)$, $\bm{\Xi}^{(i,k)}$ is a $n \times n$ diagonal
matrix whose diagonal entries are $\mathrm{Var}[\me^{(i)}_{k
  \ell}]$. 
In the special case of the COISIE model we have
$\operatorname{Var}[\me^{(i)}_{k\ell}]=\mpp^{(i)}_{k\ell}(1-\mpp^{(i)}_{k\ell})$
which yields the covariance matrix $\bm{\Upsilon}^{(k)}$ in Theorem~\ref{thm:WhatU_k-U_k->norm_part2}.

\begin{sloppypar}
	Now let $\tilde\my_{i,\ell}^{(k)}=(\bm{\Upsilon}_{\muu_c}^{(k)})^{-1/2}\my_{i,\ell}^{(k)}$ 
and set $\tilde\my_{i,\ell}^{(k,1)}=(\bm{\Upsilon}_{\muu_c}^{(k)})^{-1/2}m^{-1}\me^{(i,1)}_{k\ell}\muu_c^\top\muu^{(i)}(\mr^{(i)\top})^{-1}v_\ell$
and $\tilde\my_{i,\ell}^{(k,2)}=(\bm{\Upsilon}_{\muu_c}^{(k)})^{-1/2}m^{-1}\me^{(i,2)}_{k\ell}\muu_c^\top\muu^{(i)}(\mr^{(i)\top})^{-1}v_\ell$.
From the assumption $\sigma_{\min}(\bm{\Upsilon}_{\muu_c}^{(k)})\gtrsim
m^{-1}n^{-2}\rho_n^{-1}$, we have
$\|(\bm{\Upsilon}_{\muu_c}^{(k)})^{-1/2}\|\lesssim m^{1/2}n\rho_n^{1/2}$.
Then for any $i\in[m],\ell\in[n]$, we can bound the spectral norm of $\tilde\my_{i,\ell}^{(k,1)}$ by
\end{sloppypar}
\begin{equation}
\label{eq:||Y_ij||}
	\begin{aligned}
	\big\|\tilde\my_{i,\ell}^{(k,1)}\big\|
	&\leq \|(\bm{\Upsilon}_{\muu_c}^{(k)})^{-1/2}\| \cdot  m^{-1}|\me_{k\ell}^{(i,1)}|
	\cdot \|\muu_c^\top\muu^{(i)}\|
	\cdot \|(\mr^{(i)\top})^{-1}\|\cdot \|v_\ell\|
	\\ &\lesssim m^{1/2} n \rho_n^{1/2}	\cdot   m^{-1}\cdot 1\cdot 1 \cdot (n\rho_n)^{-1}
	\cdot d_i^{1/2}n^{-1/2} 
	\lesssim d_i^{1/2} (mn\rho_n)^{-1/2}
\end{aligned}
\end{equation}
almost surely.
For any fixed $\epsilon > 0$,
Eq.~\eqref{eq:||Y_ij||} implies that, for sufficiently large $n$, we
have 
$\bigl\|\tilde\my_{i,\ell}^{(k,1)}\bigr\|\leq \epsilon$
almost surely.

For $\tilde\my_{i,\ell}^{(k,2)}$, because $\me_{k\ell}^{(i,2)}$ is
sub-Gaussian  with
$\|\me_{k\ell}^{(i,2)}\|_{\psi_2}\lesssim \rho_n^{1/2}$, by a similar
analysis to Eq.~\eqref{eq:||Y_ij||} we have 
$\|\|\tilde\my_{i,\ell}^{(k,2)}\|\|_{\psi_2}\lesssim
d_i^{1/2}(mn)^{-1/2}$. Now, for any fixed but arbitrary $\epsilon > 0$, we have
\begin{equation*}
  \begin{split}
\mathbb{E}\Bigl[\bigl\|\tilde\my_{i,\ell}^{(k)}\bigr\|^{2}
   \cdot\mathbb{I}\bigl\{\|\tilde\my_{i,\ell}^{(k)}\|>\epsilon
   \bigr\}\Bigr] &
      \leq \mathbb{E}\Bigl[\bigl\|\tilde\my_{i,\ell}^{(k)}\bigr\|^{2}
   \mathbb{I}\bigl\{\|\tilde\my_{i,\ell}^{(k,1)}\|>\tfrac{\epsilon}{2}\bigr\}\Bigr]
      + \mathbb{E}\Bigl[\bigl\|\tilde\my_{i,\ell}^{(k)}\bigr\|^{2}
      \mathbb{I}\bigl\{\|\tilde\my_{i,\ell}^{(k,2)}\|>\tfrac{\epsilon}{2}\bigr\}\Bigr].
  \end{split}
\end{equation*}
Therefore, if $n$ is sufficiently large, we have
\begin{equation}
  \label{eq:reduction_thm2_proof}
 \mathbb{E}\Bigl[\bigl\|\tilde\my_{i,\ell}^{(k)}\bigr\|^{2}
   \cdot\mathbb{I}\bigl\{\|\tilde\my_{i,\ell}^{(k)}\|>\epsilon
   \bigr\}\Bigr] \leq \mathbb{E}\Bigl[\bigl\|\tilde\my_{i,\ell}^{(k)}\bigr\|^{2}
      \mathbb{I}\bigl\{\|\tilde\my_{i,\ell}^{(k,2)}\|>\tfrac{\epsilon}{2}\bigr\}\Bigr].
\end{equation}
      Furthermore, we also have
      \begin{equation*}
        \begin{split}
        \mathbb{E}\Bigl[\bigl\|\tilde\my_{i,\ell}^{(k)}\bigr\|^{2}
      \mathbb{I}\bigl\{\|\tilde\my_{i,\ell}^{(k,2)}\|>\tfrac{\epsilon}{2}\bigr\}\Bigr]
      & \leq 
      \mathbb{E}\Bigl[2 \bigl(\bigl\|\tilde\my_{i,\ell}^{(k,1)}\bigr\|^{2}
        + \bigl\|\tilde\my_{i,\ell}^{(k,2)}\bigr\|^{2}\bigr)
      \cdot
      \mathbb{I}\bigl\{\|\tilde\my_{i,\ell}^{(k,2)}\|>\tfrac{\epsilon}{2}\bigr\}\Bigr]
      \\ & \leq 2 \mathbb{E}\Bigl[\bigl\|\tilde\my_{i,\ell}^{(k,1)}\bigr\|^{2}\Bigr]
        \cdot \mathbb{P}(\|\tilde\my_{i,\ell}^{(k,2)}\|>\tfrac{\epsilon}{2}\bigr\})
    + 4 \epsilon^{-1} \mathbb{E}\Bigl[\bigl\|\tilde\my_{i,\ell}^{(k,2)}\bigr\|^{3}\Bigr],
        \end{split}
      \end{equation*}
      where the second inequality follows from the independence of
      $\tilde{\my}_{i,\ell}^{(k,1)}$ and $\tilde{\my}_{i,\ell}^{(k,2)}$ (as
      $\me^{(i,1)}_{k \ell}$ is independent of $\me_{k\ell}^{(i,2)}$). 
      As $\|\tilde{\my}^{(k,2)}_{i, \ell}\|$ is sub-Gaussian with
$\|\|\tilde\my_{i,\ell}^{(k,2)}\|\|_{\psi_2}\lesssim
d_i^{1/2}(mn)^{-1/2}$,
      there exists a constant $C > 0$ such that
\begin{equation*}
 	\begin{aligned}
 	\mathbb{P}\big[\big\|\tilde\my_{i,\ell}^{(k,2)}\big\| \geq
      \tfrac{\epsilon}{2} \big]
	\leq 2\exp\Big(\frac{- Cmn \epsilon^2}{4d_i}\Big), \quad
\Bigl(\mathbb{E}[\|\tilde\my_{i,\ell}^{(k,2)}\|^3]\Bigr)^{1/3} \leq
C d_i^{1/2} (m n)^{-1/2}.
 \end{aligned}
 \end{equation*}
 See Eq.~(2.14) and Eq.~(2.15) in \cite{vershynin2018high} for
 more details on the above bounds. Combining the above bounds and Eq.~\eqref{eq:||Y_ij||}, we therefore have
 \begin{equation}
   \label{eq:reduction2_theorem2_proof}
   \mathbb{E}\Bigl[\bigl\|\tilde\my_{i,\ell}^{(k)}\bigr\|^{2}
   \mathbb{I}\bigl\{\|\tilde\my_{i,\ell}^{(k,2)}\|>\tfrac{\epsilon}{2}\bigr\}\Bigr]
   \lesssim \frac{d_i}{m n \rho_n} \exp\Big(\frac{- Cmn
     \epsilon^2}{4d_i}\Big) + \epsilon^{-1} d_i^{3/2} (mn)^{-3/2}.
 \end{equation}
 Substituting Eq.~\eqref{eq:reduction2_theorem2_proof} into
 Eq.~\eqref{eq:reduction_thm2_proof} and then summing over $i \in [m]$
 and $\ell \in [n]$ we obtain
\begin{equation*}
 \begin{aligned}
\lim_{n \rightarrow \infty} \sum_{i=1}^{m}
		\sum_{\ell=1}^{n}
\mathbb{E}\Bigl[\bigl\|\tilde\my_{i,\ell}^{(k)}\bigr\|^{2}
   \cdot\mathbb{I}\bigl\{\|\tilde\my_{i,\ell}^{(k)}\|>\epsilon
   \bigr\}\Bigr] &\lesssim \lim_{n \rightarrow \infty}  d_i (n\rho_n)^{-1} \Big[n\exp\Big(\frac{- Cmn
     \epsilon^2}{4d_i}\Big)\Big] + \epsilon^{-1} d_i^{3/2} (m n)^{-1/2} = 0.
	\end{aligned}
\end{equation*}
As $\epsilon > 0$ is fixed but arbitrary, 
the collection $\{\tilde{\my}^{(k)}_{i \ell}\}$ satisfies 
the condition of the Lindeberg-Feller central limit theorem (see e.g.,
Proposition~2.27 in \cite{van2000asymptotic}) and hence
\begin{equation}
\label{eq:thm2_proof3}
(\bm{\Upsilon}^{(k)}_{\muu_c})^{-1/2}\sum_{i=1}^m\sum_{\ell=1}^n\my_{i,\ell}^{(k)}
=\sum_{i=1}^m\sum_{\ell=1}^n\tilde\my_{i,\ell}^{(k)}
\leadsto \mathcal{N}\big(\mathbf{0},\mathbf{I}_{d_{0,\muu}}\big)
\end{equation}
as $n \rightarrow \infty$. For the second term on the right hand side of Eq.~\eqref{eq:thm2_proof1} we have 
$$
\begin{aligned}
	\|(\bm{\Upsilon}^{(k)}_{\muu_c})^{-1/2}q_{\muu_c,k}\|
	&\lesssim \|(\bm{\Upsilon}_{\muu_c}^{(k)})^{-1/2}\|\cdot\|q_{\muu_c,k}\|\\
	&\lesssim m^{1/2} n\rho_n^{1/2}\cdot d_{\max}^{1/2}n^{-1/2}(n\rho_n)^{-1} t
	\lesssim d_{\max}^{1/2}m^{1/2}(n\rho_n)^{-1/2} t
\end{aligned}
$$
with probability at least $1 - n^{-c}- O(me^{-t})$ for any $c>0$. If $m \log^2 m= o(n \rho_n)$ then we can choose  
$t$ depending on $n$ such that $me^{-t} \rightarrow 0$ and $d_{\max}^{1/2} m^{1/2} (n \rho_n)^{-1/2} t \rightarrow 0$ as $n \rightarrow \infty$. In other words we have
\begin{equation}
\label{eq:thm2_proof2}
(\bm{\Upsilon}^{(k)}_{\muu_c})^{-1/2}{q_{\muu_c,k}}\stackrel{p}{\longrightarrow} \bm{0}
\end{equation}
as $n \rightarrow \infty$. 
Combining Eq.~\eqref{eq:thm2_proof1}, Eq.~\eqref{eq:thm2_proof3} and Eq.~\eqref{eq:thm2_proof2}, and applying Slutsky's theorem, we obtain
$$
(\bm{\Upsilon}^{(k)}_{\muu_c})^{-1/2}(\mw_{\muu_c}^\top \hat{u}_{c,k}- u_{c,k})
\leadsto \mathcal{N}\big(\mathbf{0},\mathbf{I}_{d_{0,\muu}}\big)
$$
as $n \rightarrow \infty.$
\hspace*{\fill} \qedsymbol

\subsection{Formal statements of some theoretical results for the COSIE model}
\label{sec:formal COSIE}

\begin{definition1}[Common subspace independent edge graphs] 
\label{CSIEgraph} 
For each $i \in [m]$, let $\mathbf{R}^{(i)}$ be a $d \times d$ matrix, and let $\mathbf{U}$ and $\mathbf{V}$ be $n \times d$ orthonormal matrices representing the shared subspaces across all $i$, such that $u_t^\top \mathbf{R}^{(i)} v_k \in [0,1]$ for all $t, k \in [n]$ and $i \in [m]$,  
where $u_t$ and $v_k$ denote the $t$th and $k$th rows of $\mathbf{U}$ and $\mathbf{V}$, respectively.
We say that the random adjacency matrices $\{\mathbf{A}^{(i)}\}_{i=1}^m$ are jointly distributed according to the common subspaces independent edge graphs model with $\mathbf{U}$, $\mathbf{V}$, $\{\mathbf{R}^{(i)}\}_{i=1}^m$, if, for each $i \in [m]$, $\mathbf{A}^{(i)}$ is an $n \times n$ random matrix whose entries $\{\mathbf{A}_{tk}^{(i)}\}$ are independent Bernoulli random variables with  
$
\mathbb{P}[\mathbf{A}_{tk}^{(i)} = 1] = u_t^{\top} \mathbf{R}^{(i)} v_k.
$
In other words,  
$$
\mathbb{P}\big(\mathbf{A}^{(i)} \mid \mathbf{U}, \mathbf{V}, \mathbf{R}^{(i)}\big) 
= \prod_{t\in[n]} \prod_{k\in[n]} \big(u_t^{\top} \mathbf{R}^{(i)} v_k\big)^{\mathbf{A}_{tk}^{(i)}} \big(1 - u_t^{\top} \mathbf{R}^{(i)} v_k\big)^{1-\mathbf{A}_{tk}^{(i)}}.
$$ 
We denote the multiple networks by $\left(\mathbf{A}^{(1)}, \ldots, \mathbf{A}^{(m)}\right) \sim \operatorname{COSIE}(\mathbf{U},\mathbf{V},\{\mathbf{R}^{(i)}\}_{i=1}^m)$, and write  
\begin{equation*}
    \mathbf{P}^{(i)} = \mathbf{U}\mathbf{R}^{(i)} \mathbf{V}^{\top}  \end{equation*}
to represent the (unobserved) edge probabilities matrix for each network $\mathbf{A}^{(i)}$.
\end{definition1}

\begin{algorithm}[htbp!]
\caption{Estimation of COSIE parameters}	
\label{Alg}
\begin{algorithmic}
\REQUIRE Adjacency matrices $\mathbf{A}^{(1)}, \dots, \mathbf{A}^{(m)}$, embedding dimension $d_1, \dots, d_m$, a final embedding dimension $d$. 
\begin{enumerate}
	\item For each $i \in [m]$, obtain $\hat{\mathbf{U}}^{(i)}$ and $\hat{\mathbf{V}}^{(i)}$ as the $n \times d_i$ matrices whose columns are the $d_i$ leading left and right singular vectors of $\mathbf{A}^{(i)}$, respectively.
	
	\item Compute $\hat{\mathbf{U}}$ as the $n \times d$ matrix whose columns are the leading left singular vectors of $[\hat{\mathbf{U}}^{(1)} \mid \cdots \mid \hat{\mathbf{U}}^{(m)}]$, and compute $\hat{\mathbf{V}}$ as the $n \times d$ matrix whose columns are the leading left singular vectors of $[\hat{\mathbf{V}}^{(1)} \mid \cdots \mid \hat{\mathbf{V}}^{(m)}]$.

	\item For each $i \in [m]$, compute $\hat{\mathbf{R}}^{(i)} = \hat{\mathbf{U}}^{\top} \mathbf{A}^{(i)} \hat{\mathbf{V}}$.
\end{enumerate} 
\ENSURE $\hat{\mathbf{U}}, \hat{\mathbf{V}}, \{\hat{\mathbf{R}}^{(i)}\}_{i=1}^m$.
\end{algorithmic}
\end{algorithm}

\begin{assumption}
  \label{ass:main}
  The following conditions hold for sufficiently large $n$. 
  \begin{itemize}
\item	The matrices $\muu$ and $\mv$ are $n \times d$ matrices with bounded coherence, i.e.,
      $$\|\muu\|_{2 \to
        \infty} \lesssim d^{1/2}n^{-1/2} \quad \text{and} \quad \|\mv\|_{2 \to \infty}
      \lesssim d^{1/2}n^{-1/2}.$$
\item There exists a factor $\rho_n \in [0,1]$ depending on $n$
        such that for each $i \in [m]$, $\mathbf{R}^{(i)}$ is a
        $d \times d$ matrix with $\|\mathbf{R}^{(i)}\| =
        \Theta(n \rho_n)$ where $n\rho_n\geq C\log n$ for some sufficiently large but finite constant $C>0$. 
        We interpret $n \rho_n$ as the growth rate for the average degree
        of the graphs $\ma^{(i)}$ generated from $\mpp^{(i)}$. 
\item The matrices $\{\mr^{(i)}\}_{i=1}^{m}$ have bounded condition numbers, i.e., there
        exists a finite constant $M$ such that
        $$\max_{i\in[m]}
        \frac{\sigma_1(\mathbf{R}^{(i)})}{\sigma_d(\mathbf{R}^{(i)})}
          \leq M,$$
       where $\sigma_{1}(\mathbf{R}^{(i)})$ and $\sigma_d(\mr^{(i)})$ denote the
       largest and smallest singular values of $\mathbf{R}^{(i)}$, respectively.
      \end{itemize}
\end{assumption}

\begin{theorem1}
  \label{thm:UhatW-U=EVR^{-1}+xxx}
   Consider $\left(\mathbf{A}^{(1)}, \ldots, \mathbf{A}^{(m)}\right) \sim \operatorname{COSIE}(\mathbf{U}, \mathbf{V}, \{\mathbf{R}^{(i)}\}_{i=1}^m)$ under the conditions in Assumption~\ref{ass:main}. Let $\hat{\mathbf{U}}$ be the estimate of $\mathbf{U}$ obtained by Algorithm~\ref{Alg}, and let $\mathbf{W}_{\mathbf{U}}$ be the minimizer of $\|\hat{\mathbf{U}} \mathbf{O} - \mathbf{U}\|_{F}$ over all $d \times d$ orthogonal matrices $\mathbf{O}$.
Then
\begin{equation*}
  \hat{\mathbf{U}} \mathbf{W}_{\mathbf{U}} -\mathbf{U} =
  \frac{1}{m}\sum_{i=1}^m \me^{(i)} \mv(\mr^{(i)})^{-1}+\mq_{\muu},
\end{equation*}
where $\me^{(i)}=\ma^{(i)}-\mpp^{(i)}$ and $\mq_{\muu}$ is a random matrix satisfying
\begin{equation*}
  \begin{aligned}
 &\|\mq_{\muu}\| \lesssim (n \rho_n)^{-1}\max\{1,d^{1/2}\rho_n^{1/2}\log^{1/2} n\}, \\
&\|\mq_{\muu}\|_{2\to\infty} \lesssim d^{1/2}n^{-1/2}(n\rho_n)^{-1}\log n
  \end{aligned}
\end{equation*}
with high probability. 
And for any $k \in [n]$, the $k$th row $q_{\muu,k}$ of $\mq_{\muu}$ satisfies $$\|q_{\muu,k}\| \lesssim d^{1/2} n^{-1/2} (n \rho_n)^{-1} t$$ with probability at least $1 - n^{-c}- O(me^{-t})$ for any $c>0$. 

The estimate $\hat{\mathbf{V}}$ has an analogous expansion, with $\mathbf{E}^{(i)}$, $\mv$, $\mathbf{R}^{(i)}$, and $\mq_{\muu}$ replaced by $\mathbf{E}^{(i)\top}$, $\muu$, $\mathbf{R}^{(i)\top}$, and $\mq_{\mv}$, respectively.
\end{theorem1}

\begin{proposition1}
  \label{thm:Vhat-VW}
  Consider the setting in Theorem~\ref{thm:UhatW-U=EVR^{-1}+xxx} and
furthermore assume that $\{\ma^{(i)}\}_{i=1}^{m}$ are mutually independent.
We then have
\begin{equation*}
  \begin{aligned}
  	  &\|\hat{\mathbf{U}}\mathbf{W}_{\mathbf{U}} -\mathbf{U}\|_{2 \to
    \infty}  \lesssim  d^{1/2} (mn)^{-1/2} (n \rho_n)^{-1/2} \log^{1/2}n+d^{1/2}n^{-1/2}(n\rho_n)^{-1}\log n,\\
   & \|\hat{\mathbf{U}}\mw_{\muu}-\mathbf{U}\|_{F} 
\lesssim d^{1/2}m^{-1/2}(n\rho_n)^{-1/2}
 +d^{1/2}(n\rho_n)^{-1}\max\{1, (d\rho_n\log n)^{1/2}\}
  \end{aligned}
\end{equation*}
with high probability.
Similar results hold for $\hat{\mathbf{V}}$.
\end{proposition1}

\begin{theorem1}
  \label{thm:WhatU_k-U_k->norm}
  Consider the setting in Theorem~\ref{thm:UhatW-U=EVR^{-1}+xxx} and
furthermore assume that $\{\ma^{(i)}\}_{i=1}^{m}$ are mutually independent.
For each $i\in[m]$ and $k\in[n]$, let $\mathbf{\Xi}^{(i,k)}$ be a $n\times n$ diagonal matrix whose diagonal elements are of the form
$$
\mathbf{\Xi}^{(i,k)}_{\ell\ell}=\mpp_{k  \ell}^{(i)}(1-\mpp_{k    \ell}^{(i)}).
$$
Define $\bm{\Upsilon}^{(k)}_\muu$ as the $d \times d$ symmetric matrix 
$$
\bm{\Upsilon}^{(k)}_\muu=\frac{1}{m^2}\sum_{i=1}^{m}(\mr^{(i)\top})^{-1}\mv^\top\mathbf{\Xi}^{(k,i)}\mv(\mr^{(i)})^{-1}.
$$
Note that $\|\bm{\Upsilon}^{(k)}_{\muu}\|\lesssim (mn^2\rho_n)^{-1}$. 
Further suppose $\sigma_{\min}(\bm{\Upsilon}^{(k)}_{\muu})\gtrsim
(mn^2\rho_n)^{-1}$. 
Then for the $k$th rows $\hat{u}_{k}$ and ${u}_{k}$ of $\hat\muu$ and $\muu$, we have
\begin{equation*}
  \begin{aligned}
  	(\bm{\Upsilon}^{(k)}_{\muu})^{-1/2} \big(\mw_{\muu}^\top \hat{u}_{k}-{u}_{k}\big)
\rightsquigarrow \mathcal{N}\big(\mathbf{0},\mi_{d}\big)
  \end{aligned}
\end{equation*}
as $n \rightarrow \infty$.
Similar results hold for $\hat{\mathbf{V}}$ and its rows $\hat v_{k}$ with $\mpp^{(i)}$, $\mathbf{R}^{(i)}$, and $\mv$ replaced by $\mpp^{(i)\top}$, $\mathbf{R}^{(i)\top}$, and $\muu$ respectively.
\end{theorem1}

\subsection{Proof of Theorem~\ref{thm:What_U Rhat What_v^T-R->norm}}
\label{proof:thm3_3}

We begin with the statement of several lemmas that we use in the following proof. 
We first define the matrices
  $$
  \begin{aligned}
  	&\mathbf{M}^{(i)} = \muu^{\top} \me^{(i)} \mv \text{ for any }i\in[n], \\
  &\mathbf{N}^{(ij)} = \muu^\top
  \mathbf{E}^{(i)}\mathbf{E}^{(j)\top} \muu,
  \tilde{\mathbf{N}}^{(ij)} =  \mv^\top \mathbf{E}^{(i)\top}
  	  \mathbf{E}^{(j)}\mv \text{ for any }i,j\in[n],
  \end{aligned}
  $$ 
and let $\vartheta_n = \max\{1,d^{1/2}\rho_n^{1/2}(\log n)^{1/2}\}$.

\begin{lemma1}
  \label{lemma:R V^T Q_v,R^T U^T Q_u}
  Consider the setting in Theorem~\ref{thm:UhatW-U=EVR^{-1}+xxx}.  
      We then have
  $$
  \begin{aligned}
  	  \mv^\top \mq_\mv
  	  =&-\frac{1}{m} \sum_{j=1}^{m} \mathbf{M}^{(j)\top}(\mathbf{R}^{(j)\top})^{-1}
  	  -\frac{1}{2m^2}\sum_{j=1}^m\sum_{k=1}^m
  	  \bigl(\mr^{(j)}\bigr)^{-1}\mathbf{N}^{(jk)}
      \bigl(\mr^{(k)\top}\bigr)^{-1} +O_p((n\rho_n)^{-3/2}\vartheta_n),\\
  	  \muu^\top \mq_\muu
  	  =&-\frac{1}{m} \sum_{j=1}^{m} \mathbf{M}^{(j)} (\mathbf{R}^{(j)})^{-1} 
	-\frac{1}{2m^2}\sum_{j=1}^m\sum_{k=1}^m
  	 \bigl(\mr^{(j)\top}\bigr)^{-1}\tilde{\mathbf{N}}^{(jk)}\bigl(\mr^{(k)}\bigr)^{-1} +O_p((n\rho_n)^{-3/2}\vartheta_n).
  \end{aligned}
  $$
\end{lemma1}

\begin{lemma1}
  \label{lemma:U^T EE^T U R^-1->}
  Consider the setting in Theorem~\ref{thm:UhatW-U=EVR^{-1}+xxx}. For any $i \in [m]$, let
 $\mathbf{F}^{(i)}$ be the $d \times d$ matrix defined by
  $$
\begin{aligned}
\mathbf{F}^{(i)} &=\frac{1}{m}\sum_{j=1}^m \mathbf{N}^{(ij)} (\mr^{(j)\top})^{-1}
+\frac{1}{m}\sum_{j=1}^m(\mr^{(j)\top})^{-1}\tilde{\mathbf{N}}^{(ji)}
\\ &-\frac{1}{2m^2} \sum_{j=1}^{m} \sum_{k=1}^{m} \mathbf{R}^{(i)}(\mathbf{R}^{(j)})^{-1} {\mathbf{N}}^{(jk)} (\mathbf{R}^{(k) \top})^{-1}
-\frac{1}{2m^2} \sum_{j=1}^{m} \sum_{k=1}^{m} (\mathbf{R}^{(j) \top})^{-1} \tilde{\mathbf{N}}^{(jk)}(\mathbf{R}^{(k)})^{-1}\mathbf{R}^{(i)}.
\end{aligned}
$$
  We then have, for any $i\in[m]$, 
  $$\rho_n^{-1/2}\bigl(\mathrm{vec}(\mathbf{F}^{(i)}) -
  \bm{\mu}^{(i)}\bigr) \stackrel{p}{\longrightarrow} \bm{0}$$
  as $n\rightarrow \infty$, where $\bm{\mu}^{(i)}$ is defined in the statement of Theorem~\ref{thm:What_U Rhat
    What_v^T-R->norm}.
\end{lemma1}

\begin{lemma1}
  \label{lemma:U^T E V->norm}
  Consider the setting in Theorem~\ref{thm:UhatW-U=EVR^{-1}+xxx}. Then for any $i \in [m]$, we have
  $$
  (\boldsymbol{\Sigma}^{(i)})^{-1/2}\operatorname{vec}\big(\mathbf{U}^{\top} \mathbf{E}^{(i)} \mathbf{V}\big) 
  \leadsto \mathcal{N}\big(\bm{0},\mi\big)
  $$
  as $n\rightarrow \infty$, where $\boldsymbol{\Sigma}^{(i)}$ is defined in the statement of Theorem~\ref{thm:What_U Rhat
    What_v^T-R->norm}.
\end{lemma1}

The proofs of Lemma~\ref{lemma:R V^T Q_v,R^T U^T Q_u} through
Lemma~\ref{lemma:U^T E V->norm} are presented in
Section~\ref{sec:technical_lemmas2} in the supplementary material. We now proceed with
the proof of Theorem~\ref{thm:What_U Rhat What_v^T-R->norm}.
Recall that $\hat\mr^{(i)}=\hat\muu^\top\ma^{(i)}\hat\mv$ and let $\zeta^{\star} = \mw_\muu^\top\hat\mr^{(i)}\mw_\mv$. 
Then, by Theorem~\ref{thm:UhatW-U=EVR^{-1}+xxx}, we have with high probability the following decomposition
for $\zeta^\star$
\begin{equation}
  \label{eq:What_U Rhat What_v^T-R}
\begin{split}
\zeta^{\star} &=\mw_\muu^\top\hat\muu^\top\ma^{(i)}\hat\mv\mw_\mv\\
&=(\mw_\muu^\top\hat\muu^\top - \muu^\top + \muu^\top)\ma^{(i)}(\hat\mv\mw_\mv -
\mv + \mv)\\
&=\muu^\top \ma^{(i)}\mv
+\muu^\top \ma^{(i)}\frac{1}{m}\sum_{k=1}^m\me^{(k)\top}\muu(\mr^{(k)\top})^{-1}
+\muu^\top \ma^{(i)}\mq_\mv
\\ &
+\frac{1}{m}\sum_{j=1}^m(\mr^{(j)\top})^{-1}\mv^\top\me^{(j)\top}\ma^{(i)}\mv
+\frac{1}{m}\sum_{j=1}^m(\mr^{(j)\top})^{-1}\mv^\top\me^{(j)\top}\ma^{(i)}\mq_\mv \\
&+\frac{1}{m}\sum_{j=1}^m(\mr^{(j)\top})^{-1}\mv^\top\me^{(j)\top}\ma^{(i)}\frac{1}{m}\sum_{k=1}^m\me^{(k)\top}\muu(\mr^{(k)\top})^{-1}\\
&+\mq_\muu^\top \ma^{(i)}\mv
+\mq_\muu^\top \ma^{(i)}\frac{1}{m}\sum_{k=1}^m\me^{(k)\top}\muu(\mr^{(k)\top})^{-1}
+\mq_\muu^\top \ma^{(i)}\mq_\mv.
\end{split}
\end{equation}
We now analyze each of the nine terms on the right
hand side of Eq. (\ref{eq:What_U Rhat What_v^T-R}). Note that we
always expand $\ma^{(i)}$ as $\ma^{(i)} = \mpp^{(i)} + \me^{(i)}$.
In the following proof, for any matrix $\mm$, we write $\mm=O_p(a_n)$ to denote $\|\mm\|=O_p(a_n)$.

Let $\zeta_1=\muu^\top \ma^{(i)}\mv$. We then have
\begin{equation}
  \label{eq:What_U Rhat What_v^T-R term1}
\begin{split}
\zeta_1 = \mr^{(i)}+
\muu^\top \me^{(i)} \mv.
\end{split}
\end{equation}

Let $\zeta_2 = \muu^\top
\ma^{(i)}\frac{1}{m}\sum_{k=1}^m\me^{(k)\top}\muu(\mr^{(k)\top})^{-1}$. We
then have
\begin{equation}
  \label{eq:What_U Rhat What_v^T-R term2}
\begin{split}
\zeta_2 =\frac{1}{m}\sum_{k=1}^m \mr^{(i)} \mathbf{M}^{(k)\top}
(\mr^{(k)\top})^{-1} +\frac{1}{m}\sum_{k=1}^m \mathbf{N}^{(ik)}(\mr^{(k)\top})^{-1}.
\end{split}
\end{equation}

Let $\zeta_ 3 = \muu^\top \ma^{(i)}\mq_\mv = \muu^\top
\big(\mpp^{(i)}+\me^{(i)}\big) \mq_\mv$. 
Using Lemma~\ref{lemma:R V^T Q_v,R^T U^T Q_u}, we obtain
\begin{equation}
  \label{eq:What_U Rhat What_v^T-R term3}
\begin{split}
\zeta_3 &=
\mr^{(i)}\mv^\top \mq_\mv
+\muu^\top \me^{(i)} \mq_\mv\\
&= -\frac{1}{m} \sum_{j=1}^{m} \mathbf{R}^{(i)} \mathbf{M}^{(j)\top}
(\mathbf{R}^{(j) \top})^{-1} - \frac{1}{2m^2} \sum_{j=1}^{m} \sum_{k=1}^{m} \mathbf{R}^{(i)}(\mathbf{R}^{(j)})^{-1} \mathbf{N}^{(jk)} (\mathbf{R}^{(k) \top})^{-1} 
+O_p((n\rho_n)^{-1/2}\vartheta_n),
\end{split}
\end{equation}
where the last equality follows from combining Lemma~\ref{lemma:|E|_2,|UEV|F} and Theorem~\ref{thm:UhatW-U=EVR^{-1}+xxx} to bound
$$
\begin{aligned}
    &\|\mr^{(i)}\| \times O_p((n\rho_n)^{-3/2}\vartheta_n)
    \lesssim  (n\rho_n)^{-1/2}\vartheta_n,
    \\
	&\|\muu^\top \me^{(i)} \mq_\mv\|
	\leq \|\me^{(i)}\|\cdot\|\mq_\mv\|
	\lesssim (n\rho_n)^{-1/2}\vartheta_n
  \end{aligned}
$$
with high probability.

Next let $\zeta_4 =
\tfrac{1}{m}\sum_{j=1}^m(\mr^{(j)\top})^{-1}\mv^\top\me^{(j)\top}\ma^{(i)}\mv$. 
We then have
\begin{equation}
  \label{eq:What_U Rhat What_v^T-R term4}
\begin{split}
  \zeta_4  
	=&\frac{1}{m}\sum_{j=1}^m(\mr^{(j)\top})^{-1}\mathbf{M}^{(j)\top}\mr^{(i)}
	+\frac{1}{m}\sum_{j=1}^m(\mr^{(j)\top})^{-1}\tilde{\mathbf{N}}^{(ji)}.
\end{split}
\end{equation}

Now let $\zeta_5 =
\tfrac{1}{m}\sum_{j=1}^m(\mr^{(j)\top})^{-1}\mv^\top\me^{(j)\top}\ma^{(i)}\mq_\mv$. We
then have
\begin{equation}
  \label{eq:What_U Rhat What_v^T-R term5}
\begin{aligned}
	\zeta_5 &= \frac{1}{m}\sum_{j=1}^m(\mr^{(j)\top})^{-1}\mathbf{M}^{(j)\top}\mr^{(i)}\mv^\top\mq_\mv
	+\frac{1}{m}\sum_{j=1}^m(\mr^{(j)\top})^{-1}\mv^\top\me^{(j)\top}\me^{(i)}\mq_\mv
	\\ &
	=O_p((n\rho_n)^{-1}\vartheta_n^2),
\end{aligned}
\end{equation}
where the final bound in Eq.~\eqref{eq:What_U Rhat What_v^T-R term5}
follows from
Lemma~\ref{lemma:|E|_2,|UEV|F} and Theorem~\ref{thm:UhatW-U=EVR^{-1}+xxx}, i.e.,
$$
\begin{aligned}
	&\|(\mr^{(j)\top})^{-1} \mathbf{M}^{(j)\top} \mr^{(i)} \mv^{\top} \mq_{\mv}\| 
	\leq \|(\mr^{(j)})^{-1}\| \cdot \|\mathbf{M}^{(j)}\| \cdot \|\mr^{(i)}\| \cdot \|\mq_{\mv}\| 
	\lesssim d^{1/2}n^{-1/2}(n\rho_n)^{-1/2} (\log n)^{1/2}\vartheta_n, \\
&\|(\mr^{(j)\top})^{-1}\mv^\top\me^{(j)\top}\me^{(i)}\mq_\mv\| \leq
\|(\mr^{(j)})^{-1}\| \cdot \|\me^{(j)}\| \cdot \|\me^{(i)}\| \cdot
\|\mq_{\mv}\| \lesssim (n \rho_n)^{-1}\vartheta_n \end{aligned}
$$
with high probability.

Let $\zeta_6 =
\tfrac{1}{m}\sum_{j=1}^m(\mr^{(j)\top})^{-1}\mv^\top\me^{(j)\top}\ma^{(i)}\frac{1}{m}\sum_{k=1}^m\me^{(k)\top}\muu(\mr^{(k)\top})^{-1}$. We
then have
\begin{equation}
  \label{eq:What_U Rhat What_v^T-R term6}
\begin{aligned}
\zeta_6
&=\frac{1}{m^2}\sum_{j=1}^m\sum_{k=1}^m(\mr^{(j)\top})^{-1}\mathbf{M}^{(j)\top}\mr^{(i)}\mathbf{M}^{(k)\top}(\mr^{(k)\top})^{-1}
+\frac{1}{m^2}\sum_{j=1}^m\sum_{k=1}^m(\mr^{(j)\top})^{-1}\mv^\top\me^{(j)\top}\me^{(i)}\me^{(k)\top}\muu(\mr^{(k)\top})^{-1}
\\&
=O_p((n\rho_n)^{-1/2}),
\end{aligned}
\end{equation}
where the final bound in Eq.~\eqref{eq:What_U Rhat What_v^T-R term6} follows 
from Lemma~\ref{lemma:|E|_2,|UEV|F}, i.e.,
$$
\begin{aligned}
     \|(\mr^{(j)\top})^{-1} \mv^{\top} \me^{(j)\top} \me^{(i)} \me^{(k)\top} (\mr^{(k)\top})^{-1}\| 
     &\leq \|(\mr^{(j)})^{-1}\|  \cdot \|\me^{(j)}\| \cdot \|\me^{(i)}\| \cdot \|\me^{(k)}\| \cdot \|(\mr^{(k)})^{-1}\|
     \\ & \lesssim (n \rho_n)^{-1/2}
\end{aligned}
$$
with high probability.

Let $\zeta_7 = \mq_\muu^\top \ma^{(i)}\mv$. From Lemma \ref{lemma:R
  V^T Q_v,R^T U^T Q_u} we have
\begin{equation}
  \label{eq:What_U Rhat What_v^T-R term7}
\begin{split}
  \zeta_7 &=\mq_\muu^\top \muu\mr^{(i)}+\mq_\muu^\top \me^{(i)}\mv\\
	&=-\frac{1}{m}\sum_{j=1}^m (\mathbf{R}^{(j)\top})^{-1} \mathbf{M}^{(j)\top} \mr^{(i)} 
	-\frac{1}{2m^2} \sum_{j=1}^{m} \sum_{k=1}^{m} (\mathbf{R}^{(j) \top})^{-1} \tilde{\mathbf{N}}^{(jk)} (\mathbf{R}^{(k)})^{-1}\mathbf{R}^{(i)}
	 +O_p((n\rho_n)^{-1/2}\vartheta_n),
\end{split}
\end{equation}
where the last equality follows from Lemma~\ref{lemma:|E|_2,|UEV|F} and Theorem~\ref{thm:UhatW-U=EVR^{-1}+xxx}, i.e., 
$$
\begin{aligned}
&\|\mr^{(i)}\| \times O_p((n\rho_n)^{-3/2}\vartheta_n)
\lesssim (n\rho_n)^{-1/2}\vartheta_n,
\\
	&\|\mq_\muu^\top \me^{(i)}\mv\|
	\leq \|\mq_\muu\|\cdot \|\me^{(i)}\|
	\lesssim (n\rho_n)^{-1/2} \vartheta_n
\end{aligned}
$$
with high probability.

Now let $\zeta_8 = \mq_\muu^\top
\ma^{(i)}\frac{1}{m}\sum_{k=1}^m\me^{(k)\top}\muu(\mr^{(k)\top})^{-1}$. We
then have
\begin{equation}
  \label{eq:What_U Rhat What_v^T-R term8}
\begin{aligned}
\zeta_8 &=\frac{1}{m}\sum_{k=1}^m\mq_\muu^\top \muu\mr^{(i)} \mathbf{M}^{(k)\top} (\mr^{(k)\top})^{-1}
	+\frac{1}{m}\sum_{k=1}^m\mq_\muu^\top \me^{(i)}\me^{(k)\top}\muu(\mr^{(k)\top})^{-1}
	\\ &=O_p((n\rho_n)^{-1}\vartheta_n^2),
\end{aligned}
\end{equation}
where the last bound follows from Lemma~\ref{lemma:|E|_2,|UEV|F} and Theorem~\ref{thm:UhatW-U=EVR^{-1}+xxx}, i.e.,
$$
\begin{aligned}
	&\|\mq_\muu^\top \muu\mr^{(i)} \mathbf{M}^{(k)\top}
    (\mr^{(k)\top})^{-1}\| \leq \|\mq_{\muu}\| \cdot \|\mr^{(i)}\|
    \cdot \|\mathbf{M}^{(k)}\| \cdot \|(\mr^{(k)})^{-1}\| \lesssim
    d^{1/2}(n\rho_n)^{-1} (\rho_n\log n)^{1/2}\vartheta_n, \\
&\|\mq_\muu^\top \me^{(i)}\me^{(k)\top}\muu(\mr^{(k)\top})^{-1}\| \leq
\|\mq_{\muu}\| \cdot \|\me^{(i)}\| \cdot \|\me^{(k)}\| \cdot
\|(\mr^{(k)})^{-1} \| \lesssim (n \rho_n)^{-1}\vartheta_n \end{aligned}
$$
with high probability.

Finally, let $\zeta_9=\mq_\muu^\top \ma^{(i)}\mq_\mv$, we once again have from Lemma~\ref{lemma:|E|_2,|UEV|F} and Theorem~\ref{thm:UhatW-U=EVR^{-1}+xxx} that
\begin{equation}
  \label{eq:What_U Rhat What_v^T-R term9}
\begin{split}
\zeta_9=\mq_\muu^\top \muu \mr^{(i)}\mv^\top\mq_\mv+\mq_\muu^\top \me^{(i)}\mq_\mv=O_p((n\rho_n)^{-1}\vartheta_n^2).
\end{split}
\end{equation}

Combining Eq.~(\ref{eq:What_U Rhat What_v^T-R}) through
Eq.~(\ref{eq:What_U Rhat What_v^T-R term9}) and noting that one 
term in $\zeta_2$ cancels out another term in $\zeta_3$ while one term
in $\zeta_4$ cancels out another term in $\zeta_7$, 
we obtain 
\begin{equation}
\label{eq:W_U^T hatR W_V-R=...}
	\begin{aligned}
	 \mw_\muu^\top\hat\mr^{(i)}\mw_\mv -\mr^{(i)}
     &=\muu^\top \me^{(i)} \mv + \mathbf{F}^{(i)} + O_p((n\rho_n)^{-1/2}\vartheta_n),
\end{aligned}
\end{equation}
where $\mathbf{F}^{(i)}$ is defined in the statement of Lemma~\ref{lemma:U^T EE^T U R^-1->}. 
We then show in Lemma \ref{lemma:U^T EE^T U R^-1->} that
$$
\begin{aligned}
\rho_n^{-1/2}(\operatorname{vec}(\mathbf{F}^{(i)}) -
\bm{\mu}^{(i)}\bigr) \stackrel{p}{\longrightarrow} \bm{0}.
\end{aligned}
$$
In addition we also show in  Lemma \ref{lemma:U^T E V->norm} that
\begin{equation}
	\label{eq:Sigma d->}
	(\boldsymbol{\Sigma}^{(i)})^{-1/2}\operatorname{vec}(\mathbf{U}^{\top} \mathbf{E}^{(i)} \mathbf{V}) 
    \leadsto \mathcal{N}\big(\bm{0},\mi\big).\\
\end{equation}
From the assumption $\sigma_{\min}(\mSigma^{(i)})\gtrsim \rho_n$, we have $\|(\boldsymbol{\Sigma}^{(i)})^{-1/2}\|\lesssim \rho_n^{-1/2}$,
hence
\begin{gather}
    \label{eq:F p->}
    (\boldsymbol{\Sigma}^{(i)})^{-1/2}\bigl(\operatorname{vec}(\mathbf{F}^{(i)})
  	  - \bm{\mu}^{(i)}\bigr) \stackrel{p}{\longrightarrow} \bm{0}.
\end{gather}
Finally, because we assume $n\rho_n=\omega(n^{1/2})$, we have
\begin{equation}
	\label{eq:Op p->}
  (\boldsymbol{\Sigma}^{(i)})^{-1/2} O_p((n\rho_n)^{-1/2}\vartheta_n)
  \stackrel{p}{\longrightarrow} \bm{0}.
\end{equation}
Combining Eq.~\eqref{eq:W_U^T hatR W_V-R=...} through Eq.~\eqref{eq:Op p->}, and applying
Slutsky's theorem, we obtain
$$
\begin{aligned}
	\bigl(\boldsymbol{\Sigma}^{(i)}\bigr)^{-1/2}\Bigl(\operatorname{vec}\big(\mw_\muu^\top\hat\mr^{(i)}\mw_\mv-\mr^{(i)}\big)
    - \bm{\mu}^{(i)} \Bigr)     \leadsto \mathcal{N}\big(\bm{0},\mi\big)
\end{aligned}
$$
as $n \rightarrow \infty$. 
Finally $\me^{(i)}$ is independent of $\me^{(j)}$ for $i
\not = j$, and hence 
$\mw_\muu^\top\hat\mr^{(i)}\mw_\mv$ and
$\mw_\muu^\top\hat\mr^{(j)}\mw_\mv$
are asymptotically independent for any $i \not = j.$
\hspace*{\fill} \qedsymbol

\subsection{Proof of Theorem~\ref{thm:HT}}
\label{sec:proof_thmHT}
We first consider $\mathbb{H}_0 \colon \mr^{(i)} =
\mr^{(j)}$ versus $\mathbb{H}_{A} \colon \mr^{(i)} \not =
\mr^{(j)}$ for some $i \not = j$. Define 
\begin{equation*}
\zeta_{ij} = \mathrm{vec}^\top\big( \hat{\mathbf{R}}^{(i)}-
\hat{\mathbf{R}}^{(j)} \big)
(\mw_\mv\otimes\mw_\muu)(\mathbf{\Sigma}^{(i)}+\mathbf{\Sigma}^{(j)})^{-1}(\mw_\mv\otimes\mw_\muu)^\top\operatorname{vec}\big(
\hat{\mathbf{R}}^{(i)}- \hat{\mathbf{R}}^{(j)} \big).
\end{equation*}
Now suppose $\mr^{(i)} =
\mr^{(j)}$. Then $\zeta_{ij}
\leadsto \chi^{2}_{d^2}$; see
Eq.~\eqref{eq:HT:initial}. As $d$ is finite, we conclude that $\zeta_{ij}$ is bounded in
probability. 
Now $\|\mSigma^{(i)}+\mSigma^{(j)}\|\leq
\|\mSigma^{(i)}\|+\|\mSigma^{(j)}\|\lesssim \rho_n$, and hence, by the assumption $\sigma_{\min}(\mathbf{\Sigma}^{(i)}+\mathbf{\Sigma}^{(j)})\asymp \rho_n$, we have $\sigma_r\big((\mathbf{\Sigma}^{(i)}+\mathbf{\Sigma}^{(j)})^{-1}\big)\asymp
\rho_n^{-1}$ for any $r\in[d^2]$. We thus have
%
$\zeta_{ij} \asymp \rho_n^{-1} \| \hat{\mathbf{R}}^{(i)}-
\hat{\mathbf{R}}^{(j)}\|_{F}^2$, 
i.e., $\rho_n^{-1/2}\|\hat{\mathbf{R}}^{(i)} - \hat{\mathbf{R}}^{(j)}\|_{F}$ is
bounded in probability.

Let $\mathbf{W}_* = \mathbf{W}_{\mv} \otimes \mathbf{W}_{\muu}$.
Then by Lemma~\ref{lemma:|sigma inverse-sigma|,order of sigma}, we have
$$
\big\|\mw_*(\mSigma^{(i)}+\mSigma^{(j)})^{-1}\mw_*^\top
		-\big(\hat\mSigma^{(i)}+\hat\mSigma^{(j)}\big)^{-1}\big\|
		\lesssim d (n\rho_n)^{-1/2}(\log n)^{1/2}\times\rho_n^{-1}
$$
with high probability. Now recall the definition of 
$T_{ij}$ in Theorem~\ref{thm:HT}. Under the assumption $n\rho_n=\omega(\log n)$, we then have
\begin{equation}
\begin{aligned}
  \label{eq:zeta_T}
	|\zeta_{ij} - T_{ij}| &\leq \big\|\mw_*(\mSigma^{(i)}+\mSigma^{(j)})^{-1}\mw_*^\top
		-\big(\hat\mSigma^{(i)}+\hat\mSigma^{(j)}\big)^{-1}\big\|
		\cdot \|\hat{\mathbf{R}}^{(i)}- \hat{\mathbf{R}}^{(j)}\|_F^2\\
	&\lesssim  (d (n\rho_n)^{-1/2}(\log n)^{1/2}) \cdot (\rho_n^{-1}\| \hat{\mathbf{R}}^{(i)}- \hat{\mathbf{R}}^{(j)} \|_F^2) \stackrel{p}{\rightarrow}0.
\end{aligned}
\end{equation}
Therefore, by Slutsky's theorem, we have $T_{ij}\leadsto\chi^2_{d^2}$
under $\mathbb{H}_0$. 

We now consider the case where $\mr^{(i)} \not = \mr^{(j)}$ satisfies a
local alternative hypothesis, i.e., 
\begin{gather}
\label{eq:HT:local condition2}
	 \operatorname{vec}^\top(\mathbf{R}^{(i)}-\mathbf{R}^{(j)})
	    (\mathbf{\Sigma}^{(i)}+\mathbf{\Sigma}^{(j)})^{-1}
	    \operatorname{vec}(\mathbf{R}^{(i)}-\mathbf{R}^{(j)})
	    {\rightarrow}\eta
\end{gather}
for some finite $\eta>0$. 
As $\|\big(\mSigma^{(i)}+\mSigma^{(j)}\big)^{-1}\|\asymp \rho_n^{-1}$, 
Eq.~\eqref{eq:HT:local condition2} implies $\rho_n^{-1/2}\|\mathbf{R}^{(i)}-\mathbf{R}^{(j)}\|$
is bounded and hence $n^{2}\rho_n^{3/2}
\|(\mathbf{R}^{(i)})^{-1}-(\mathbf{R}^{(j)})^{-1}\|$ is also bounded.
Indeed, by Assumption~\ref{ass:main} we have 
$ \|(\mr^{(i)})^{-1}\|\asymp (n\rho_n)^{-1}$ for all $i$ and thus
\begin{equation*}
  \label{eq:HT:tmp2}
	\begin{aligned}
		n^2\rho_n^{3/2} \|(\mathbf{R}^{(i)})^{-1}-(\mathbf{R}^{(j)})^{-1}\|
		\leq &(n\rho_n)^2\rho_n^{-1/2} \|(\mr^{(i)})^{-1}\|\cdot\|\mr^{(i)} - \mr^{(j)}\|\cdot\|(\mr^{(j)})^{-1}\|\\
		\lesssim & \rho_n^{-1/2} \|\mr^{(i)} - \mr^{(j)}\|.
	\end{aligned}
\end{equation*}
Now recall the expression for $\bm\mu^{(i)}$ and $\bm{\mu}^{(j)}$
given in Theorem~\ref{thm:What_U Rhat What_v^T-R->norm}. Then 
by Lemma~\ref{lemma:mui-muj}
we have
$$
\begin{aligned}
	\|\bm\mu^{(i)}-\bm\mu^{(j)}\|
	\lesssim d^{1/2}m^{-1}\big(n\rho_n\|(\mr^{(i)})^{-1}-(\mr^{(j)})^{-1}\|+d(n\rho_n)^{-1}\|\mr^{(i)}-\mr^{(j)}\|\big).
\end{aligned}
$$
We therefore have $n\rho_n^{1/2}\|\bm\mu^{(i)}-\bm\mu^{(j)}\|$ is bounded. 
Next define $\xi_{ij}$ and $\tilde{\xi}_{ij}$ by
\begin{equation}
\label{eq:HT:local7}
  \begin{aligned}
  \xi_{ij} &= (\mathbf{\Sigma}^{(i)}+\mathbf{\Sigma}^{(j)})^{-1/2} \, \mathrm{vec}(\mathbf{R}^{(i)}-\mathbf{R}^{(j)}), \\
\tilde{\xi}_{ij} &= (\mathbf{\Sigma}^{(i)}+\mathbf{\Sigma}^{(j)})^{-1/2} \bigl(\mathrm{vec}(\mathbf{R}^{(i)}-\mathbf{R}^{(j)})+\boldsymbol{\mu}^{(i)}-\boldsymbol{\mu}^{(j)}\bigr).
    \end{aligned}
\end{equation}
We then have 
\begin{equation*}
	\begin{aligned}
      \|\xi_{ij} - \tilde{\xi}_{ij}\|
	&\lesssim  \rho_n^{-1/2} \|\bm{\mu}^{(i)} - \bm{\mu}^{(j)}\| = (n \rho_n)^{-1} n \rho_n^{1/2} \|\bm\mu^{(i)} - \bm\mu^{(j)}\|
    {\rightarrow}0.
\end{aligned}
\end{equation*}
Since $\|\xi_{ij}\|^2 
{\rightarrow} \eta$, we have $\|\tilde{\xi}_{ij}\|^2 
{\rightarrow} \eta$.
Now recall Theorem~\ref{thm:What_U Rhat What_v^T-R->norm}. In
particular we have
\begin{equation*}
\begin{aligned}
    (\bm{\Sigma}^{(i)} + \bm{\Sigma}^{(j)})^{-1/2}\mw_*^\top\mathrm{vec}(\hat\mr^{(i)}-\hat\mr^{(j)})-\tilde{\xi}_{ij}
    \rightsquigarrow \mathcal{N}\bigl(\bm{0}, \mathbf{I}\bigr).
\end{aligned}
\end{equation*}
We conclude that
$\zeta_{ij} \rightsquigarrow \chi^{2}_{d^2}(\eta)$, where $\zeta_{ij}$ is defined at the beginning of the current proof. As $\eta$ is
finite, $(\bm{\Sigma}^{(i)} +
\bm{\Sigma}^{(j)})^{-1/2}\mw_*^\top \mathrm{vec}(\hat\mr^{(i)}-\hat\mr^{(j)})$
is also bounded in
probability. Finally, using the same argument as that for
deriving Eq.~\eqref{eq:zeta_T} under $\mathbb{H}_0$, we also have
$\zeta_{ij} - T_{ij} \stackrel{p}{\rightarrow}0$ under the local
alternative in Eq.~\eqref{eq:HT:local condition2} and hence 
$T_{ij} \rightsquigarrow \chi^2_{d^2}(\eta)$ as desired. 
We next consider $\mathbb{H}_0 \colon \mr^{(1)}=\cdots=\mr^{(m)}$ versus
$\mathbb{H}_A \colon \mr^{(i)} \not = \mr^{(j)}$ for some $i
\not =j$. Define
\begin{equation*}
  \begin{split}
\zeta &= \sum_{i=1}^m\mathrm{vec}^\top\big( \hat{\mathbf{R}}^{(i)}-
\bar{\hat{\mathbf{R}}} \big)
(\mw_\mv\otimes\mw_\muu)
{\bar\mSigma}^{-1}
(\mw_\mv\otimes\mw_\muu)^\top\operatorname{vec}\big(
        \hat{\mathbf{R}}^{(i)}- \bar{\hat{\mathbf{R}}} \big) \\
    &= \sum_{i=1}^m \mathrm{vec}^\top\big(  \mw_{\muu}^{\top} (\hat{\mathbf{R}}^{(i)}-
\bar{\hat{\mathbf{R}}} ) \mw_{\mv} \big)
\bar{\bm{\Sigma}}^{-1/2}\bar{\bm{\Sigma}}^{-1/2}\operatorname{vec}\big(  \mw_{\muu}^{\top}
      (\hat{\mathbf{R}}^{(i)}- \bar{\hat{\mathbf{R}}}) \mw_{\mv} \big).
  \end{split}
\end{equation*}
Now suppose $\mr^{(1)} = \dots = \mr^{(m)}$. 
Then $\bm\mu^{(i)}=\cdots=\bm\mu^{(m)}$ and $\bm{\Sigma}^{(1)}
= \dots = \bm{\Sigma}^{(m)}$.
Hence, by Theorem~\ref{thm:What_U Rhat What_v^T-R->norm}, 
$(\bm{\Sigma}^{(i)})^{-1/2} \bigl(\operatorname{vec}(\mw_\muu^\top
\hat\mr^{(i)}\mw_\mv-\mr^{(i)}) - \bm{\mu}^{(i)}\bigr)
\rightsquigarrow \mathcal{N}(\mathbf{0}, \mi)$ for all $i \in [m]$ and furthermore
$
\mw_\muu^\top \hat\mr^{(1)}\mw_\mv$,
$\cdots$ , $
\mw_\muu^\top \hat\mr^{(m)}\mw_\mv$ are asymptotically independent.
We let $\mathbf{Y}$ be the $d^2 \times m$ matrix with
{\em columns} $\{(\bm{\Sigma}^{(i)})^{-1/2} \bigl(\operatorname{vec}(\mw_\muu^\top
\hat\mr^{(i)}\mw_\mv-\mr^{(i)}) - \bm{\mu}^{(i)}\bigr)\}$. As $n
\rightarrow \infty$, $\mathbf{Y}$ converges in distribution to a $d^{2} \times
m$ matrix whose entries are iid $\mathcal{N}(0,1)$ random
variables. 
We therefore have
\begin{equation*}
  \begin{split}
\zeta &= \mathrm{tr}\,\Big[ \Bigl(\mathbf{I} - \frac{\bm{1}
      \bm{1}^{\top}}{m}\Bigr)\mathbf{Y}^{\top} \mathbf{Y} \Bigl(\mathbf{I} - \frac{\bm{1}
      \bm{1}^{\top}}{m}\Bigr)\Big] = \mathrm{tr}\,\Big[ \mathbf{Y} \Bigl(\mathbf{I} - \frac{\bm{1}
        \bm{1}^{\top}}{m}\Bigr) \mathbf{Y}^{\top}\Big].
  \end{split}
\end{equation*}
Then by Corollary~3 in \cite{matrix_quadratic}, 
$\mathbf{Y} (\mathbf{I} - \bm{1} \bm{1}^{\top}/m) \mathbf{Y}^{\top}$
converges in distribution to a $d^2 \times d^2$ Wishart random matrix with $m-1$
degrees of freedom and scale matrix $\mathbf{I}$. Therefore, 
$$\zeta = \mathrm{tr}\,\Big[ \mathbf{Y} \Bigl(\mathbf{I} - \frac{\bm{1} \bm{1}^{\top}}{m}\Bigr) \mathbf{Y}^{\top}\Big] \rightsquigarrow
    \mathrm{tr} \,\Bigl[\sum_{i=1}^{m-1} \bm{g}_i \bm{g}_i^{\top} \Bigr]
    =     \sum_{i=1}^{m-1} \bm{g}_i^{\top} \bm{g}_i \sim \chi^2_{(m-1)d^2},$$
    where $\{\bm{g}_i\}$ are iid $\mathcal{N}(\bm{0}, \mathbf{I})$.
    Under the assumption $\sigma_{\min}(\mSigma^{(i)})\asymp \rho_n$ for all $i\in[m]$, by Weyl's inequality we have $\sigma_{r}(\bar\mSigma^{-1})\asymp \rho_n^{-1}$ for all $r\in[d^2]$. Then we can now follow the same argument as that for deriving Eq.~\eqref{eq:zeta_T}
and show that $T - \zeta \stackrel{p}{\rightarrow} 0$ 
 and hence 
$T \rightsquigarrow \chi^2_{(m-1)d^2}$ under $\mathbb{H}_0$. 
We now consider the case where $\mr^{(i)} \not = \mr^{(j)}$  for some
$i\neq j$. Suppose $\{\mr^{(i)}\}$ satisfy
\begin{equation*}
	\begin{aligned}
	 &\sum_{i=1}^m\operatorname{vec}^\top(\mathbf{R}^{(i)}-\bar{\mathbf{R}})
	    (\bar{\mathbf{\Sigma}})^{-1}
	    \operatorname{vec}(\mathbf{R}^{(i)}-\bar{\mathbf{R}})
	    {\rightarrow}\eta
	\end{aligned}
\end{equation*}
for some finite $\eta>0$. 
As $\sigma_{r}(\bar\mSigma^{-1})\asymp \rho_n^{-1}$ for all $r\in[d^2]$, $\max_{i \in [m]} \rho_n^{-1/2}\|\mathbf{R}^{(i)}-\bar{\mathbf{R}}\|$
is bounded in probability. By Theorem~\ref{thm:What_U Rhat What_v^T-R->norm},  
 $\max_{i \in [m]} \rho_n^{-1/2}\|\mathrm{vec}\bigl(\mathbf{W}_{\mathbf{U}}^{\top}
        \hat{\mathbf{R}}^{(i)} \mathbf{W}_{\mathbf{V}} -
        \mr^{(i)}\bigr) - \bm{\mu}^{(i)}\|$ is also bounded in probability.
Define
$$
\begin{aligned}
	\theta_i=(\bm{\Sigma}^{(i)})^{-1/2} \Bigl(\mathrm{vec}\bigl(\mathbf{W}_{\mathbf{U}}^{\top}
        \hat{\mathbf{R}}^{(i)} \mathbf{W}_{\mathbf{V}} -
        \mr^{(i)}\bigr) - \bm{\mu}^{(i)} \Bigr),\\
	\tilde\theta_i=(\bar{\bm{\Sigma}})^{-1/2} \Bigl(\mathrm{vec}\bigl(\mathbf{W}_{\mathbf{U}}^{\top}
        \hat{\mathbf{R}}^{(i)} \mathbf{W}_{\mathbf{V}} -
        \mr^{(i)}\bigr) - \bm{\mu}^{(i)} \Bigr).
\end{aligned}
$$
By an similar argument to that for deriving
Lemma~\ref{lemma:mui-muj}, we have 
\begin{gather*}
 \|\bm\mSigma^{(i)}-\bar{\bm\mSigma}\|
	\lesssim dn^{-1}\|\mr^{(i)}-\bar\mr\|, \\
 	\|(\bm{\mSigma}^{(i)})^{-1}-(\bar{\bm\mSigma})^{-1}\|
	\lesssim \|(\bm\mSigma^{(i)})^{-1}\|\cdot\|\bm\mSigma^{(i)}-\bar{\bm\mSigma}\|\cdot\|(\bar{\bm\mSigma})^{-1}\|
	\lesssim dn^{-1}\rho_n^{-2}\|\mr^{(i)}-\bar\mr\|.
\end{gather*}
As $\sigma_{r}(\mSigma^{(i)})\asymp \rho_n$ and $\sigma_{r}(\bar\mSigma)\asymp \rho_n$ for all $r\in[d^2]$, by Weyl's inequality we have $\sigma_{\min}((\bm{\mSigma}^{(i)})^{-1/2}+(\bar{\bm\mSigma})^{-1/2})\asymp \rho_n^{-1/2}$,
and we therefore have (see e.g., Problem X.5.5 in \cite{bhatia2013matrix}),
$$	\|(\bm{\mSigma}^{(i)})^{-1/2}-(\bar{\bm\mSigma})^{-1/2}\| \lesssim
\rho_n^{1/2}
\|(\bm{\mSigma}^{(i)})^{-1}-(\bar{\bm\mSigma})^{-1}\|
\lesssim dn^{-1} \rho_n^{-3/2} \|\mr^{(i)}-\bar\mr\|,$$
and hence
$$
\begin{aligned}
	\|\theta_i-\tilde\theta_i \|
    \leq &\|(\bm{\mSigma}^{(i)})^{-1/2}-(\bar{\bm\mSigma})^{-1/2}\|
    \cdot \|\mathrm{vec}\bigl(\mathbf{W}_{\mathbf{U}}^{\top}
        \hat{\mathbf{R}}^{(i)} \mathbf{W}_{\mathbf{V}} -
        \mr^{(i)}\bigr) - \bm{\mu}^{(i)}\|\\
     \lesssim & d n^{-1/2} (n \rho_n)^{-1/2}
     \cdot (\rho_n^{-1/2}\|\mr^{(i)}-\bar\mr\|)
    \cdot (\rho_n^{-1/2}\|\mathrm{vec}\bigl(\mathbf{W}_{\mathbf{U}}^{\top}
        \hat{\mathbf{R}}^{(i)} \mathbf{W}_{\mathbf{V}} -
        \mr^{(i)}\bigr) - \bm{\mu}^{(i)}\|)
      \stackrel{p}{\rightarrow}   0.
 \end{aligned}
$$
From Theorem~\ref{thm:What_U Rhat What_v^T-R->norm} we have $\theta_i\rightsquigarrow
        \mathcal{N}\bigl(\bm{0}, \mathbf{I} \bigr)$ and hence, by
        Slutsky's theorem, 
$\tilde\theta_i\rightsquigarrow
        \mathcal{N}\bigl(\bm{0}, \mathbf{I} \bigr)$.
Now define
\begin{gather*}
\tilde{\xi}_{i} = (\bar{\mathbf{\Sigma}})^{-1/2}
	    \big(\operatorname{vec}(\mathbf{R}^{(i)}-\bar{\mathbf{R}}
	   ) +\bm\mu^{(i)}-\bar{\bm\mu}\big),
\end{gather*}
where $\bar{\bm\mu}=m^{-1}\sum_{i=1}^m\bm\mu^{(i)}$.
Then, using the same argument as that for controlling
the quantities $\tilde{\xi}_{ij}$ in Eq.~\eqref{eq:HT:local7}, 
we have $\sum_{i=1}^m
	    \|\tilde{\xi}_{i}\|^2
	    {\rightarrow}\eta$. Finally, note that $\zeta$ can be written as
$$\zeta = \mathrm{tr}\:\Big[ (\tilde{\bm{\Theta}} + \tilde{\bm{\Xi}})
\Bigl(\mathbf{I} - \frac{\bm{1}{\bm{1}}^{\top}}{m}\Bigr) (\tilde{\bm{\Theta}} +
\tilde{\bm{\Xi}})^\top\Big],$$
where $\tilde{\bm{\Theta}}$ and $\tilde{\bm{\Xi}}$ are 
$d^2 \times
m$ matrices with columns $\{\tilde{\theta}_i\}$ and $\{(\bar{\mathbf{\Sigma}})^{-1/2}
	    \big(\operatorname{vec}(\mathbf{R}^{(i)}
	   ) +\bm\mu^{(i)}\big)\}$, respectively. 
As $\mathrm{tr}\big[ \tilde{\bm{\Xi}}
\big(\mathbf{I} - \frac{\bm{1}{\bm{1}}^{\top}}{m}\big)
\tilde{\bm{\Xi}}^\top\big]=\sum_{i=1}^m
	    \|\tilde{\xi}_{i}\|^2
	    {\rightarrow}\eta$,  we have by Corollary~3 in \cite{matrix_quadratic}
that $\zeta \rightsquigarrow \chi^2_{(m-1)d^2}(\eta)$ as
$n \rightarrow \infty$. Once again, using the same derivations as that
        for Eq.~\eqref{eq:zeta_T}, we obtain $|\zeta - T| \stackrel{p}{\rightarrow}
        0$ 
        and thus $T \rightsquigarrow
        \chi^2_{(m-1)d^2}(\eta)$, as desired.
\hspace*{\fill} \qedsymbol

\subsection{Proof of Theorem~\ref{thm:PCA_UhatW-U=(I-UU^T)EVR^{-1}+xxx_previous}}
\label{Appendix:A_thm:PCA_UhatW-U=(I-UU^T)EVR^{-1}+xxx_previous}
The proof follows a similar argument to that presented in the proof of Theorem~\ref{thm:UhatW-U=EVR^{-1}+xxx_part2}. We begin with the statement of several important bounds that we use throughout the following derivations.
\begin{lemma1}
  \label{lemma:PCA_|E|_previous}
  Consider the setting in
  Theorem~\ref{thm:PCA_UhatW-U=(I-UU^T)EVR^{-1}+xxx_previous}. For
  $i\in[m]$ 
  let $\me^{(i)}=\hat\mSigma^{(i)}-\mSigma^{(i)}$. Let 
$$r_i = \frac{\mathrm{tr} (\bm{\Sigma}^{(i)})}{\lambda_1} = \frac{1}{\lambda^{(i)}_1}
\Bigl(\sum_{k=1}^{d_i}\lambda^{(i)}_k+(D-d_i)\sigma_i^2\Bigr)
\asymp
D^{1-\gamma}$$
be the effective rank of $\bm{\Sigma}^{(i)}$. We then have
 \begin{equation*}
 \begin{aligned}
& 	\|\me^{(i)}\|
\lesssim D^\gamma \varphi,
\quad \|\me^{(i)}\muu^{(i)}\|_{2\to\infty}
\lesssim d_i^{1/2} D^{\gamma/2} \tilde\varphi,
\quad \|\me^{(i)}\|_\infty
\lesssim D \tilde{\varphi} 
 \end{aligned}
\end{equation*}
with high probability.
And for $i\neq j$, we have
$$
\|\muu^{(i)}(\mi-\muu^{(j)}\muu^{(j)\top})\mathbf{E}^{(j)}\muu^{(j)}\|
\lesssim  D^{\gamma/2} \tilde{\varphi}$$
with high probability.
Here we define
\begin{gather*}
\varphi =
\Bigl(\frac{\max\{r, \log D\}}{n}\Bigr)^{1/2}, 
\quad \tilde\varphi =\Bigl(\frac{\log D}{n}\Bigr)^{1/2}.
\end{gather*} 
Note that $\varphi \leq r^{1/2} \tilde{\varphi} \asymp D^{(1 -
  \gamma)/2} \tilde{\varphi}$. 
Furthermore, under the assumption $n = \omega(\max\{D^{1 - \gamma}, \log
D\})$ in Theorem~\ref{thm:PCA_UhatW-U=(I-UU^T)EVR^{-1}+xxx_previous}
we have $\varphi = o(1)$ and $\tilde{\varphi} =
o(1)$.
\end{lemma1}

We next state an important technical lemma for bounding the error of $\hat{\mathbf{U}}^{(i)}$ as an estimate for the true $\mathbf{U}^{(i)}$ for each $i \in [m]$.

\begin{lemma1}
  \label{lemma:PCA_UhatW-U_previous}
Consider the setting in Theorem~\ref{thm:PCA_UhatW-U=(I-UU^T)EVR^{-1}+xxx_previous}. Fix an $i \in [m]$ and write the eigendecomposition of $\hat{\bm{\Sigma}}^{(i)}$ as $\hat{\bm{\Sigma}}^{(i)} = \hat{\mathbf{U}}^{(i)}\hat{\bm{\Lambda}}^{(i)}(\hat{\mathbf{U}}^{(i)})^\top+ \hat{\mathbf{U}}^{(i)}_\perp\hat{\bm{\Lambda}}^{(i)}_\perp(\hat{\mathbf{U}}^{(i)}_\perp)^\top$. Next define $\mathbf{W}^{(i)}$ as a minimizer of $\|\hat{\mathbf{U}}^{(i)} \mathbf{O} - \mathbf{U}^{(i)}\|_{F}$ over all $d_i \times d_i$ orthogonal matrix $\mathbf{O}$. We then have
\begin{equation*} \label{Eq:PCA_UW-U for each i_previous}
\hat{\mathbf{U}}^{(i)}\mathbf{W}^{(i)}-\mathbf{U}^{(i)}=(\mathbf{I}-\mathbf{U}^{(i)}\mathbf{U}^{(i)\top})(\hat{\bm{\Sigma}}^{(i)} - \bm{\Sigma}^{(i)}) \mathbf{U}^{(i)}(\bm{\Lambda}^{(i)})^{-1}+\mathbf{T}^{(i)},
\end{equation*}
where the residual matrix $\mathbf{T}^{(i)}$ satisfies
$$
\|\mathbf{T}^{(i)}\| 
\lesssim D^{-\gamma} \varphi + \varphi^2
$$
with high probability. Furthermore, if $n=\omega(D^{2-2\gamma}\log D)$, we have
\begin{equation*}
  \begin{split}
\|\mathbf{T}^{(i)}\|_{2 \to \infty} 
 &\lesssim  d_i^{1/2}D^{-3\gamma/2}\tilde{\varphi} (1 + D \tilde{\varphi})
\end{split}
\end{equation*}
with high probability.
\end{lemma1}

The proofs of Lemma~\ref{lemma:PCA_|E|_previous} and
Lemma~\ref{lemma:PCA_UhatW-U_previous} are provided in
Section~\ref{Appendix:proof of lemma:PCA_|E|_previous and
  lemma:PCA_UhatW-U_previous}. 
We now complete the proof of
Theorem~\ref{thm:PCA_UhatW-U=(I-UU^T)EVR^{-1}+xxx_previous}. 
Suppose that the bounds in
Lemma~\ref{lemma:PCA_|E|_previous} and Lemma~\ref{lemma:PCA_UhatW-U_previous} hold.    
We then invoke Theorem~\ref{thm3_general}. 
More specifically, for each $i\in[m]$, by Lemma~\ref{lemma:PCA_UhatW-U_previous} we have
\begin{equation*}
  \hat\muu^{(i)}\mw^{(i)}-\muu^{(i)}=\mt^{(i)}_0+\mt^{(i)},
\end{equation*}
where
$\mt^{(i)}_0=(\mi-\muu^{(i)}\muu^{(i)\top})(\hat{\mSigma}^{(i)} -
\mSigma^{(i)}) \muu^{(i)}(\mLambda^{(i)})^{-1}$ and $\mw^{(i)} \in \mathcal{O}_{d_i}$.
Recall $\me^{(i)}=\hat{\mSigma}^{(i)} -
\mSigma^{(i)}$. Then by Lemma~\ref{lemma:PCA_|E|_previous} we have
$$
\begin{aligned}
	&\|\mt^{(i)}_0\|
\leq\|\me^{(i)}\|\cdot \|(\mLambda^{(i)})^{-1}\|
\lesssim D^\gamma \varphi\cdot (D^\gamma )^{-1}
\lesssim \varphi,\\
&\|\mt^{(i)}_0\|_{2\to\infty}
\leq\|\me^{(i)}\muu^{(i)}\|_{2\to\infty}\cdot \|(\mLambda^{(i)})^{-1}\|
\lesssim d_i^{1/2} D^{\gamma/2} \tilde\varphi\cdot (D^\gamma )^{-1}
\lesssim d_i^{1/2} D^{-\gamma/2} \tilde\varphi
\end{aligned}
$$
with high probability. Thus given the condition $\varphi=o(1)$ and $D=O(1)$ we have
$$
\begin{aligned}
	\max_{i\in[m]}\big(2\|\mt^{(i)}_0\|
	+2\|\mt^{(i)}\|
	+\|\mt^{(i)}_0+\mt^{(i)}\|^2\big)
	\lesssim \varphi
\end{aligned}
$$
with high probability and thus $\max_{i\in[m]}\big(2\|\mt^{(i)}_0\|
	+2\|\mt^{(i)}\|
	+\|\mt^{(i)}_0+\mt^{(i)}\|^2\big)=o_p(1)$.
Under the assumption that
$\|\mathbf\Pi_s\|=\|m^{-1}\sum_{i=1}^m\muu_s^{(i)}\muu_s^{(i)\top}\|\leq 1-c_s$ for some constant $0<c_s\leq 1$, we have $\frac{1}{2}(1-\|\mathbf\Pi_s\|)\geq \frac{c_s}{2}$. Then for large enough $n$ we have 
$$
\begin{aligned}
	\max_{i\in[m]}\big(2\|\mt^{(i)}_0\|
	+2\|\mt^{(i)}\|
	+\|\mt^{(i)}_0+\mt^{(i)}\|^2\big)
	\leq c(1-\|\mathbf\Pi_s\|)
	<\frac{1}{2}(1-\|\mathbf\Pi_s\|)
\end{aligned}
$$
with high probability for any constant $c<\frac{1}{2}$.
Now we have
\begin{equation}
\label{eq:ep_t,ze_t_3}
	\begin{aligned}
		&\epsilon_{\mt_0}
		=\max_{i\in[m]}\|\mt^{(i)}_0\|
		\lesssim \varphi,\\
		&\zeta_{\mt_0}
		=\max_{i\in[m]}\|\mt^{(i)}_0\|_{2\to\infty}
		\lesssim d_{\max}^{1/2} D^{-\gamma/2} \tilde\varphi,\\
		&\epsilon_{\mt}
		=\max_{i\in[m]}\|\mt^{(i)}\|
		\lesssim D^{-\gamma} \varphi + \varphi^2,\\
		&\zeta_{\mt}
		=\max_{i\in[m]}\|\mt^{(i)}\|_{2\to\infty}
		\lesssim d_{\max}^{1/2}D^{-3\gamma/2}\tilde\varphi (1 + D \tilde{\varphi})
	\end{aligned}
\end{equation}
with high probability, where $d_{\max}=\max_{i\in[m]}d_i$. Notice the bound of $\zeta_{\mt}$ holds when $n=\omega(D^{2-2\gamma}\log D)$.
By the assumption about $\muu$, we have
\begin{equation}
\label{eq:ze_u_3}
	\begin{aligned}
		\zeta_\muu=\max_{i\in[m]}\|\muu^{(i)}\|_{2\to\infty}
		\lesssim d_{\max}^{1/2}D^{-1/2}.
	\end{aligned}
\end{equation}
And for $\epsilon_{\star}$, first we have
$$
\epsilon_{\star}
		=\max_{i,j\in[m]}\|\muu^{(i)\top}\mt^{(j)}_0\|
		=\max_{i\neq j}\|\muu^{(i)\top}(\mi-\muu^{(j)}\muu^{(j)\top})(\hat{\mSigma}^{(j)} -
\mSigma^{(j)}) \muu^{(j)}(\mLambda^{(j)})^{-1}\|.
$$
Then when there are no individual components $\{\hat\muu_s^{(i)}\}$, we have $\epsilon_{\star}=0$, while when there exist individual components, by Lemma~\ref{lemma:PCA_|E|_previous} we have
\begin{equation}
\label{eq:ep_cs_3}
	\begin{aligned}
		\epsilon_{\star}
		\leq \max_{i\neq j}\|\muu^{(i)\top} (\mi-\muu^{(j)}\muu^{(j)\top})(\hat{\mSigma}^{(j)} -
\mSigma^{(j)}) \muu^{(j)}\|\cdot\|(\mLambda^{(j)})^{-1}\|
\lesssim D^{-\gamma/2} \tilde{\varphi}
	\end{aligned}
\end{equation}
with high probability.
Therefore by Theorem~\ref{thm3_general}, for the estimation of $\muu_c$ we have
$$
\begin{aligned}
	\hat\muu_c\mw_{\muu_c}-\muu_c
	&=\frac{1}{m}\sum_{i=1}^m\mt_0^{(i)}\muu^{(i)\top}\muu_c+\mq_{\muu_c}
	\\
	&=\frac{1}{m}\sum_{i=1}^m(\mi-\muu^{(i)}\muu^{(i)\top})(\hat{\mSigma}^{(i)} -\mSigma^{(i)}) \muu_c(\mLambda_c^{(i)})^{-1}
+\mq_{\muu_c},
\end{aligned}
$$
where $\mw_{\muu_c}$ is a minimizer of $\|\hat\muu_c\mo-\muu_c\|_{F}$ over all $\mo \in \mathcal{O}_{d_0}$, 
and by Eq.~\eqref{eq:ep_t,ze_t_3}, Eq.~\eqref{eq:ze_u_3} and Eq.~\eqref{eq:ep_cs_3}, 
 when there are no individual components,
$\mq$ satisfies
$$
\begin{aligned}
	\|\mq_{\muu_c}\|
	&\lesssim \epsilon_{\star} + \epsilon_{\mt_0}^2+ \epsilon_{\mt} 
     \lesssim \varphi^2 +(D^{-\gamma} \varphi + \varphi^2)
	\lesssim D^{-\gamma} \varphi + \varphi^2,\\
	\|\mq_{\muu_c}\|_{2\to\infty}
	&\lesssim \zeta_\muu(\epsilon_{\star}+\epsilon_{\mt_0}^2+\epsilon_{\mt})
       +\zeta_{\mt_0}(\epsilon_{\star}+\epsilon_{\mt_0}+\epsilon_{\mt})
       +\zeta_{\mt}\\
    &\lesssim d_{\max}^{1/2}D^{-1/2}\cdot(D^{-\gamma} \varphi + \varphi^2)
    +d_{\max}^{1/2} D^{-\gamma/2} \tilde\varphi\cdot \varphi 
    +d_{\max}^{1/2}D^{-3\gamma/2}\tilde\varphi (1 + D \tilde{\varphi})\\
    &\lesssim d_{\max}^{1/2}D^{-3\gamma/2}\tilde\varphi (1 + D \tilde{\varphi})
\end{aligned}
$$
with high probability. Notice the bound of $\|\mq_{\muu_c}\|_{2\to\infty}$ holds when $n=\omega(D^{2-2\gamma}\log D)$.
When there are individual components,
$\mq$ satisfies
$$
\begin{aligned}
	\|\mq_{\muu_c}\|
	&\lesssim \epsilon_{\star} + \epsilon_{\mt_0}^2+ \epsilon_{\mt} 
     \lesssim D^{-\gamma/2} \tilde{\varphi}+\varphi^2 +(D^{-\gamma} \varphi + \varphi^2)
	\lesssim D^{-\gamma/2} \tilde{\varphi}+D^{-\gamma} \varphi + \varphi^2,\\
	\|\mq_{\muu_c}\|_{2\to\infty}
	&\lesssim \zeta_\muu(\epsilon_{\star}+\epsilon_{\mt_0}^2+\epsilon_{\mt})
       +\zeta_{\mt_0}(\epsilon_{\star}+\epsilon_{\mt_0}+\epsilon_{\mt})
       +\zeta_{\mt}\\
    &\lesssim d_{\max}^{1/2}D^{-1/2}\cdot(D^{-\gamma/2} \tilde{\varphi}+D^{-\gamma} \varphi + \varphi^2)
    +d_{\max}^{1/2} D^{-\gamma/2} \tilde\varphi\cdot \varphi 
    +d_{\max}^{1/2}D^{-3\gamma/2}\tilde\varphi (1 + D \tilde{\varphi})\\
    &\lesssim d_{\max}^{1/2}D^{-(\gamma+1)/2} \tilde{\varphi}
    +d_{\max}^{1/2}D^{-3\gamma/2}\tilde\varphi (1 + D \tilde{\varphi})
\end{aligned}
$$
with high probability.
And for each $i\in[m]$, the estimation for the possibly distinct subspace $\muu_s^{(i)}$ has the expansion
$$
\begin{aligned}
	\hat\muu_s^{(i)}\mw_{\muu_s}^{(i)}-\muu_s^{(i)}
	&=\mt_0^{(i)}\muu^{(i)\top}\muu_s^{(i)}+\mq_{\muu_s}^{(i)}
	\\
	&=(\mi-\muu^{(i)}\muu^{(i)\top})(\hat{\mSigma}^{(i)} -\mSigma^{(i)}) \muu_s^{(i)}(\mLambda_s^{(i)})^{-1}
+\mq_{\muu_s}^{(i)},
\end{aligned}
$$
where $\mw_{\muu_s}^{(i)}$ is a minimizer of $\|\hat\muu_s^{(i)}\mo-\muu_s^{(i)}\|_{F}$ over all $\mo \in \mathcal{O}_{d_i-d_0}$, and 
$\mq_{\muu_s}^{(i)}$ satisfies the same upper bounds as those for $\mq_{\muu_c}$.
\hspace*{\fill} \qedsymbol

\subsection{Proof of Proposition~\ref{prop:PCA_||hatUW-U||_previous}}
\begin{sloppypar}
Let $\overline{\bm{\Pi}}_{\muu}^{(i)} = \mi - \muu ^{(i)}\muu^{(i)\top}$. 
By Theorem~\ref{thm:PCA_UhatW-U=(I-UU^T)EVR^{-1}+xxx_previous} we have
\begin{equation}\label{Eq:PCA_prop_total_previous}
  \hat\muu_c\mw_{\muu_c}-\muu_c
=\frac{1}{m}\sum_{i=1}^m \overline{\bm{\Pi}}_{\muu}^{(i)} (\hat\mSigma^{(i)}-\mSigma^{(i)})\muu_c (\mLambda^{(i)})_c^{-1} +\mq_{\muu_c},
\end{equation}
where 
$
\|\mq_{\muu_c}\|_F
\leq d_0^{1/2}\|\mq_{\muu_c}\|
\lesssim d_0^{1/2} D^{-\gamma/2} \tilde{\varphi}+d_0^{1/2}D^{-\gamma}\varphi+d_0^{1/2}\varphi^2
$
with high probability.
We now expand
\begin{equation}\label{Eq:PCA_prop_decomposition_previous}
	\begin{aligned}
	\Big\|\frac{1}{m}\sum_{i=1}^m \overline{\bm{\Pi}}_{\muu}^{(i)} (\hat\mSigma^{(i)}-\mSigma^{(i)})\muu^{(i)}_c (\mLambda^{(i)}_c)^{-1}\Big\|^2_F
	&=\frac{1}{m^2}\sum_{i=1}^m\|\piu^{(i)}\me^{(i)}\muu_c(\mLambda^{(i)}_c)^{-1}\|^2_F\\&
	+\frac{1}{m^2}\sum_{i\neq j}\operatorname{tr}[
(\mLambda^{(i)}_c)^{-1}\muu^{(i)\top}_c\me^{(i)}\piu^{(i)}\piu^{(j)}\me^{(j)}\muu_c(\mLambda^{(i)}_c)^{-1}].
\end{aligned}
\end{equation}
For the first term on the right hand side of Eq.~\eqref{Eq:PCA_prop_decomposition_previous}, by Eq.~\eqref{Eq:PCA_||E||,||EU||_previous} we have
\begin{equation}\label{Eq:PCA_prop_decomposition_part1_previous}
	\begin{aligned}
	\|\piu^{(i)}\me^{(i)}\muu_c(\mLambda^{(i)}_c)^{-1}\|_F
	\leq d_0^{1/2}\|\me^{(i)}\|\cdot \|(\mLambda^{(i)})^{-1}\|
	\lesssim d_0^{1/2}\varphi
\end{aligned}
\end{equation}
with high probability.
For the second term on the right hand side of
Eq.~\eqref{Eq:PCA_prop_decomposition_previous}, for $i \not = j$ we have
$$\mathbb{E}[(\mLambda^{(i)}_c)^{-1}\muu_c^{\top}\me^{(i)}\piu^{(i)}\piu^{(j)}\me^{(j)}\muu_c(\mLambda_c^{(i)})^{-1}]
= \mathbf{0}.$$
We now consider the variance for the entries of it. 
For each $k\in[d_0]$, let $\zeta_k$ be the $k$th diagonal entry of $(\mLambda^{(i)}_c)^{-1}\muu_c^{\top}\me^{(i)}\piu^{(i)}\piu^{(j)}\me^{(j)}\muu_c(\mLambda_c^{(i)})^{-1}$
and let $\tilde{u}_{c,k}$ denote the $k$th column of $\muu_c$. 
Then we have
$$
\begin{aligned}
\zeta_k 
=&\frac{1}{(\mLambda_c)^{2}_{kk}}\tilde{u}_{c,k}^\top
\hat{\mSigma}^{(i)}\piu^{(i)}\piu^{(j)}\hat{\mSigma}^{(j)} \tilde{u}_{c,k},
\end{aligned}
$$
where we had used the fact that $\muu^{(i)\top} \bm{\Sigma}^{(i)} \piu^{(i)} =
\bm{0}$. 
Then by Lemma~4 and Lemma~9 in \cite{neudecker1986symmetry} we have
$$
\begin{aligned}
\mathrm{Var}[\zeta_k] &=\frac{1}{(\mLambda_c)^{4}_{kk}}\operatorname{Var} \Big(
\mathbb{E}\big[ \tilde{u}_{c,k}^\top
\hat{\mSigma}^{(i)}\piu^{(i)}\piu^{(j)} \hat{\mSigma}^{(j)} \tilde{u}_{c,k}
\big|\hat{\mSigma}^{(j)}\big]\Big)
\\& +\frac{1}{(\mLambda_c)^{4}_{kk}} \mathbb{E}\Big(
\operatorname{Var}\big[ \tilde{u}_{c,k}^\top
\hat{\mSigma}^{(i)}\piu^{(i)}\piu^{(j)} \hat{\mSigma}^{(j)} \tilde{u}_{c,k}
\big|\hat{\mSigma}^{(j)}\big]\Big)\\
&=0+\frac{1}{(\mLambda_c)^{4}_{kk}} \mathbb{E}\Big(
\operatorname{Var}\big[ 
(\tilde{u}_{c,k}^{\top} \hat{\mSigma}^{(j)} \piu^{(j)}\piu^{(i)} \otimes \tilde{u}_{c,k}^{\top})\operatorname{vec}(\hat{\mSigma}^{(i)})
\big|\hat{\mSigma}^{(j)}\big]\Big)\\
&=\frac{1}{n(\mLambda_c)^{4}_{kk}} 
\mathbb{E}\Big(
(\tilde{u}_{c,k}^{\top} \hat{\mSigma}^{(j)} \piu^{(j)}\piu^{(i)} \otimes
\tilde{u}_{c,k}^{\top})(\mSigma^{(i)}\otimes\mSigma^{(i)})(\mi_{D^2}
+\mathcal{K}_{D}) 
(\piu^{(i)}\piu^{(j)}
\hat{\mSigma}^{(j)}\tilde{u}_{c,k}\otimes \tilde{u}_{c,k})\Big)\\
&=\frac{1}{n(\mLambda_c)^{4}_{kk}} \mathbb{E}\Big(
\tilde{u}_{c,k}^{\top} \hat{\mSigma}^{(j)} \piu^{(j)}\piu^{(i)}
\mSigma^{(i)}
\piu^{(i)}\piu^{(j)}
\hat{\mSigma}^{(j)}\tilde{u}_{c,k}
 \cdot \tilde{u}_{c,k}^{\top} \mSigma^{(i)} \tilde{u}_{c,k} 
+ \bigl(\tilde{u}_{c,k}^{\top} \hat{\mSigma}^{(j)} \piu^{(j)}\piu^{(i)} \mSigma^{(i)} \tilde{u}_{c,k}\bigr)^2 \Big)\\
&=\frac{\sigma_i^2}{n(\mLambda_c)^{3}_{kk}} \mathbb{E}\Big(
\tilde{u}_{c,k}^{\top} \hat{\mSigma}^{(j)} \piu^{(j)}\piu^{(i)}\piu^{(i)}\piu^{(j)} \hat{\mSigma}^{(j)} \tilde{u}_{c,k}\Big),
\end{aligned}
$$
where $\mathcal{K}_D$ is the $D^2 \times D^2$ commutation matrix. Now since $\mathbb{E}[\tilde{u}_{c,k}^{\top} \hat{\mSigma}^{(j)}
\piu^{(j)} \piu^{(i)}]= \tilde{u}_{c,k}^{\top} \bm{\Sigma}^{(j)} \piu^{(j)}\piu^{(i)} = \bm{0}$, we
have
\begin{equation*}
  \begin{split}
\mathrm{Var}[\zeta_k] & 
= \frac{\sigma_i^2}{n (\mLambda_c)_{kk}^3} \mathrm{tr} \,
\mathrm{Var}\Bigl[(\tilde{u}_{c,k}^{\top}  \otimes \piu^{(i)}\piu^{(j)})
\hat{\bm{\Sigma}}^{(j)}\Bigr] \\ 
&= \frac{\sigma_i^2}{n^2 (\mLambda_c)_{kk}^3}
\mathrm{tr} (\tilde{u}_{c,k}^{\top}  \otimes \piu^{(i)}\piu^{(j)}) (\mSigma^{(j)}
\otimes \mSigma^{(j)}) (\mathbf{I}_{D^2} + \mathcal{K}_D) (\tilde{u}_{c,k} \otimes \piu^{(j)}\piu^{(i)}) 
\\ &= \frac{\sigma_i^2}{n^2 (\mLambda_c)_{kk}^3} \tilde{u}_{c,k}^{\top} \mSigma^{(j)} \tilde{u}_{c,k} \cdot\mathrm{tr} \,
(\piu^{(i)}\piu^{(j)} \mSigma^{(j)}
\piu^{(j)}\piu^{(i)} )
\\ &\leq \frac{\sigma_i^2}{n^2 (\mLambda_c)_{kk}^2} \mathrm{tr} \,
(\piu^{(j)} \mSigma^{(j)}
\piu^{(j)})\cdot \|\piu^{(i)}\|^2
= \frac{\sigma_i^2\sigma_j^2 D}{n^2 (\mLambda_c)_{kk}^2} \lesssim n^{-1} D^{-\gamma} \varphi^{2}.
\end{split}
\end{equation*}
Hence, by Chebyshev inequality, for $i\neq j, k\in[d_0]$ we have
\begin{equation}\label{Eq:PCA_prop_decomposition_part2_previous}
	\begin{aligned}
		\big[\muu_c^\top\muu^{(i)}(\mLambda^{(i)})^{-1}\muu^{(i)\top}\me^{(i)}\piu^{(i)}\piu^{(j)}\me^{(j)}\muu^{(i)}(\mLambda^{(i)})^{-1}\muu^{(i)\top}\muu_c\big]_{kk}
		\lesssim n^{-1/2}D^{-\gamma/2}\varphi
	\end{aligned}
\end{equation}
with probability converging to one. 
Combining Eq.~\eqref{Eq:PCA_prop_decomposition_previous}, Eq.~\eqref{Eq:PCA_prop_decomposition_part1_previous} and Eq.~\eqref{Eq:PCA_prop_decomposition_part2_previous}, we therefore have
$$
\begin{aligned}
	\Big\|\frac{1}{m}\sum_{i=1}^m \overline{\bm{\Pi}}_{\muu}^{(i)} (\hat\mSigma^{(i)}-\mSigma^{(i)})\muu_c (\mLambda^{(i)}_c)^{-1} \Big\|^2_F
	&\lesssim m^{-1} (d_0^{1/2}\varphi)^2
	+d_0 n^{-1/2}D^{-\gamma/2}\varphi
	\lesssim d_0m^{-1} \varphi^2
\end{aligned}
$$
with high probability. Recalling Eq.~\eqref{Eq:PCA_prop_total_previous} we have
$$
\|\hat\muu_c\mw_{\muu_c}-\muu_c\|_F\lesssim d_0^{1/2}m^{-1/2} \varphi 
+d_0^{1/2} D^{-\gamma/2} \tilde{\varphi}+d_0^{1/2}D^{-\gamma}\varphi+d_0^{1/2}\varphi^2
$$
with high probability, as desired. The analysis for the bound of $\|\hat\muu_s^{(i)}\mw_{\muu_s}^{(i)}-\muu_s^{(i)}\|_F$ and the analysis for the case with no individual subspaces follows similar arguments.
\hspace*{\fill} \qedsymbol
\end{sloppypar}

\subsection{Proof of Theorem~\ref{thm:PCA_(UhatW-U)_k->normal_previous}}
We now derive the normal approximation for $\hat u_{c,k}$. The result for $\hat u_{s,k}^{(i)}$ follows from similar arguments.
By Theorem~\ref{thm:PCA_UhatW-U=(I-UU^T)EVR^{-1}+xxx_previous} and
$\muu^{(i)\top}\mSigma^{(i)}(\mi-\muu^{(i)}\muu^{(i)\top})=\mathbf{0}$ we have
\begin{equation}\label{Eq:PCA_normal_decomposition_previous}
	\begin{aligned}
		\mw_{\muu_c}^\top\hat u_{c,k}-u_{c,k}
&=\frac{1}{m}\sum_{i=1}^m(\mLambda^{(i)}_c)^{-1}\muu_c^{\top}(\hat\mSigma^{(i)}-\mSigma^{(i)})(\mi-\muu^{(i)}\muu^{(i)\top})e_k+q_{\muu_c,k}\\
&=\sum_{i=1}^m\sum_{j=1}^n\my^{(k)}_{ij}
        + q_{\muu_c,k},
	\end{aligned}
\end{equation}
where $e_k$ is the $k$th basis vector, $q_{\muu_c,k}$ denotes the $k$th row of $\mq_{\muu_c}$, and we define
$$
\my^{(k)}_{ij}=\frac{1}{mn}(\mLambda^{(i)}_c)^{-1}\muu_c^{\top} X^{(i)}_jX^{(i)\top}_j(\mi-\muu^{(i)}\muu^{(i)\top})e_k.
$$
Note that $\{\my^{(k)}_{ij}\}_{i\in[m],j\in[n]}$ are independent mean $\bm{0}$ random vectors.
Let $\zeta_{i,k} := (\mi-\muu^{(i)}\muu^{(i)\top}) e_k$. Then 
for any $i\in[m],k\in[n]$, by Lemma~4 and Lemma~9 in \cite{neudecker1986symmetry}, the variance of $\my^{(k)}_{ij}$ is
\begin{equation}
  \label{eq:var_pca_1}
\begin{aligned}
\operatorname{Var}\big[\my^{(k)}_{ij}\big]
&=\frac{1}{m^2n^2}(\zeta_{i,k}^\top\otimes\mLambda_c^{(i)-1}\muu_c^{\top})
(\mSigma^{(i)}\otimes\mSigma^{(i)}) \times(\mi_{D^2}+\mathcal{K} _D)
(\zeta_{i,k}\otimes \muu_c\mLambda_c^{(i)-1})\\
&=\frac{1}{m^2n^2}(\zeta_{i,k}^\top\otimes \mLambda_c^{(i)-1}\muu_c^{\top})
(\mSigma^{(i)}\otimes\mSigma^{(i)})
 \times (\zeta_{i,k}\otimes \muu_c\mLambda_c^{(i)-1}+\muu_c\mLambda_c^{(i)-1}\otimes \zeta_{i,k})\\
&=\frac{1}{m^2n^2}\zeta_{i,k}^\top\mSigma^{(i)}\zeta_{i,k}\otimes (\mLambda^{(i)}_c)^{-1}
\\&=\frac{\sigma_i^2(1-\|u^{(i)}_k\|^2)}{m^2n^2} (\mLambda^{(i)}_c)^{-1},
\end{aligned}
\end{equation}
where $\mathcal{K}_D$
denotes the $D^2 \times D^2$ commutation matrix. See
Theorem~3.1 in \cite{magnus1979commutation} for a summary of some
simple but widely used relationships between commutation matrices and Kronecker products. 
As $\|\muu^{(i)}\|_{2 \to \infty} \lesssim
d_i^{1/2} D^{-1/2}$, we have $\|u^{(i)}_k\|^2 = o(1)$ for all $k$, and hence for each $i\in[m]$, 
$$
\sum_{j=1}^n\operatorname{Var}\big[\my^{(k)}_{ij}\big]
=(1+o(1))\bm{\Upsilon}_{\muu_c}^{(i)},
$$
where we define $\bm{\Upsilon}_{\muu_c}^{(i)}
    :=\frac{1}{Nm} \sigma_i^2(\mLambda^{(i)}_c)^{-1}$. Note that $\bm{\Upsilon}_{\muu_c}=\sum_{i=1}^m\bm{\Upsilon}_{\muu_c}^{(i)}$, 
where $\bm{\Upsilon}_{\muu_c}$ is defined in the statement of Theorem~\ref{thm:PCA_(UhatW-U)_k->normal_previous}.
As $\{\my_{ij}^{(k)}\}_{j \in [n]}$ are iid,
by the (multivariate) central limit theorem we have
$$
(\bm{\Upsilon}_{\muu_c}^{(i)})^{-1/2}\sum_{j=1}^n \my^{(k)}_{ij} 
    \leadsto \mathcal{N}\big(\mathbf{0},\mathbf{I}\big)
$$
as $n,D \rightarrow \infty$.
Then as $\big\{\sum_{j=1}^n \my^{(k)}_{ij}\big\}_{i\in[m]}$ are independent, we have
\begin{equation}\label{Eq:PCA_normal_decomposition1_previous}
	\begin{aligned}
	\bm{\Upsilon}_{\muu_c}^{-1/2}\sum_{i=1}^m\sum_{j=1}^n\my^{(k)}_{ij}
     \leadsto \mathcal{N}\big(\mathbf{0},\mathbf{I}\big).
	\end{aligned}
\end{equation}
as $n,D \rightarrow \infty$.

For the second term on the right hand side of Eq.~\eqref{Eq:PCA_normal_decomposition_previous}, from Theorem~\ref{thm:PCA_UhatW-U=(I-UU^T)EVR^{-1}+xxx_previous} we have
$$
\begin{aligned}
	\|\bm{\Upsilon}_{\muu_c}^{-1/2}q_{\muu_c,k}\|
	\leq &\|\bm{\Upsilon}_{\muu_v}^{-1/2}\|\cdot \|\mq_{\muu_c}\|_{2\to\infty}\\
	\lesssim & m^{1/2}n^{1/2}D^{\gamma/2}\cdot(d_{\max}^{1/2}D^{-(\gamma+1)/2}n^{-1/2}\log^{1/2}D
	+d_{\max}^{1/2}D^{-3\gamma/2}n^{-1/2}\log^{1/2}D\\&+d_{\max}^{1/2}D^{1-3\gamma/2}n^{-1}\log D)\\
    \lesssim &m^{1/2}d_{\max}^{1/2}\Bigl(\frac{\log^{1/2}D}{D^{1/2}}
  +\frac{\log^{1/2} D}{D^{\gamma}}
    +\frac{D^{1-\gamma}\log D}{n^{1/2}}\Bigr)
    \end{aligned}
$$
with high probability.
We then  have 
\begin{equation}
\label{Eq:PCA_normal_decomposition2_previous}
\bm{\Upsilon}_{\muu_c}^{-1/2}q_{\muu_c,k}\stackrel{p}{\longrightarrow} \bm{0}
\end{equation}
as $n,D \rightarrow \infty$, 
provided that $m=o(D^{2\gamma}/\log D)$ and $m=o(n/(D^{2-2\gamma}\log^2 D))$ 
as
assumed in the statement of Theorem~\ref{thm:PCA_(UhatW-U)_k->normal_previous}.
Combining Eq.~\eqref{Eq:PCA_normal_decomposition_previous}, Eq.~\eqref{Eq:PCA_normal_decomposition1_previous} and Eq.~\eqref{Eq:PCA_normal_decomposition2_previous}, and applying Slutsky's theorem, we obtain
$$
\bm{\Upsilon}_{\muu_c}^{-1/2}\big(\mw_{\muu_c}^\top\hat u_{c,k}-u_{c,k}\big)
\leadsto \mathcal{N}\big(\mathbf{0},\mathbf{I}\big)
$$
as $n,D \rightarrow \infty$.
\hspace*{\fill} \qedsymbol

\subsection{Proof of Theorem~\ref{thm:PCA_UhatW-U=(I-UU^T)EVR^{-1}+xxx}}
We begin with the statement of several basic bounds that are used
frequently in the subsequent derivations; these bounds are
reformulations of Theorem~6 and Theorem~9 in \cite{yan2021inference}
to the setting of the current paper. For ease of reference we
will use the same notations as that in \cite{yan2021inference}. Define   
$$
\begin{aligned}
	\mm^{(i)}
	=n^{-1/2}\mx^{(i)},\text{ }
	\mm^{\natural(i)}
	=\mathbb{E}[\mm^{(i)}|\mathbf{F}^{(i)}]
	=n^{-1/2}\my^{(i)}, \text{ }
	\me^{(i)}=\mm^{(i)}-\mm^{\natural(i)}=n^{-1/2}\mz^{(i)},
\end{aligned}
$$ 
and let the singular value decomposition of $\mm^{\natural(i)}$ be
$\mm^{\natural(i)}=\muu^{\natural(i)}\mSigma^{\natural(i)}\mv^{\natural(i)\top}$. We
note that if $n \geq d_i$ then, almost surely, there exists a $d_i \times d_i$ orthogonal
matrix $\mw^{\natural(i)}$ such that $\muu = \muu^{\natural(i)}
\mw^{\natural(i)}$.
\begin{lemma1}
  \label{lemma:PCA_|E|_basis}
  Consider the setting in
  Theorem~\ref{thm:PCA_UhatW-U=(I-UU^T)EVR^{-1}+xxx} and suppose
  $\tfrac{\log(n+D)}{n}\lesssim 1$. We then have
  \begin{gather*}
  \|\muu^{\natural (i)}\|_{2\to\infty} \lesssim d_i^{1/2}D^{-1/2},
  \quad\mSigma^{\natural(i)}_{rr}\asymp D^{\gamma/2}\quad \text{for
    any }r\in[d_i], \\
	\max_{k\in[D], \ell \in[n]}|\me^{(i)}_{k \ell}|
		\lesssim 
		n^{-1/2}\log^{1/2}(n+D),\quad
    \|\mv^{\natural(i)}\|_{2\to\infty}
	\lesssim d_i^{1/2}n^{-1/2}\log^{1/2}(n+D)
    \end{gather*}
    with probability at least $1 - O((n+D)^{-10})$.
    Here $\mSigma^{\natural(i)}_{rr}$ denote the $r$th largest
    singular value of $\mathbf{M}^{\natural(i)}$. 
\end{lemma1}

\begin{lemma1}
  \label{lemma:PCA_|E|}
  Consider the setting in
  Theorem~\ref{thm:PCA_UhatW-U=(I-UU^T)EVR^{-1}+xxx} and suppose
  $\tfrac{\log^2(n+D)}{n}\lesssim 1$. 
  We then have
 \begin{equation*}
 \begin{aligned}
&\|\me^{(i)}\|
\lesssim 
\Bigl(1+\frac{D}{n}\Bigr)^{1/2}, \quad 
\|\me^{(i)}\mv^{\natural(i)}\|_{2\to\infty}
\lesssim 
 d_i^{1/2}n^{-1/2}\log(n+D), \\
&\|\muu^{(i)\top}\me^{(j)}\mv^{\natural(j)}\|_F
\lesssim d_j^{1/2}n^{-1/2}\log(n+D)
 \end{aligned}
\end{equation*}
with probability at least $1-O((n+D)^{-10})$.
\end{lemma1}

Finally we state a technical lemma for the error of
$\hat\muu^{(i)}$ as an estimate for the true $\muu$. 
\begin{lemma1}
  \label{lemma:PCA_UhatW-U}
Consider the setting in
Theorem~\ref{thm:PCA_UhatW-U=(I-UU^T)EVR^{-1}+xxx}.
Define $$\phi=\frac{(n+D)\log(n+D)}{nD^\gamma} = \frac{\log(n+D)}{D^{\gamma}} \Bigl(1 +
\frac{D}{n}\Bigr).$$
Suppose
$\tfrac{\log^3(n+D)}{\min\{n,D\}}\lesssim 1$ and $\phi \ll 1$. 
Fix an $i\in[m]$ and let $\mw^{(i)}$ be a minimizer of
$\|\hat{\muu}^{(i)} \mathbf{O} - \muu^{(i)}\|_{F}$ over all
$d_i \times d_i$ orthogonal matrix $\mo$. 
Then conditional on $\mathbf{F}^{(i)}$ we have
$$
\hat\muu^{(i)}\mw^{(i)}-\muu^{(i)}=\me^{(i)}\mv^{\natural(i)}(\mSigma^{\natural(i)})^{-1}\mw^{\natural(i)}+\mt^{(i)},
$$
where $\mw^{\natural(i)}$ is such that $\muu^{(i)}=\muu^{\natural(i)}\mw^{\natural(i)}$. 
The residual matrix $\mt^{(i)}$ satisfies
\begin{equation}
  \label{eq:mt_1_pca_2}
 \begin{aligned}
	\|\mt^{(i)}\|_{2\to\infty}
    	&\lesssim \frac{d_i^{1/2}\phi}{(n+D)^{1/2}}
    	+\frac{d_i^{1/2} \phi}{D^{1/2}\log(n+D)}+\frac{d_i\phi^{1/2}}{(n+D)^{1/2}D^{1/2}}
\end{aligned}
\end{equation}
with probability as least $1-O((n+D)^{-10})$.
\end{lemma1}
The proofs of Lemma~\ref{lemma:PCA_|E|_basis} through Lemma~\ref{lemma:PCA_UhatW-U} are presented in Section~\ref{Appendix:proof of lemma:PCA_chen}.
We now complete the proof of
Theorem~\ref{thm:PCA_UhatW-U=(I-UU^T)EVR^{-1}+xxx} by invoking Theorem~\ref{thm3_general}. 
More specifically, for each $i\in[m]$, by Lemma~\ref{lemma:PCA_UhatW-U} we have the expansion 
$i \in [m]$
\begin{equation*}
  \hat\muu^{(i)}\mw^{(i)}-\muu^{(i)}=\mt^{(i)}_0+\mt^{(i)}
\end{equation*}
for some orthogonal matrix $\mw^{(i)}$, 
where
$\mt^{(i)}_0=\me^{(i)}\mv^{\natural(i)}(\mSigma^{\natural(i)})^{-1}\mw^{\natural(i)}$.
By Lemma~\ref{lemma:PCA_|E|_basis} and Lemma~\ref{lemma:PCA_|E|} we have
$$
\begin{aligned}
	&\|\mt^{(i)}_0\|
\leq\|\me^{(i)}\|\cdot \|(\mSigma^{\natural(i)})^{-1}\|
\lesssim\Bigl(1+\frac{D}{n}\Bigr)^{1/2} \cdot (D^{\gamma/2})^{-1}
\lesssim \Bigl(\frac{n+D}{nD^\gamma}\Bigr)^{1/2},\\
&\|\mt^{(i)}_0\|_{2\to\infty}
\leq\|\me^{(i)}\mv^{\natural(i)}\|_{2\to\infty}\cdot \|(\mSigma^{\natural(i)})^{-1}\|
\lesssim \frac{d_i^{1/2}\log(n+D)}{n^{1/2}D^{\gamma/2}}
\end{aligned}
$$
with probability at least $1-O((n+D)^{-10})$. 
Notice $\|\mt^{(i)}\|\leq D^{1/2}\|\mt^{(i)}\|_{2\to\infty}$. Then under the condition $\phi\ll 1$ and $\frac{\log(n+D)}{n+D}\lesssim 1$, we have
$$
\begin{aligned}
	\max_{i\in[m]}\big(2\|\mt^{(i)}_0\|
	+2\|\mt^{(i)}\|
	+\|\mt^{(i)}_0+\mt^{(i)}\|^2\big)
	\lesssim \Bigl(\frac{n+D}{nD^\gamma}\Bigr)^{1/2}
\end{aligned}
$$
with probability at least $1-O((n+D)^{-10})$, and thus under the assumption $\phi\ll 1$, we have $\max_{i\in[m]}\big(2\|\mt^{(i)}_0\|
	+2\|\mt^{(i)}\|
	+\|\mt^{(i)}_0+\mt^{(i)}\|^2\big)\ll \log^{-1/2}(n+D)$.
Under the assumption that $\|\mathbf\Pi_s\|=\|m^{-1}\sum_{i=1}^m\muu_s^{(i)}\muu_s^{(i)\top}\|=1-c_s$ for some constant $0<c_s\leq 1$, we have $\frac{1}{2}(1-\|\mathbf\Pi_s\|)\geq \frac{c_s}{2}$.
Then for large enough $n$ and $D$, under our assumption we hae
$$
\begin{aligned}
	\max_{i\in[m]}\big(2\|\mt^{(i)}_0\|
	+2\|\mt^{(i)}\|
	+\|\mt^{(i)}_0+\mt^{(i)}\|^2\big)
	\leq c(1-\|\mathbf\Pi_s\|)
	<\frac{1}{2}(1-\|\mathbf\Pi_s\|)
\end{aligned}
$$
with probability at least $1-O((n+D)^{-10})$ for any constant $c<\frac{1}{2}$.
Now we have
\begin{equation}
\label{eq:ep_t,ze_t_4}
	\begin{aligned}
		&\epsilon_{\mt_0}
		=\max_{i\in[m]}\|\mt^{(i)}_0\|
		\lesssim \Bigl(\frac{n+D}{nD^\gamma}\Bigr)^{1/2},\\
		&\zeta_{\mt_0}
		=\max_{i\in[m]}\|\mt^{(i)}_0\|_{2\to\infty}
		\lesssim \Bigl(\frac{d_{\max}\log^2(n+D)}{nD^{\gamma}}\Bigr)^{1/2},\\
		&\epsilon_{\mt}
		=\max_{i\in[m]}\|\mt^{(i)}\|
		\lesssim \frac{d_{\max}^{1/2}D^{1/2}\phi}{(n+D)^{1/2}}
    	+\frac{d_{\max}^{1/2} \phi}{\log(n+D)}+\frac{d_{\max}\phi^{1/2}}{(n+D)^{1/2}},\\
		&\zeta_{\mt}
		=\max_{i\in[m]}\|\mt^{(i)}\|_{2\to\infty}
		\lesssim \frac{d_{\max}^{1/2}\phi}{(n+D)^{1/2}}
    	+\frac{d_{\max}^{1/2} \phi}{D^{1/2}\log(n+D)}+\frac{d_{\max}\phi^{1/2}}{(n+D)^{1/2}D^{1/2}}
	\end{aligned}
\end{equation}
with probability at least $1-O((n+D)^{-10})$. 
By the assumption about $\muu$, we have
\begin{equation}
\label{eq:ze_u_4}
	\begin{aligned}
		\zeta_\muu=\max_{i\in[m]}\|\muu^{(i)}\|_{2\to\infty}
		\lesssim d_{\max}^{1/2}D^{-1/2}.
	\end{aligned}
\end{equation}
And by Lemma~\ref{lemma:PCA_|E|_basis} and Lemma~\ref{lemma:PCA_|E|} we have
\begin{equation}
\label{eq:ep_cs_4}
	\begin{aligned}
		\epsilon_{\star}
		&=\max_{i,j\in[m]}\|\muu^{(i)\top}\me^{(j)}\mv^{\natural(j)}(\mSigma^{\natural(j)})^{-1}\mw^{\natural(j)}\|
\leq \max_{i\in[m]}\|\muu^{(i)\top}\me^{(j)}\mv^{\natural(j)}\|	\cdot \|(\mSigma^{\natural(j)})^{-1}\|
		\\&
		\lesssim d_{\max}^{1/2}n^{-1/2}\log(n+D)\cdot (D^{\gamma/2})^{-1}
		\lesssim \Big(\frac{d_{\max}\log^2(n+D)}{n D^{\gamma}}\Big)^{1/2}
		\end{aligned}
\end{equation}
with probability at least $1-O((n+D)^{-10})$. 
Therefore by Theorem~\ref{thm3_general}, we have
$$
\begin{aligned}
	\hat\muu_c\mw_{\muu_c}-\muu_c
	&=\frac{1}{m}\sum_{i=1}^m\mt_0^{(i)}\muu^{(i)\top}\muu_c+\mq_{\muu_c}
	=\frac{1}{m}\sum_{i=1}^m\me^{(i)}\mv^{\natural(i)}(\mSigma^{\natural(i)})^{-1}\mw^{\natural(i)}\muu^{(i)\top}\muu_c+\mq_{\muu_c}\\
	&=\frac{1}{m} \sum_{i=1}^{m} \mz^{(i)} (\my^{(i)})^{\dagger} \muu_c + \mq_{\muu_c},
\end{aligned}
$$
where $\mw_{\muu_c}$ is a minimizer of $\|\hat\muu_c\mo-\muu_c\|_{F}$ over all orthogonal matrix $\mo$,
and by Eq.~\eqref{eq:ep_t,ze_t_4}, Eq.~\eqref{eq:ze_u_4} and Eq.~\eqref{eq:ep_cs_4}, $\mq$ satisfies
$$
\begin{aligned}
	\|\mq_{\muu_c}\|_{2\to\infty}
	&\lesssim \zeta_\muu(\epsilon_{\star}+\epsilon_{\mt_0}^2+\epsilon_{\mt})
       +\zeta_{\mt_0}(\epsilon_{\star}+\epsilon_{\mt_0}+\epsilon_{\mt})
       +\zeta_{\mt}\\
    &\lesssim \frac{d_{\max}(n+D)^{1/2}\log(n+D)}{nD^\gamma}
     	+\frac{d_{\max}(n+D)}{nD^{1/2+\gamma}}
     	+\frac{ d_{\max}(n+D)^{1/2}D^{1/2}\log^2(n+D)}{n^{3/2}D^{3\gamma/2}}
     	\\&+\frac{ d_{\max} \log(n+D)}{n^{1/2}D^{(1+\gamma)/2}}
\end{aligned}
$$
with probability at least $1-O((n+D)^{-10})$. 
And for each $i\in[m]$, the estimation for $\muu_s^{(i)}$ has the expansion
$$
\begin{aligned}
	\hat\muu_s^{(i)}\mw_{\muu_s}^{(i)}-\muu_s^{(i)}
	&=\mt_0^{(i)}\muu^{(i)\top}\muu_s^{(i)}+\mq_{\muu_s}^{(i)}
	=\me^{(i)}\mv^{\natural(i)}(\mSigma^{\natural(i)})^{-1}\mw^{\natural(i)}\muu^{(i)\top}\muu_s^{(i)}+\mq_{\muu_s}^{(i)}\\
	&=\mz^{(i)} (\my^{(i)})^{\dagger} \muu_s^{(i)} + \mq_{\muu_s}^{(i)},
\end{aligned}
$$
where $\mw_{\muu_s}^{(i)}$ is a minimizer of $\|\hat\muu_s^{(i)}\mo-\muu_s^{(i)}\|_{F}$ over all orthogonal matrix $\mo$,
and $\mq_{\muu_s}^{(i)}$ satisfies the same upper bounds as that for $\mq_{\muu_c}$.
\hspace*{\fill} \qedsymbol

\subsection{Proof of Theorem~\ref{thm:PCA_(UhatW-U)_k->normal}}
We now derive the normal approximation for $\hat u_{c,k}$. The result for $\hat u_{s,k}^{(i)}$ follows from similar arguments.
By Theorem~\ref{thm:PCA_UhatW-U=(I-UU^T)EVR^{-1}+xxx} we have
\begin{equation}\label{Eq:PCA_normal_decomposition}
	\begin{aligned}
		\mw_{\muu_c}^\top\hat u_{c,k}-u_{c,k}
        &=\frac{1}{m}\sum_{i=1}^m \muu_c^\top (\my^{(i)})^{\dagger\top} \mz^{(i)\top}e_k+q_{\muu_c,k}
        \\ &=\frac{1}{m}
		\sum_{i=1}^m
		\sum_{\ell=1}^n 
		\mz^{(i)}_{k \ell}\muu_c^\top(\my^{ (i)})^{\dagger}_{\ell}
		+q_{\muu_c,k},
	\end{aligned}
\end{equation}
where $e_k$ is the $k$th basis vector, $(\my^{ (i)})^{\dagger}_{\ell}$ denotes the $\ell$th row of $(\my^{ (i)})^{\dagger}$, 
and $q_{\muu_c,k}$ denotes the $k$th row
of $\mq_{\muu_c}$.

We now follow the arguments used in the proof of Lemma~9 in \cite{yan2021inference}.
We first derive the limiting distribution of the first term on the right hand side of 
Eq.~\eqref{Eq:PCA_normal_decomposition}.
This term is, conditional on $\{\mathbf{F}^{(i)}\}$, the sum of independent mean $\bm{0}$ random vectors $\{\xi^{(k)}_{il}\}_{i\in[m],l\in[n]}$, where
$$
\xi^{(k)}_{i \ell}=\frac{1}{m}\mz^{(i)}_{k \ell}\muu_c^\top(\my^{ (i)})^{\dagger}_{\ell}
$$
and $(\my^{ (i)})^{\dagger}_{\ell}$ is the $l$th row of $(\my^{ (i)})^{\dagger}$.
Let
$\tilde{\bm{\Upsilon}}=\sum_{i=1}^m\sum_{\ell=1}^n\operatorname{Var}\big[\xi^{(k)}_{i
  \ell}|\mathbf{F}^{(i)}\big]$ and
$\bm{\Upsilon}=\bm{\Upsilon}_{\muu_c}$. Recall the definition of $\bm{\Upsilon}_{\muu_c}$ in the statement of
Theorem~\ref{thm:PCA_(UhatW-U)_k->normal}.
Let $\mathcal{E}_{\mathrm{good}}^{(i)}$ denote the event defined in Lemma~6 of \cite{yan2021inference} where $\mathcal{E}_{\mathrm{good}}^{(i)}$ is measurable with respect to the sigma-algebra generated by $\mathbf{F}^{(i)}$ and
$\mathbb{P}(\mathcal{E}_{\text{good}}^{(i)})\geq 1-O((n+D)^{-10})$.
Now let $\mathcal{E}_{\mathrm{good}} = \cap_{i=1}^{m} \mathcal{E}_{\mathrm{good}}^{(i)}$
and note that $\mathbb{P}(\mathcal{E}_{\text{good}}) \geq 1 - O(m(n+D)^{-10})$.
Next assume (unless stated otherwise) that the event $\mathcal{E}_{\mathrm{good}}$ occurs and
$\frac{\log^3D}{n} = o(1)$.
Then by 
Lemma~8 in \cite{yan2021inference} and Weyl's inequality, we have
\begin{equation}\label{Eq:PCA_normal_||Upsilon||}
\begin{aligned}
	&\|\tilde{\bm{\Upsilon}}-{\bm{\Upsilon}}\|
	\lesssim \frac{d_{\max}^{1/2}\log^{3/2}(n+D)}{mn^{3/2} D^\gamma},\\
	&\lambda_{i}({\bm{\Upsilon}})\asymp \frac{1}{mn D^\gamma},\quad
	\lambda_{i}(\tilde{\bm{\Upsilon}})\asymp \frac{1}{mn D^\gamma},\quad \text{for any }i\in[d_0].
\end{aligned}
\end{equation}
Because $\xi^{(k)}_{i \ell}=\frac{1}{m}\mz^{(i)}_{k
  \ell}\muu_c^\top(\my^{
  (i)})^{\dagger}_{\ell}=\frac{1}{m}\me^{(i)}_{k\ell}\mw^{\natural
  (i)\top}(\mSigma^{\natural (i)})^{-1}v^{\natural (i)}_\ell$ where
$v^{\natural (i)}_\ell$ is the
$\ell$th row of $\mv^{\natural (i)}_\ell$, by Lemma~\ref{lemma:PCA_|E|_basis} the spectral norm of
$\tilde{\bm{\Upsilon}}^{-1/2}\xi^{(k)}_{i \ell}$ can be bounded as
\begin{equation}\label{Eq:PCA_normal_||Yij||}
	\begin{aligned}
		\|\tilde{\bm{\Upsilon}}^{-1/2}\xi^{(k)}_{i \ell}\|
		&\leq 
		\|\tilde{\bm{\Upsilon}}^{-1/2}\|
		\cdot m^{-1}|\me^{(i)}_{k \ell}|
		\cdot \|\mv^{\natural (i)}\|_{2\to\infty}
		\cdot \|(\mSigma^{\natural (i)})^{-1}\|\\
		&\lesssim m^{1/2}n^{1/2}D^{\gamma/2}
		\cdot \frac{\log^{1/2}(n+D)}{m n^{1/2}}
		\cdot \frac{d_{\max}^{1/2}\log^{1/2}(n+D)}{n^{1/2}}
		\cdot D^{-\gamma/2}
		\\ &
		\lesssim \frac{d_{\max}^{1/2}\log(n+D)}{m^{1/2}n^{1/2}}.
	\end{aligned}
\end{equation}
Now fix an arbitrary $\epsilon>0$. Then under the assumption
$\frac{\log^2(n+D)}{n}=o(1)$,
we have from Eq,~\eqref{Eq:PCA_normal_||Yij||} that for sufficiently large
$n$ and $D$, 
$
\|\tilde{\bm{\Upsilon}}^{-1/2}\xi^{(k)}_{i \ell}\| \leq \epsilon
$
for all $i\in[m], \ell\in[n]$. We thus have
$$
\begin{aligned}
		\sum_{i=1}^{m}\sum_{\ell=1}^{n}
\mathbb{E}\Bigl[\bigl\|\tilde{\bm{\Upsilon}}^{-1/2}\xi^{(k)}_{i \ell}\bigr\|^{2}
\cdot\mathbb{I}\bigl\{\|\tilde{\bm{\Upsilon}}^{-1/2}\xi^{(k)}_{i \ell}\|>\epsilon \bigr\}\Bigr] \longrightarrow 0.
	\end{aligned}
$$
Therefore, by the Lindeberg-Feller central limit theorem (see e.g., Proposition~2.27
in \cite{van2000asymptotic}), we have
\begin{equation}\label{Eq:PCA_normal_clt1}
	\begin{aligned}
	\tilde{\bm{\Upsilon}}^{-1/2}
		\sum_{i=1}^m
		\sum_{\ell=1}^n 
		\xi^{(k)}_{i \ell}
		\leadsto \mathcal{N}(\mathbf{0},\mathbf{I})
	\end{aligned}
\end{equation}
as $(n+D) \rightarrow \infty$. Next we have
\begin{equation*}
  \begin{split}
 \Bigl\|\bm{\Upsilon}^{-1/2} \sum_{i=1}^m \sum_{\ell=1}^n 
		\xi^{(k)}_{i \ell} - \tilde{\bm{\Upsilon}}^{-1/2} \sum_{i=1}^m \sum_{\ell=1}^n 
		\xi^{(k)}_{i \ell}\Bigr\| &= \Bigl\|\bm{\Upsilon}^{-1/2}
                             (\tilde{\bm{\Upsilon}}^{1/2} - \bm{\Upsilon}^{1/2}) \tilde{\bm{\Upsilon}}^{-1/2} 
\sum_{i=1}^m \sum_{\ell=1}^n 
                             \xi^{(k)}_{i \ell}\Bigr\| \\
    & \leq \|\bm{\Upsilon}^{-1/2}(\tilde{\bm{\Upsilon}}^{1/2} -
      \bm{\Upsilon}^{1/2})\| \cdot \Bigl\|\tilde{\bm{\Upsilon}}^{-1/2} \sum_{i=1}^m \sum_{\ell=1}^n 
                             \xi^{(k)}_{i \ell}\Bigr\|.
  \end{split}
\end{equation*}
Eq.~\eqref{Eq:PCA_normal_||Upsilon||} then implies (see e.g., Problem X.5.5. in \cite{bhatia2013matrix})
\begin{equation}
  \label{eq:theorem8_technica3}
  \begin{split}
\|\bm{\Upsilon}^{-1/2}(\tilde{\bm{\Upsilon}}^{1/2} -
      \bm{\Upsilon}^{1/2})\| &\leq \|\bm{\Upsilon}^{-1/2}\| \cdot \|\tilde{\bm{\Upsilon}}^{1/2} -
      \bm{\Upsilon}^{1/2}\| \lesssim \frac{d_{\max}^{1/2} \log^{3/2}(n + D)}{n^{1/2}}.
  \end{split}
\end{equation}
Combining Eq.~\eqref{Eq:PCA_normal_clt1} and
Eq.~\eqref{eq:theorem8_technica3}, under the assumption $\frac{\log^2(n+D)}{n}=o(1)$ we obtain
\begin{equation*}
  \Bigl\|\bm{\Upsilon}^{-1/2} \sum_{i=1}^m \sum_{\ell=1}^n 
		\xi^{(k)}_{i \ell} - \tilde{\bm{\Upsilon}}^{-1/2} \sum_{i=1}^m \sum_{\ell=1}^n 
		\xi^{(k)}_{i \ell}\Bigr\| \overset{p}{\longrightarrow}  0  
\end{equation*}
as $n \rightarrow \infty$, and hence, by Slutsky's theorem
\begin{equation}\label{Eq:PCA_normal_clt3}
	\begin{aligned}
	{\bm{\Upsilon}}^{-1/2}
		\sum_{i=1}^m
		\sum_{\ell=1}^n 
		\xi^{(k)}_{i \ell}
		\leadsto \mathcal{N}(\mathbf{0},\mathbf{I})
	\end{aligned}
\end{equation}
as $(n+D) \rightarrow \infty$; we emphasize that
Eq.~\eqref{Eq:PCA_normal_clt3} is
conditional on $\mathcal{E}_{\mathrm{good}}$ and $\{\mathbf{F}^{(i)}\}$ so that the only source of
randomness is in $\{\mathbf{Z}^{(i)}\}$.

\begin{sloppypar}
	We now remove the conditioning on $\mathcal{E}_{\mathrm{good}}$ and $\{\mathbf{F}^{(i)}\}$. Let $\mathcal{Y} = \bm{\Upsilon}^{-1/2}
		\sum_{i=1}^m
		\sum_{\ell=1}^n 
		\xi^{(k)}_{i \ell}$ and $\mathcal{Z} \sim \mathcal{N}(\mathbf{0},\mathbf{I})$. Then
for any convex set $\mathcal{B}$ in $\mathbb{R}^{d}$, we have
\begin{equation}\label{Eq:PCA_normal_clt4}
\begin{aligned}
	\Big|\mathbb{P}\Big(
	\mathcal{Y}
		\in \mathcal{B}
	\Big)
	-\mathbb{P}\Big(
	\mathcal{Z}
		\in \mathcal{B}
	\Big)\Big|
	=&\Big|\mathbb{E} \Big[\Big[\mathbb{P}\Big(
	\mathcal{Y}\in \mathcal{B}
		|\{\mathbf{F}^{(i)}\}\Big)
		-\mathbb{P}\Big(
	\mathcal{Z}
		\in \mathcal{B}
	\Big)\Big]
	\mathbb{I}_{\mathcal{E}_{\text {good }}}\Big]
	\Big|\\
	+&\Big|\mathbb{E} \Big[\Big[\mathbb{P}\Big(
	\mathcal{Y}\in \mathcal{B}
		|\{\mathbf{F}^{(i)}\}\Big)
		-\mathbb{P}\Big(
	\mathcal{Z}
		\in \mathcal{B}
	\Big)\Big]
	\mathbb{I}_{\mathcal{E}^c_{\text {good }}}\Big]
	\Big|\\
	\leq &\Big|\mathbb{E} \Big[\Big[\mathbb{P}\Big(
	\mathcal{Y}\in \mathcal{B}
		|\{\mathbf{F}^{(i)}\}\Big)
		-\mathbb{P}\Big(
	\mathcal{Z}
		\in \mathcal{B}
	\Big)\Big]
	\mathbb{I}_{\mathcal{E}_{\text {good }}}\Big]
	\Big|
	+2\mathbb{P}(\mathcal{E}^c_{\text {good }}).
\end{aligned}
\end{equation}
Combining Eq.~\eqref{Eq:PCA_normal_clt3}, Eq.~\eqref{Eq:PCA_normal_clt4}, and $\mathbb{P}(\mathcal{E}^c_{\text {good }})\leq O(m(n+D)^{-10})$, we obtain the {\em unconditional} limit result
${\bm{\Upsilon}}^{-1/2}
		\sum_{i=1}^m
		\sum_{\ell=1}^n 
		\xi^{(k)}_{i \ell}
		\leadsto \mathcal{N}(\mathbf{0},\mathbf{I})
$ as $(n+D) \rightarrow \infty$, i.e.,
\begin{equation}\label{Eq:PCA_normal_decomposition1}
	\begin{aligned}
	{\bm{\Upsilon}}^{-1/2}
		\frac{1}{m}\sum_{i=1}^m \muu_c^\top (\my^{(i)})^{\dagger\top} \mz^{(i)\top}e_k
		\leadsto \mathcal{N}(\mathbf{0},\mathbf{I})
	\end{aligned}
\end{equation}
as $(n+D) \rightarrow \infty$.
For the term involving $q_{\muu_c,k}$ in  
Eq.~\eqref{Eq:PCA_normal_decomposition}, from Theorem~\ref{thm:PCA_UhatW-U=(I-UU^T)EVR^{-1}+xxx} we have
$$
\begin{aligned}
	\|\bm{\Upsilon}^{-1/2}q_{\muu_c,k}\|
	\leq &\|\bm{\Upsilon}^{-1/2}\|\cdot\|\mq_{\muu_c}\|_{2\to\infty}\\
    \lesssim &\frac{m^{1/2}d_{\max}(n+D)^{1/2}\log(n+D)}{n^{1/2}D^{\gamma/2}}
    	+\frac{m^{1/2}d_{\max}(n+D)}{n^{1/2}D^{1/2+\gamma/2}}
    	\\&\!+ \frac{m^{1/2} d_{\max}(n+D)^{1/2}D^{1/2}\log^2(n+D)}{nD^{\gamma}}
    	+\frac{m^{1/2} d_{\max} \log(n+D)}{D^{1/2}}
\end{aligned}
$$
with probability as least $1-O((n+D)^{-10})$. We then have
\begin{equation}\label{Eq:PCA_normal_decomposition2}
	\begin{aligned}
		\bm{\Upsilon}^{-1/2}q_{\muu_c,k}\stackrel{p}{\longrightarrow} \bm{0}
	\end{aligned}
\end{equation}
as $(n+D) \rightarrow \infty$, provided the following conditions hold
$$
m=o\Big(\frac{nD^\gamma}{(n+D)\log^2(n+D)}\big), \quad
m=o\Big(D^{1+\gamma}/n\Big).$$
Combining Eq.~\eqref{Eq:PCA_normal_decomposition}, Eq.~\eqref{Eq:PCA_normal_decomposition1} and Eq.~\eqref{Eq:PCA_normal_decomposition2}, and applying Slutsky's theorem, we have
$$
\bm{\Upsilon}^{-1/2}_{\muu_c} \big(\mw^\top\hat{u}_{c,k}-u_{c,k}\big)
\leadsto \mathcal{N}(\mathbf{0},\mathbf{I})
$$
as $(n+D) \rightarrow \infty$. 
\hspace*{\fill} \qedsymbol
\end{sloppypar}

\section{Important Technical Lemmas}
\label{Appendix:B}
\subsection{Proof of Lemma~\ref{lemma:|E|_2,|UEV|F}}
\label{Appendix:proof of lemma:|E|_2,|UEV|F}

For ease of exposition, we fix values of $i, j$ and omit the indices $i, j$ from $\muu^{(i)}$, $\me^{(i)}/\me^{(j)}$, $\mv^{(i)}/\mv^{(j)}$, and $d_i/d_j$. Specifically, when no ambiguity arises, we write $\muu$, $\me$, $\me^{(1)}$, $\me^{(2)}$, $\mv$, and $d$ in place of $\muu^{(i)}$, $\me^{(i)}/\me^{(j)}$, $\me^{(i,1)}/\me^{(j,1)}$, $\me^{(i,2)}/\me^{(j,2)}$, $\mv^{(i)}/\mv^{(j)}$, and $d_i/d_j$, respectively.

We first bound $\|\me^{(1)}\|, \|\muu^{\top} \me^{(1)} \mv\|,
\|\me^{(1)} \mv\|_{2 \to \infty}$, and $\|\me^{(1)\top} \muu\|_{2
  \to \infty}$. For ease of exposition in our subsequent
derivations we will let $C$ 
denote a {\em universal} constant that can change from line to line,
i.e., $C$  can depend on $\{C_1, C_2, C_3\}$ but does not depend
on $m,n$ or $\rho_n$.

For $\|\mathbf{E}^{(1)}\|$, according to Remark~3.13 of
\cite{bandeira2016sharp}, there exists for any $0<\varepsilon\leq 1/2$
a universal constant $\tilde{c}_{\varepsilon}$ such that for every
$t\geq 0$
\begin{equation*}
\begin{aligned}
	\mathbb{P}\Big(\|\me^{(1)}\| \geq(1+\varepsilon) 2 \sqrt{2} \tilde{\sigma}+t\Big) \leq n e^{-t^{2} / \tilde{c}_{\varepsilon} \tilde{\sigma}_{*}^{2}},
\end{aligned}
\end{equation*}
where $\tilde{\sigma}_{*} =\max
_{k,\ell\in[n]}\|\mathbf{E}^{(1)}_{k\ell}\|_{\infty} \leq C_1$ almost surely and
$$
\begin{aligned}
	\tilde{\sigma}^2
	=&\max\Bigl\{\max_{k\in[n]}
    \sum_{\ell=1}^n\mathrm{Var}[\me^{(1)}_{k \ell}], \max_{\ell\in[n]}
       \sum_{k=1}^n\mathrm{Var}[\me^{(1)}_{k \ell}]\Bigr\}
	\leq C_2 n\rho_n.
\end{aligned}
$$
Let $t=C(n\rho_n)^{1/2}$ for some sufficiently large constant $C$. We then have
$$
\begin{aligned}
	\mathbb{P}\Big(\|\me^{(1)}\| \geq(1+\varepsilon) 2 \sqrt{2} \tilde{\sigma}+C(n\rho_n)^{1/2}\Big)
	&\leq n e^{-C^2(n\rho_n) / \tilde{c}_{\varepsilon} \tilde{\sigma}_{*}^{2}}.
\end{aligned}
$$
From the assumption $n\rho_n=\Omega(\log n)$, we have
$\|\mathbf{E}^{(1)}\|\lesssim (n\rho_n)^{1/2}$ with high probability.

For $\muu^\top\me^{(1)}\mv$ we follow the argument for Claim~S.4 in
\cite{zhang2022perturbation}. Let
$\mz^{(i;k,\ell)}
    =\me^{(1)}_{k\ell} u_k v_\ell^\top,
$
where $u_k$ denotes the $k$th row of $\muu$ and $v_\ell$ denotes
the $\ell$th row of $\mv$. Then $\muu^\top\me^{(1)}\mv$
is the sum of  independent mean $\mathbf{0}$ random matrices
$\{\mz^{(k,\ell)}\}_{k,\ell\in[n]}$ where, for any
$\mz^{(k,\ell)}$, we have
    $$
    \begin{aligned}
    	\|\mz^{(k,\ell)}\|
    	&\leq |\me^{(1)}_{k\ell}|
    	\cdot \|\muu\|_{2\to\infty}
    	\cdot \|\mv\|_{2\to\infty}
    	\lesssim C_1 \cdot
    	d^{1/2}n^{-1/2}\cdot d^{1/2}n^{-1/2}
    	\lesssim d n^{-1}
    \end{aligned}
    $$
    almost surely.
    Now $\mz^{(k,\ell)}(\mz^{(k,\ell)})^{\top}
    =(\me^{(1)}_{k\ell})^2 \|v_\ell\|^2u_k  u_k ^\top$ and hence, by Weyl's inequality, we have
    $$
    \begin{aligned}
    	\Big\|\sum_{k=1}^n\sum_{\ell=1}^n \mathbb{E}[\mz^{(k,\ell)}(\mz^{(k,\ell)})^{\top}]\Big\|
    	&\leq \max_{k,\ell\in[n]} \mathbb{E}[(\me^{(1)}_{k\ell})^2] 
    	\cdot \sum_{\ell=1}^n \|v_\ell\|^2
    	\cdot \Big\|\sum_{k=1}^n u_k  u_k ^\top\Big\|\\
       	&\lesssim  \rho_n
    	\cdot n\|\mv\|_{2\to\infty}^2
    	\cdot \|\muu^{\top}\muu\|
    	\lesssim d\rho_n.
    \end{aligned}
    $$
    Similarly, we also have
    $$
    \Big\|\sum_{k=1}^n\sum_{\ell=1}^n \mathbb{E}[(\mz^{(k,\ell)})^{\top}\mz^{(k,\ell)}]\Big\|
    	\lesssim d\rho_n.
    $$
    Therefore, by Theorem~1.6 in \cite{tropp2012user}, there exists
    a $C > 0$ such that for all $t>0$ 
    we have
    $$
    \begin{aligned}
    	\mathbb{P}\Big(\|\muu^{\top}\me^{(1)}\mv\|\geq t\Big)
    	&\leq 2d\cdot \exp \Big(\frac{-Ct^2}{d\rho_n+dn^{-1}
    	t/3}\Big),
    \end{aligned}
    $$
    and hence, with $t \asymp d^{1/2}(\rho_n \log n)^{1/2}$, we obtain
    $$
    \begin{aligned}
    	\|\muu^{\top}\me^{(1)}\mv\|
    	\lesssim 
    	d^{1/2}(\rho_n \log n)^{1/2}
    \end{aligned}
    $$
    with high probability.

For $\me^{(1)}\mv$, its $k$th row is
$\sum_{\ell=1}^n \me^{(1)}_{k \ell} v_{\ell}$
where $v_{\ell}$ represents the $\ell$th row of
$\mv$. Once again, by Theorem~1.6 in
\cite{tropp2012user}, we
have $$\Bigl\|\sum_{\ell=1}^{n} \me^{(1)}_{k \ell}
v_{\ell}\Bigr\|\lesssim d^{1/2}(\rho_n \log
n)^{1/2}$$ with high probability.
Taking a union over $k \in [n]$ we obtain
$ \|\me^{(1)}\mv \|_{2\to\infty}\lesssim d^{1/2}(\rho_n \log n)^{1/2}$
with high probability.
The proof for $\me^{(1)\top}\muu$ is identical
and is thus omitted. 
 If we further assume $\{\me^{(i,1)}\}$ are independent, by almost identical proof we have $\|\frac{1}{m}\sum_{i=1}^m\me^{(i,1)}\mv (\mr^{(i)})^{-1}\|_{2\to\infty}\lesssim d^{1/2}(mn)^{-1/2}(n\rho_n)^{-1/2}\log^{1/2} n$
with high probability.

We now bound $\|\me^{(2)}\|, \|\muu^{\top} \me^{(2)} \mv\|,
\|\me^{(2)} \mv\|_{2 \to \infty}$, and $\|\me^{(2)\top} \muu\|_{2
  \to \infty}$. 
The matrix $\rho_n^{-1/2}
\me^{(2)}$ contains independent mean-zero sub-Gaussian random
variables whose Orlicz-$2$ norms are bounded from above by
$C_3$. Therefore, by a standard $\epsilon$-net argument (see e.g.,
Theorem~4.4.5 in \cite{vershynin2018high}), we have
$$\|\rho_n^{-1/2} \me^{(2)}\| \lesssim C_3(n^{1/2} + \log^{1/2}{n})$$
with high probability. We thus obtain $\|\me^{(2)}\| \lesssim (n
\rho_n)^{1/2}$ with high probability.

Next, for $\|\muu^{\top} \me^{(2)} \mv\|$ we have
$$\|\muu^{\top} \me^{(2)} \mv\| = \sup_{\bm{x}, \bm{y}}
\Bigl|\bm{x}^{\top} \muu^{\top} \me^{(2)} \mv \bm{y}\Bigr|$$
where the supremum is over all $\bm{x} \in \mathbb{R}^{d}, \bm{y} \in
\mathbb{R}^{d}, \|\bm{x}\| = \|\bm{y}\| = 1$. Fix vectors $\bm{x}$
and $\bm{y}$ of unit norms and let $\bm{\xi} = \muu
\bm{x} \in \mathbb{R}^{n}$ and $\bm{\zeta} = \mv \bm{y} \in
\mathbb{R}^{n}$. Then
$$\bm{x}^{\top} \muu^{\top} \me^{(2)} \mv \bm{y} =
\mathrm{vec}(\me^{(2)})^{\top} \mathrm{vec}(\bm{\xi}
\bm{\zeta}^{\top}) = \sum_{k=1}^{n} \sum_{\ell=1}^{n} \me^{(2)}_{k \ell} \xi_k \zeta_{\ell}$$
is a sum of independent mean-zero sub-gaussian random variables.  
Hence, by the general form of Hoeffding's inequality
(see e.g., Theorem~2.6.3 in \cite{vershynin2018high}), there exists a $C>0$ such that for all $t>0$ we have
$$\mathbb{P}\Bigl(\Bigl|\sum_{k=1}^n \sum_{\ell=1}^n \me^{(2)}_{k \ell} \xi_k \zeta_{\ell}\Bigr| \geq t\Bigr) \leq 2
\exp\Bigl(\frac{-Ct^2}{C_3^2 \rho_n \|\mathrm{vec}(\bm{\xi} \bm{\zeta}^{\top})\|^2}\Bigr).$$
Now $ \|\muu \bm{x}\| = \|\bm{x}\| =  1 = \|\bm{y}\| = \|\mv \bm{y}\|$
and hence $\|\mathrm{vec}(\bm{\xi} \bm{\zeta}^{\top})\|^2 = \|\bm{\xi}\|^2 \cdot
\|\bm{\zeta}\|^2 = 1$.
Then with $t \asymp (\rho_n \log n)^{1/2}$ we obtain
$$
\Big|\bm{x}^{\top} \muu^{\top} \me^{(2)} \mv \bm{y}\Big|
\lesssim (\rho_n \log n)^{1/2}
$$
with high probability.
Let $\mathcal{M}$ be a $\epsilon$-net of the unit sphere in
$\mathbb{R}^{d}$, and set $\epsilon=1/3$. Then the cardinality of $\mathcal{M}$ is bounded by $|\mathcal{M}|\leq 18^d$. As $d$ is fixed we have
$$
\max_{\bm{x} \in
  \mathcal{M}, \bm{y} \in \mathcal{M}}\Big|\bm{x}^{\top} \muu^{\top} \me^{(2)} \mv \bm{y}\Big|
\lesssim (\rho_n \log n)^{1/2}
$$
with high probability.
By a standard $\epsilon$-net argument, we have
$$
\|\muu^{\top} \me^{(2)} \mv\|
\leq 
\frac{1}{1-\epsilon^2-2\epsilon}\max_{\bm{x} \in
  \mathcal{M}, \bm{y} \in \mathcal{M}}\Big|\bm{x}^{\top} \muu^{\top} \me^{(2)} \mv \bm{y}\Big|
  \lesssim \frac{9}{2}(\rho_n \log n)^{1/2}
  \lesssim (\rho_n \log n)^{1/2}
$$
with high probability.

For $\me^{(2)} \mv$, its $k$th row is of the form
$\sum_{\ell=1}^n \me^{(2)}_{k \ell} v_{\ell}$. As $\me^{(2)}_{k
  \ell}$ is mean-zero sub-Gaussian, $\me^{(2)}_{k \ell} v_{\ell}$ is a
mean-zero sub-Gaussian random vector
, i.e.,
\begin{equation*}
  \begin{split}
 \Bigl(\mathbb{E}\bigl[\|\me^{(2)}_{k \ell} v_{\ell}\|^{p}\bigr]\Bigr)^{1/p}
= \bigl(\mathbb{E}[|\me^{(2)}|^{p} \|v_{\ell}\|^{p}]\bigr)^{1/p} 
    \leq
\bigl(\mathbb{E}[|\me^{(2)}|^{p}])^{1/p} \times \|\mv\|_{2 \to
  \infty} \lesssim d^{1/2} n^{-1/2} \|\me^{(2)}_{k \ell}\|_{\psi_2} p^{1/2}.
  \end{split}
\end{equation*}
Therefore, by Lemma~2 and Corollary~7 in \cite{short_note_concentration}, we have
$$\Bigl\|\sum_{\ell=1}^n \me^{(2)}_{k \ell} v_{\ell}\Bigr\| \lesssim
\Bigl(\sum_{\ell=1}^n (d^{1/2} n^{-1/2} \|\me^{(2)}_{k
  \ell}\|_{\psi_2})^{2} (\log d + \log n)\Bigr)^{1/2}  \lesssim
d^{1/2} (\rho_n \log n)^{1/2}$$
with high probability. A union bound over all $k \in [n]$ yields
$\|\me^{(2)} \mv\| \lesssim d^{1/2}(\rho_n \log n)^{1/2}$ with high
probability. 
The bound for $\|\me^{(2)\top} \muu\|_{2 \to \infty}$ is identical
and is once again omitted.    
If we further assume $\{\me^{(i,2)}\}$ are independent, by almost identical proof we have $\|\frac{1}{m}\sum_{i=1}^m\me^{(i,2)}\mv (\mr^{(i)})^{-1}\|_{2\to\infty}\lesssim d^{1/2}(mn)^{-1/2}(n\rho_n)^{-1/2}\log^{1/2} n$
with high probability.

Combining the above bounds about $\me^{(1)}$ and $\me^{(2)}$, the bounds for $\me=\me^{(1)}+\me^{(2)}$ in Lemma~\ref{lemma:|E|_2,|UEV|F} can be derived.
\hspace*{\fill} \qedsymbol  

\subsection{Proof of Lemma~\ref{lemma:Uhat-UW}}
\label{Appendix:proof of lemma:Uhat-UW}
  We only prove the result for $\hat{\muu}^{(i)}
  \mw_{\muu}^{(i)} - \muu^{(i)}$ as the proof for $\hat{\mv}^{(i)}
  \mw_{\mv}^{(i)} - \mv^{(i)}$ is identical. For ease of exposition, we fix a value of $i$ and thereby 
drop the index $i$ from our matrices and quantities.

  First consider 
  the singular value decomposition of
  $\mpp$ as $\mpp=\muu^*\mSigma\mv^{*\top}$.
  Since $\muu^{*}$ spans the same invariant subspace as $\muu$, we
  have $\mathbf{U} \mathbf{U}^{\top}=\muu^{*}\muu^{*\top}$. Similarly,
  we also have $\mv \mv^{\top} = \mv^{*} \mv^{*\top}$.
  There thus exists $d \times d$ orthogonal matrices $\mw_1$ and
  $\mw_2$ such that $\muu^*=\muu\mw_1, \mv^*=\mv\mw_2$ and $\mr=\mw_1
  \mSigma\mw_2^\top$. We emphasize that $\mw_1$ and $\mw_2$ can
  depend on $i$. Indeed, while $\muu$ and $\mv$ are pre-specified and
  does not depend on the choice of $i$,
  $\muu^{*}$ and $\mv^{*}$ are defined via the singular value
  decomposition of $\mpp^{(i)}$. 
  
  Note that 
  $$
  \begin{aligned}
  	\hat\muu&
  	=\ma\hat\mv\hat\mSigma^{-1}
  	=\mpp\hat\mv\hat\mSigma^{-1}+\me\hat\mv\hat\mSigma^{-1}
  	=\muu\mr\mv^\top\hat\mv\hat\mSigma^{-1}+\me\hat\mv\hat\mSigma^{-1}\\
  	&=\muu\muu^\top\hat\muu+\muu\mr(\mv^\top\hat\mv\hat\mSigma^{-1}-\mr^{-1}\muu^\top\hat\muu)+\me\hat\mv\hat\mSigma^{-1}.
  \end{aligned}
  $$
  Hence for any $d \times d$ orthogonal matrices $\mw$ and $\tilde{\mw}$, we have
  \begin{equation}\label{eq:expansion_UhatW-U}
  	\begin{aligned}
  		\hat{\muu}\mw-\muu
  		&=\me\mv\mr^{-1}
  		+\underbrace{\muu(\muu^\top\hat\muu-\mw^\top)\mw}_{\mt_1}
  		+\underbrace{\muu\mr(\mv^\top\hat\mv\hat\mSigma^{-1}-\mr^{-1}\muu^\top\hat\muu)\mw}_{\mt_2}\\
  		&+\underbrace{\me\mv(\tilde\mw^\top\hat\mSigma^{-1}\mw-\mr^{-1})}_{\mt_3}
  		+\underbrace{\me(\hat\mv\tilde\mw-\mv)\tilde\mw^\top\hat\mSigma^{-1}\mw}_{\mt_4}.
  	\end{aligned}
  \end{equation}
 
    Now let $\mw_{\muu}$ and $\mw_{\mv}$ minimize $\|\hat{\muu} \mo -
    \muu\|_{F}$ and $\|\hat{\mv} \mo - \mv\|_{F}$ over all $d \times
    d$ orthogonal matrices $\mo$, respectively. 
	By Lemma \ref{lemma:T1}, Lemma~\ref{lemma:T2}, Lemma~\ref{lemma:T3} and Lemma~\ref{lemma:T4} we have, for these choices of
    $\mw = \mw_{\muu}$ and $\tilde{\mw} = \mw_{\mv}$, that
	$$
	\begin{aligned}
	\|\sum_{r=1}^4\mt_r\|	
	&\lesssim \|\mt_1\|+\|\mt_2\|+\|\mt_3\|+\|\mt_4\|
		\lesssim (n\rho_n)^{-1}\max\{1,d^{1/2}\rho_n^{1/2}(\log n)^{1/2}\},\\
	\|\sum_{r=1}^4\mt_r\|_{2\to\infty}	
	&\lesssim \|\mt_1\|_{2\to\infty}	+\|\mt_2\|_{2\to\infty}	+\|\mt_3\|_{2\to\infty}+\|\mt_4\|_{2\to\infty}	
	\lesssim d^{1/2}n^{-1/2}(n\rho_n)^{-1}\log n
	\end{aligned}
	$$
    with high probability. The proof is completed by defining $\mt = \mt_1 + \mt_2 + \mt_3+\mt_4$.
\hspace*{\fill} \qedsymbol

\subsection{Technical lemmas for $\mt_4$ in Lemma~\ref{lemma:Uhat-UW}}
\label{sec:mt_4}
We now present technical lemmas for bounding the term $\mt_4$ 
used in the above proof of Lemma~\ref{lemma:Uhat-UW}. Technical
lemmas for $\mt_1,\mt_2$ and $\mt_3$ are presented in
Section~\ref{sec:mt_1-mt_3}.  For ease of exposition we include the
index $i$ in the statement of these lemmas but we will generally drop
this index in the proofs.

Our bound for $\mt_4$ is based on a series of technical lemmas
with the most important being Lemma~\ref{lemma:||E(hatVtildeW-V)||2toinfty} which provides a
high-probability bound for
$\|\me(\hat\mv\tilde\mw-\mv)\|_{2\to\infty}$. Lemma~\ref{lemma:||E(hatVtildeW-V)||2toinfty}
is an adaptation of the leave-one-out analysis presented in Theorem~3.2 of \cite{xie2021entrywise}.
Leave-one-out arguments provide a
simple and elegant approach for handling the (often times) complicated 
dependencies between the rows of $\hat{\muu}$. See
\cite{abbe2020entrywise,chen2021spectral,javanmard_montanari,zhong_boumal,lei2019unified}
for other examples of leave-one-out analysis in the context of
random graphs inference, linear regression using lasso, and phase
synchronization. We can also prove Lemma~\ref{lemma:||E(hatVtildeW-V)||2toinfty} using the techniques in
\cite{cape2019signal,mao_sarkar} but this require a slightly stronger
assumption of $n \rho_n = \omega(\log^{c}n)$ for some
$c > 1$ as opposed to $n \rho_n = \Omega(\log n)$ in
the current paper. 

We first introduce some notations. Let $\ma = \ma^{(i)}$ be an observed
adjacency matrix and define the following collection of auxiliary matrices
$\ma^{[1]},\dots,\ma^{[n]}$ generated from $\ma$.
For each row index $h\in[n]$, the matrix
$\ma^{[h]}=(\ma^{[h]}_{k \ell})_{n\times n}$ is obtained by replacing the
entries in the $h$th row of $\ma$ with their expected values, i.e.,
$$
\ma^{[h]}_{k \ell}= \begin{cases}\ma_{k \ell}, & \text { if } k \neq h, \\ 
\mpp_{k \ell}, & \text { if } k=h .\end{cases}
$$
Denote the SVD of $\ma$ and $\ma^{[h]}$ as
\begin{gather*}
  \ma=\hat\muu\hat\mSigma\hat\mv^{\top}
  +\hat\muu_\perp\hat\mSigma_\perp\hat\mv_\perp^{\top}, \\
  \ma^{[h]}=\hat\muu^{[h]}\hat\mSigma^{[h]}\hat\mv^{[h]\top}
  +\hat\muu^{[h]}_\perp\hat\mSigma^{[h]}_\perp \hat\mv^{[h]\top}_\perp.
\end{gather*}
\begin{lemma1}
   \label{lemma:||hatV||2toinfty}
   Consider the setting in Lemma~\ref{lemma:Uhat-UW} for some fixed
   $i$ where, for ease of exposition, we will drop the index $i$ in
   all matrices. We then have
   $$
   \begin{aligned}
   	  & \|\hat\muu\|_{2\to\infty}
   \lesssim d^{1/2}n^{-1/2},
   \quad \|\hat\mv\|_{2\to\infty}
   \lesssim d^{1/2}n^{-1/2},\\
   &\|\hat\muu^{[h]}\|_{2\to\infty}
   \lesssim d^{1/2}n^{-1/2},
   \quad \|\hat\mv^{[h]}\|_{2\to\infty}
   \lesssim d^{1/2}n^{-1/2}
   \end{aligned}
   $$
   with high probability.
\end{lemma1}
\begin{proof}
Consider the Hermitian dilations 
$$
\begin{aligned}
\mpp'&=\begin{bmatrix} \bm{0} & \mpp \\
  \mpp^\top & \bm{0} \end{bmatrix}
  =\muu'\mSigma'\muu'^\top,\quad
  \ma'=\begin{bmatrix} \bm{0} & \ma \\
  \ma^\top & \bm{0} \end{bmatrix}
  =\hat\muu'\hat\mSigma'\hat\muu'^\top
  +\hat\muu'_\perp\hat\mSigma'_\perp\hat\muu'^\top_\perp,
\end{aligned}
$$
where we define
$$
\begin{aligned}
  &\muu'=\frac{1}{\sqrt{2}}\begin{bmatrix} \muu^* & \muu^* \\
  \mv^*& -\mv^* \end{bmatrix},\quad
 \hat\muu'=\frac{1}{\sqrt{2}}\begin{bmatrix} \hat\muu & \hat{\muu} \\
  \hat{\mv} & -\hat\mv \end{bmatrix},\quad
  \hat\muu'_\perp=\frac{1}{\sqrt{2}}\begin{bmatrix} \hat\muu_\perp & \hat{\muu}_\perp \\
  \hat{\mv}_{\perp} & -\hat\mv_\perp \end{bmatrix} ,\\
 &\mSigma'=\begin{bmatrix} \mSigma & \bm{0} \\
  \bm{0} & -\mSigma \end{bmatrix}, \quad
  \hat\mSigma'=\begin{bmatrix} \hat\mSigma & \bm{0} \\
  \bm{0} & -\hat{\mSigma} \end{bmatrix},\quad
 \hat\mSigma'_\perp=\begin{bmatrix} \hat\mSigma_\perp & \bm{0} \\
  \bm{0} & -\hat{\mSigma}_{\perp} \end{bmatrix}.
\end{aligned}
$$
Then from Lemma~B.5 in \cite{xie2021entrywise} (see also Theorem~2.1
in \cite{abbe2020entrywise}), we have
$$
\begin{aligned}
	\max\{\|\hat\muu\|_{2\to\infty},\|\hat\mv\|_{2\to\infty}\}
	= \|\hat\muu'\|_{2\to\infty}
	&\lesssim \|\muu'\|_{2\to\infty}\\
	&\lesssim \max\{\|\muu^*\|_{2\to\infty},\|\mv^*\|_{2\to\infty}\}\\
	&\lesssim \max\{\|\muu\|_{2\to\infty},\|\mv\|_{2\to\infty}\}
	\lesssim d^{1/2}n^{-1/2}
\end{aligned}
$$
with high probability.
The analysis of $\|\hat\muu^{[h]}\|_{2\to\infty}$ and
$\|\hat\mv^{[h]}\|_{2\to\infty}$ follows the same argument and is thus omitted.
\end{proof}

\begin{lemma1}
   \label{lemma:|| sinTheta(hatV[h],hatV)||}
   Consider the setting in Lemma~\ref{lemma:||hatV||2toinfty}. We then have
   $$
   \begin{aligned}
   &\|\sin\Theta(\hat\mv^{[h]},\hat\mv)\|
\lesssim d^{1/2} n^{-1/2}(n \rho_n)^{-1/2} \log^{1/2}{n}
   \end{aligned}
   $$
   with high probability.
\end{lemma1}

\begin{proof}
From Eq.~\eqref{eq:Weyl} we have $\sigma_{d+1}(\ma)\lesssim
E_n$ with high probability.
By the construction of $\ma^{[h]}$ and Lemma~\ref{lemma:|E|_2,|UEV|F}, it follows that
\begin{equation*} \label{eq:||A[h]-A||}
	\|\ma^{[h]}-\ma\|
\leq \big(\sum_{\ell=1}^n\me_{h\ell}^2\big)^{1/2}
\leq \|\me\|_{2\to\infty}
\leq \|\me\|
\lesssim (n \rho_n)^{1/2}
\end{equation*}
with high probability. We thus obtain
\begin{equation*}
	\|\ma^{[h]}-\mpp\|
	\leq \|\ma-\ma^{[h]}\|+\|\me\|
	\lesssim (n \rho_n)^{1/2}
\end{equation*}
with high probability. Therefore, by Weyl's inequality for singular
values (see e.g., Problem~III.6.13 in \cite{bhatia2013matrix}), we have
\begin{equation*}
	\max_{k\in[n]}|\sigma_k(\ma^{[h]})-\sigma_k(\mpp)|
	\leq \|\ma^{[h]}-\mpp\|
	\lesssim (n \rho_n)^{1/2} 
\end{equation*}
with high probability. As $\sigma_k(\mpp)=\sigma_k(\mr)\asymp n\rho_n$ for all $k\leq d$
and $\sigma_{k}(\mpp)=0$ otherwise, we have with high probability that
\begin{equation}
  \label{eq:eigen_Ah}
  \begin{split}
  &\sigma_k(\ma^{[h]})\asymp n \rho_n \quad \text{for all $1 \leq k \leq
  d$}, \\ & \sigma_k(\ma^{[h]})\lesssim (n \rho_n)^{1/2} \quad \text{for all
  $k\geq d+1$}.
  \end{split}
\end{equation}
Therefore, by Wedin's $\sin\Theta$ Theorem (see e.g., Theorem~4.4 in \cite{stewart_sun}), 
\begin{equation}\label{eq:sinTheta_decomposition}
	\begin{aligned}
	\|\sin\Theta(\hat\mv^{[h]},\hat\mv)\|
	&\leq \frac{\max\{\|(\ma^{[h]}-\ma)\hat\mv^{[h]}\|,\|\hat\muu^{[h]\top}(\ma^{[h]}-\ma)\|\}}{\sigma_d(\ma^{[h]})-\sigma_{d+1}(\ma)}\\
	&\lesssim
      \frac{\max\{\|(\ma^{[h]}-\ma)\hat\mv^{[h]}\|_F,\|\hat\muu^{[h]\top}(\ma^{[h]}-\ma)\|_F\}}{n
      \rho_n}
\end{aligned}
\end{equation}
with high probability.

From Lemma~\ref{lemma:|E|_2,|UEV|F} and Lemma~\ref{lemma:||hatV||2toinfty}, we have
\begin{equation}\label{eq:sinTheta_decomposition1}
	\begin{aligned}
		\|\hat\muu^{[h]\top}(\ma^{[h]}-\ma)\|_F
		=\Big(\sum_{\ell=1}^n\sum_{r=1}^d(\me_{h \ell}\hat\muu^{[h]}_{hr})^2\Big)^{1/2}
		\leq \|\me\|_{2\to\infty}\cdot \|\hat\muu^{[h]}\|_{2\to\infty}
		\leq \|\me\|\cdot \|\hat\muu^{[h]}\|_{2\to\infty}
		\lesssim d^{1/2} \rho_n^{1/2}
	\end{aligned}
\end{equation}
with high probability. We now consider $\|(\ma^{[h]}-\ma)\hat\mv^{[h]}\|_F$. Write
\begin{equation}\label{eq:sinTheta_decomposition2}
	\begin{aligned}
		\|(\ma^{[h]}-\ma)\hat\mv^{[h]}\|_F
		=\Big\|\sum_{\ell=1}^n\me_{h \ell} \hat v^{[h]}_{\ell}\Big\| =
      \Big\|\sum_{\ell=1}^n(\me^{(1)}_{h \ell} + \me^{(2)}_{h \ell}) \hat v^{[h]}_{\ell}\Big\| 
      ,
	\end{aligned}
\end{equation}
where $\hat{v}^{[h]}_{\ell}$ represents the $\ell$th row of
$\hat\mv^{[h]}$ and $\me^{(1)}_{h \ell}$ and $\me^{(2)}_{h \ell}$
denote the $h \ell$th element of $\me^{(i,1)}$ and $\me^{(i,2)}$; recall that we had fixed an $i \in [m]$ and use $\me$ to denote
$\me^{(i)} = \me^{(i,1)} + \me^{(i,2)}$. 
For any $t\geq 1$, define the events
$$
\begin{aligned}
&\mathcal{E}_1=\Big\{\mathbf{A}:\Big\|\sum_{l=1}^n\me^{(1)}_{h\ell} \hat v^{[h]}_{\ell}\Big\|
\leq C (t^2 \|\hat\mv^{[h]}\|_{2 \rightarrow \infty}+  
  \rho_{n}^{1/2} t \|\hat\mv^{[h]}\|_{F})\Big\}, \\
&\mathcal{E}_2=\Big\{\mathbf{A}:\Big\|\sum_{l=1}^n\me^{(2)}_{h\ell} \hat v^{[h]}_{\ell}\Big\|
\leq C \rho_n^{1/2} t\|\hat\mv^{[h]}\|_{F}\Big\}.
\end{aligned}
$$
Now the $h$th row of $\me$ is independent of $\hat\mv^{[h]}$ and
hence, by Lemma~B.1 and Lemma~B.2 in \cite{xie2021entrywise}, there exists some finite constant $C>0$ that can depend on $\{C_1,
C_2, C_3\}$ in Assumption~\ref{ass:general} but does not depend on $n,
m$ and $\rho_n$, and
 for any $t\geq 1$ we have
\begin{gather*}
\mathbb{P}(\mathcal{E}_1)=\sum_{\mathbf{A}^{[h]}}
\mathbb{P}(\mathcal{E}_1 \mid \mathbf{A}^{[h]})
\mathbb{P}(\mathbf{A}^{[h]}) \geq \sum_{\mathbf{A}^{[h]}}(1-28
e^{-t^2}) \mathbb{P}(\mathbf{A}^{[h]})=1-28 e^{-t^2}, \\
\mathbb{P}(\mathcal{E}_2) = \sum_{\mathbf{A}^{[h]}}
\mathbb{P}(\mathcal{E}_2 \mid \mathbf{A}^{[h]})
 \mathbb{P}(\mathbf{A}^{[h]}) \geq 1 - 2(d+1)e^{-t^2}.
\end{gather*}
Thus with sufficiently large $t \asymp (\log n)^{1/2}$, by Lemma~\ref{lemma:||hatV||2toinfty} and the assumption $n\rho_n=\Omega(\log n)$ we have
\begin{equation*}
	\begin{aligned}
		\Big\|\sum_{l=1}^n\me^{(1)}_{h\ell} \hat v^{[h]}_{\ell}\Big\|
		&\lesssim \log n\|\hat\mv^{[h]}\|_{2 \rightarrow \infty}+  
  (\rho_{n}\log n)^{1/2}  \|\hat\mv^{[h]}\|_{F}
  \\&\lesssim \log n\|\hat\mv^{[h]}\|_{2 \rightarrow \infty}+  
  (n\rho_{n}\log n)^{1/2}  \|\hat\mv^{[h]}\|_{2 \rightarrow \infty}
  \lesssim d^{1/2}(\rho_n \log n)^{1/2},\\
  \Big\|\sum_{l=1}^n\me^{(2)}_{h\ell} \hat v^{[h]}_{\ell}\Big\|
  &\lesssim (\rho_n \log n)^{1/2}\|\hat\mv^{[h]}\|_{F}
  \lesssim d^{1/2}(\rho_n \log n)^{1/2}
  	\end{aligned}
\end{equation*}
with high probability.
We therefore have
\begin{equation}\label{eq:sinTheta_decomposition2a}
\begin{aligned}
\Big\|\sum_{\ell=1}^n\me_{h\ell} \hat v^{[h]}_{\ell}\Big\|
\leq \Big\|\sum_{l=1}^n\me^{(1)}_{h\ell} \hat v^{[h]}_{\ell}\Big\|
+\Big\|\sum_{l=1}^n\me^{(2)}_{h\ell} \hat v^{[h]}_{\ell}\Big\|
\lesssim d^{1/2}(\rho_n \log n)^{1/2}
\end{aligned}
\end{equation}
with high probability. Combining Eq.~\eqref{eq:sinTheta_decomposition}, Eq.~\eqref{eq:sinTheta_decomposition1}, Eq.~\eqref{eq:sinTheta_decomposition2} and Eq.~\eqref{eq:sinTheta_decomposition2a}, we obtain
$$
\|\sin\Theta(\hat\mv^{[h]},\hat\mv)\|
\lesssim d^{1/2} n^{-1/2}(n \rho_n)^{-1/2} \log^{1/2}{n}
$$
with high probability as desired.
\end{proof}

\begin{lemma1}
   \label{lemma:||e_h^top E (hat V[h]hat V[h]^topV-V)||}
   Consider the setting in Lemma~\ref{lemma:||hatV||2toinfty}. We then have
   $$
   \|e^\top_h\me(\hat\mv^{[h]}\hat\mv^{[h]\top}\mv-\mv)\|
   \lesssim d^{1/2}n^{-1/2}\log {n}
   $$
   with high probability.
\end{lemma1}

\begin{proof}
From the proof of Lemma~\ref{lemma:|| sinTheta(hatV[h],hatV)||} (see Eq.~\eqref{eq:eigen_Ah}) we
have
\begin{equation*}
  \begin{split}
  &\sigma_k(\ma^{[h]})\asymp n \rho_n \quad \text{for all $1 \leq k \leq
  d$}, \\ & \sigma_k(\ma^{[h]})\lesssim (n \rho_n)^{1/2} \quad \text{for all
  $k\geq d+1$}.
  \end{split}
\end{equation*}
Let $\mz^{[h]}=\hat\mv^{[h]}\hat\mv^{[h]\top}\mv-\mv$.
Thus by
Wedin's $\sin\Theta$ Theorem, we have
\begin{equation}
\label{eq:Z1}
	\begin{aligned}
	\|\mz^{[h]}\|
	 =\|(\hat\mv^{[h]}\hat\mv^{[h]\top}-\mi)\mv\|=
	\|\sin\Theta(\hat\mv^{[h]},\mv)\|
	=\|\sin\Theta(\hat\mv^{[h]},\mv^*)\|
	\leq \frac{\|\ma^{[h]}-\mpp\|}{\sigma_d(\ma^{[h]})-\sigma_{d+1}(\mpp)}
	\lesssim (n \rho_n)^{-1/2}
\end{aligned}
\end{equation}
with high probability. 
Let $\mw^{[h]}$ be
orthogonal Procrustes problem between $\hat\mv^{[h]}$ and $\mv$. Then we have
\begin{equation*}
	\begin{aligned}
		\|\hat\mv^{[h]\top}\mv-\mw^{[h]}\| 
		\leq \|\sin\Theta(\hat\mv^{[h]},\mv)\|^2
		\lesssim (n \rho_n)^{-1}
	\end{aligned}
\end{equation*}
with high probability. Finally, by Lemma~\ref{lemma:||hatV||2toinfty}, we have
\begin{equation}\label{eq:Z2}
\begin{aligned}
	\|\mz^{[h]}\|_{2\to\infty}
	&\leq \|\hat\mv^{[h]}\|_{2\to\infty} + \|\mv\|_{2\to\infty} 
	\lesssim d^{1/2}n^{-1/2}
\end{aligned}
\end{equation}
with high probability. We now follow the same argument as that
for deriving Eq.~\eqref{eq:sinTheta_decomposition2a}.
First define the events
$$
\begin{aligned}
&\mathcal{E}_1=\Big\{\mathbf{A}:\|e^\top_h\me^{(1)}\mz^{[h]}\|
\leq C (t^2 \|\mz^{[h]}\|_{2 \rightarrow \infty}+  
  \rho_{n}^{1/2} t\|\mz^{[h]}\|_{F})\Big\}. \\
 &\mathcal{E}_2=\Big\{\mathbf{A}:\|e^\top_h\me^{(2)}\mz^{[h]}\|
\leq C \rho_n^{1/2} t \|\mz^{[h]}\|_{F}\Big\},
\end{aligned}
$$
By the definition of $\mz^{(h)}$, $e_h^{\top} \me$
and $\mz^{(h)}$ are independent. 
Once again by Lemma~B.1 and
Lemma~B.2 in
\cite{xie2021entrywise}, there exists some finite constant $C>0$ that can depend on $\{C_1,
C_2, C_3\}$ in Assumption~\ref{ass:general} but does not depend on $n,
m$ and $\rho_n$, such that
for any $t\geq 1$ we have
$$
\mathbb{P}(\mathcal{E}_1) \geq 1-28 e^{-t^2}, \quad
\mathbb{P}(\mathcal{E}_2) \geq 1 - 2(d+1)e^{-t^2}.
$$
Thus with sufficiently large $t \asymp (\log n)^{1/2}$, by Eq.~\eqref{eq:Z1}, Eq.~\eqref{eq:Z2} and the assumption $n\rho_n=\Omega(\log n)$ we have
\begin{equation*}
	\begin{aligned}
		\Big\|e^\top_h\me^{(1)}\mz^{[h]}\Big\|
		&\lesssim \log n\|\hat\mz^{[h]}\|_{2 \rightarrow \infty}+  
  (\rho_{n}\log n)^{1/2}  \|\hat\mz^{[h]}\|_{F}\\
  &\lesssim \log n\|\hat\mz^{[h]}\|_{2 \rightarrow \infty}+  
  (d\rho_{n}\log n)^{1/2}  \|\hat\mz^{[h]}\|
  \\ & \lesssim d^{1/2} n^{-1/2} \log n + d^{1/2} (\rho_n \log
       n)^{1/2} (n \rho_n)^{-1/2}
\\
  \Big\|e^\top_h\me^{(2)}\mz^{[h]}\Big\|
  &\lesssim d^{1/2} (\rho_n \log n)^{1/2}\|\hat\mz^{[h]}\|_{F}
  \lesssim d^{1/2}(\rho_n \log n)^{1/2} (n \rho_n)^{-1/2}
  	\end{aligned}
\end{equation*}
with high probability. Adding the above two bounds we obtain
\begin{equation}\label{eq:sinTheta_decomposition22}
\begin{aligned}
\|e^\top_h\me\mz^{[h]}\|
\lesssim d^{1/2} n^{-1/2} \log n + d^{1/2} (\rho_n \log n)^{1/2} (n
\rho_n)^{-1/2} \lesssim d^{1/2} n^{-1/2} \log n
\end{aligned}
\end{equation}
with high probability.
\end{proof}

\begin{lemma1}
   \label{lemma:||E(hatVtildeW-V)||2toinfty}
   Consider the setting in Lemma~\ref{lemma:Uhat-UW}. We then have
   $$
   \|\me (\hat\mv\mw_\mv-\mv)\|_{2\to\infty}
   \lesssim d^{1/2}n^{-1/2}\log {n}
   $$
   with high probability.
\end{lemma1}

\begin{proof}
We will drop the dependency on the index $i$ from our matrices. 
First we have
\begin{equation}\label{eq:eh E(hatVW-V)_decomposition}
	\begin{aligned}
		\|e_h^\top\me(\hat\mv\mw_\mv-\mv)\|
		\leq \|e_h^\top\me\hat\mv(\mw_\mv-\hat\mv^\top\mv)\|
		+\|e_h^\top\me(\hat\mv\hat\mv^\top-\hat\mv^{[h]}\hat\mv^{[h]\top})\mv\|
		+\|e_h^\top\me(\hat\mv^{[h]}\hat\mv^{[h]\top}\mv-\mv)\|
	\end{aligned}
\end{equation}
for each row
$h \in [n]$. 
We now bound each term in the right hand side of the above display.
For the first term we have
\begin{equation}
  \label{eq:lemma_t3_technical1}
	\begin{aligned}
		\|e_h^\top\me\hat\mv(\mw_\mv-\hat\mv^\top\mv)\|
		&\leq \|e_h^{\top} \me \hat\mv\| \cdot
        \|\mw_\mv-\hat\mv^\top\mv\|.
	\end{aligned}
\end{equation}
Now, by Lemma~2 in \cite{abbe2020entrywise}, we know that 
$\hat{\mv}^{\top} \mv$ is invertible and $\|(\hat{\mv}^{\top} \mv)^{-1}\| \leq 2$ with high probability. 
Then for $e_h^{\top} \me \hat{\mv}$ we have
\begin{equation*}
  \begin{split}
    \|e_h^{\top} \me \hat{\mv}\| &=
    \|e_h^{\top}  \me (\hat{\mv} \hat{\mv}^{\top} - \hat{\mv}^{[h]}
                                   \hat{\mv}^{[h]\top} + \hat{\mv}^{[h]}
                                   \hat{\mv}^{[h]\top}) \mv
                                   (\hat{\mv}^{\top} \mv)^{-1} \| \\
    & \leq \|e_h^{\top}  \me (\hat{\mv} \hat{\mv}^{\top} - \hat{\mv}^{[h]}
      \hat{\mv}^{[h]\top}) \mv (\hat{\mv}^{\top} \mv)^{-1}\| +
      \|e_h^{\top} \me\hat{\mv}^{[h]}  \hat{\mv}^{[h]\top} \mv
                                   (\hat{\mv}^{\top} \mv)^{-1} \|
\\ & \leq \bigl(\|\me\| \cdot \|\hat{\mv} \hat{\mv}^{\top} - \hat{\mv}^{[h]}
           \hat{\mv}^{[h]\top}\| +
           \|e_h^{\top} \me\hat{\mv}^{[h]}\|\bigr)\cdot \|(\hat{\mv}^{\top} \mv)^{-1} \|.
        \end{split}
\end{equation*}
In Eq.~\eqref{eq:sinTheta_decomposition2a} we have $
  	\|e_h^{\top} \me \hat{\mv}^{[h]}\|
\lesssim
d^{1/2} (\rho_n \log
n)^{1/2} 
$
with high probability. Combining this bound, Lemma~\ref{lemma:|E|_2,|UEV|F} and
Lemma~\ref{lemma:|| sinTheta(hatV[h],hatV)||} we obtain
\begin{equation*}
  \begin{split}
     \|e_h^{\top} \me \hat{\mv}\| & \leq 2 \bigl(\|e_h^{\top} \me \hat{\mv}^{[h]}\| + 
        \|\me\| \cdot \|\sin \Theta(\hat{\mv}^{[h]}, \hat{\mv})\|\bigr) 
         \lesssim d^{1/2} (\rho_n \log n)^{1/2}
\end{split}
\end{equation*}
with high probability. Substituting Eq.~\eqref{eq:sin_thetav} and the above bound into
Eq.~\eqref{eq:lemma_t3_technical1} yields
\begin{equation*}
\begin{split}
  \|e_h^\top\me\hat\mv(\mw_\mv-\hat\mv^\top\mv)\|
		& \lesssim d^{1/2} n^{-1/2} (n \rho_n)^{-1/2} \log^{1/2}{n}
\end{split}
        \end{equation*}
 with high probability. 
 For the second term, by Lemma~\ref{lemma:|E|_2,|UEV|F} and Lemma~\ref{lemma:|| sinTheta(hatV[h],hatV)||} we have
\begin{equation*}
	\begin{aligned}
		\|e_h^\top\me(\hat\mv\hat\mv^\top-\hat\mv^{[h]}\hat\mv^{[h]\top})\mv\|
		&\leq 2 \|\me\|\cdot \|\sin\Theta (\hat\mv^{[h]},\hat\mv)\|
\lesssim  d^{1/2} n^{-1/2} \log^{1/2}{n} 
	\end{aligned}
\end{equation*}
with high probability. For the third term, by Lemma~\ref{lemma:||e_h^top E (hat V[h]hat V[h]^topV-V)||} we have
\begin{equation*}
	\|e^\top_h\me(\hat\mv^{[h]}\hat\mv^{[h]\top}\mv-\mv)\|
   \lesssim d^{1/2}n^{-1/2}\log{n}
\end{equation*}
with high probability. Combining the above bounds for the terms on the
right hand side of Eq.~\eqref{eq:eh E(hatVW-V)_decomposition}, we obtain
the bound for $\|\me(\hat\mv\mw_\mv-\mv)\|_{2\to\infty}$ as claimed.
\end{proof}

\begin{lemma1}
  \label{lemma:T4}
  Consider the setting of Lemma \ref{lemma:Uhat-UW}. Define
  $$\mt_4 = \me (\hat\mv \mw_\mv-\mv)\mw_\mv^\top(\hat\mSigma )^{-1}\mw_\muu.$$
  We then have
  $$
  \begin{aligned}
  	  \|\mt_4 \|\lesssim (n \rho_n)^{-1}, \quad \text{and}
    \quad \|\mt_4   \|_{2\to\infty}\lesssim d^{1/2}n^{-1/2} (n
    \rho_n)^{-1} \log n
  \end{aligned}
  $$
  with high probability.
\end{lemma1}
 
\begin{proof}
By Lemma~\ref{lemma:|E|_2,|UEV|F}, Eq.~\eqref{eq:Weyl} and Eq.~\eqref{eq:T2_parta}, we have
$$
\|\mt_4\|
\leq \|\me\|\cdot \|\hat\mv\mw_\mv-\mv\|\cdot \|\hat\mSigma^{-1}\|
\lesssim (n \rho_n)^{-1}
$$
with high probability.
By Lemma~\ref{lemma:||E(hatVtildeW-V)||2toinfty} and Eq.~\eqref{eq:T2_parta}, we have
$$
\|\mt_4\|_{2\to\infty}
\leq \|\me(\hat\mv\mw_\mv-\mv)\|_{2\to\infty}\cdot \|\hat\mSigma^{-1}\|
\lesssim d^{1/2} n^{-1/2} (n \rho_n)^{-1} \log{n}
$$
with high probability.
\end{proof}

\section{Remaining Technical Lemmas}
\label{Appendix:C}

\subsection{Technical lemmas for $\mt_1,\mt_2$ and $\mt_3$ in Lemma~\ref{lemma:Uhat-UW}}
\label{sec:mt_1-mt_3}
We now present upper bounds for $\mt_1,
\mt_2$ and $\mt_3$ as used in the proof of
Lemma~\ref{lemma:Uhat-UW}; an upper bound for $\mt_4$
was given in Section~\ref{sec:mt_4}.
For ease of exposition, we drop the index $i$ from our matrices and quantities. 
in the proofs. 
\begin{lemma1}
  \label{lemma:T1}
  Consider the setting of Lemma \ref{lemma:Uhat-UW}. Define
  $$\mt_1  = \muu (\muu^{\top} \hat{\muu}  -
  \mw_{\muu}^{ \top}) \mw_{\muu} .$$
  We then have
  $$
  \begin{aligned}
  	  \|\mt_1 \|\lesssim (n \rho_n)^{-1}, \quad \text{and}
       \quad \|\mt_1  \|_{2\to\infty}\lesssim d_i^{1/2} n^{-1/2} (n \rho_n)^{-1}
  \end{aligned}
  $$
  with high probability.
\end{lemma1}
 
\begin{proof}
First by Lemma~\ref{lemma:|E|_2,|UEV|F}, we have
$\|\me\|\lesssim(n\rho_n)^{1/2}$ with high probability,
hence by applying perturbation theorem for singular values (see Problem~III.6.13 in \cite{bhatia2013matrix}) we have 
\begin{equation}
  \label{eq:Weyl}
\max_{1 \leq j \leq n} |\sigma_{j}(\ma)-\sigma_{j}(\mpp)|\leq
\left\|\me\right\|\lesssim (n \rho_n)^{1/2}
\end{equation}
with high probability. Since $\sigma_{k}(\mpp)=\sigma_{k}(\mr)\asymp
n\rho_n$ for all $k \leq d$ and $\sigma_{k}(\mpp) = 0$ otherwise,
we have that, with high probability, $\sigma_{k}(\ma)\asymp S_n$ for all $k \leq d$
and $\sigma_{k}(\ma) \lesssim (n \rho_n)^{1/2}$ for all $k \geq d +
1$. 
Then by Wedin's $\sin\Theta$ Theorem (see e.g., Theorem~4.4 in \cite{stewart_sun}), we have
\begin{equation}
  \label{eq:wedin_sintheta}
\begin{aligned}
	\max\{\|\sin\Theta(\hat\muu,\muu)\|,\|\sin\Theta(\hat\mv,\mv)\|\}
	&
	=
    \max\{\|\sin\Theta(\hat\muu,\muu^*)\|,\|\sin\Theta(\hat\mv,\mv^{*})\|\}
    \\&
	\leq\frac{\|\me\|}{\sigma_d(\ma) - \sigma_{d+1}(\mpp)} 
	\lesssim (n \rho_n)^{-1/2}
    \end{aligned}
\end{equation}
with high probability. 
Now recall that $\mw_{\muu}$ is the solution of orthogonal Procrustes
problem between $\hat{\muu}$ and $\muu$, i.e., $\mw_{\muu} = \mo_2
\mo_1^{\top}$ where $\mo_1 \cos\Theta(\muu, \hat{\muu}) \mo_2^{\top}$ is the
singular value decomposition of $\muu^{\top} \hat{\muu}$. We therefore have
\begin{equation}
  \label{eq:wedin_sintheta2}
\begin{aligned}
	\|\muu^{\top} \hat{\muu} - \mw_{\muu}^\top\|
	&
	= \|\cos \Theta(\muu, \hat{\muu}) - \mi\| 
\\&
	= \max_{1 \leq j \leq d} 1 - \sigma_j(\muu^{\top} \hat{\muu})
	\leq \max_{1 \leq j \leq d} 1-\sigma_{j}^2(\muu^{\top} \hat{\muu}) 
	= \|\sin\Theta(\hat\muu,\muu)\|^2 \lesssim (n \rho_n)^{-1}
    \end{aligned}
\end{equation}
with high probability. We therefore obtain
\begin{equation*}
  \begin{aligned}
	&\|\mt_1\| \leq \|\muu^\top\hat\muu-\mw_{\muu}^\top\|
	\lesssim (n \rho_n)^{-1} \\
	&\|\mt_1\|_{2\to\infty} \leq \|\muu\|_{2\to\infty}\cdot\|\muu^\top\hat\muu-\mw_{\muu}^\top\|
	\lesssim d^{1/2}n^{-1/2} (n \rho_n)^{-1}
    \end{aligned}
\end{equation*}
with high probability.
\end{proof}

\begin{lemma1}
  \label{lemma:T2}
  Consider the setting of Lemma \ref{lemma:Uhat-UW}. Define
  \begin{gather*}
    \begin{aligned}
   	\mt_2  &= \muu\mr \big(\mv^\top\hat\mv \hat\mSigma ^{-1}-\mr^{-1}\muu^\top\hat\muu \big)\mw_\muu .
   	      \end{aligned}
    \end{gather*}
    Let $\vartheta_{n} = \max\{1,d^{1/2}\rho_n^{1/2}(\log n)^{1/2}\}$.
 We then have
 $$
 \begin{aligned}
 	&\|\mt_2 \|\lesssim (n\rho_n)^{-1} \vartheta_n, \quad
   &\|\mt_2  \|_{2 \to \infty} \lesssim d^{1/2}n^{-1/2}
   (n\rho_n)^{-1} \vartheta_n
 \end{aligned}
 $$
with high probability.
\end{lemma1}

\begin{proof}
  
  Let
  $\tilde{\mt}_{2}=\mv^{*\top}\hat\mv\hat\mSigma^{-1}-\mSigma^{-1}\muu^{*\top}\hat\muu$
  and note that $\mv^\top\hat\mv\hat\mSigma^{-1}-\mr^{-1}\muu^\top\hat\muu
    = \mw_{2} \tilde\mt_2$. We then have
	$$
	\begin{aligned}
		\mSigma\tilde{\mt}_{2}\hat\mSigma
		=\mSigma\mv^{*\top}\hat\mv-\muu^{*\top}\hat\muu\hat\mSigma
		=\muu^{*\top}\mpp\hat\mv-\muu^{*\top}\ma\hat\mv
		=-\muu^{*\top}\me(\hat\mv\mw_\mv-\mv)\mw_\mv^\top
		-\muu^{*\top}\me\mv\mw_\mv^\top.
	\end{aligned}
	$$
   We now bound each term in the right hand side of the above
   display. First note that, by Lemma~\ref{lemma:|E|_2,|UEV|F} we have
   \begin{equation}
     \label{eq:T2_part0}
	\|\muu^{*\top}\me\mv\mw_\mv^\top\|
	\leq \|\muu^\top\me\mv\|
	\lesssim d^{1/2}\rho_n^{1/2}(\log n)^{1/2}
\end{equation}
with high probability. Next, by Eq.~\eqref{eq:wedin_sintheta}, we have
  $ \|\sin\Theta(\hat\mv,\mv)\| \lesssim (n \rho_n)^{-1/2}$
  with high probability and hence, using the same argument as that for deriving
  Eq.~\eqref{eq:wedin_sintheta2}, we have
  \begin{equation}
    \label{eq:sin_thetav}
    \|\hat\mv^{\top} {\mv} - \mw_{\mv}\|\lesssim (n\rho_n)^{-1}
    \end{equation}
with high probability. We therefore have
\begin{equation}
  \label{eq:T2_parta}
\begin{aligned}
	\|\hat\mv\mw_\mv-\mv\|
	&
	\leq \|(\mi-\mv\mv^{\top})\hat\mv\|+\|\mv\|\cdot\|\hat\mv^{\top}\mv-\mw_{\mv}\|
	\\&
	\leq \|\sin\Theta(\hat\mv,\mv)\|+\|\mv\|\cdot\|\hat\mv^{\top}\mv-\mw_{\mv}\|
    \lesssim (n\rho_n)^{-1/2}
\end{aligned}
\end{equation}
with high probability. Lemma~\ref{lemma:|E|_2,|UEV|F} and
Eq.~\eqref{eq:T2_parta} then imply
\begin{equation}
  \label{eq:T2_partb0}
	\|\muu^{*\top}\me(\hat\mv\mw_\mv-\mv)\mw_\mv^\top\|
	\leq \|\me\|\cdot\|\hat\mv\mw_\mv-\mv\|
	\lesssim 1
  \end{equation}
with high probability.

Combining Eq.~\eqref{eq:T2_part0} and Eq.~\eqref{eq:T2_partb0} we have $\|\mSigma\tilde{\mt}_{2}\hat\mSigma\|\lesssim \vartheta_n$ with high probability, and hence
$$
\begin{aligned}
		\|\tilde{\mt}_{2}\|
		\leq \|\mSigma\tilde{\mt}_{2}\hat\mSigma\|
		\cdot \|\mSigma^{-1}\|\cdot \|\hat\mSigma^{-1}\|
		\lesssim (n\rho_n)^{-2}\vartheta_n\end{aligned}
$$
with high probability. In summary we obtain
$$
\begin{aligned}
	&\|\mt_2\|
	\leq \|\mr\|\cdot\|\tilde{\mt}_{2}\|
	\lesssim (n\rho_n)^{-1}\vartheta_n\\
	&\|\mt_2\|_{2\to\infty}
	\leq \|\muu\|_{2\to\infty}\cdot\|\mr\|\cdot\|\tilde{\mt}_{2}\|
	\lesssim d^{1/2}n^{-1/2} (n\rho_n)^{-1}\vartheta_n
\end{aligned}
$$
with high probability.
\end{proof}

\begin{lemma1}
   \label{lemma:T3}
   Consider the setting of Lemma \ref{lemma:Uhat-UW}. Define 
   $$\mt_3  = \me  \mv \bigl(\mw_{\mv}^{ \top}
   \hat{\mSigma} ^{-1} \mw_{\muu}  - \mr ^{-1}\bigr)$$
   Let $\vartheta_{n} = \max\{1,d^{1/2}\rho_n^{1/2}(\log n)^{1/2}\}$.
   We then have
   $$
	\begin{aligned}
		&\|\mt_3 \|\lesssim  (n\rho_n)^{-3/2} \vartheta_n, 
         \quad &\|\mt_3 \|_{2\to\infty}\lesssim
        d^{1/2}n^{-1/2} (n\rho_n)^{-3/2}(\log n)^{1/2} \vartheta_n
	\end{aligned}
	$$
    with high probability.
\end{lemma1}

\begin{proof}
  Let
  $\tilde{\mt}_{3}=\mw_2^\top\mw_{\mv}^\top\hat\mSigma^{-1}-\mSigma^{-1}\mw_1^\top\mw_{\muu}^\top$
  where $\mw_1$ and $\mw_2$ are defined in the proof of Lemma~\ref{lemma:Uhat-UW}. Note that $\mw_{\mv}^{\top} \hat{\mSigma}^{-1} \mw_{\muu} - \mr^{-1}
  = \mw_{2} \tilde{\mt}_{3} \mw_{\muu}$. We then have
	$$
	\begin{aligned}
		\mSigma\tilde{\mt}_{3}\hat\mSigma
		&=\mSigma\mw_2^\top\mw_{\mv}^\top-\mw_1^\top\mw_{\muu}^\top\hat\mSigma\\
		&=\mSigma\mw_2^\top(\mw_{\mv}^\top-\mv^{\top}\hat\mv)+(\mSigma\mv^{*\top}\hat\mv-\muu^{*\top}\hat\muu\hat\mSigma)+\mw_1^\top(\muu^{\top}\hat\muu-\mw_{\muu}^\top)\hat\mSigma.
	\end{aligned}
	$$
    We now bound each term in the right hand side of the above display.
    First recall Eq.~\eqref{eq:sin_thetav}. We then have
    \begin{equation}
      \label{eq:term3_a}
	\begin{aligned}
		\|\mSigma\mw_2^\top(\mw_{\mv}^\top-\mv^{\top}\hat\mv)\|
		&\leq\|\mSigma\|\cdot\|\mw_{\mv}^{\top}-\mv^{\top}\hat\mv\|
		\lesssim n\rho_n\cdot (n \rho_{n})^{-1}
		\lesssim 1
	\end{aligned}
    \end{equation}
    with high probability. 
	For the second term, we have
	$$
	\begin{aligned}
	\mSigma\mv^{*\top}\hat\mv-\muu^{*\top}\hat\muu\hat\mSigma
	=\muu^{*\top}\mpp\hat\mv-\muu^{*\top}\ma\mv
	=-\muu^{*\top}\me\hat\mv
	=-\mw_1^\top\muu^\top\me\mv\mv^\top\hat\mv-\mw_1^\top\muu^\top\me(\mi-\mv\mv^\top)\hat\mv,
	\end{aligned}
	$$
and hence, by Lemma~\ref{lemma:|E|_2,|UEV|F} and Eq.~\eqref{eq:wedin_sintheta}, we have 
\begin{equation}
  \label{eq:T3_term2}
	\begin{aligned}
		\|\mSigma\mv^{*\top}\hat\mv-\muu^{*\top}\hat\muu\hat\mSigma\|
		&\leq
        \|\muu^\top\me\mv\|+\|\me\|\cdot\|(\mi-\mv\mv^\top)\hat\mv\|
        \\ &
		\lesssim d^{1/2}\rho_n^{1/2}(\log n)^{1/2}+(n\rho_n)^{1/2}\cdot (n\rho_n)^{-1/2}
		\lesssim \vartheta_n 
	\end{aligned}
    \end{equation}
    with high probability. For the third term, Eq.~\eqref{eq:Weyl}
    and Eq.~\eqref{eq:wedin_sintheta2} together imply
    \begin{equation}
      \label{eq:term3_c}
	\begin{aligned}
		\|\mw_1^\top(\muu^{\top}\hat\muu-\mw^\top_\muu)\hat\mSigma\|
		&\leq \|\hat\mSigma\|\cdot\|\muu^{\top}\hat\muu-\mw^\top_\muu\|
		\lesssim n\rho_n \cdot (n \rho_{n})^{-1}
		\lesssim 1.
	\end{aligned}
    \end{equation}
    with high probability.

	 Combining Eq.~\eqref{eq:term3_a}, Eq.~\eqref{eq:T3_term2} and
     Eq.~\eqref{eq:term3_c} we have
	 $\|\bm{\Sigma}\tilde{\mt}_{3}\hat\mSigma\|\lesssim \vartheta_n$
     with high probability, and hence
	 $$
	\begin{aligned}
		\|\tilde{\mt}_{3}\|
		&\leq \|\mSigma\tilde{\mt}_{3}\hat\mSigma\|\cdot
        \|\bm{\Sigma}^{-1}\| \cdot \|\hat{\bm{\Sigma}}^{-1}\|
		\lesssim (n\rho_n)^{-2}\vartheta_n \end{aligned}
	$$
    with high probability.
    In summary we obtain
	$$
	\begin{aligned}
		&\left\|\mt_3\right\|
		\leq \|\me\|\cdot\|\tilde{\mt}_3\|
		\lesssim  (n\rho_n)^{-3/2}\vartheta_n, \\
		&\|\mt_3\|_{2\to\infty}
		\leq \left\|\me\mv\right\|_{2\to\infty}\cdot\|\tilde{\mt}_3\|
		\lesssim d^{1/2}n^{-1/2} (n\rho_n)^{-3/2}(\log n)^{1/2}\vartheta_n
	\end{aligned}
	$$
    with high probability.
\end{proof}

\subsection{Technical lemmas for Theorem~\ref{thm:What_U Rhat What_v^T-R->norm}}
\label{sec:technical_lemmas2}

\begin{lemma1}
  \label{lemma:(U^T Uhat-W_Q)W_Q^T=...}
  Consider the setting in Theorem~\ref{thm:UhatW-U=EVR^{-1}+xxx}. 
  Let $\vartheta_{n} = \max\{1,d^{1/2}\rho_n^{1/2}(\log n)^{1/2}\}$. 
  We then have
  \begin{equation*}
  \begin{aligned}
  	  \muu^\top \hat\muu\mw_{\muu}- \mathbf{I} &= -\frac{1}{2m^2}\sum_{j=1}^m\sum_{k=1}^m
  	  (\mathbf{R}^{(j)\top})^{-1}\mv^\top \mathbf{E}^{(j)\top}
  	  \mathbf{E}^{(k)}\mv(\mathbf{R}^{(k)})^{-1} +O_p((n\rho_n)^{-3/2}\vartheta_n). 
  \end{aligned}
  \end{equation*}
\end{lemma1}

\begin{proof}
  First recall the statement of Theorem~\ref{thm:UhatW-U=EVR^{-1}+xxx}, i.e.,
  $$\hat\muu \mw_\muu-\muu = \frac{1}{m} \sum_{j=1}^{m} \mathbf{E}^{(j)} \mathbf{V}(\mathbf{R}^{(j)})^{-1}+\mathbf{Q}_\muu$$
  with $\mq_\muu$ satisfying $\|\mq_\muu\|\lesssim (n\rho_n)^{-1} \vartheta_{n}$.
  Now let $\mathbf{E}^* = \muu^\top \hat\muu\hat\muu^\top\muu-\mi$. We
  then have
\begin{equation}
\label{eq:U^T Uhat Uhat^T U-I}
	\begin{aligned}
      \mathbf{E}^* =
      &-(\hat\muu \mw_\muu-\muu)^\top(\hat\muu
      \mw_\muu-\muu)+\muu^\top
      (\hat\muu \mw_\muu-\muu)(\hat\muu \mw_\muu-\muu)^\top\muu\\
	=&-(\hat\muu \mw_\muu-\muu)^\top(\hat\muu \mw_\muu-\muu)+O_{p}((n\rho_n)^{-2})\\
=&-\frac{1}{m^2}\sum_{j=1}^m\sum_{k=1}^m
  	  (\mathbf{R}^{(j)\top})^{-1}\mv^\top \mathbf{E}^{(j)\top}
  	  \mathbf{E}^{(k)}\mv(\mathbf{R}^{(k)})^{-1}
  	  +O_{p}((n\rho_n)^{-3/2}\vartheta_{n}),
\end{aligned}
\end{equation}
where the second equality in the above display follows from
Eq.~\eqref{eq:wedin_sintheta2}, i.e.,
$$
\begin{aligned}
	\|\muu^\top (\hat\muu \mw_\muu-\muu)\|
	=\|(\muu^\top \hat\muu -\mw_\muu^\top)\mw_\muu\|
	=\|\muu^\top \hat\muu -\mw_\muu^\top\|
	\lesssim (n \rho_{n})^{-1}
\end{aligned}
$$
with high probability. Eq.~\eqref{eq:U^T Uhat Uhat^T U-I} also implies
$\left\|\me^*\right\|=O_{p}((n\rho_n)^{-1})$  with high probability.

Denote the singular value decomposition of $\muu^\top \hat\muu$ by
$\muu'\mSigma'\mv'^\top$. Recall that $\mw_\muu$ is the solution of
orthogonal Procrustes problem between $\hat\muu$ and $\muu$, i.e.,
$\mw_\muu=\mv'\muu'^\top$. We thus have
$$
\begin{aligned}
	\muu^\top \hat\muu\mw_\muu
	=\muu'\mSigma'\muu'^\top
	=\bigl((\muu'\mSigma'\mv'^\top)(\mv'\mSigma'\muu'^\top)\bigr)^{1/2} 
	=\bigl(\mathbf{I}+\me^*\bigr)^{1/2}.
\end{aligned}
$$
Then by applying Theorem~2.1 in \cite{carlsson2018perturbation}, we obtain
$$
\begin{aligned}
	\muu^\top \hat\muu\mw_\muu
	&=\mathbf{I} + \frac{1}{2}\mathbf{E}^{*} + O(\|\mathbf{E}^{*}\|^2) \\
  	 &=\mi-\frac{1}{2m^2}\sum_{j=1}^m\sum_{k=1}^m
  	  (\mathbf{R}^{(j)\top})^{-1}\mv^\top \mathbf{E}^{(j)\top}
  	  \mathbf{E}^{(k)}\mv(\mathbf{R}^{(k)})^{-1}+O_{p}((n\rho_n)^{-3/2}\vartheta_{n})
\end{aligned}
$$
as desired.
\end{proof}

\begin{lemma1}
	\label{lemma:R^T (U^T hatU hatLambda^(-1)-U^T hatU)W_U)}
	Consider the setting in Theorem~\ref{thm:UhatW-U=EVR^{-1}+xxx}. Let $\vartheta_n =
    \max\{1,d^{1/2}\rho_n^{1/2}(\log n)^{1/2}\}$. We
    then have
$$
\begin{aligned}
  \mathbf{U}^{\top} \hat{\mathbf{U}}
    (\hat{\mathbf{\Lambda}}^{-1}- \mathbf{I}) \mw_\muu
	=&-\frac{1}{m} 
    \sum_{j=1}^{m} \bigl(\muu^{\top}\mathbf{E}^{(j)} \mathbf{V}(\mathbf{R}^{(j)})^{-1} 
	+ (\mathbf{R}^{(j) \top})^{-1} \mathbf{V}^{\top} \mathbf{E}^{(j)
      \top}\muu\bigr) \\
	&-
    \muu^{\top}\mathbf{L}\muu-\frac{1}{m} \muu^{\top}
    \tilde{\mathbf{E}} \sum_{k=1}^{m} \mathbf{E}^{(k)}
    \mathbf{V}(\mathbf{R}^{(k)})^{-1} +O_{p}((n\rho_n)^{-3/2}\vartheta_n),
\end{aligned}
$$	
where the matrix $\tilde\me$ and $\ml$ are defined in Eq.~\eqref{eq:def_e_l_thm3}
and $\hat{\bm{\Lambda}}$ is the matrix containing the $d$ largest
eigenvalues of $\sum_{i=1}^{m} \hat{\muu}^{(i)}
\hat{\muu}^{(i)\top}$, i.e.,
\begin{equation}
  \label{eq:UU^T}
\begin{split}
\hat\muu\hat\mLambda \hat\muu^\top+\hat\muu_\perp\hat\mLambda_\perp \hat\muu_\perp^\top
=\frac{1}{m}\sum_{i=1}^m\hat\muu^{(i)}(\hat\muu^{(i)})^\top
=\muu\muu^{\top}+\tilde\me.
\end{split}
\end{equation}
\end{lemma1}
\begin{proof}
We first bound $\tilde\me$ and $\ml$ for the setting in
Theorem~\ref{thm:UhatW-U=EVR^{-1}+xxx}. By plugging Eq.~\eqref{eq:ep_t,ze_t} into Eq.~\eqref{eq:epsilon_e_thm3}, we have
\begin{equation}
\label{eq:COSIE_tildeE_L}
	\begin{aligned}
	    &\|\ml\|
	    =\epsilon_\ml
	    \lesssim \epsilon_{\mt_0}^2+\epsilon_{\mt}
	    \lesssim [(n\rho_n)^{-1/2}]^2+(n \rho_n)^{-1}\vartheta_n
	    \lesssim (n \rho_n)^{-1}\vartheta_n,\\
		&\|\tilde{\mathbf{E}}\|
		=\epsilon_{\tilde{\mathbf{E}}}
		\lesssim \epsilon_{\mt_0}+\epsilon_{\mt}
		\lesssim (n\rho_n)^{-1/2}+(n \rho_n)^{-1}\vartheta_n
		\lesssim (n\rho_n)^{-1/2}
	\end{aligned}
\end{equation}
with high probability.

We note that
\begin{equation}
\label{eq:decomposition of R^T (U^T hatU hatLambda^(-1)-U^T hatU)W_U)}
\begin{split}
	\mathbf{U}^{\top} \hat{\mathbf{U}}
    \bigl(\hat{\mathbf{\Lambda}}^{-1}- \mathbf{I}\bigr) \mw_\muu
	&=\mathbf{U}^{\top} \hat{\mathbf{U}}\bigl(\mathbf{I} -
    \hat{\bm{\Lambda}}\big) \hat{\mathbf{\Lambda}}^{-1}\mw_\muu 
    = -\muu^{\top} \tilde{\me}
    \hat{\muu} \mw_\muu
    \mw_\muu^\top\hat{\mathbf{\Lambda}}^{-1}\mw_\muu,
    \end{split}
\end{equation}
where the last equality follows from Eq.~\eqref{eq:UU^T}.
Let $\mt_0^{(k)} = \me^{(k)} \mv (\mr^{(k)})^{-1}$. 
Using the definition of $\tilde{\me}$ and the expansion for $(\muu-\hat\muu \mw_\muu)$
in Theorem~\ref{thm:UhatW-U=EVR^{-1}+xxx}, we have
\begin{equation}
\label{eq:R^T(U^T hatU - U^T hatU hatLambda)W_U}
	\begin{aligned}
	\muu^{\top} \tilde{\mathbf{E}} \hat\muu \mw_\muu
	=& \muu^{\top} \tilde{\mathbf{E}} \muu+ \muu^{\top} \tilde{\mathbf{E}}(\hat\muu \mw_\muu-\muu)\\
	=& \muu^{\top} \tilde{\mathbf{E}} \muu
	+ \muu^{\top} \tilde{\mathbf{E}}\Big[ \frac{1}{m}\sum_{k=1}^{m} \mathbf{E}^{(k)} \mathbf{V}(\mathbf{R}^{(k)})^{-1}+\mathbf{Q}_\muu\Big]\\
	=& \muu^{\top} \tilde{\mathbf{E}} \muu
	+ \frac{1}{m}\muu^{\top} \tilde{\mathbf{E}} \sum_{k=1}^{m} \mathbf{E}^{(k)} \mathbf{V}(\mathbf{R}^{(k)})^{-1}
	+O_{p}((n\rho_n)^{-3/2} \vartheta_n)\\
	=&\frac{1}{m} \sum_{j=1}^{m} \bigl[\muu^{\top}\mathbf{E}^{(j)} \mathbf{V}(\mathbf{R}^{(j)})^{-1} 
	+ (\mathbf{R}^{(j) \top})^{-1}
    \mathbf{V}^{\top} \mathbf{E}^{(j) \top}\muu\bigr] \\
     +&
    \muu^{\top}\mathbf{L}\muu
	+ \frac{1}{m} \muu^{\top} \tilde{\mathbf{E}} \sum_{k=1}^{m} \mathbf{E}^{(k)} \mathbf{V}(\mathbf{R}^{(k)})^{-1}
	+O_p((n\rho_n)^{-3/2}\vartheta_n),
\end{aligned}
\end{equation}
where the third equality follows from Eq.~\eqref{eq:COSIE_tildeE_L} and Theorem~\ref{thm:UhatW-U=EVR^{-1}+xxx}. i.e.,
\begin{equation*}
\begin{aligned}
	\| \muu^{\top} \tilde{\mathbf{E}}\mq_\muu\|
	\leq \|\tilde{\mathbf{E}}\|\cdot\|\mq_\muu\|
	\lesssim (n\rho_n)^{-3/2} \vartheta_n
\end{aligned}
\end{equation*}
with high probability. Eq.~\eqref{eq:COSIE_tildeE_L} and
Lemma~\ref{lemma:|E|_2,|UEV|F} then imply
\begin{equation}
  \label{eq:utilde_me_uhat}
\begin{aligned}
	\|\mathbf{U}^{\top} \tilde{\me} \hat{\mathbf{U}} \mw_\muu\|
	&\lesssim \frac{1}{m} \sum_{j=1}^{m} \|\muu^{\top}\mathbf{E}^{(j)} \mathbf{V}\|\cdot \|(\mathbf{R}^{(j)})^{-1}\|+\|\ml\|
	+\frac{1}{m} \|\tilde{\mathbf{E}}\| \sum_{j=1}^{m}\|\mathbf{E}^{(k)}\|\cdot \|(\mathbf{R}^{(k)})^{-1}\|+(n\rho_n)^{-3/2}\vartheta_n\\
	&\lesssim d^{1/2} n^{-1} \rho_n^{-1/2}  (\log n)^{1/2}
	+(n \rho_n)^{-1}\vartheta_n
	+ (n \rho_n)^{-1}+(n\rho_n)^{-3/2}\vartheta_n\\
	&\lesssim (n \rho_n)^{-1}\vartheta_n
\end{aligned}
\end{equation}
with high probability.

Now for the diagonal matrix $\hat\mLambda$, we have for any $j \in
[d]$ that
$$
\begin{aligned}
	\hat\mLambda_{jj}^{-1}-1
	=\frac{1}{1-(1-\hat\mLambda_{jj})}-1
	=\sum_{k\geq 1}\big(1-\hat\mLambda_{jj}\big)^k
	=O_{p}((n\rho_n)^{-1/2})
\end{aligned}
$$
where the last equality follows from Eq.~\eqref{eq:hatLambda_separate_2} and Eq.~\eqref{eq:COSIE_tildeE_L}. We
therefore have
\begin{equation}
\label{eq:hatLambda^(-1)}
	\hat{\bm{\Lambda}}^{-1}=\mi+O_p((n\rho_n)^{-1/2}).
\end{equation}
Combining Eq.~\eqref{eq:decomposition of R^T (U^T hatU
  hatLambda^(-1)-U^T hatU)W_U)}, Eq.~\eqref{eq:utilde_me_uhat}, and Eq.~\eqref{eq:hatLambda^(-1)}, we obtain
\begin{equation*}
\begin{aligned}
	\mathbf{U}^{\top} \hat{\mathbf{U}}
    (\hat{\mathbf{\Lambda}}^{-1}- \mi) \mw_\muu
	=&-\Bigl[\muu^{\top}\tilde{\me} \hat{\muu} \mw_\muu\Big]\mw_\muu^\top\Big[\mi+O_p((n\rho_n)^{-1/2})\Bigr]\mw_\muu\\
	=&-\muu^{\top} \tilde{\me} \hat{\muu} \mw_{\muu} +O_{p}((n\rho_n)^{-3/2}\vartheta_n).
\end{aligned}
\end{equation*}
We complete the proof by substituting Eq.~\eqref{eq:R^T(U^T hatU - U^T hatU
  hatLambda)W_U} into the above display. 
\end{proof}

\begin{lemma1}
  \label{lemma:U^T Uhat lambdahat^-2 W_Q^T= I + ...}
  Consider the setting in Theorem~\ref{thm:UhatW-U=EVR^{-1}+xxx}. We
  then have
  \begin{gather*}
  \begin{aligned}
  	  \mathbf{U}^{\top} \hat{\mathbf{U}} \hat{\mathbf{\Lambda}}^{-2} \mw_\muu
  	  =\mi+O_p((n\rho_n)^{-1/2}), \quad
  	  \mathbf{U}^{\top} \hat{\mathbf{U}} \hat{\mathbf{\Lambda}}^{-3} \mw_\muu
  	  =\mi+O_p((n\rho_n)^{-1/2}).
  \end{aligned}
  \end{gather*}
\end{lemma1}

\begin{proof}
  We only derive the result for $\mathbf{U}^{\top} \hat{\mathbf{U}}
  \hat{\mathbf{\Lambda}}^{-2} \mw_\muu$ as the result for $\mathbf{U}^{\top} \hat{\mathbf{U}}
  \hat{\mathbf{\Lambda}}^{-3} \mw_\muu$ follows an almost identical argument. 
  First recall Eq.~\eqref{eq:hatLambda_separate_2}
  We then have, for any $j \in [d]$, 
  \begin{equation*}
    \hat{\bm{\Lambda}}_{jj}^{-2} - 1 = \sum_{k \geq 1} (1 -
    \hat{\bm{\Lambda}}_{jj}^2)^{k} = O_{p}((n \rho_n)^{-1/2}).
  \end{equation*}
  and hence $\|\hat{\bm{\Lambda}}^{-2} - \mathbf{I}\| =
  O_{p}((n\rho_n)^{-1/2})$. We therefore have
  \begin{equation*}
\begin{aligned}
	\mathbf{U}^{\top} \hat{\mathbf{U}} \hat{\mathbf{\Lambda}}^{-2} \mw_\muu
	&=\mathbf{U}^{\top} \hat{\mathbf{U}} \mw_\muu + O_{p}((n
    \rho_n)^{-1/2}) \\ &= \mi
	+(\mathbf{U}^{\top} \hat{\mathbf{U}}  -\mw_\muu^\top)\mw_\muu + O_{p}((n
    \rho_n)^{-1/2})
	=\mi
	+O_{p}((n\rho_n)^{-1/2}),
\end{aligned}
\end{equation*}
where the last equality follows the bounds in Eq.~\eqref{eq:wedin_sintheta2}, i.e.,
$$
\begin{aligned}
	&\|(\mathbf{U}^{\top} \hat{\mathbf{U}}  -\mw_\muu^\top)\mw_\muu\|
	\leq \|\mathbf{U}^{\top} \hat{\mathbf{U}}  -\mw_\muu^\top\|
	\lesssim (n\rho_n)^{-1}
\end{aligned}
$$
with high probability.
\end{proof}

\begin{proof}[Proof of Lemma~\ref{lemma:R V^T Q_v,R^T U^T Q_u}]
We will only prove the result for $\muu^\top \mq_\muu$ as
the proof for $\mv^\top \mq_\mv$ follows an almost identical argument.
Recall Eq.~\eqref{eq:Q1-Q5} and let
$\mq_\muu=\mq_{\muu,1}+\mq_{\muu,2}+\mq_{\muu,3}+\mq_{\muu,4}+\mq_{\muu,5}$.
We now analyze each of the terms $\muu^{\top} \mq_{\muu,1}$ through
$\muu^{\top} \mq_{\muu,5}$. 
For $\muu^\top \mq_{\muu,1} $ we have
$$
\begin{aligned}
	\muu^\top \mq_{\muu,1} 
	=&\mathbf{U}^{\top} \hat{\mathbf{U}} \hat{\mathbf{\Lambda}}^{-1}
    \mw_{\muu} - \mi = \mathbf{U}^{\top} \hat{\mathbf{U}}
    \bigl(\hat{\mathbf{\Lambda}}^{-1}- \mi\bigr) \mw_\muu +
    \bigl(\muu^{\top} \hat{\muu} \mw_{\muu} - \mi\bigr).
  \end{aligned}
$$
Therefore, by Lemma~\ref{lemma:(U^T
  Uhat-W_Q)W_Q^T=...} and Lemma~\ref{lemma:R^T (U^T
  hatU hatLambda^(-1)-U^T hatU)W_U)}, we have
\begin{equation*}
  \begin{aligned}
	\muu^\top \mq_{\muu,1} 
  	=&-\frac{1}{m} \sum_{j=1}^{m} \bigl(\mathbf{M}^{(j)}(\mathbf{R}^{(j)})^{-1} + (\mathbf{R}^{(j) \top})^{-1}
    \mathbf{M}^{(j)\top}\bigr) - \muu^{\top}\mathbf{L}\muu\\ &
	-\frac{1}{m}  \muu^{\top} \tilde{\mathbf{E}}\sum_{k=1}^{m} \mathbf{E}^{(k)} \mathbf{V}(\mathbf{R}^{(k)})^{-1}
	-\frac{1}{2m^2}\sum_{j=1}^m\sum_{k=1}^m (\mr^{(j)\top})^{-1}
    \tilde{\mathbf{N}}^{(jk)} (\mr^{(k)})^{-1}
	+O_p((n\rho_n)^{-3/2}\vartheta_n).
\end{aligned}
\end{equation*}
We next consider $\muu^\top \mq_{\muu,2} $. We have
$$
\begin{aligned}
	\muu^\top \mq_{\muu,2}
	&=\frac{1}{m} \sum_{j=1}^{m}\mathbf{M}^{(j)}(\mathbf{R}^{(j)})^{-1}\big(\mathbf{U}^{\top}
    \hat{\mathbf{U}} \hat{\boldsymbol{\Lambda}}^{-2}\mathbf{W}_\muu - \mi\big)
    =O_p(d^{1/2} n^{-1/2} (n \rho_n)^{-1}(\log n)^{1/2}),
\end{aligned}
$$
where the final equality follows from Lemma~\ref{lemma:|E|_2,|UEV|F}
and Lemma~\ref{lemma:U^T Uhat lambdahat^-2 W_Q^T= I + ...} , i.e., 
$$
\begin{aligned}
	\bigl\|\mathbf{M}^{(j)}(\mathbf{R}^{(j)})^{-1} \big(\mathbf{U}^{\top}
    \hat{\mathbf{U}}
    \hat{\boldsymbol{\Lambda}}^{-2}\mathbf{W}_\muu - \mi \big)\bigr\| &\leq 
	\|\mathbf{M}^{(j)}\| \cdot \|(\mathbf{R}^{(j)})^{-1}\| \cdot \|
\muu^{\top} \hat{\muu}
\hat{\boldsymbol{\Lambda}}^{-2}\mathbf{W}_\muu - \mi\|
\\ &\lesssim d^{1/2}\rho_n^{1/2}(\log n)^{1/2}\cdot (n\rho_n)^{-1}\cdot (n\rho_n)^{-1/2}
\\ &\lesssim d^{1/2}n^{-1/2}(n \rho_n)^{-1} (\log n)^{1/2}
\end{aligned}
$$
with high probability. 
For $\muu^\top \mq_{\muu,3}$, we once again use Lemma~\ref{lemma:|E|_2,|UEV|F} and Lemma~\ref{lemma:U^T
  Uhat lambdahat^-2 W_Q^T= I + ...} to obtain
$$
\begin{aligned}
	\muu^\top \mq_{\muu,3}
	&=\frac{1}{m} \sum_{j=1}^{m}(\mathbf{R}^{(j) \top})^{-1} \mathbf{M}^{(j)\top}\mathbf{U}^{\top} \hat{\mathbf{U}} \hat{\mathbf{\Lambda}}^{-2} \mw_{\muu}\\
	&=\frac{1}{m} \sum_{j=1}^{m}(\mathbf{R}^{(j) \top})^{-1}
    \mathbf{M}^{(j)\top}+O_p(d^{1/2}n^{-1/2}(n\rho_n)^{-1}(\log n)^{1/2}).
  \end{aligned}
$$
For $\muu^\top \mq_{\muu,4}$, we have from Lemma~\ref{lemma:U^T
  Uhat lambdahat^-2 W_Q^T= I + ...} and Eq.~\eqref{eq:COSIE_tildeE_L}
that
$$
\begin{aligned}
    \muu^\top \mq_{\muu,4}
    &=\muu^\top\mathbf{L} \mathbf{U} \mathbf{U}^{\top} \hat{\mathbf{U}} \hat{\boldsymbol{\Lambda}}^{-2}\mw_{\muu}
    =\muu^\top\mathbf{L} \mathbf{U}+O_p((n\rho_n)^{-3/2}\vartheta_n).
\end{aligned}
$$

Finally, for $\muu^\top \mq_{\muu,5} $, we have
$$
\begin{aligned}
	\muu^\top \mq_{\muu,5} 
	&=\muu^\top  \tilde{\mathbf{E}}^{2} \mathbf{U} \mathbf{U}^{\top} \hat{\mathbf{U}} \hat{\mathbf{\Lambda}}^{-3} \mw_\muu
	+\sum_{k=3}^{\infty}\muu^\top  \tilde{\mathbf{E}}^{k} \mathbf{U} \mathbf{U}^{\top} \hat{\mathbf{U}} \hat{\mathbf{\Lambda}}^{-(k+1)} \mw_\muu\\
	&=\muu^\top  \tilde{\mathbf{E}}^{2} \mathbf{U}+O_p((n\rho_n)^{-3/2}),
\end{aligned}
$$
where the last equality follows from Lemma~\ref{lemma:U^T
  Uhat lambdahat^-2 W_Q^T= I + ...} and Eq.~\eqref{eq:COSIE_tildeE_L}, e.g., 
$$
\begin{aligned}
&\|\muu^\top  \tilde{\mathbf{E}}^{2} \mathbf{U}(\mathbf{U}^{\top} \hat{\mathbf{U}} \hat{\mathbf{\Lambda}}^{-3} \mw_\muu-\mi)\|
\leq \|\tilde{\mathbf{E}}\|^2\cdot \|\mathbf{U}^{\top} \hat{\mathbf{U}} \hat{\mathbf{\Lambda}}^{-3} \mw_\muu-\mi\|\lesssim (n\rho_n)^{-3/2},
\\
	&\|\sum_{k=3}^{\infty}\muu^\top  \tilde{\mathbf{E}}^{k} \mathbf{U} \mathbf{U}^{\top} \hat{\mathbf{U}} \hat{\mathbf{\Lambda}}^{-(k+1)} \mw_\muu\|
	\leq \sum_{k=3}^{\infty}  \|\tilde{\mathbf{E}}\|^{k}
	\lesssim \sum_{k=3}^{\infty} (n\rho_n)^{-k/2} \lesssim (n\rho_n)^{-3/2}
\end{aligned}
$$
with high probability.

Combining the bounds for $\muu^{\top}\mq_{\muu,1}$ through $\muu^{\top}\mq_{\muu,5}$, 
and noting that 
$\muu^{\top} \ml \muu$ appeared in both $\muu^{\top} \mq_{\muu,1}$ and
$\muu^{\top} \mq_{\muu,4}$ but with different signs while 
$\frac{1}{m}\sum_{j=1}^m (\mathbf{R}^{(j) \top})^{-1} \mathbf{V}^{\top}
\mathbf{E}^{(j) \top}\muu$ appeared in both $\muu^{\top}
\mq_{\muu,1}$ and $\muu^{\top}\mq_{\muu,3}$ but with different signs,
we obtain
$$
\begin{aligned}
	\muu^\top \mq_\muu
	=&-\frac{1}{m} \sum_{j=1}^{m} \mathbf{M}^{(j)}(\mathbf{R}^{(j)})^{-1} 
	-\frac{1}{2m^2}\sum_{j=1}^m\sum_{k=1}^m
  	  (\mathbf{R}^{(j)\top})^{-1}\tilde{\mathbf{N}}^{(jk)}(\mathbf{R}^{(k)})^{-1}
	\\
	+&\muu^\top  \tilde{\mathbf{E}} \Bigl(\tilde{\mathbf{E}}\mathbf{U}-\frac{1}{m} \sum_{k=1}^{m} \mathbf{E}^{(k)} \mathbf{V}(\mathbf{R}^{(k)})^{-1}\Bigr)+O_p((n\rho_n)^{-3/2}\vartheta_n) \\
  	 =&-\frac{1}{m} \sum_{j=1}^{m} \mathbf{M}^{(j)}(\mathbf{R}^{(j)})^{-1} 
	-\frac{1}{2m^2}\sum_{j=1}^m\sum_{k=1}^m
  	  (\mathbf{R}^{(j)\top})^{-1}\tilde{\mathbf{N}}^{(jk)}(\mathbf{R}^{(k)})^{-1}
  	  +O_p((n\rho_n)^{-3/2}\vartheta_n)
\end{aligned}
$$
where the last equality follows from
 and Lemma~\ref{lemma:|E|_2,|UEV|F}, i.e.,
$$
\begin{aligned}
	\Big\|\muu^\top  \tilde{\mathbf{E}} \Bigl(\tilde{\mathbf{E}} \mathbf{U}- \frac{1}{m} \sum_{k=1}^{m} \mathbf{E}^{(k)} \mathbf{V}(\mathbf{R}^{(k)})^{-1}\Bigr)\Big\|
	=&\Big\|\muu^\top  \tilde{\mathbf{E}}\Bigl(\frac{1}{m} \sum_{k=1}^{m} \mathbf{U}(\mathbf{R}^{(k) \top})^{-1} \mathbf{V}^{\top} \mathbf{E}^{(k) \top}\muu+\mathbf{L}\muu\Bigr)\Big\|\\
	\lesssim & \|\tilde{\mathbf{E}}\|\Big(\|(\mathbf{R}^{(k) })^{-1}\| \cdot \|\muu^{\top} \mathbf{E}^{(k) }\mv\|_F+\|\ml\|\Big)\\
	\lesssim & (n\rho_n)^{-3/2}\vartheta_n
\end{aligned}
$$
with high probability.
\end{proof}

\begin{proof}[Proof of Lemma~\ref{lemma:U^T EE^T U R^-1->}]
\begin{sloppypar}
Recall the definition of $\mathbf{F}^{(i)}$ in the statement of Lemma~\ref{lemma:U^T EE^T U R^-1->} as
  $$
  \begin{aligned}
    \mathbf{F}^{(i)} &=\frac{1}{m}\sum_{j=1}^m \muu^\top \me^{(i)}
    \me^{(j)\top} \muu  (\mr^{(j)\top})^{-1}
    +\frac{1}{m}\sum_{j=1}^m(\mr^{(j)\top})^{-1} \mv^\top\me^{(j)\top}\me^{(i)}\mv
    \\ &-\frac{1}{2m^2} \sum_{j=1}^{m} \sum_{k=1}^{m}
    \mathbf{R}^{(i)}(\mathbf{R}^{(j)})^{-1}
\mathbf{U}^{\top} \mathbf{E}^{(j)} \mathbf{E}^{(k) \top} \mathbf{U} (\mathbf{R}^{(k) \top})^{-1}
    -\frac{1}{2m^2} \sum_{j=1}^{m} \sum_{k=1}^{m} (\mathbf{R}^{(j)
     \top})^{-1}
\mathbf{V}^{\top} \mathbf{E}^{(j)
        \top} \mathbf{E}^{(k)}
      \mathbf{V} (\mathbf{R}^{(k)})^{-1}\mathbf{R}^{(i)},
  \end{aligned}
  $$
  and recall from the statement of Theorem~\ref{thm:What_U Rhat
    What_v^T-R->norm} that
  $\tilde{\mathbf{D}}^{(i)}$ is a $n \times n$ diagonal matrix with
  $$\tilde{\mathbf{D}}^{(i)}_{kk} = \sum_{\ell=1}^n \mathbf{P}^{(i)}_{k \ell} (1 -
  \mathbf{P}^{(i}_{k \ell}). 
  $$
  We now prove that the elements of $\rho_n^{-1/2}\sum_{j=1}^m\muu^\top
\me^{(i)} \me^{(j)\top}\muu(\mr^{(j)\top})^{-1}$ converge in
probability to the elements of $\rho_n^{-1/2}\mathbf{U}^{\top}
\tilde{\mathbf{D}}^{(i)} \muu \bigl(\mr^{(i)\top}\bigr)^{-1}$. 
The convergence of the remaining terms in
$\mathbf{F}^{(i)}$ to their corresponding terms in $\bm{\mu}^{(i)}$
follows the same idea and is thus omitted.
\end{sloppypar}

Define $\zeta_{st}^{(ij)}$ for $i \in [m], j \in [m], s \in [n]$
and $t \in [n]$ as the $st$th element of $\muu^\top \me^{(i)}
\me^{(j)\top}\muu(\mr^{(j)\top})^{-1}$. 
We then have
$$
\begin{aligned}
  \zeta_{st}^{(ij)} =\sum_{k_1=1}^n\sum_{k_2=1}^n\sum_{k_3=1}^n\sum_{\ell=1}^d \muu_{k_1 s}\muu_{k_3 \ell}\bigl((\mr^{(i)})^{-1}\bigr)_{t \ell}\me^{(i)}_{k_1 k_2}\me^{(j)}_{k_3k_2}.
\end{aligned}
$$
We will compute the mean and variance for
$\zeta_{st}^{(ij)}$ when $i\neq j$ and when $i=j$ separately.
First suppose that $i\neq j$. It is then obvious that $\mathbb{E}[\zeta_{st}^{(ij)}] = 0$. 
We now consider the variance. Note that even though some of
$\big\{\me^{(i)}_{k_1 k_2}\me^{(j)}_{k_3k_2}\big\}_{k_1,k_2,k_3\in[n]}$ are dependent,
such as $\me^{(i)}_{12}\me^{(j)}_{32}$ and
$\me^{(i)}_{12}\me^{(j)}_{42}$, their covariances are always $0$, e.g., 
$$
\begin{aligned}
    \operatorname{Cov}\big(\me^{(i)}_{12}\me^{(j)}_{32},\me^{(i)}_{12}\me^{(j)}_{42}\big)	
    &=\mathbb{E}\big[(\me^{(i)}_{12}\me^{(j)}_{32}-\mathbb{E}[\me^{(i)}_{12}\me^{(j)}_{32}])(\me^{(i)}_{12}\me^{(j)}_{42}-\mathbb{E}[\me^{(i)}_{12}\me^{(j)}_{42}])\big]\\
    &=\mathbb{E}\Big[\mathbb{E}\big[(\me^{(i)}_{12}\me^{(j)}_{32}-\mathbb{E}[\me^{(i)}_{12}\me^{(j)}_{32}])(\me^{(i)}_{12}\me^{(j)}_{42}-\mathbb{E}[\me^{(i)}_{12}\me^{(j)}_{42}])\big]\big|\me^{(i)}_{12}\Big]\\
    &=\mathbb{E}\Big[\me^{(i)2}_{12}\mathbb{E}\big[(\me^{(j)}_{32}-\mathbb{E}[\me^{(j)}_{32}])(\me^{(j)}_{42}-\mathbb{E}[\me^{(j)}_{42}])\big]\big|\me^{(i)}_{12}\Big]\\
    &=\mathbb{E}\Big[\me^{(i)2}_{12}
    \mathbb{E}(\me^{(j)}_{32}-\mathbb{E}[\me^{(j)}_{32}])\mathbb{E}(\me^{(j)}_{42}-\mathbb{E}[\me^{(j)}_{42}])\big|\me^{(i)}_{12}\Big]
    =0.
\end{aligned}
$$
Thus $\mathrm{Var}[\zeta_{st}^{(ij)}]$ can be written as the
sum of variances of $\big\{\me^{(i)}_{k_1
  k_2}\me^{(j)}_{k_3k_2}\big\}_{k_1,k_2,k_3\in[n]}$. Define
$$
\begin{aligned}
	\mathrm{Var}\big[\me^{(i)}_{k_1 k_2}\me^{(j)}_{k_3k_2}\big]
	&=\mathbb{E}\big[(\me^{(i)}_{k_1 k_2})^2(\me^{(j)}_{k_3k_2})^2\big]-\mathbb{E}\big[\me^{(i)}_{k_1 k_2}\me^{(j)}_{k_3k_2}\big]^2\\
	&=\mathbb{E}\big[(\me^{(i)}_{k_1 k_2})^2\big]\mathbb{E}\big[(\me^{(j)}_{k_3k_2})^2\big]-\mathbb{E}\big[\me^{(i)}_{k_1 k_2}\big]^2\mathbb{E}\big[\me^{(j)}_{k_3k_2}\big]^2
	\\ &=\mpp^{(i)}_{k_1 k_2}(1-\mpp^{(i)}_{k_1 k_2})\mpp^{(j)}_{k_3k_2}(1-\mpp^{(j)}_{k_3k_2}).
\end{aligned}
$$
We therefore have
$$
\begin{aligned}
	\mathrm{Var}[\zeta_{st}^{(ij)}] 	&=
    \sum_{k_1=1}^n\sum_{k_2=1}^n\sum_{k_3=1}^n\sum_{\ell=1}^d \muu_{k_1
      s}^2\muu_{k_1
      \ell}^2\bigl(\bigl(\mr^{(i)}\bigr)^{-1}\bigr)_{t\ell}^2\mathrm{Var}\big[\me^{(i)}_{k_1
      k_2}\me^{(j)}_{k_3k_2}\big] \\
	&\lesssim n^3d\cdot d^2n^{-2}\cdot (n\rho_n)^{-2}\cdot\rho_n^2
	\lesssim d^3n^{-1}.
\end{aligned}
$$
Next suppose that $i=j$. We then have
$$
\begin{aligned}
	\mathbb{E}[\zeta_{st}^{(ii)}]	&=\sum_{k_1=1}^n\sum_{k_2=1}^n\sum_{\ell=1}^d \muu_{k_1 s}\muu_{k_1
      \ell}\bigl(\bigl(\mr^{(i)}\bigr)^{-1}\bigr)_{t \ell}\mathbb{E}\big[(\me^{(i)}_{k_1 k_2})^2\big]\\
	&=\sum_{k_1=1}^n\sum_{k_2=1}^n\sum_{\ell=1}^d \muu_{k_1 s}\muu_{k_1 \ell}\bigl(\bigl(\mr^{(i)}\bigr)^{-1}\bigr)_{t \ell}\mpp^{(i)}_{k_1 k_2}(1-\mpp^{(i)}_{k_1 k_2}).
\end{aligned}
$$
Now for $\mathrm{Var}[\zeta_{st}^{(ii)}]$, similarly to the case
$i\neq j$, the covariances of the $\big\{\me^{(i)}_{k_1
  k_2}\me^{(i)}_{k_3k_2}\big\}_{k_1,k_2,k_3\in[n]}$ are all equal to
$0$. Define
$$
\begin{aligned}
	&\mathrm{Var}\big[(\me^{(i)}_{k_1 k_2})^2\big]
	=\mathbb{E}\big[(\me^{(i)}_{k_1 k_2})^4\big]-\mathbb{E}\big[(\me^{(i)}_{k_1 k_2})^2\big]^2
	=\mpp^{(i)}_{k_1 k_2}(1-\mpp^{(i)}_{k_1 k_2})(1-2\mpp^{(i)}_{k_1 k_2})^2,\\
	&\mathrm{Var}\big[\me^{(i)}_{k_1 k_2}\me^{(i)}_{k_3 k_2}\big]
    =\mpp^{(i)}_{k_1 k_2}(1-\mpp^{(i)}_{k_1
      k_2})\mpp^{(i)}_{k_3k_2}(1-\mpp^{(i)}_{k_3k_2}) \quad \text{if
      $k_3 \not = k_1$.}
\end{aligned}
$$
We therefore have
\begin{equation*}
\begin{aligned}
	\mathrm{Var}[\zeta_{st}^{(ii)}]
    &=\sum_{k_1=1}^n\sum_{k_2=1}^n\sum_{\ell=1}^d \muu_{k_1
      s}^2\muu_{k_1 \ell}^2(\mr^{(i)-1})_{t \ell}^2\mathrm{Var}\big[(\me^{(i)}_{k_1 k_2})^2\big]
	\\&
	+\sum_{k_1=1}^n\sum_{k_2=1}^n\sum_{k_3\neq k_1}\sum_{\ell=1}^d
    \muu_{k_1 s}^2\muu_{k_1 \ell}^2\bigl(\bigl(\mr^{(i)}\bigr)^{-1}\bigr)_{t
      \ell}^2\mathrm{Var}\big[\me^{(i)}_{k_1 k_2}\me^{(i)}_{k_3
      k_2}\big] \\
 &\lesssim  n^2d\cdot d^2n^{-2}\cdot (n\rho_n)^{-2}\cdot\rho_n\cdot 1^2
	+d^3n^{-1}\lesssim d^3n^{-1}.
\end{aligned}
\end{equation*}
Therefore, by Chebyshev inequality, we have
\begin{equation*}
\begin{split}
    \rho_n^{-1/2}\Bigl(\sum_{j=1}^m \zeta_{st}^{(ij)}\Bigr) - \rho_n^{-1/2}\mathbb{E}[\zeta_{st}^{(ii)}]
  	  \stackrel{p}{\longrightarrow} 0.
   \end{split}
\end{equation*}
We conclude the proof by noting that $\mathbb{E}[\zeta_{st}^{(ii)}]$
can also be written as
$$ 	  \sum_{k_1=1}^n\sum_{k_2=1}^n\sum_{\ell=1}^d \muu_{k_1
        s}\muu_{k_1 \ell}\bigl(\bigl(\mr^{(i)}\bigr)^{-1}\bigr)_{t \ell}\mpp^{(i)}_{k_1
        k_2}(1-\mpp^{(i)}_{k_1 k_2}) = \bm{u}_s^{\top}
      \tilde{\mathbf{D}}^{(i)} \bm{z}_{t},$$
      where $\bm{u}_s$ is the $s$th column of $\muu$ and $\bm{z}_t$ is
      the $t$th column of
      $\muu\bigl(\mr^{(i){\top}}\bigr)^{-1}$. Collecting all the
      terms $\mathbb{E}[\zeta_{st}^{(ii)}]$ into a matrix yields the
      desired claim. 
\end{proof}

\begin{proof}[Proof of Lemma~\ref{lemma:U^T E V->norm}]
We observe that $\operatorname{vec}\big(\mathbf{U}^{\top}
\mathbf{E}^{(i)} \mathbf{V}\big)$ is a sum of independent random
vectors. More specifically, let $\mz =
(\mv\otimes\muu)^{\top}\in\mathbb{R}^{d^2 \times n^2}$ and let $z_k$
denote the $k$th column of $\mz$. Next let $\my^{(i)}_{k_1 k_2} \in
\mathbb{R}^{d^2}$ be the random vector
$$
\begin{aligned}
	\my_{k_1,k_2}^{(i)}=\mathbf{E}_{k_1 k_2}^{(i)}z_{k_1+(k_2-1)n}.
\end{aligned}
$$
For a fixed $i$ and varying $k_1 \in [n]$ and $k_2 \in [n]$, the
collection $\{\my_{k_1,k_2}^{(i)}\}$ are mutually independent mean $0$ random vectors.
We then have
\begin{equation*}
\begin{aligned}
	\operatorname{vec}\big(\mathbf{U}^{\top} \mathbf{E}^{(i)} \mathbf{V}\big)
	=\big(\mv\otimes\muu\big)^\top\operatorname{vec}\big(\mathbf{E}^{(i)}\big)
		=\sum_{k_1=1}^{n} \sum_{k_2=1}^{n} \mathbf{E}_{k_1 k_2}^{(i)}z_{k_1+(k_2-1)n}
	=\sum_{k_1=1}^{n} \sum_{k_2=1}^{n}\my_{k_1,k_2}^{(i)}.
\end{aligned}
\end{equation*}
Next we observe that, for any $k_1,k_2\in[n]$, 
$$
\begin{aligned}
	\operatorname{Var}\big[\my_{k_1,k_2}^{(i)}\big]
	=\mpp^{(i)}_{k_1k_2}(1-\mpp^{(i)}_{k_1k_2})z_{k_1+(k_2-1)n}z_{k_1+(k_2-1)n}^\top.
\end{aligned}
$$
Then we have
\begin{equation*}\begin{split}
	\sum_{k_1=1}^{n} \sum_{k_2=1}^{n}\operatorname{Var}\big[\my_{k_1,k_2}^{(i)}\big]
	&=\sum_{k_1=1}^{n} \sum_{k_2=1}^{n}\mpp^{(i)}_{k_1k_2}(1-\mpp^{(i)}_{k_1k_2})z_{k_1+(k_2-1)n}z_{k_1+(k_2-1)n}^\top \\&= (\mv \otimes \muu)^{\top} \mathbf{D}^{(i)} (\mv \otimes \muu)=\mSigma^{(i)},
\end{split}
\end{equation*}
where $\mSigma^{(i)}$ is defined in the statement of Theorem~\ref{thm:What_U Rhat What_v^T-R->norm}.

Let $\tilde{\my}^{(i)}_{k_1,k_2} = (\bm{\Sigma}^{(i)})^{-1/2} \my_{k_1,k_2}^{(i)}$.
For any $i\in[m]$, 
we assume $\sigma_{\min}(\mSigma^{(i)})\gtrsim \rho_n$, thus $\|(\mSigma^{(i)})^{-1/2}\|\lesssim \rho_n^{-1/2}$.
For any $k_1,k_2\in[n]$, by the definition of $z_{k_1+n(k_2-1)}$ and our assumption of $\muu$ and $\mv$, we have $\|z_{k_1+n(k_2-1)}\|\lesssim d^2n^{-1}$. 
Then for any
$k_1, k_2 \in [n]$, we can bound the spectral norm of $\tilde\my_{k_1,k_2}^{(i)}$ by
\begin{equation}
\label{eq:||Y_k1k2||}
	\begin{aligned}
	\|\tilde \my_{k_1,k_2}^{(i)}\|
	\leq \|(\mSigma^{(i)})^{-1/2}\|\cdot |\mathbf{E}_{k_1 k_2}^{(i)}|\cdot \|z_{k_1+n(k_2-1)}\|
	\lesssim \rho_n^{-1/2}\cdot 1\cdot d^2n^{-1}\lesssim d^2n^{-1/2}(n\rho_n)^{-1/2}.\end{aligned}
\end{equation}
For any fixed but arbitrary $\epsilon>0$, Eq.~\eqref{eq:||Y_k1k2||} implies that, for sufficiently large $n$, we have
\begin{equation*}
\label{eq:Y_k1k2<epsilon}
	\begin{aligned}
	\max_{k_1,k_2} \|\tilde\my_{k_1,k_2}^{(i)}\|\leq \epsilon.
    \end{aligned}
\end{equation*}
We therefore have
$$
\begin{aligned}
	\sum_{k_1=1}^{n} \sum_{k_2=1}^{n}\mathbb{E}\Big[\|\tilde\my_{k_1,k_2}^{(i)}\|^2\cdot\mathbb{I}\big\{\|\tilde\my_{k_1,k_2}^{(i)}\|>\epsilon\big\}\Big]
	\longrightarrow 0.
\end{aligned}
$$
as $n\rightarrow \infty$.
Applying the
Lindeberg-Feller central limit theorem, see e.g. Proposition 2.27 in
\cite{van2000asymptotic}, we finally have
$$
  (\boldsymbol{\Sigma}^{(i)})^{-1/2}\operatorname{vec}\big(\mathbf{U}^{\top} \mathbf{E}^{(i)} \mathbf{V}\big) 
  \leadsto \mathcal{N}\big(\bm{0},\mi\big)
$$
as $n\rightarrow \infty$.
\end{proof}

\subsection{Technical lemmas for Theorem~\ref{thm:HT}}
\label{sec:technical_HT}
\begin{lemma1}
\label{lemma:sigmahat-sigma}
	Consider the setting of Theorem~\ref{thm:What_U Rhat
      What_v^T-R->norm}. Then for any $i \in [m]$ we have
	$$\begin{aligned}
		\big\|(\mw_\mv\otimes\mw_\muu)\mathbf{\Sigma}^{(i)}(\mw_\mv\otimes\mw_\muu)^\top
		-\hat\mSigma^{(i)}\big\| 
		\lesssim dn^{-1}(n\rho_n)^{1/2}(\log n)^{1/2}
	\end{aligned}
	$$
	with high probability.
\end{lemma1}
\begin{proof}
  We first recall Theorem~\ref{thm:UhatW-U=EVR^{-1}+xxx} and Eq.~\eqref{eq:W_U^T hatR W_V-R=...}. In particular we have
	\begin{equation}
  \label{eq:technical_uhat-uw_2inf}
	\begin{aligned}
		&
      \|\hat\muu\mw_\muu-\muu\|_{2\to\infty}\lesssim d^{1/2}n^{-1/2}(n\rho_n)^{-1/2}(\log n)^{1/2}, \\
      &\|\hat\mv\mw_\mv-\mv\|_{2\to\infty}\lesssim d^{1/2}n^{-1/2}(n\rho_n)^{-1/2}(\log n)^{1/2},\\
      &\|\mw_\muu^\top\hat\mr^{(i)}\mw_\mv-\mr^{(i)}\|
      \lesssim \vartheta_n
	\end{aligned}
	\end{equation}
	with high probability, where $\vartheta_n=\max\{1,d\rho_n^{1/2}(\log n)^{1/2}\}$.
	Then under the assumption $n\rho_n=\Omega(\log n)$, we have the bound of $\|\hat{\muu}\|_{2 \to \infty}$, $\|\hat{\mv}\|_{2 \to \infty}$ and $\|\hat\mr^{(i)}\|$ as
	\begin{equation} \label{eq:technical_uhat_2inf}
		\begin{aligned}
		\|\hat{\muu}\|_{2 \to \infty}
		&\leq \|\muu\|_{2\to\infty}
		+\|\hat\muu\mw_\muu-\muu\|_{2\to\infty}
		\lesssim d^{1/2}n^{-1/2},\\
		\|\hat{\mv}\|_{2 \to \infty}
		&\leq \|\mv\|_{2\to\infty}
		+\|\hat\mv\mw_\muu-\mv\|_{2\to\infty}
		\lesssim d^{1/2}n^{-1/2},\\
		\|\hat\mr^{(i)}\|
		&\leq \|\mr^{(i)}\|+\|\mw_\muu^\top\hat\mr^{(i)}\mw_\mv-\mr^{(i)}\|
		\lesssim n\rho_n
	\end{aligned}
	\end{equation}
	with high probability. 
	Next recall that $\mpp^{(i)} =\muu\mr^{(i)}\mv^\top$ and
	$\hat\mpp^{(i)} = \hat{\muu} \hat{\mr}^{(i)} \hat{\mv}^{\top}$.
  We thus have
  \begin{equation*}
  \begin{split}
  \bigl\|\hat{\mpp}^{(i)} - \mpp^{(i)}\bigr\|_{\max} 
  &
  \leq
  \bigl\|(\muu - \hat{\muu} \mw_{\muu})\mr^{(i)}\mv^{\top}\|_{\max} 
  + \bigl\|\hat{\muu} \mw_{\muu}(\mr^{(i)} - \mw_{\muu}^{\top}\hat{\mr}^{(i)}\mw_{\mv})\mv^{\top}\bigr\|_{\max} 
  \\&
  + 
  \bigl\|\hat{\muu}\hat{\mr}^{(i)}(\mw_{\mv} \mv^{\top} - \hat{\mv}^{\top})\bigr\|_{\max} .
  \end{split}
  \end{equation*}
  Now for any two matrices $\mathbf{A}$ and $\mathbf{B}$ whose product $\mathbf{A}\mathbf{B}^{\top}$ is well defined, we have
  $$\|\ma \mb^{\top}\|_{\max} \leq \|\ma\|_{2 \to \infty} \cdot \|\mb\|_{2 \to \infty}.$$
  Thus, by Eq.~\eqref{eq:technical_uhat-uw_2inf} and Eq.~\eqref{eq:technical_uhat_2inf}, we have
  \begin{equation*}
  \begin{split}
  \bigl\|(\muu - \hat{\muu} \mw_{\muu})\mr^{(i)}\mv^{\top}\|_{\max}
   &\leq \|\hat{\muu} \mw_{\muu}-\muu\|_{2 \to \infty} \cdot \|\mv\mr^{(i)\top}\|_{2 \to \infty} \\ 
  &\leq  \|\hat{\muu} \mw_{\muu}-\muu \|_{2 \to \infty} \cdot \|\mv\|_{2 \to \infty} \cdot \|\mr^{(i)}\| 
  \lesssim dn^{-1}(n \rho_n)^{1/2}(\log n)^{1/2},\\
  \bigl\|\hat{\muu} \mw_{\muu}(\mr^{(i)} - \mw_{\muu}^{\top}\hat{\mr}^{(i)}\mw_{\mv})\mv^{\top}\bigr\|_{\max}
  &\leq \|\hat{\muu} \mw_{\muu}\|_{2 \to \infty}\cdot \|\mv(\mr^{(i)} - \mw_{\muu}^{\top}\hat{\mr}^{(i)}\mw_{\mv})^\top\|_{2\to\infty}\\
  &\leq \|\hat{\muu}\|_{2 \to \infty}\cdot\|\mv\|_{2 \to \infty}\cdot \|\mr^{(i)} - \mw_{\muu}^{\top}\hat{\mr}^{(i)}\mw_{\mv}\|
  \lesssim d n^{-1}\vartheta_n,\\
  \bigl\|\hat{\muu}\hat{\mr}^{(i)}(\mw_{\mv} \mv^{\top} - \hat{\mv}^{\top})\bigr\|_{\max}
  & \leq \|\hat{\muu}\|_{2\to\infty}\cdot \|(\mv\mw_\mv^\top-\hat\mv)\hat{\mr}^{(i)\top}\|_{2\to\infty}\\
  &\leq \|\hat{\muu}\|_{2\to\infty}\cdot \|\hat\mv\mw_\mv-\mv\|_{2\to\infty}\|\hat{\mr}^{(i)}\|
  \lesssim dn^{-1}(n \rho_n)^{1/2}(\log n)^{1/2}
  \end{split}
  \end{equation*}
	with high probability. We thus have 
	$$
	\bigl\|\hat{\mpp}^{(i)} - \mpp^{(i)}\bigr\|_{\max} 
	\lesssim dn^{-1}(n \rho_n)^{1/2}(\log n)^{1/2}
	$$
	with high probability. Hence
	\begin{equation}\label{eq:technical_Dhat-D_max}
		\begin{aligned}
			\bigl\|\hat{\md}^{(i)} - \md^{(i)}\bigr\|=\bigl\|\hat{\md}^{(i)} - \md^{(i)}\bigr\| _{\max}
	\lesssim dn^{-1}(n \rho_n)^{1/2}(\log n)^{1/2}
		\end{aligned}
	\end{equation}
	with high probability. The diagonal matrices $\hat{\md}^{(i)}$ and
    $\md^{(i)}$ are defined in Eq.~\eqref{eq:HT:hatsigma} and Theorem~\ref{thm:What_U Rhat What_v^T-R->norm}, respectively.
	
	Now recall the definitions of $\hat\mSigma^{(i)}$ and $\mSigma^{(i)}$. We then have
	$$
	\begin{aligned}
		\big\|(\mw_\mv\otimes\mw_\muu)\mathbf{\Sigma}^{(i)}(\mw_\mv\otimes\mw_\muu)^\top
		-\hat\mSigma^{(i)}\big\|
		&\leq  \|(\mv\mw_\mv^\top\otimes\muu\mw_\muu^\top-\hat\mv\otimes\hat\muu)^\top
		\md^{(i)}
		(\mv\mw_\mv^\top\otimes\muu\mw_\muu^\top)\|\\
		&+\|(\hat\mv\otimes\hat\muu)^\top(\md^{(i)}-\hat\md^{(i)})(\mv\mw_\mv^\top\otimes\muu\mw_\muu^\top)\|\\
		&+\|(\hat\mv\otimes\hat\muu)^\top\hat\md^{(i)}(\mv\mw_\mv^\top\otimes\muu\mw_\muu^\top-\hat\mv\otimes\hat\muu)\|.
	\end{aligned}
	$$
    From Eq.~\eqref{eq:|(I-UcUc^)hat Uc|} we have
    \begin{equation*}
      \|\muu\mw_\muu^\top-\hat\muu\| \lesssim (n \rho_n)^{-1/2}, \quad
      \|\mv\mw_\mv^\top-\hat\mv\| \lesssim (n \rho_n)^{-1/2}
    \end{equation*}
    and hence
    \begin{equation*}
	\begin{aligned}
	\|\mv\mw_\mv^\top\otimes\muu\mw_\muu^\top-\hat\mv\otimes\hat\muu\|
	&\leq\|(\mv\mw_\mv^\top-\hat\mv)\otimes \muu\mw_\muu^\top\|
	+\|\hat\mv\otimes(\muu\mw_\muu^\top-\hat\muu)\|\\
	&\leq\|\mv\mw_\mv^\top-\hat\mv\|
	+\|\muu\mw_\muu^\top-\hat\muu\|
	\lesssim (n\rho_n)^{-1/2}
	\end{aligned}
    \end{equation*}
	with high probability. Next, as we assume 
    $\mathbf{P}^{(i)}_{k_1k_2} \lesssim \rho_n$ for all $k_1 \in [n]$
    and $k_2 \in [n]$, we have $\|\mathbf{D}^{(i)}\| \lesssim
    \rho_n$ and hence, by Eq.~\eqref{eq:technical_Dhat-D_max},
    $\|\hat{\mathbf{D}}^{(i)}\| \lesssim \rho_n$ with high probability.
    We therefore have
	$$
	\begin{aligned}
		\|(\mv\mw_\mv^\top\otimes\muu\mw_\muu^\top-\hat\mv\otimes\hat\muu)^\top
		\md^{(i)}
		(\mv\mw_\mv^\top\otimes\muu\mw_\muu^\top)\|
		&\leq \|\mv\mw_\mv^\top\otimes\muu\mw_\muu^\top-\hat\mv\otimes\hat\muu\|\cdot \|\md^{(i)}\|
		\lesssim n^{-1}(n \rho_n)^{1/2},\\
		\|(\hat\mv\otimes\hat\muu)^\top(\md^{(i)}-\hat\md^{(i)})(\mv\mw_\mv^\top\otimes\muu\mw_\muu^\top)\|
		&\leq \|\md^{(i)}-\hat\md^{(i)}\|
		\lesssim dn^{-1}(n \rho_n)^{1/2}(\log n)^{1/2},\\
		\|(\hat\mv\otimes\hat\muu)^\top\hat\md^{(i)}(\mv\mw_\mv^\top\otimes\muu\mw_\muu^\top-\hat\mv\otimes\hat\muu)\|
		&\leq  \|\hat\md^{(i)}\|
		\cdot \|\mv\mw_\mv^\top\otimes\muu\mw_\muu^\top-\hat\mv\otimes\hat\muu\|
		\lesssim n^{-1}(n \rho_n)^{1/2}
	\end{aligned}
	$$
	with high probability. 
    In summary we obtain
	$$
	\big\|(\mw_\mv\otimes\mw_\muu)\mathbf{\Sigma}^{(i)}(\mw_\mv\otimes\mw_\muu)^\top
		-\hat\mSigma^{(i)}\big\|
		\lesssim dn^{-1}(n \rho_n)^{1/2}(\log n)^{1/2}
	$$
	with high probability.
	\end{proof}

\begin{proof}[Proof of Lemma~\ref{lemma:|sigma inverse-sigma|,order of sigma}]
    Now recall Lemma~\ref{lemma:sigmahat-sigma}, i.e.,
    \begin{equation}
    \label{eq:HT:A-B}
    \begin{aligned}
    	\big\|(\mw_\mv\otimes\mw_\muu)(\mSigma^{(i)}+\mSigma^{(j)})(\mw_\mv\otimes\mw_\muu)^\top
		-\big(\hat\mSigma^{(i)}+\hat\mSigma^{(j)}\big)\big\|
		\lesssim & d n^{-1/2} \rho_n^{1/2}(\log n)^{1/2}
    \end{aligned}
    \end{equation}
    with high probability.
    Applying Weyl's inequality, with the assumption $\sigma_{\min}(\mSigma^{(i)}+\mSigma^{(j)})\asymp \rho_n$ we have that 
    \begin{equation}
    \label{eq:HT:min singular2}
    \begin{aligned}
		\sigma_{\min}\big(\hat\mSigma^{(i)}+\hat\mSigma^{(j)}\big)
		\asymp \rho_n
    \end{aligned}
    \end{equation}
    with high probability. 
    From the assumption, Eq.~\eqref{eq:HT:A-B} and Eq.~\eqref{eq:HT:min
      singular2} we obtain
    \begin{equation}
    \label{eq:HT:A^-1, B^-1}
    \begin{aligned}
    	&\|(\mw_\mv\otimes\mw_\muu)(\mSigma^{(i)}+\mSigma^{(j)})^{-1}(\mw_\mv\otimes\mw_\muu)^\top\|
    	\asymp
        \rho_n^{-1}, \\ 
		&\|\big(\hat\mSigma^{(i)}+\hat\mSigma^{(j)}\big)^{-1}\|
		\asymp 
		\rho_n^{-1}
    \end{aligned}
    \end{equation}
    with high probability. 
    Now since $
    \|\ma^{-1}-\mb^{-1}\| \leq \|\ma^{-1}\| \cdot \|\ma-\mb\| \cdot \|\mb^{-1}\|$ for any invertible matrices $\ma$ and
    $\mb$, we have by Eq.~\eqref{eq:HT:A-B} and Eq.~\eqref{eq:HT:A^-1,
      B^-1} that
    $$
    \begin{aligned}
    	\rho_n \big\|(\mw_\mv\otimes\mw_\muu)(\mSigma^{(i)}+\mSigma^{(j)})^{-1}(\mw_\mv\otimes\mw_\muu)^\top
		-\big(\hat\mSigma^{(i)}+\hat\mSigma^{(j)}\big)^{-1}\big\| 
		\lesssim d (n\rho_n)^{-1/2}(\log n)^{1/2}
    \end{aligned}
    $$
    with high probability.
\end{proof}

\begin{lemma1}
	\label{lemma:mui-muj}
	Consider the setting of Theorem~\ref{thm:What_U Rhat
      What_v^T-R->norm}. Recall the expression for $\bm\mu^{(i)}$ given in Theorem~\ref{thm:What_U Rhat What_v^T-R->norm}. 
Then 
we have
$$
\begin{aligned}
	\|\bm\mu^{(i)}-\bm\mu^{(j)}\|
	\lesssim d^{1/2}m^{-1}\big(n\rho_n\|(\mr^{(i)})^{-1}-(\mr^{(j)})^{-1}\|+d(n\rho_n)^{-1}\|\mr^{(i)}-\mr^{(j)}\|\big).
\end{aligned}
$$
\end{lemma1}

\begin{proof}
We first bound $\|\tilde\md^{(i)}\|$ and $\|\tilde\md^{(i)}-\tilde\md^{(j)}\|$. For $\tilde\md^{(i)}$ we have
\begin{equation}
\label{eq:mui-muj_||D||}
	\begin{aligned}
		\|\tilde\md^{(i)}\|=\max_{s\in[n]}|\tilde\md^{(k)}_{ss}|
		=\max_{s\in[n]} \sum_{t=1}^n\mpp^{(i)}_{st}(1-\mpp^{(i)}_{st})
		\lesssim n\cdot\rho_n\cdot 1\lesssim n\rho_n.
	\end{aligned}
\end{equation}
For $\tilde\md^{(i)}-\tilde\md^{(j)}$, we have
$$
\begin{aligned}
	\tilde\md^{(i)}_{ss}-\tilde\md^{(j)}_{ss}
	&=\sum_{t=1}^n\mpp^{(i)}_{st}(1-\mpp^{(i)}_{st})-\mpp^{(j)}_{st}(1-\mpp^{(j)}_{st})
	=\sum_{t=1}^n(\mpp^{(i)}_{st}-\mpp^{(j)}_{st})(1-\mpp^{(i)}_{st})
	+\mpp^{(j)}_{st}(\mpp^{(j)}_{st}-\mpp^{(i)}_{st}).
\end{aligned}
$$
Now $\mpp^{(j)}_{st}\in[0,1]$ for all $\{s,t\}$ and hence
\begin{equation}
\label{eq:mui-muj_||Di-Dj||}
\begin{aligned}
	\big\|\tilde\md^{(i)}-\tilde\md^{(j)}\big\|
	&\leq n\|\mpp^{(i)}-\mpp^{(j)}\|_{\max}\\
	&\leq 2n\|\muu(\mr^{(i)}-\mr^{(j)})\mv^\top\|_{\max}\\
	&\leq 2n\|\muu\|_{2\to\infty}\cdot \|\mv\|_{2\to\infty}\cdot\|\mr^{(i)}-\mr^{(j)}\|
	\lesssim d\|\mr^{(i)}-\mr^{(j)}\|.
\end{aligned}
\end{equation}
Recall the expression for $\boldsymbol\mu^{(i)}$ in
Theorem~\ref{thm:What_U Rhat What_v^T-R->norm}. We now bound the terms
appearing in $\boldsymbol\mu^{(i)}-\boldsymbol\mu^{(j)}$. 
For $\frac{1}{m}\muu^\top\tilde\md^{(i)}\muu(\mr^{(i)\top})^{-1}
	-\frac{1}{m}\muu^\top\tilde\md^{(j)}\muu(\mr^{(j)\top})^{-1}$, by applying Eq.~\eqref{eq:mui-muj_||D||} and Eq.~\eqref{eq:mui-muj_||Di-Dj||}, we have
\begin{equation}
\label{eq:mui-muj_part1}
\begin{aligned}
	&\Big\|\frac{1}{m}\muu^\top\tilde\md^{(i)}\muu(\mr^{(i)\top})^{-1}
	-\frac{1}{m}\muu^\top\tilde\md^{(j)}\muu(\mr^{(j)\top})^{-1}\Big\|\\
	\leq &\Big\|\frac{1}{m}\muu^\top\tilde\md^{(i)}\muu\big((\mr^{(i)})^{-1}-(\mr^{(j)})^{-1}\big)^{\top}\Big\|
	+\Big\|\frac{1}{m}\muu^\top(\tilde\md^{(i)}-\tilde\md^{(j)})\muu(\mr^{(j)\top})^{-1}\Big\|\\
	\leq & m^{-1}\|\tilde\md^{(i)}\|\cdot \|(\mr^{(i)})^{-1}-(\mr^{(j)})^{-1}\|
	+m^{-1}\|\tilde\md^{(i)}-\tilde\md^{(j)}\|\cdot \|(\mr^{(j)})^{-1}\|\\
	\lesssim & m^{-1}\cdot n\rho_n\cdot \|(\mr^{(i)})^{-1}-(\mr^{(j)})^{-1}\|
	+m^{-1}\cdot d\|\mr^{(i)}-\mr^{(j)}\|\cdot (n\rho_n)^{-1}\\
	\lesssim & m^{-1}n\rho_n\|(\mr^{(i)})^{-1}-(\mr^{(j)})^{-1}\|
	+dm^{-1}(n\rho_n)^{-1}\|\mr^{(i)}-\mr^{(j)}\|.
\end{aligned}
\end{equation}
Similarly, we have
\begin{equation}
\label{eq:mui-muj_part2}
\begin{aligned}
	&
	\Big\|\frac{1}{m}(\mr^{(i)\top})^{-1}\mv^\top\breve\md^{(i)}\mv
	-\frac{1}{m}(\mr^{(j)\top})^{-1}\mv^\top\breve\md^{(j)}\mv\Big\|\\
	\lesssim &
	 m^{-1}n\rho_n\|(\mr^{(i)})^{-1}-(\mr^{(j)})^{-1}\|
	+dm^{-1}(n\rho_n)^{-1}\|\mr^{(i)}-\mr^{(j)}\|.
\end{aligned}
\end{equation}
For $\frac{1}{2m^2}\sum_{k=1}^m(\mr^{(i)}-\mr^{(j)})(\mr^{(k)})^{-1}\muu^\top\tilde\md^{(k)}\muu(\mr^{(k)\top})^{-1}$, we have
\begin{equation}
\label{eq:mui-muj_part3}
\begin{aligned}
	\Big\|\frac{1}{2m^2}\sum_{k=1}^m(\mr^{(i)}-\mr^{(j)})(\mr^{(k)})^{-1}\muu^\top\tilde\md^{(k)}\muu(\mr^{(k)\top})^{-1}\Big\|
	\lesssim &
	 m^{-2}\sum_{k=1}^m\|\tilde\md^{(k)}\|\cdot\|(\mr^{(k)})^{-1}\|^2\cdot \|\mr^{(i)}-\mr^{(j)}\|\\
	\lesssim & m^{-2}\cdot m\cdot (n\rho_n)\cdot (n\rho_n)^{-2}\|\mr^{(i)}-\mr^{(j)}\|\\
	\lesssim & m^{-1}(n\rho_n)^{-1}\|\mr^{(i)}-\mr^{(j)}\|.
\end{aligned}
\end{equation}
Similarly, we have
\begin{equation}
\label{eq:mui-muj_part4}
\begin{aligned}
	\Big\|\frac{1}{2m^2}\sum_{k=1}^m(\mr^{(k)\top})^{-1}\mv^\top\breve\md^{(k)}\mv(\mr^{(k)})^{-1}(\mr^{(i)}-\mr^{(j)})\Big\|
	\lesssim m^{-1}(n\rho_n)^{-1}\|\mr^{(i)}-\mr^{(j)}\|.
\end{aligned}
\end{equation}
Combining Eq.~\eqref{eq:mui-muj_part1} through Eq.~\eqref{eq:mui-muj_part4}, we obtain
$$
\begin{aligned}
	\|\boldsymbol\mu^{(i)}-\boldsymbol\mu^{(j)}\|
	\lesssim &d^{1/2}\big[m^{-1}n\rho_n\|(\mr^{(i)})^{-1}-(\mr^{(j)})^{-1}\|
	+dm^{-1}(n\rho_n)^{-1}\|\mr^{(i)}-\mr^{(j)}\|
	+m^{-1}(n\rho_n)^{-1}\|\mr^{(i)}-\mr^{(j)}\|\big]\\
	\lesssim & d^{1/2}m^{-1}(n\rho_n\|(\mr^{(i)})^{-1}-(\mr^{(j)})^{-1}\|
	+d(n\rho_n)^{-1}\|\mr^{(i)}-\mr^{(j)}\|).
\end{aligned}
$$
as claimed.
\end{proof}

\subsection{Proof of technical lemmas for Theorem~\ref{thm:PCA_UhatW-U=(I-UU^T)EVR^{-1}+xxx_previous}}
\label{Appendix:proof of lemma:PCA_|E|_previous and lemma:PCA_UhatW-U_previous}


\begin{proof}[Proof of Lemma~\ref{lemma:PCA_|E|_previous}]
Under the assumption $\varphi=o(1)$ and $\lambda^{(i)}_1 \asymp \lambda^{(i)}_{d_i} \asymp D^\gamma$, we have
\begin{equation*}\label{Eq:PCA_||E||,||EU||_previous}
\|\mathbf{E}^{(i)}\|
\lesssim D^\gamma \varphi,
\quad \|\mathbf{E}^{(i)}\mathbf{U}^{(i)}\|_{2\to\infty}
\lesssim \nu(X) d_i^{1/2} D^{\gamma/2} \tilde{\varphi},
\end{equation*}
with probability at least $1-\frac{1}{3}D^{-2}$, where $\nu(X)=\max_{\ell\in[D]}\operatorname{Var}(X_{(\ell)})$ and $X_{(\ell)}$ represents the $\ell$th variate in $X$; see Eq.~(1.3) in \cite{Lounici2017} for the bound for $\|\mathbf{E}^{(i)}\|$ and see the proof of Theorem~1.1 in \cite{cape2019two} for the bounds for $\|\mathbf{E}^{(i)} \mathbf{U}^{(i)}\|_{2 \to \infty}$. We note that the bound as presented in \cite{cape2019two} is somewhat sub-optimal as it uses the factor $\varphi$ as opposed to $\tilde{\varphi}$; using the same argument but with more careful book-keeping yields the bound presented here. 
Now $\sigma^2$ is bounded and $\mathbf{U}$ has bounded coherence and hence $\nu(X)$ is also bounded. The bound for $\|\mathbf{E}^{(i)}\mathbf{U}^{(i)}\|_{2\to\infty}$ in Lemma~\ref{lemma:PCA_|E|_previous} are thereby established.

Next, by Eq.~(12) in \cite{fan2018eigenvector}, we have
\begin{equation*}\label{Eq:PCA_||E||_infty_previous}
\|\mathbf{E}^{(i)}\|_\infty
\lesssim (\sigma_i^2D+D^\gamma)\tilde{\varphi}
\lesssim D \tilde{\varphi}
\end{equation*}
with probability at least $1-D^{-1}$. We note that the notations in \cite{fan2018eigenvector} are somewhat different from the notations used in the current paper; in particular \cite{fan2018eigenvector} used $r$ to denote our $d_i$ and used $d$ to denote our $D$. 

For $\muu^{(i)\top}(\mi-\muu^{(j)}\muu^{(j)\top})\mathbf{E}^{(j)}\muu^{(j)}$,
we have
$$
\begin{aligned}
\muu^{(i)\top}(\mi-\muu^{(j)}\muu^{(j)\top})\mathbf{E}^{(j)}\muu^{(j)}
=\muu^{(i)\top}(\mi-\muu^{(j)}\muu^{(j)\top})\hat\mSigma^{(j)}\muu^{(j)}
=\frac{1}{n_j}\sum_{k=1}^{n_j}\big(\muu^{(i)\top}a_k\big)b_k^\top\end{aligned}
$$
where we define $a_k := (\mi-\muu^{(j)}\muu^{(j)\top}) X^{(j)}_k\in\mathbb R^{D},b_k := \muu^{(j)\top}X^{(j)}_k\in\mathbb R^{d_j}$, and $\{X^{(i)}_k\}_{k\in[n_i]}$ are $n_i$ iid mean-zero Gaussian random vectors in $\mathbb{R}^D$ with covariance matrix $\mSigma^{(i)}$.
Therefore $\{a_k\}_{k\in[n_j]}$ are iid mean-zero Gaussian vectors with $\|a_k\|_{\psi_2}\lesssim\sigma_j\lesssim 1$, and $\{b_k\}_{k\in[n_j]}$ are iid mean-zero Gaussian vectors with $\|b_k\|_{\psi_2}\lesssim (\lambda_1^{(j)})^{1/2}\asymp D^{\gamma/2}$. 
Notice that for each $k\in[n_j]$, we have
$\mathbb E\big[(\muu^{(i)\top}a_k)b_k^\top\big]=\muu^{(i)\top}\operatorname{Cov}(a_k,b_k)=\mathbf 0.$
For any fixed unit vectors $u\in\mathbb S^{d_i-1}$ and $v\in\mathbb S^{d_j-1}$, consider
$$u^\top\Big(\frac1{n_j}\sum_{k=1}^{n_j}(\muu^{(i)\top}a_k)b_k^\top\Big)v
=\frac1{n_j}\sum_{k=1}^{n_j}\big((\muu^{(i)}u)^\top a_k\big)\big(b_k^\top v\big).$$
For any $k\in[n_i]$, since $(\muu^{(i)}u)^\top a_k$ is sub-Gaussian with $\|(\muu^{(i)}u)^\top a_k\|_{\psi_2}\lesssim\sigma_j\lesssim 1$, and $b_k^\top v$ is sub-Gaussian with $\|b_k^\top v\|_{\psi_2}\lesssim D^{\gamma/2}$, $\big((\muu^{(i)}u)^\top a_k\big)\big(b_k^\top v\big)$ is sub-exponential with
$$\big\|\big((\muu^{(i)}u)^\top a_k\big)\big(b_k^\top v\big)\big\|_{\psi_1}
\le \big\|(\muu^{(i)}u)^\top a_k\big\|_{\psi_2}\big\|b_k^\top v\big\|_{\psi_2}
\lesssim D^{\gamma/2}.$$
Since $\{\big((\muu^{(i)}u)^\top a_k\big)\big(b_k^\top v\big)\}_{k\in[n_i]}$ are independent, mean-zero, sub-exponential random variables. Then by Theorem~2.9.1 (Sub-exponential Bernstein inequality) in \cite{vershynin2018high}, for any $t\geq 0$, we have
$$\mathbb{P}\left(\Big|\sum_{k=1}^{n_j}\big((\muu^{(i)}u)^\top a_k\big)\big(b_k^\top v\big)\Big|\geq t\right)
\le 2\exp\left(-c\min\Big\{\frac{t^2}{n_i D^{\gamma}} \frac{t}{D^{\gamma/2}}\Big\}\right)$$
for some constant $c>0$.
Therefore we have
$$
\Big|\frac1{n_j}\sum_{k=1}^{n_j}\big((\muu^{(i)}u)^\top a_k\big)\big(b_k^\top v\big)\Big|
\lesssim  n_i^{-1/2}D^{\gamma/2}\log^{1/2} D+n_i^{-1}D^{\gamma/2}\log D
\lesssim D^{\gamma/2}\tilde\varphi
$$
with high probability.
Let $\mathcal N_u$ and $\mathcal N_v$ be $\tfrac14$-nets of $\mathbb S^{d_i-1}$ and $\mathbb S^{d_j-1}$. By Corollary~4.2.11 (Covering numbers of the Euclidean ball) in \cite{vershynin2018high}, we have  $|\mathcal N_u|\le 9^{d_i}\lesssim 1$ and $|\mathcal N_v|\le 9^{d_j}\lesssim 1$. 
Then by Lemma~4.4.2 (Maximizing quadratic forms on a net) in \cite{vershynin2018high}, we have
$$\|\muu^{(i)}(\mi-\muu^{(j)}\muu^{(j)\top})\mathbf{E}^{(j)}\muu^{(j)}\|=\Big\|\frac1{n_j}\sum_{k=1}^{n_j}(\muu^{(i)\top}a_k)b_k^\top\Big\|
\leq  2\max_{u\in\mathcal N_u,\,v\in\mathcal N_v}\Big|\frac1{n_j}\sum_{k=1}^{n_j}\big((\muu^{(i)}u)^\top a_k\big)\big(b_k^\top v\big)\Big|
\lesssim D^{\gamma/2}\tilde\varphi$$
with high probability.
\end{proof}

\begin{proof}[Proof of Lemma~\ref{lemma:PCA_UhatW-U_previous}]
\begin{sloppypar}
For simplicity of notation, we will omit the superscript $``(i)"$ from the matrices such as $\mathbf{U}^{(i)},\hat{\mathbf{U}}^{(i)},\bm{\Sigma}^{(i)},\hat{\bm{\Sigma}}^{(i)},\bm{\Lambda}^{(i)},\hat{\bm{\Lambda}}^{(i)},\mathbf{W}^{(i)},\mathbf{E}^{(i)},\mathbf{T}^{(i)}$ and the subscript $i$ from notations such as $d_i, \sigma_i$ as it should cause minimal confusion.
From Lemma~\ref{lemma:PCA_|E|_previous} and Weyl's inequality, we have
$\lambda_1(\hat\mSigma)\asymp\lambda_d(\hat\mSigma)\asymp D^\gamma$ with high
probability.
Therefore, by the Davis-Kahan theorem \citep{yu2015useful,davis70}, we have
\begin{equation}\label{Eq:PCA_||Uperp^T hatU||_previous}
\begin{aligned}
	\|(\mi - \muu \muu)^\top\hat\muu\|=\|\sin\Theta(\hat\muu,\muu)\|
	&\leq \frac{C\|\me\|}{\lambda_d(\hat\mSigma)-\lambda_{d+1}(\mSigma)}
	\lesssim \varphi
\end{aligned}
\end{equation}
with high probability. As $\mw$ is the solution of orthogonal
Procrustes problem between $\hat\muu$ and $\muu$,
we have
\begin{equation}
  \label{Eq:PCA_||U^T hat U-W^T||,||hat U-UW^T||_previous}
	\begin{aligned}
	&\|\muu^\top\hat\muu-\mw^\top\|
	\leq \|\sin\Theta(\hat\muu,\muu)\|^2
	\lesssim \varphi^2, \\
	&\|\hat\muu-\muu\mw^\top\|
	\leq \|\sin\Theta(\hat\muu,\muu)\|+\|\muu^\top\hat\muu-\mw^\top\|
	\lesssim \varphi
  \end{aligned}
\end{equation}
with high probability.
\end{sloppypar}

Define the matrices $\mt_1$ through $\mt_4$ by
$$
\begin{aligned}
	&\mt_1=\muu(\muu^\top\hat\muu-\mw^\top)\mw,\\
	&\mt_2=\sigma^2(\mi - \muu \muu^{\top})\hat\muu\hat\mLambda^{-1}\mw,\\
	&\mt_3=-\muu\muu^\top\me(\hat\muu-\muu\mw^\top)\hat\mLambda^{-1}\mw,\\
	&\mt_4=-\muu\muu^\top\me\muu(\mw^\top\hat\mLambda^{-1}\mw-\mLambda^{-1}).
\end{aligned}
$$
Then for $\hat\muu\mw-\muu$, we have the decomposition
\begin{equation}\label{Eq:PCA_T123_previous}
\begin{aligned}
	\hat\muu\mw-\muu
	&=(\mi-\muu\muu^\top)\hat\mSigma\hat\muu\hat\mLambda^{-1}\mw+\mt_1\\
	&=(\mi-\muu\muu^\top)\me\hat\muu\hat\mLambda^{-1}\mw+\mt_1+\mt_2\\
	&=\me\hat\muu\hat\mLambda^{-1}\mw-\muu\muu^\top\me\muu\mLambda^{-1}+\mt_1+\mt_2+\mt_3+\mt_4.
\end{aligned}	
\end{equation}

The spectral norms of $\mt_1$, $\mt_2$ and $\mt_3$ can be bounded by
\begin{equation}\label{Eq:PCA_T123_bound_previous}
	\begin{aligned}
		&\|\mt_1\|
		\leq \|\muu^\top\hat\muu-\mw^\top\|
		\lesssim \varphi^2,  \\
		&\|\mt_2\|
		\leq \sigma^2\|(\mi - \muu \muu^{\top})\hat\muu\|\cdot \|\hat\mLambda^{-1}\|
		\lesssim 
		D^{-\gamma} \varphi, \\
		&\|\mt_3\|
		\leq \|\me\|\cdot\|\hat\muu-\muu\mw^\top\|\cdot\|\hat\mLambda^{-1}\|
		\lesssim \varphi^2 \\
	\end{aligned}
\end{equation}
with high probability. For $\mt_4$ we have
\begin{equation}
  \label{eq:mt4_pca_expand}
  \begin{split}
    \mt_4 &= -\muu\muu^\top\me\muu\mLambda^{-1}\bigl[\mLambda
    \mw^{\top} - \mw^{\top} \hat{\mLambda} \bigr]
    \hat{\mLambda}^{-1} \mw 
    \\ &=-\muu\muu^\top\me\muu\mLambda^{-1}[\mLambda(\mw^\top-\muu^\top\hat\muu)-\muu^\top\me\hat\muu+(\muu^\top\hat\muu-\mw^\top)\hat\mLambda]\hat\mLambda^{-1}\mw,\\
\end{split}
\end{equation}
which implies
\begin{equation}
  \label{eq:mt4_spectral}
  \begin{split}
		\|\mt_4\|
		 \leq\|\me\| \cdot \big((\|\mLambda^{-1}\|+\|\hat\mLambda^{-1}\|)\|\muu^\top\hat\muu-\mw^\top \|+\|\mLambda^{-1}\|\cdot\|\hat\mLambda^{-1}\|\cdot\|\me\|\big)
		\lesssim \varphi^2
         \end{split}
\end{equation}
with high probability.
We now bound the $2 \to \infty$ norms of $\mt_1$ through $\mt_4$.
Recall that, from the assumption in
Theorem~\ref{thm:PCA_UhatW-U=(I-UU^T)EVR^{-1}+xxx_previous} we have
  $\|\muu\|_{2 \to \infty} \lesssim d^{1/2} D^{-1/2}$. 
As
$\mt_1, \mt_3$ and $\mt_4$ all include $\muu$ as the first term in the
matrix products, we have
\begin{equation}\label{Eq:PCA_T123_bound2_previous}
	\begin{aligned}
	 &\|\mt_1\|_{2\to\infty}
		 \leq \|\muu\|_{2\to\infty}\cdot\|\muu^\top\hat\muu-\mw^\top\|
		 \lesssim d^{1/2} D^{-1/2} \varphi^2, \\
	&\|\mt_3\|_{2\to\infty}
		 \leq \|\muu\|_{2\to\infty}\cdot \|\me\|\cdot\|\hat\muu-\muu\mw^\top\|\cdot\|\hat\mLambda^{-1}\|
		 \lesssim d^{1/2} D^{-1/2} \varphi^2, \\
		&\|\mt_4\|_{2\to\infty}
		 \leq \|\muu\|_{2\to\infty}\cdot\|\me\| \cdot \big((\|\mLambda^{-1}\|+\|\hat\mLambda^{-1}\|)\|\muu^\top\hat\muu-\mw^\top \|+\|\mLambda^{-1}\|\cdot\|\hat\mLambda^{-1}\|\cdot\|\me\|\big)
		 \lesssim d^{1/2} D^{-1/2} \varphi^2
	\end{aligned}
\end{equation}
with high probability.
Bounding $\|\mt_2\|_{2 \to \infty}$ requires slightly more
effort. 
Let $\bm{\Pi}_{\muu} = \muu \muu^{\top}$ and
$\piu = \mi - \muu \muu^{\top}$. Then
\begin{equation*}
  \begin{split}
    \mt_2 &= \sigma^2 \piu \hat{\muu} \hat{\mLambda}^{-1} \mw
    = \sigma^2 \piu \hat{\mSigma} \hat{\muu} \hat{\mLambda}^{-2} \mw 
    = \sigma^2 \piu (\me + \mSigma) \hat{\muu} \hat{\mLambda}^{-2} \mw
\\ & = \bigl(\sigma^2 \piu \me + \sigma^4 \piu\bigr)
\hat{\muu} \hat{\mSigma}^{-2} \mw 
= \bigl(\sigma^2 \me \bm{\Pi}_{\muu} + 
\sigma^2 \me \piu - \sigma^2
\bm{\Pi}_{\muu} \me + \sigma^4 \piu\bigr) \hat{\muu} \hat{\mLambda}^{-2} \mw.
  \end{split}
\end{equation*}
We now have, by Lemma~\ref{lemma:PCA_|E|_previous} and the condition
$n = \omega(D^{2 - 2 \gamma} \log D)$ that
\begin{gather*}
\|\me \piu \hat{\muu} \hat{\mLambda}^{-2} \mw \|_{2 \to
  \infty} \leq \|\me\|_{\infty} \cdot \|\piu \hat{\muu}
\hat{\mLambda}^{-1}\|_{2 \to \infty} \cdot \|\hat{\mLambda}^{-1}\|
\lesssim D^{1 - \gamma} \tilde{\varphi} \|\mt_2\|_{2 \to \infty}
= o(\|\mt_2\|_{2 \to \infty}), \\
\|\piu \hat{\muu} \hat{\mLambda}^{-2} \mw\|_{2 \to \infty} \leq \|\piu
\hat{\muu} \hat{\mLambda}^{-1}\|_{2 \to \infty} \cdot
\|\hat{\mLambda}^{-1}\| \lesssim \|\mt_2\|_{2 \to \infty} \cdot
\|\hat{\mLambda}^{-1}\| = o(\|\mt_2\|_{2 \to \infty})
\end{gather*}
We therefore have
\begin{equation}
  \label{eq:mt2_2inf}
		\|\mt_2\|_{2\to\infty}
		 \leq (1 + o(1)) \sigma^2\bigl(\|\me\muu\|_{2\to\infty} +\|\muu\|_{2\to\infty}\cdot\|\me\|\bigr)\|\hat\mLambda^{-1}\|^2
	\lesssim d^{1/2}D^{-3\gamma/2}\tilde\varphi
  \end{equation}
  with high probability.
From Lemma~\ref{lemma:PCA_|E|_previous}, we know the spectra of $\hat\mLambda$ and $\me$ are disjoint from one another with high probability, therefore $\hat\muu$ has a von Neumann series expansion as 
\begin{equation*}
	\begin{aligned}
		\hat\muu
		=\sum_{k=0}^{\infty}\me^k\mSigma\hat\muu\hat\mLambda^{-(k+1)}
		=\sum_{k=0}^{\infty}\me^k\muu\mLambda\muu^\top\hat\muu\hat\mLambda^{-(k+1)}
		+\sigma^2\sum_{k=0}^{\infty}\me^k(\mi-\muu\muu^\top)\hat\muu\hat\mLambda^{-(k+1)}
	\end{aligned}
\end{equation*}
with high probability. 
Suppose the above series expansion for $\hat{\muu}$ holds
and define the matrices
\begin{equation*}
	\begin{aligned}
      &\mt_5=\me\muu(\mw^\top\hat\mLambda^{-1}\mw-\mLambda^{-1}) = \me\muu\mLambda^{-1}[\mLambda(\mw^\top-\muu^\top\hat\muu)-\muu^\top\me\hat\muu+(\muu^\top\hat\muu-\mw^\top)\hat\mLambda]\hat\mLambda^{-1}\mw, \\
	    &\mt_6=\me\muu(\muu^\top\hat\muu-\mw^\top)\hat\mLambda^{-1}\mw,\\
	    &\mt_7=\me\muu\mLambda(\muu^\top\hat\muu\hat\mLambda^{-1}-\mLambda^{-1}\muu^\top\hat\muu)\hat\mLambda^{-1}\mw
        = -\me\muu\muu^\top\me\hat\muu\hat\mLambda^{-2}\mw,
 \\
	    &\mt_8=\sum_{k=2}^{\infty}\me^k\muu\mLambda\muu^\top\hat\muu\hat\mLambda^{-(k+1)}\mw,\\
		&\mt_9=\sigma^2\sum_{k=1}^{\infty}\me^k(\mi-\muu\muu^\top)\hat\muu\hat\mLambda^{-(k+1)}\mw.
	\end{aligned}
\end{equation*}
Note that the second expression for $\mt_5$ is similar to that
for Eq.~\eqref{eq:mt4_pca_expand}.
We then have
\begin{equation}\label{Eq:PCA_T4-8_previous}
	\begin{aligned}
		\me\hat\muu\hat\mLambda^{-1}\mw
		=&\me\muu\mLambda^{-1}
		+\mt_5
		+\mt_6
		+\mt_7
		+\mt_8
		+\mt_9.
	\end{aligned}
\end{equation}
Using Lemma~\ref{lemma:PCA_|E|_previous}, Eq.~\eqref{Eq:PCA_||Uperp^T
  hatU||_previous} and Eq.~\eqref{Eq:PCA_||U^T hat U-W^T||,||hat
  U-UW^T||_previous}, the spectral norms of $\mt_5$ through $\mt_9$
can be bounded by
\begin{equation}\label{Eq:PCA_T4-8_bound_previous}
	\begin{aligned}
		&\|\mt_5\|
		\leq\|\me\| \cdot \big((\|\mLambda^{-1}\|+\|\hat\mLambda^{-1}\|)\|\muu^\top\hat\muu-\mw^\top \|+\|\mLambda^{-1}\|\cdot\|\hat\mLambda^{-1}\|\cdot\|\me\|\big)
		\lesssim \varphi^2, \\
		&\|\mt_6\|
		\leq \|\me\|\cdot\|\muu^\top\hat\muu-\mw^\top\|\cdot\|\hat\mLambda^{-1}\|
		\lesssim \varphi^{3},\\
		&\|\mt_7\|
		\leq \|\me\|^2\cdot \|\hat\mLambda^{-1}\|^2
		\lesssim \varphi^2, \\
		&\|\mt_8\|
		\leq \sum_{k=2}^{\infty}\|\me\|^k\cdot \|\mLambda\|\cdot\|\hat\mLambda^{-1}\|^{k+1}
		\lesssim \varphi^2 ,\\
		&\|\mt_9\|
		\leq \sigma^2 \sum_{k=1}^{\infty}\|\me\|^k \cdot\|(\mi - \muu \muu^{\top})\hat\muu\|\cdot\|\hat\mLambda^{-1}\|^{k+1}
		\lesssim D^{-\gamma}\varphi^2
	\end{aligned}
\end{equation}
with high probability. Furthermore, the $2 \to \infty$ norm for
$\mt_5$ through $\mt_9$ can be bounded by
\begin{equation}\label{Eq:PCA_T4-8_bound2_previous}
	\begin{aligned}
	 &\|\mt_5\|_{2\to\infty}
		 \leq\|\me\muu\|_{2\to\infty}\cdot \big((\|\mLambda^{-1}\|+\|\hat\mLambda^{-1}\|)\|\muu^\top\hat\muu-\mw^\top \|+\|\mLambda^{-1}\|\cdot\|\hat\mLambda^{-1}\|\cdot\|\me\|\big)
		 \lesssim d^{1/2}D^{-\gamma/2} \varphi\tilde\varphi,\\
	 &\|\mt_6\|_{2\to\infty}
		 \leq \|\me\muu\|_{2\to\infty}\cdot\|\muu^\top\hat\muu-\mw^\top\|\cdot\|\hat\mLambda^{-1}\|
		 \lesssim d^{1/2}D^{-\gamma/2}\varphi^2\tilde\varphi,\\
	 &\|\mt_7\|_{2\to\infty}
		 \leq \|\me\muu\|_{2\to\infty}\cdot\|\me\|\cdot \|\hat\mLambda^{-1}\|^2
		 \lesssim d^{1/2} D^{-\gamma/2} \varphi\tilde\varphi,\\
		&\|\mt_8\|_{2\to\infty}
		 \leq \sum_{k=2}^{\infty}\|\me\|_\infty^{k-1}\cdot \|\me\muu\|_{2\to\infty}\cdot \|\mLambda\|\cdot\|\hat\mLambda^{-1}\|^{k+1}
		 \lesssim d^{1/2}D^{1-3\gamma/2}\tilde\varphi^2, 
\\ 		&\|\mt_9\|_{2\to\infty}
 \leq\sum_{k=1}^{\infty}\|\me\|_\infty^k\cdot\|\sigma^2(\mi-\muu\muu^\top)\hat\muu\hat\mLambda^{-1}\|_{2\to\infty}\cdot\|\hat\mLambda^{-1}\|^k
 \lesssim d^{1/2}D^{1-5\gamma/2}\tilde\varphi^2
	\end{aligned}
\end{equation}
with high probability. N
ote that bounds for $\|\mathbf{T}_8\|_{2 \to \infty}$ and $\|\mathbf{T}_9\|_{2 \to \infty}$ require $n = \omega(D^{2 - 2\gamma} \log D)$; in contrast, bounds for $\|\mathbf{T}_5\|_{2 \to \infty}, \|\mathbf{T}_6\|_{2 \to \infty}$, and $\|\mathbf{T}_7\|_{2 \to \infty}$ require the weaker assumption $n = \omega(\max\{D^{1 - \gamma}, \log D\})$. Furthermore the bound for $\|\mathbf{T}_9\|_{2 \to \infty}$ also uses the bound for $\|\mathbf{T}_2\|_{2 \to \infty}$ derived earlier in the proof.

Recall Eq.~\eqref{Eq:PCA_T123_previous} and
Eq.~\eqref{Eq:PCA_T4-8_previous}, and define $\mt = \mt_1 + \mt_2 +
\dots + \mt_9$. The bounds for $\|\mt\|$ and $\|\mt\|_{2 \to \infty}$
in Lemma~\ref{lemma:PCA_UhatW-U_previous} follow directly from
Eq.~\eqref{Eq:PCA_T123_bound_previous}, Eq.~\eqref{eq:mt4_spectral}, Eq.~\eqref{Eq:PCA_T123_bound2_previous}, Eq.~\eqref{eq:mt2_2inf}, Eq.~\eqref{Eq:PCA_T4-8_bound_previous} 
and Eq.~\eqref{Eq:PCA_T4-8_bound2_previous}.
\end{proof}

\subsection{Proof of technical lemmas for Theorem~\ref{thm:PCA_UhatW-U=(I-UU^T)EVR^{-1}+xxx}}
\label{Appendix:proof of lemma:PCA_chen}

\begin{proof}[Proof of Lemma~\ref{lemma:PCA_|E|_basis}]
    Recall that $\me_{k \ell}^{(i)}$ is distributed $\mathcal{N}(0,
    \sigma_i^2/n)$ for $k \in [D], \ell \in [n]$ and $i \in [m]$. 
    By the tail bound for a Gaussian random variable, we have
	$$
	\begin{aligned}
		\max_{k\in[D], \ell \in[n]}|\me^{(i)}_{kl}|
		\lesssim \frac{\sigma_i\log^{1/2}(n+D)}{n^{1/2}}
	\end{aligned}
	$$
	with probability at least $1-O((n+D)^{-10})$.
	As $\muu^{(i)}$ and $\muu^{\natural(i)}$ represent the same column space
    for $\mathbf{X}^{(i)}$, there exists an orthogonal matrix
    $\mw^{\natural(i)}$
    such that $\muu^{(i)}=\muu^{\natural(i)}\mw^{\natural(i)}$ and hence
    $$
    \begin{aligned}
    	\|\muu^{\natural (i)}\|_{2\to\infty}
    	= \|\muu^{(i)}\|_{2\to\infty}
    	\lesssim d_i^{1/2}D^{-1/2}.
    \end{aligned}
    $$
	Finally by Lemma~6 in \cite{yan2021inference} we have, under the
    assumption $\frac{\log(n+D)}{n}\lesssim 1$, that
	$$
	\begin{aligned}
		\mSigma^{\natural(i)}_{rr}\asymp D^{\gamma/2}\quad \text{for any }r\in[d_i]
		 \quad\text{and}\quad
		 \|\mv^{\natural(i)}\|_{2\to\infty}
		 \lesssim \frac{d_i^{1/2}\log^{1/2}(n+D)}{n^{1/2}}
	\end{aligned}
	$$
	with probability at least $1-O((n+D)^{-10})$.
\end{proof}
\begin{proof}[Proof of Lemma~\ref{lemma:PCA_|E|}]
Let $c > 0$ be fixed but arbitrary. 
Then by applying Theorem~3.4 in \cite{chen2021spectral} there exists a
constant $C(c)$ depending only on $c$ such that
	$$
	\begin{aligned}
		\mathbb{P}\Big(\|\me^{(i)}\| \geq C(c) \frac{\sigma_i(n+D)^{1/2}}{n^{1/2}}+t\Big)
		\leq (n+D)\exp\Big(-\frac{ct^2n}{\sigma_i^2\log(n+D)}\Big). 
	\end{aligned}
	$$
    We can thus set $t=C\sigma_i\bigl(1 + D/n\bigr)^{1/2}$
	for some universal constant $C$ not depending on $n$
    and $D$ (provided that $n \geq \log D$) 
    such that
    $$
	 \begin{aligned}
&\|\me^{(i)}\| \lesssim \sigma_i \Big (1 + \frac{D}{n}\Big)^{1/2}
\end{aligned}
	$$
	with probability as least $1-O((n+D)^{-10})$.

	For $\muu^{(i)\top}\me^{(j)}\mv^{\natural(j)}$, we notice 
    $$
    \begin{aligned}
    \muu^{(i)\top}\me^{(j)}\mv^{\natural(j)}=\sum_{k=1}^n\sum_{\ell=1}^n\me^{(j)}_{k\ell} u_k^{(i)} v^{\natural(j)\top}_\ell
    \text{ and }
    \me^{(j)}_{k\ell} u_k^{(i)} v^{\natural(j)\top}_\ell
    =n^{1/2}\sigma_j^{-1}\me^{(j)}_{k\ell} \cdot \mb^{(i,j;k,\ell)},
    \end{aligned}
    $$
    where $u_k^{(i)}$ denotes the $k$th row of $\muu^{(i)}$, $v^{\natural(j)}_\ell$ denotes the $\ell$th row of $\mv^{\natural(j)}$, and $\mb^{(i,j;k,\ell)}=n^{-1/2}\sigma_j u_k^{(i)} v^{\natural(j)\top}_\ell$. 
    Note that $\{n^{1/2}\sigma_j^{-1}\me^{(j)}_{k\ell}\}_{k\in[D],\ell\in[n]}$ are independent standard normal random variables.
    Let $\mathcal{A}$ be the event $\{\|\mv^{\natural(j)\top}\|_{2\to\infty}\lesssim d_j^{1/2}n^{-1/2}\log^{1/2}(n+D)\}$. Then by Lemma~\ref{lemma:PCA_|E|_basis}, we have $\mathbb{P}(\mathcal{A})\geq 1-O((n+D)^{-10})$.
    Next suppose that $\mathcal{A}$ holds.
    Then from $\mb^{(i,j;k,\ell)}(\mb^{(i,j;k,\ell)})^{\top}
    =n^{-1} \sigma_j^2\|v^{\natural(j)\top}_\ell\|^2u_k^{(i)}  u_k ^{(i)\top}$ and Weyl's inequality, we have
    $$
    \begin{aligned}
    	\Big\|\sum_{k=1}^D\sum_{\ell=1}^n \mb^{(i,j;k,\ell)}(\mb^{(i,j;k,\ell)})^{\top}\Big\|
    	&\leq n^{-1} \sigma_j^2\sum_{\ell=1}^n\|v^{\natural(j)\top}_\ell\|^2\Big\|\sum_{k=1}^D u_k^{(i)}  u_k ^{(i)\top}\Big\|\\
       	&\leq n^{-1}\sigma_j^2
    	\cdot n\|\mv^{\natural(j)\top}\|_{2\to\infty}^2
    	\cdot \|\muu^{(i)\top}\muu^{(i)}\|\\
    	&\lesssim n^{-1}\sigma_j^2
    	\cdot n(d_j^{1/2}n^{-1/2}\log^{1/2}(n+D))^2 
    	\cdot 1
    	\lesssim \sigma_j^2d_jn^{-1}\log(n+D).
    \end{aligned}
    $$
    Similarly, we also have
    $$
    \Big\|\sum_{k=1}^n\sum_{\ell=1}^n (\mb^{(i,j;k,\ell)})^{\top}\mb^{(i,j;k,\ell)}\Big\|
    	\lesssim \sigma_j^2d_in^{-1}.
    $$
    Hence, by Theorem~1.5 in \cite{tropp2012user}, there exist a constant $C>0$ such that for all $t>0$ 
    we have
    $$
    \begin{aligned}
    	\mathbb{P}\Big\{\|\muu^{(i)\top}\me^{(j)}\mv^{\natural(j)}\|\geq t\Big\}
    	&\leq (n+D)\cdot \exp \Big(\frac{-t^2/2}{2C\max\{\sigma_j^2d_jn^{-1}\log(n+D),\sigma_j^2d_in^{-1}\}}\Big),
    \end{aligned}
    $$
    from which we obtain
    $$
    \begin{aligned}
    	\|\muu^{(i)\top}\me^{(j)}\mv^{\natural(j)\top}\|
    	\lesssim \sigma_j d_j^{1/2}n^{-1/2}\log(n+D)
    \end{aligned}
    $$
    with probability as least $1-O((n+D)^{-10})$. Finally we unconditioned on the event $\mathcal{A}$ to obtain the desired upper bound for $\|\muu^{(i)\top}\me^{(j)}\mv^{\natural(j)}\|$.

    For $\me^{(i)}\mv^{\natural(i)}$, we notice
$
\|\me^{(i)}\mv^{\natural(i)}\|_{2\to\infty}	=\max_{k\in[n]} \|\big(\me^{(i)}\mv^{\natural(i)}\big)_k\|,
$
where $\big(\me^{(i)}\mv^{\natural(i)}\big)_k$ represents the $k$th row of
$\big(\me^{(i)}\mv^{\natural(i)}\big)$. Similarly, by Theorem~1.5 in \cite{tropp2012user}, we
have 
$\|\big(\me^{(i)}\mv^{\natural(i)}\big)_k\|\lesssim \sigma d^{1/2}n^{-1/2}\log (n+D)$ with probability as least $1-O((n+D)^{-10})$. In summary we have
$
	\|\me^{(i)}\mv^{\natural(i)}\|_{2\to\infty}\lesssim \sigma d^{1/2}n^{-1/2}\log (n+D)
$
with probability as least $1-O((n+D)^{-10})$.
\end{proof}

\begin{proof}[Proof of Lemma~\ref{lemma:PCA_UhatW-U}]
Recall Lemma~\ref{lemma:PCA_|E|_basis}. In particular we have
    $$
    \|\muu^{\natural (i)}\|_{2\to\infty}
    \lesssim \sqrt{\frac{\mu^\natural d}{D}}, \quad
    \|\mv^{\natural (i)}\|_{2\to\infty}
    \lesssim \sqrt{\frac{\mu^\natural d}{n}}
    $$
    where $\mu^\natural =1+\log(n+D)$, and furthermore the
    $\me_{k \ell}^{(i)}$ are independent random variables with
    $$
    \mathbb{E}(\me^{(i)}_{kl})=0,\quad
	\max_{k \ell} \operatorname{Var}(\me^{(i)}_{k \ell})\leq \tilde\sigma^2,\\
	\quad |\me^{(i)}_{kl}| \lesssim B
    $$
    with probablity at least $1 - O((n+D)^{-10})$; here $\tilde\sigma_i^2=\frac{\sigma_i^2}{ n},B=\sqrt{\frac{\sigma_i^2\log(n+D)}{n}}$.
    
    Then, by Theorem~9 in \cite{yan2021inference}, we have
    \begin{equation*}
    \hat\muu^{(i)}\mw^{(i)}-\muu^{(i)}=\me^{(i)}\mv^{\natural(i)}(\mSigma^{\natural(i)})^{-1}\mw^{\natural(i)}+\mt^{(i)}
    \end{equation*}
    where $\mt^{(i)}$ satisfies
    $$
    \begin{aligned}
    	\|\mt^{(i)}\|_{2\to\infty}
    	&\lesssim \frac{\sigma_i^2d_i^{1/2}(n+D)^{1/2}\log(n+D)}{nD^\gamma}
    	+\frac{\sigma_i^2d^{1/2}(n+D)}{nD^\gamma D^{1/2}}+\frac{\sigma_i
          d_i\log^{1/2}(n+D)}{n^{1/2}D^{(1 + \gamma)/2}}
    \end{aligned}
    $$
    with probability at least $1-O((n+D)^{-10})$, provided that
    \begin{equation}
    \begin{aligned}
      \label{eq:condition_pca_chen}
    	&\frac{\sigma_i \log^{1/2}(n+D)}{n^{1/2}}
    	\lesssim \sigma_i \sqrt{\frac{\min\{n,D\}}{n(1+\log(n+D))\log(\max\{n,D\})}},\\
    	&\sigma_i \sqrt{\frac{\max\{n,D\}\log(\max\{n,D\})}{n}}
    	 \ll D^{\gamma/2}.
    \end{aligned}
    \end{equation}
    The conditions in Eq.~\eqref{eq:condition_pca_chen} follows from
    the conditions
    $$ \begin{aligned}
    	\frac{\log^3(n+D)}{\min\{n,D\}}\lesssim 1
    	\quad\text{and}\quad
    	\frac{(n+D)\log(n+D)}{nD^\gamma}\ll 1.
    \end{aligned}
    $$
    stated in Theorem~\ref{thm:PCA_UhatW-U=(I-UU^T)EVR^{-1}+xxx}.
\end{proof}

\subsection{Additional discussion for Section~\ref{sec:related_works}}

\label{sec:related_works_technical}

{\color{black}
We now discuss the results of \cite{ma2026optimal} on the ``aggregate-then-estimate'' approach in the presence of individual subspaces. \cite{ma2026optimal} considers a matrix denoising framework $\mathbf{A}^{(i)} = \mathbf{P}^{(i)} + \mathbf{E}^{(i)}$, where $\mathbf{E}^{(i)}$ has zero-mean sub-Gaussian entries and each signal matrix has a partially shared left singular subspace, i.e., $\mathbf{P}^{(i)} = [\muu_c \mid \muu_s^{(i)}]\mSigma^{(i)}\mv^{(i)\top}$ with $\mSigma^{(i)}$ a diagonal matrix of singular values. Thus, apart from the different noise model, this can be viewed as a special case of the COISIE model. \cite{ma2026optimal} studies the stacked-SVD estimator, which takes the leading left singular subspace of the stacked matrix $[\ma^{(1)} \mid \cdots \mid \ma^{(m)}]$; note that, for symmetric networks, this coincides with the leading eigenvectors of $\sum_{i=1}^m (\ma^{(i)})^2$.
\cite{ma2026optimal} shows that stacked-SVD achieves minimax rate-optimality when the singular subspaces are fully shared, i.e., when only $\muu_c$ is present and there are no individual subspaces $\{\muu_s^{(i)}\}$. Once individual subspaces are present, however, stacked-SVD can easily fail; for example, it fails when the unshared signals are strong. 
Moreover, even when stacked-SVD is effective, the analysis of \cite{ma2026optimal} relies on the relatively restrictive assumption that the individual subspaces $\{\muu_s^{(i)}\}$ are all mutually orthogonal, which is not required by the COISIE model.
}

{\color{black}
As mentioned in Section~\ref{sec:related_works}, for degree-corrected multiple networks with layer-specific degree heterogeneity, the ``estimate-then-aggregate'' approach can accommodate the degree heterogeneity within each individual estimate, whereas ``aggregate-then-estimate'' approaches fail. We now provide a detailed discussion of this setting, focusing on the degree-corrected SBM. For degree-corrected SBMs, theoretical results similar to our results were derived in \cite{agterberg_dcsbm}. In particular, Theorem~A.3 in \cite{agterberg_dcsbm} is closely related to our Theorem~\ref{thm:UhatW-U=EVR^{-1}+xxx}. The main difference is that the notations $\hat{\muu}$ and $\muu$ in Theorem~A.3 of \cite{agterberg_dcsbm} refer to matrices that are related to but distinct from those in our paper.

The algorithm studied in \cite{agterberg_dcsbm} is also an ``estimate-then-aggregate'' algorithm. More specifically, \cite{agterberg_dcsbm} parametrizes their degree-corrected SBMs as $\mpp^{(i)} = \bm{\Theta}^{(i)} \mathbf{Z} \mathbf{B}^{(i)} \mathbf{Z}^{\top} \bm{\Theta}^{(i)}$. Let $\hat{\my}^{(i)}$ denote the matrix whose rows are the {  normalized} rows of $\hat{\muu}^{(i)}$, where $\hat{\muu}^{(i)}$ contains the leading eigenvectors of $\ma^{(i)}$. Then their $\hat{\muu}$ is the matrix of leading eigenvectors of $\sum_{i} \hat{\my}^{(i)} \hat{\my}^{(i)\top}$; see Algorithm~1 in \cite{agterberg_dcsbm} for details. Similarly, let $\my^{(i)}$ denote the matrix whose rows are the {  normalized} rows of $\muu^{(i)}$, where $\muu^{(i)}$ contains the leading eigenvectors of $\mpp^{(i)}$. Their $\muu$ is then the matrix of leading eigenvectors of $\sum_{i} \my^{(i)} \my^{(i)\top}$. Proposition~2.1 of \cite{agterberg_dcsbm} then showed that $\muu = \mathbf{Z} \mathbf{M}$ for some matrix $\mathbf{M}$ with {  distinct} rows. This allows one to use $\muu$ to recover the common community assignments induced by $\mathbf{Z}$, and hence, as $\hat{\muu}$ is a consistent estimate of $\muu$, this property also extends to $\hat{\muu}$. Finally, given an estimate of $\mathbf{Z}$, one can then estimate the matrices $\mathbf{B}^{(i)}$ for all $i \in [m]$ under suitable identifiability conditions.

We further note that the layer-specific $\bm\Theta^{(i)}$ in the above degree-corrected SBM setting render ``aggregate-then-estimate'' approaches ineffective, as unlike the COSIE model with $\mpp^{(i)} = \muu \mr^{(i)} \muu^{\top}$, the leading subspace of each $\mpp^{(i)}$ is now the column space of $\bm\Theta^{(i)}\muu$, which varies across layers with $\bm\Theta^{(i)}$, and hence aggregating $\{\ma^{(i)}\}$ prior to estimation, for instance via the leading eigenvectors of $\sum_{i=1}^{m}(\ma^{(i)})^2$, conflates these layer-specific subspaces and fails to recover the shared $\muu$. In contrast, because the ``estimate-then-aggregate'' approach processes each layer separately, it can incorporate a per-layer step that removes the degree correction before aggregation, as exemplified by Algorithm~1 in \cite{agterberg_dcsbm}.
}

We now compare our theoretical results with existing works on multilayer SBMs.
Under this regime of multilayer SBM, by combining our result in Theorem~\ref{thm:WhatU_k-U_k->norm} with the same argument as that for showing exact recovery in a single SBM (see, e.g., Theorem~2.6 in \cite{lyzinski2014perfect} or Theorem~5.2 in \cite{lei2019unified}), one can also show that $K$-means or $K$-medians clustering on the rows of $\hat{\muu}$ will, asymptotically almost surely, exactly recover the community assignments $\tau$ in a multilayer SBM provided that $n\rho_n = \omega(\log n)$ as $n \rightarrow \infty$. The condition $n \rho_n = \omega(\log n)$ almost matches the lower bound $n \rho_n = \Omega(\log n)$ for exact recovery in single-layer SBMs in \cite{abbe2015exact,mossel2015consistency,abbe2020entrywise}. Note, however, that \cite{abbe2015exact,mossel2015consistency,abbe2020entrywise} only consider the case of balanced SBMs where the block probabilities $\mathbf{B}$ satisfy $\mathbf{B}_{kk} \equiv p$ and $\mathbf{B}_{k \ell} \equiv q$ for all $k \neq \ell$.
Some existing works provide Frobenius norm estimation errors of $\hat{\muu}$ which only guarantee weak recovery of the community assignment $\tau$.
For example, \cite{paul2020spectral} studies community detection using two different procedures, namely a linked matrix factorization procedure (as suggested in \cite{tang2009clustering}) and a co-regularized spectral clustering procedure (as suggested in \cite{kumar2011co}), and they show that if $m n \rho_n = \omega(\log n)$ then the estimation error bounds of $\muu$ for these two procedures are
\begin{gather*}
\min_{\mw\in \mathcal{O}_d}
\|\hat{\mathbf{U}}\mw-\mathbf{U} \|_{F}
\lesssim
    d^{1/2} m^{-1/8}(\log m)^{1/4} (n\rho_n)^{-1/4}
    \log^{1+\epsilon/2}n,  \\
\min_{\mw\in \mathcal{O}_d}
\|\hat{\mathbf{U}}\mw-\mathbf{U} \|_{F}
\lesssim d^{1/2}m^{-1/4}(n\rho_n)^{-1/4}\log^{1/4+\epsilon} n
\end{gather*}
with high probability, where $\epsilon > 0$ is an arbitrary but fixed constant. See the proofs of Theorem~2 and Theorem~3 in \cite{paul2020spectral} for more details.
As another example, \cite{jing2021community} proposes a tensor-based algorithm for estimating $\mz$ in a mixture multilayer SBM model and shows that if $m n \rho_n = \omega(\log^{4}n)$ then
$$\min_{\mw\in \mathcal{O}_d} \|\hat{\mathbf{U}}\mw-\mathbf{U} \|_{F} \lesssim d^{1/2}m^{-1/2}(n\rho_n)^{-1/2}\log^{1/2}n$$
with high probability; see the condition in Corollary~1 and the proof of Theorem~5.2 in \cite{jing2021community} for more details.
If $m$ is bounded by a finite constant not depending on $n$ (as assumed in the setting of our paper), the bound in Proposition~\ref{thm:Vhat-VW} is $d^{1/2}m^{-1/2}(n\rho_n)^{-1/2}$ and is thus either equivalent to or quantitatively better than those cited above while our assumption $n \rho_n = \Omega(\log n)$ is also the same or weaker than those cited above. As discussed above regarding the differences between the two types of methods, if $m$ grows with $n$ then the above cited results allow for possibly smaller thresholds of $n \rho_n$ while still guaranteeing consistency.
Finally, \cite{lei2020bias} considers the sparse regime with $n \rho_n \leq C_0$ for some constant $C_0 > 0$ not depending on $m$ and $n$, and proposes estimating $\muu$ using the leading eigenvectors of $\sum_{i=1}^m (\mathbf{A}^{(i)})^2 - \mathbf{D}^{(i)}$ where, for each $i \in [m]$, $\mathbf{D}^{(i)}$ denotes the diagonal matrix whose diagonal entries are the vertex degrees in $\mathbf{A}^{(i)}$; the subtraction of $\mathbf{D}^{(i)}$ corresponds to a bias-removal step and is essential as the diagonal entries of $\sum_{i=1}^m (\ma^{(i)})^2$ are heavily biased when the graphs are extremely sparse. Let $\hat{\muu}_{b}$ denote the matrix containing these eigenvectors. Theorem~1 in \cite{lei2020bias} shows that if $m^{1/2}n\rho_n \gg \log^{1/2}(m+n)$ then
$$
\min_{\mathbf{W} \in \mathcal{O}(d)} \|\hat{\mathbf{U}}_{b}\mathbf{W} -\mathbf{U}\|_{F}
\lesssim d^{1/2} \Bigl[m^{-1/2}(n\rho_n)^{-1}\log^{1/2}(m+n)+n^{-1}\Bigr]
$$
with high probability. The above results for Frobenius norm estimation errors of either $\hat{\muu}$ or $\hat{\muu}_{b}$ only guarantee weak recovery of the community assignment $\tau$. More refined error bounds in $2 \to \infty$ norm for estimating $\muu$ in the "aggregate-then-estimate" setting, which also lead to exact recovery of $\tau$, are discussed in the following.

For $2 \to \infty$ norm bounds for estimating $\muu$ using ``aggregate-then-estimate'' approaches, \cite{cai2021subspace} studies subspace estimation for unbalanced matrices such as $\mathbf{A}_* = [\mathbf{A}^{(1)}\mid \dots \mid \mathbf{A}^{(m)}]$ by using the $d$ leading eigenvectors of $\mathcal{P}_{\mathrm{off\_diag}}(\mathbf{A}_* (\mathbf{A}_*)^{\top}) = \sum_{i=1}^{m} \mathcal{P}_{\mathrm{off\_diag}}((\mathbf{A}^{(i)})^2)$ where $\mathcal{P}_{\mathrm{off\_diag}}(\cdot)$ zeros out the diagonal entries of a matrix and thus serves the same purpose as the subtraction of $\mathbf{D}^{(i)}$ in \cite{lei2020bias}. Let $\tilde{\muu}_{b}$ denote the resulting leading eigenvectors. Now suppose $n \rho_n = O(1)$ and $m^{1/2}n\rho_n \gg \log (mn)$. Then by Theorem~1 in \cite{cai2021subspace} we have
\begin{equation}
  \label{eq:cai_inference}
    \min_{\mw \in \mathcal{O}_d} \|\tilde{\mathbf{U}}_b \mathbf{W} -\mathbf{U}\|_{2\to\infty}
    \lesssim d^{1/2}n^{-1/2}
    [m^{-1/2}(n\rho_n)^{-1}\log(mn)
    +dn^{-1}]
\end{equation}
    with high probability. See Section~4.3 in \cite{cai2021subspace} and subsection~\ref{sec:justification1} below for more details; note that the discussion in Section~4.3 of \cite{cai2021subspace} assumes that $\mathbf{A}_*$ is the adjacency matrix for a bipartite graph but the same argument generalizes to the multilayer SBM setting. Eq.~\eqref{eq:cai_inference} implies that clustering the rows of $\tilde{\muu}_{b}$ achieves exact recovery of $\tau$.

    The above $2 \to \infty$ norm bound can be further refined using
    results in \cite{yan2021inference} wherein the diagonal entries of
    $\sum_{i=1}^{m} (\mathbf{A}^{(i)})^2$ are iteratively imputed while computing its
    truncated eigendecompositions. In particular, let
    $\mathbf{G}^{(0)} = \sum_{i=1}^{m}
    \mathcal{P}_{\mathrm{off\_diag}}((\mathbf{A}^{(i)})^2)$ and let
    $t_{\max} \geq 0$ be a non-negative integer. Then, for $0 \leq t
    < t_{\max}$, set $\mathbf{G}^{(t+1)} =
    \mathcal{P}_{\mathrm{off\_diag}}(\mathbf{G}^{(t)}) +
    \mathcal{P}_{\mathrm{diag}}(\mathbf{G}^{(t)}_{d})$ where $\mathbf{G}^{(t)}_{d}$
    is the best rank-$d$ approximation to $\mathbf{G}^{(t)}$, and $\mathcal{P}_{\mathrm{diag}}(\cdot)$ denotes the operation which zeros out the {\em off-diagonal} entries of a matrix.
Let $\tilde{\muu}_b^{(t_{\max})}$ denote the leading eigenvectors of
    $\mathbf{G}^{(t_{\max})}$ (the estimate $\tilde{\muu}_{b}$ in Eq.~\eqref{eq:cai_inference} 
    corresponds to the case $t_{\max} = 0$). Also let
    $\mathbf{U}^{\natural} \bm{\Sigma}^{\natural} \mathbf{V}^{\natural}$ denote
    the SVD of $\mathbf{P}_* =
    [\mathbf{P}^{(1)} \mid \dots \mid
    \mathbf{P}^{(m)}]$ and denote $\mathbf{E}_* = [\mathbf{E}^{(1)} \mid
    \dots \mid \mathbf{E}^{(m)}]$. Once again suppose
    $n \rho_n = O(1)$, $m^{1/2} n \rho_n \gg \log (mn)$, and choose
    $t_{\max} \gg\log(mn\rho_n)$. Then by 
    Theorem~10 in \cite{yan2021inference}, there exists 
    $\mw_{\muu} \in \mathcal{O}_d$ such that
    \begin{equation}
      \label{eq:yan_inference}
     \tilde{\muu}^{(t_{\max})}_{b} \mw_{\muu} - \muu = \mathbf{E}_* \mv^{\natural}
    (\bm{\Sigma}^{\natural})^{-1} +
    \mathcal{P}_{\mathrm{off\_diag}}(\mathbf{E}_* \mathbf{E}_*^{\top})
    \mathbf{U}^{\natural} (\bm{\Sigma}^{\natural})^{-2} + \mathbf{Q}_b,
    \end{equation}
    where $\mathbf{Q}_b$ satisfies
    $$\|\mathbf{Q}_b\|_{2 \to \infty} \lesssim dn^{-1} m^{-1/2}(n\rho_n)^{-1}\log(mn)
+ d^{1/2}n^{-1/2}m^{-1}(n\rho_n)^{-2}\log^2(mn)
$$
with high probability; see subsection~\ref{sec:justification2} below for more details.
Eq.~\eqref{eq:yan_inference} also yields a normal approximation for the
    rows of $\tilde{\muu}^{(t_{\max})}_b$ but with more complicated
    covariance matrices than those given in
    Theorem~\ref{thm:WhatU_k-U_k->norm}; we leave the 
    precise form of these covariance matrices to the interested
    reader.

 \subsubsection{Technical details for Eq.~\eqref{eq:cai_inference}}
 \label{sec:justification1}
 We can take the matrix $\ma$
and $\ma_*$ in Section~3
of \cite{cai2021subspace} as
$$\ma = [\ma^{(1)} \mid \ma^{(2)} \mid
\dots \mid \ma^{(m)}], \quad \ma^\star = [\mpp^{(1)} \mid \mpp^{(2)} \mid
\dots \mid \mpp^{(m)}].$$
The dimensions $d_1$ and $d_2$ of $\ma$ and $\ma^\star$ are then $d_1 = n$
and $d_2 = mn$ where $m$ and $n$ denote the number of graphs and number of
vertices (as used in this paper). The rank of $\ma^\star$ in
\cite{cai2021subspace} is denoted by $r$ and corresponds to the
notation $d$ in this paper (note that $d$ in \cite{cai2021subspace} denotes $\max\{d_1,d_2\}$ and
corresponds to $mn$ in this paper). 
Now suppose that $m$ increases with
$n$ in such a way that
\begin{equation}
  \label{eq:cai2021_condition_grow}
  m^{1/2} n \rho_n= \Omega(\log (mn)).
\end{equation}
Then $\ma$ and $\ma^\star$ satisfy Assumption~1 and Assumption~2 in
\cite{cai2021subspace} with $p = 1$, $\sigma = \rho_n^{1/2}$ and 
$\|\mathbf{N}\|_{\max} = \|\ma - \ma^\star\|_{\max} \leq R$ almost surely
where $R = 1$. In particular
Eq.\eqref{eq:cai2021_condition_grow} above implies
Eq.~(10) in \cite{cai2021subspace}. Note that while \cite{cai2021subspace} assumes the entries
of $\me$ are to be mutually independent, their results
still hold for the setting considered here
where, due to the symmetric of the $\ma^{(i)}$ for each $i
\in [m]$, any two rows of $\ma - \ma^\star$ share one entry in common. 

Now let $\sigma^\star_j$ denote the $j$th largest singular
values of $\ma^\star$. Then $(\sigma_r^\star)^2$ is the smallest non-zero
eigenvalue of $\ma^\star (\ma^\star)^{\top} = \sum_{i=1}^{m} \mathbf{U} (\mathbf{R}^{(i)})^2
\mathbf{U}^{\top}$ and thus under mild conditions on $\sum_{i=1}^m
(\mathbf{R}^{(i)})^2$, we have $\sigma_j^\star \asymp m^{1/2} (n
\rho_n)$ for all $j\in[r]$, and furthermore $\ma^\star$ has bounded condition number (which is denoted by $\kappa$ in \cite{cai2021subspace}).
$\ma^\star$ also has bounded incoherence parameter (which is denoted by $\mu$ in \cite{cai2021subspace}).
Under the assumption $m^{1/2} n \rho_n\gg \log (mn)$ the above
quantities $p, \sigma, \sigma_r^\star,\kappa,\mu,d_1,d_2,d,r$ satisfy Eq.~(15) in
\cite{cai2021subspace}.
Now let $n \rho_n = O(1)$,
i.e., each $\mathbf{A}^{(i)}$ has bounded average degree. 
The quantity $\mathcal{E}_{\mathrm{general}}$ in Eq.~(17) of
\cite{cai2021subspace} is then
$$\mathcal{E}_{\mathrm{general}} \asymp \frac{\rho_n}{m (n \rho_n)^2}
\times (m^{1/2} n \log(mn)) + \frac{\rho_n}{m^{1/2} n \rho_n}\times
  (n \log (mn))^{1/2} + \frac{d}{n} \asymp \frac{\log(mn)}{m^{1/2} (n
    \rho_n)} + \frac{d}{n}.$$
   Define $\mw$ as a
  minimizer of $\|\hat{\muu}_b \mathbf{O} - \muu\|_{F}$ over all orthogonal matrix $\mo$.
  Therefore, by Eq.(16b) of Theorem~1 in \cite{cai2021subspace}, 
  there exists an
  orthogonal $\mw$ such that
  $$\|\hat{\muu}_b \mw - \muu\|_{2 \to \infty} \lesssim
d^{1/2}n^{-1/2} [m^{-1/2} (n
  \rho_n)^{-1} \log(mn) + dn^{-1}]
  $$
  with high probability, which is the bound in
  Eq.~\eqref{eq:cai_inference}. Note that
  $\muu$ and $\muu^\star$ in \cite{cai2021subspace} correspond to
  $\hat{\muu}_b$ and $\muu$ in
  this paper, respectively.  

  \subsubsection{Technical details for Eq.~\eqref{eq:yan_inference}}
  \label{sec:justification2}
 Using the notations in
Section~6.2 of \cite{yan2021inference}, we can take
$\mathbf{M}^{\natural} = [\mathbf{P}^{(1))} \mid \mathbf{P}^{(2)} \dots
\mid \mathbf{P}^{(m)}]$, $\mathbf{M} = [\mathbf{A}^{(1)} \mid
\mathbf{A}^{(2)} \mid \dots \mathbf{A}^{(m)}]$, $n_1 = n$, $n_2 =
mn$ where $m$ and $n$ denote the number of graphs and number of
vertices in this paper (note that $n$ in \cite{yan2021inference} denotes $\max\{n_1,n_2\}$ and corresponds to $mn$ in this paper). The rank of $\mathbf{M}^{\natural}$ in \cite{yan2021inference} is denoted by $r$ and corresponds to the notation $d$ in this paper. 
Once again suppose that $m$ increases with
$n$ in such a way that Eq.~\eqref{eq:cai2021_condition_grow} is satisfied. 
Then $\mathbf{M}^{\natural}$ and $\mathbf{E} = \mathbf{M} -
\mathbf{M}^{\natural}$ satisfy the conditions in Assumption~4 and
Assumption~5 of \cite{yan2021inference} with $\sigma = \rho_n^{1/2}$ and $B
= 1$; once again, while \cite{yan2021inference} also assumes that the
entries of $\me$ are independent, their results still hold for the
setting discussed here where the $\ma^{(i)}$ are symmetric matrices. 

Now let $\sigma^{\natural}_j$ denote the $j$th largest singular
values of $\mathbf{M}^{\natural}$. Similar to the above discussion for
the singular values $\sigma_r^\star$ in 
\cite{cai2021subspace}, we also have $(\sigma^{\natural}_j) \asymp m^{1/2} (n
\rho_n)$ for all $j \in [r]$, and furthermore $\mathbf{M}^{\natural}$ has bounded condition number (which is denoted by $\kappa^{\natural}$ in \cite{yan2021inference}).
$\mathbf{M}^{\natural}$ also has bounded incoherence parameter (which is denoted by $\mu^{\natural}$ in \cite{yan2021inference}).
Let $n \rho_n = O(1)$,
i.e., each $\mathbf{A}^{(i)}$ has bounded average degree. 
The quantity $\zeta_{\mathrm{op}}$ in
Eq.~(6.16) of \cite{yan2021inference} is then 
$$\zeta_{\mathrm{op}} \asymp \rho_n m^{1/2} n \log {(mn)} + \rho_n^{1/2}
m^{1/2} (n \rho_n) (n \log (mn))^{1/2} \asymp m^{1/2} (n  \rho_n) \log
(mn),$$
and furthermore $\zeta_{\mathrm{op}}$ satisfies the condition in
Eq.(6.17) of \cite{yan2021inference} under the assumption $m^{1/2} n \rho_n\gg \log (mn)$.
In particular, $m^{1/2} n \rho_n\gg \log (mn)$ implies
$$\frac{(\sigma_r^{\natural})^2}{\zeta_{\mathrm{op}}} \asymp \frac{m^{1/2} n
  \rho_n}{\log (mn)} \gg 1.$$
Letting $t_{\max} \geq \log ((\sigma_1^{\natural})^2/\zeta_{\mathrm{op}})$
we have $$
\begin{aligned}
	t_{\max} 
	&\geq \log \Big(\frac{\sigma_1^{2}(\sum_{i=1}^m(\mr^{(i)})^2)}{\max(\mathbb{E}[\me_{st}^2])m^{1/2}n\log(mn)+\max^{1/2}(\mathbb{E}[\me_{st}^2])\sigma_1(\sum_{i=1}^m(\mr^{(i)})^2)n^{1/2}\log^{1/2}(mn)}\Big)\\
	&\gtrsim \log(mn\rho_n)- \log(\log(mn)).
\end{aligned}
$$ The conditions in Theorem~10 of
\cite{yan2021inference} are satisfied under the assumptions $n \rho_n = O(1)$, $n\gtrsim \log^2(mn)$ and $m^{1/2} n \rho_n\gg \log (mn)$, and we thus obtain the
expansion for $\hat{\muu}^{(t_{\max})}_b$ as given in
Eq.\eqref{eq:yan_inference}. This corresponds to Eq.~(6.19a) of
\cite{yan2021inference} where their $\muu$ is our
$\hat{\muu}_b^{(t_{\max})}$, their $\muu^\natural$ is our $\muu$, and their
$\bm{\Psi}$ is our $\mq_{b}$. The bound for $\mq_{b}$ in
Section~\ref{sec:related_works} is then given by Eq.(6.19b) of
\cite{yan2021inference}, i.e.,
$$
\begin{aligned}
	\|\mq_{b}\|_{2 \to \infty} 
&\lesssim \frac{d}{n} \times 
\frac{\zeta_{\mathrm{op}}}{m(n\rho_n)^2} +
\frac{\zeta_{\mathrm{op}}^2}{m^2(n\rho_n)^4} \times
\frac{d^{1/2}}{n^{1/2}} \\
&\lesssim dn^{-1} m^{-1/2}(n\rho_n)^{-1}\log(mn)
+ d^{1/2}n^{-1/2}m^{-1}(n\rho_n)^{-2}\log^2(mn)
\end{aligned}
$$
with high probability.

{\color{black}

\subsection{Additional discussion for Section~\ref{sec:related_works_pca}}

\label{sec:related_works_technical_pca}

We now discuss the connections between our row-wise CLTs for distributed PCA and the eigenvector/eigenvalue central limit theorems for spiked covariance models under the $p/n \to c$ asymptotic regime.
For eigenvector/eigenvalue CLTs under the spiked covariance model, we are aware of several works such as \cite{paul2007,wang_fan_2017,fan_li_xia_zheng}. In particular, suppose we have $n$ iid mean zero random vectors $X_1, X_2, \dots, X_n$ in $\mathbb{R}^{p}$ with covariance matrix $\bm{\Sigma}$. Note that $p$ here corresponds to the dimension $D$ in the notation of our paper. Suppose also that $\bm{\Sigma}$ has eigenvalues $\lambda_1 \geq \lambda_2 \geq \dots \geq \lambda_p$ where $\lambda_1, \dots, \lambda_r$ are the {\em spiked} eigenvalues. \cite{paul2007} considers the setting where the $X_i$ are multivariate Gaussians and $p/n \rightarrow c \in (0,1)$, \cite{wang_fan_2017} considers the setting of sub-Gaussian random vectors and $p/(n \lambda_j)$ being bounded, and \cite{fan_li_xia_zheng} considers the setting where $\mathbb{E}[\|X_i\|^4] < \infty$ and $p/n \rightarrow c \in (0, \infty)$ with spiked eigenvalues satisfying $\lambda_j \asymp p$ for $j = 1,2,\dots,r$. Now assume, for simplicity, that all the spiked eigenvalues are {\em well-separated}, that is to say, there exists a constant $c_0 > 0$ such that $\lambda_j/\lambda_{j+1} \geq 1 + c_0$ for $j = 1,2,\dots,r$. Then, under appropriate regularity conditions, all of the above cited papers yield a CLT of the form
$$\sqrt{n}(\hat{\lambda}_j/\lambda_j - 1) \rightarrow \mathcal{N}(0, \sigma_j^2),\quad j = 1,2,\dots,r,$$
where the asymptotic variances $\{\sigma_j^2\}$ are known.

For the eigenvectors, \cite{paul2007,wang_fan_2017,fan_li_xia_zheng} first transform the $X_i$ to $\bm{\Gamma}^{\top} X_i$ where $\bm{\Gamma} \bm{\Lambda} \bm{\Gamma}^{\top}$ is the {\em eigendecomposition} of $\bm{\Sigma}$. Now let $\hat{\bm{\xi}}_j$ denote the eigenvector of the sample covariance matrix $\bm{\Gamma}^{\top} \hat{\bm{\Sigma}} \bm{\Gamma}$. Then \cite{paul2007,wang_fan_2017,fan_li_xia_zheng} show that
\[ \sqrt{n}\Bigl( \frac{\hat{\bm{\xi}}_{j,\leq r}}{\|\hat{\bm{\xi}}_{j,\leq r}\|} - \bm{e}_{j} \Bigr) \rightarrow \mathcal{N}(\bm 0, \bm{\Upsilon}_j), \quad j = 1,2,\dots,r. \]
Here $\hat{\bm{\xi}}_{j,\leq r} \in \mathbb{R}^{r}$ is the vector formed by keeping only the first $r$ elements of $\hat{\bm{\xi}}_j$, $\bm{e}_j$ is the $j$th elementary basis vector of $\mathbb{R}^{r}$, and the $r \times r$ covariance matrices $\{\bm{\Upsilon}_j\}$ are known. In addition, \cite{paul2007,wang_fan_2017,fan_li_xia_zheng} also show that $\hat{\bm{\xi}}_{j,>r}/\|\hat{\bm{\xi}}_{j,>r}\|$ is {\em uniformly} distributed on the unit sphere in $\mathbb{R}^{p - r}$.

From the above discussion we see that the eigenvalue CLTs are complementary to but distinct from the {\em row-wise} eigenvector CLT presented in the paper. The eigenvector CLTs of \cite{paul2007,wang_fan_2017,fan_li_xia_zheng} are also related to the row-wise eigenvector CLT in our paper, but there is an important difference in that the results of \cite{paul2007,wang_fan_2017,fan_li_xia_zheng} are for {\em all of the entries} of the {\em transformed} eigenvectors $\hat{\bm{\xi}}_j$. In other words, the eigenvectors $\hat{\bm{u}}_1, \dots, \hat{\bm{u}}_r$ (using the notation in our paper) are given by $\hat{\bm{u}}_j = \muu \hat{\bm{\xi}}_{j,\leq r} + \muu_{\perp} \hat{\bm{\xi}}_{j, > r}$ for $j = 1,2,\dots,r$ (see Eq.~(3.7) in \cite{wang_fan_2017}). Thus \cite{paul2007,wang_fan_2017,fan_li_xia_zheng} characterize the distribution of $\hat{\bm{u}}_j$ in the {\em transformed (eigenbasis) coordinates} $\hat{\bm{\xi}}_j$, corresponding to how $\hat{\bm{u}}_j$ decomposes along the signal and noise subspaces. In contrast, our CLTs characterize the {\em entrywise} behavior of $\hat{\bm{u}}_j$ in the original coordinates. We surmise that one can also extend the theoretical results and analysis in \cite{wang_fan_2017} to obtain a {\em row-wise CLT} similar to that obtained in our paper, although this would likely require additional analysis. Indeed, the results in \cite{paul2007,wang_fan_2017,fan_li_xia_zheng} are for traditional and not distributed PCA.

Finally, regarding whether our row-wise CLTs can be adapted to the $p/n \to c$ regime, we note that the row-wise CLT in Theorem~\ref{thm:PCA_(UhatW-U)_k->normal} allows the dimension $D$ (which, recall, corresponds to $p$ in the above discussion) to be comparable to the sample size $n$. In particular, Theorem~\ref{thm:PCA_(UhatW-U)_k->normal} permits $D^{1-\gamma} \ll n \ll D^{1+\gamma}$, so that when $D \asymp n$ this already covers the $p/n \to c$ regime. We also note that our results are established in the regime where the spiked eigenvalues grow with the dimension, i.e., $\lambda_j \asymp D^{\gamma}$ for some $\gamma \in (0,1]$ (see Theorem~\ref{thm:PCA_UhatW-U=(I-UU^T)EVR^{-1}+xxx_previous} to Theorem~\ref{thm:PCA_(UhatW-U)_k->normal}). 
In other regimes, such as the $p/n \to c$ regime with bounded spiked eigenvalues (i.e., $\lambda_j \asymp 1$), our current results would not apply directly; establishing analogous row-wise limit results in that regime would require a different analysis, which we leave as an interesting direction for future work.
}

\subsection{Equivalence between
Theorem~\ref{thm:PCA_UhatW-U=(I-UU^T)EVR^{-1}+xxx_previous} and 
               Theorem~\ref{thm:PCA_UhatW-U=(I-UU^T)EVR^{-1}+xxx}}
\label{sec:equivalence_thm7_thm9}

We now show that the leading terms in Theorem~\ref{thm:PCA_UhatW-U=(I-UU^T)EVR^{-1}+xxx_previous} and Theorem~\ref{thm:PCA_UhatW-U=(I-UU^T)EVR^{-1}+xxx} are
equivalent, and thus the main difference between Theorem~\ref{thm:PCA_UhatW-U=(I-UU^T)EVR^{-1}+xxx_previous} and Theorem~\ref{thm:PCA_UhatW-U=(I-UU^T)EVR^{-1}+xxx} is in bounding the residual terms, i.e., Theorem~\ref{thm:PCA_UhatW-U=(I-UU^T)EVR^{-1}+xxx_previous} analyzes the leading eigenvectors
               of $\hat{\bm{\Sigma}}^{(i)} = \mathbf{X}^{(i)}
               (\mathbf{X}^{(i)})^{\top}$ whose entries are {\em
               dependent} while Theorem~\ref{thm:PCA_UhatW-U=(I-UU^T)EVR^{-1}+xxx} analyzes the leading left
               singular vectors of $\mathbf{X}^{(i)}$ whose entries
               are {\em independent}.

\begin{sloppypar}
	For Theorem~\ref{thm:PCA_UhatW-U=(I-UU^T)EVR^{-1}+xxx_previous}, the leading order term for $\hat\muu_c\mw_{\muu_c}-\muu_c$  can be simplified as
    \begin{equation*}
        \begin{split}
        \frac{1}{m} \sum_{i=1}^m (\mi - \muu^{(i)} \muu^{(i)\top}) (\hat{\mSigma}^{(i)} - \mSigma^{(i)}) \muu_c^{(i)} (\mLambda_c^{(i)})^{-1}
        &= \frac{1}{m} \sum_{i=1}^m(\mi - \muu^{(i)} \muu^{(i)\top}) \hat{\mSigma}^{(i)} \muu_c(\mLambda_c^{(i)})^{-1} \\ 
        &= \frac{1}{mn} \sum_{i=1}^m (\mi - \muu^{(i)} \muu^{(i)\top})   (\my^{(i)} + \mz^{(i)}) (\my^{(i)} + \mz^{(i)})^{\top} \muu_c (\mLambda_c^{(i)})^{-1} \\
        &= \frac{1}{mn} \sum_{i=1}^m (\mi - \muu^{(i)} \muu^{(i)\top}) \mz^{(i)} (\my^{(i)} + \mz^{(i)})^{\top} \muu_c (\mLambda_c^{(i)})^{-1},
        \end{split}
    \end{equation*}
    where the last equality is because $\my^{(i)} = \muu^{(i)} (\bm{\Lambda}^{(i)} - \sigma^2 \mi)^{1/2} \mathbf{F}^{(i)}$ and $(\mi - \muu^{(i)} \muu^{(i)\top})\muu^{(i)}=\mathbf{0}$.
    Now, $(\mi - \muu^{(i)} \muu^{(i)\top}) \mz^{(i)} (\mz^{(i)})^{\top} \muu_c$ is a $D \times d_0$ matrix whose $rs$th entry, which we denote as $\zeta_{rk}$, is of the form
    $$\sum_{j=1}^{n} \bm{n}_{r}^{\top} Z_j^{(i)} (Z_j^{(i)})^{\top} \bm{u}_{c,k},$$
    where $\bm{n}_{r}$ is the $r$th row of $(\mi - \muu^{(i)} \muu^{(i)\top})$, $\bm{u}_{c,s}$ is the $k$th column of $\muu_c$, and $Z_j^{(i)}$ are iid $\mathcal{N}(\mathbf{0}, \sigma^2 \mi)$ random vectors. Since $$\mathbb{E}[\bm{n}_{r}^{\top} Z_j^{(i)} (Z_j^{(i)})^{\top} \bm{u}_{c,k}] = \sigma_i^2\bm{n}_r^{\top} \bm{u}_{c,k} = 0,$$ $\zeta_{rk}$ is a sum of independent mean $0$ random variables. Furthermore, as $\|\bm{n}_r\| \leq 1$ and $\|\bm{u}_{c,k}\| \leq 1$,
    $\bm{n}_{r}^{\top} Z_j^{(i)} (Z_j^{(i)})^{\top} \bm{u}_{k}$ is a sub-exponential random variable with Orlicz-1 norm bounded by $\sigma^2$ (see Lemma~2.7.7 in \cite{vershynin2018high}). We therefore have, by a standard application of Bernstein's inequality (see e.g., Theorem~2.8.1 of \cite{vershynin2018high}), that
    $|\zeta_{rk}| \precsim (n \log n)^{1/2}$
    with high probability. Therefore
    $$\Bigl\|\frac{1}{mn} \sum_{i=1}^{m}(\mi - \muu^{(i)} \muu^{(i)\top}) \mz^{(i)} \mz^{(i)\top} \muu_c (\mLambda_c^{(i)})^{-1}\Bigr\|_{2 \to \infty} \precsim \frac{d_0^{1/2} \log^{1/2}{n}}{n^{1/2}} \times \|(\mLambda_c^{(i)})^{-1}\|\lesssim d_0^{1/2}n^{-1/2} D^{-\gamma} \log^{1/2}{n}$$
    with high probability,
    and will thus be negligible as $n, D$ increase. 
    Next, we also have
    \begin{equation*}
    \begin{split}
        \|\muu^{(i)} \muu^{(i)\top} \mz^{(i)} \my^{(i)\top} \muu _c(\mLambda_c^{(i)})^{-1}\|_{2 \to \infty} 
        &\leq \|\muu^{(i)}\|_{2 \to \infty} \times \|\muu^{(i)\top} \mz^{(i)} \my^{(i)\top} \muu_c (\mLambda_c^{(i)})^{-1}\| \\ & \leq  \|\muu^{(i)}\|_{2 \to \infty} \times \|\muu ^{(i)\top}\mz^{(i)} \mathbf{F}^{(i)\top}\| \times \|(\mLambda_c^{(i)})^{-1/2}\|
        \\ & \precsim \|\muu^{(i)}\|_{2 \to \infty} \times d_i(n \log n)^{1/2} \times \|(\mLambda_c^{(i)})^{-1/2}\|
        \end{split}
    \end{equation*}
    with high probability; the final inequality in the above display follows from the fact that $\muu^{(i)\top} \mz^{(i)} (\mathbf{F}^{(i)})^{\top}$ is a $d_i \times d_i$ matrix whose $rk$th entries are of the form $(\xi_r^{(i)})^\top f_k^{(i)}$ where $\xi_r^{(i)}$ and $f_k^{(i)}$ are random vectors in $\mathbb{R}^{n}$ and their entries are independent with  bounded Orlicz-2 norms. We thus have
\begin{equation*}
      \begin{split}
      \Bigl\|\frac{1}{mn} \sum_{i=1}^{m} \muu^{(i)} \muu^{(i)\top} \mz^{(i)} \my^{(i)\top} \muu _c(\mLambda_c^{(i)})^{-1}\Bigr\|_{2 \to \infty} 
    \lesssim d_{\max}^{3/2}n^{-1/2}D^{-(1+\gamma)/2}\log^{1/2}n
      \end{split}
    \end{equation*}
    with high probability,
    which is also negligible as $n, D$ increase.
    In summary, the above chain of derivations yield the approximation
    \begin{equation}
    \label{eq:equivalence_PCA1}
        \begin{split}
     \frac{1}{m} \sum_{i=1}^m (\mi - \muu^{(i)} \muu^{(i)\top}) (\hat{\mSigma}^{(i)} - \mSigma^{(i)}) \muu^{(i)}_c (\mLambda_c^{(i)})^{-1}
     &= \frac{1}{mn} \sum_{i=1}^m \mz^{(i)} (\my^{(i)})^{\top} \muu_c (\mLambda_c^{(i)})^{-1} + \tilde{\mr} \\ &= \frac{1}{mn} \sum_{i=1}^m \mz^{(i)} (\mathbf{F}^{(i)})^{\top} (\mLambda^{(i)} - \sigma_i^2 \mi)^{1/2} \muu^{(i)\top}\muu_c(\mLambda_c^{(i)})^{-1} + \tilde{\mr},
                \end{split}
    \end{equation}  where $\tilde{\mr}$ is a $D \times d_0$ random matrix with negligible spectral and $2 \to \infty$ norms.
    \end{sloppypar}
    
    For the leading order term in Theorem~\ref{thm:PCA_UhatW-U=(I-UU^T)EVR^{-1}+xxx}, using the form for $\my^{(i)}$, we also have
    $$(\my^{(i)})^{\dagger} = (\mathbf{F}^{(i)})^{\dagger} (\bm{\Lambda}^{(i)} - \sigma_i^2 \mi)^{-1/2} \muu^{(i)\top} = \mathbf{F}^{(i)\top} (\mathbf{F}^{(i)} \mathbf{F}^{(i)\top})^{-1} (\bm{\Lambda}^{(i)} - \sigma^2 \mi)^{-1/2} \muu^{(i)\top}$$
    almost surely, provided that $n \geq d_i$. As $\mathbf{F}^{(i)}
               \mathbf{F}^{(i)\top}$ is $d_i \times d_i$ Wishart
               matrix, by Eq.~(5.11) in \cite{cai2022non} we can show that
               $n(\mathbf{F}^{(i)} \mathbf{F}^{(i)\top})^{-1} = 
                \mi + \tilde{\mr}_2$ where $\|\tilde{\mr}_2\| \lesssim n^{-1/2}\log^{1/2}n $ with high probability. We therefore have
    \begin{equation}
     \label{eq:equivalence_PCA2}
        \frac{1}{m} \sum_{i=1}^m \mz^{(i)} (\my^{(i)})^{\dagger} \muu_c
               =  \frac{1}{mn} \sum_{i=1}^m \mz^{(i)}
               \mathbf{F}^{(i)\top} (\mi + \tilde{\mr}_2) (\bm{\Lambda}^{(i)} - \sigma_i^2
               \mi)^{-1/2}\muu^{(i)\top}\muu_c
             \end{equation}
             almost surely.
    Now $\mz^{(i)} \mathbf{F}^{(i)\top}$ is a $D \times d_i$ matrix
    whose $rs$th entry are of the form
    $(z^{(i)}_{r})^\top  f^{(i)}_{ k }$ where $z^{(i)}_{r}$
               and $f^{(i)}_{ k }$ are random
               vectors in $\mathbb{R}^n$ and their entries are independent with bounded Orlicz-2 norms. We thus have $\|\mz^{(i)}
               \mathbf{F}^{(i)\top}\|_{2 \to \infty} \lesssim
               d_i^{1/2} n^{1/2} \log^{1/2}{n}$ with high probability,
               so that
               \begin{equation*}
                 \Bigl\|\frac{1}{mn} \sum_{i=1}^m \mz^{(i)}
               \mathbf{F}^{(i)\top} \tilde{\mr}_2 (\bm{\Lambda}^{(i)} - \sigma^2
               \mi)^{-1/2}\muu^{(i)\top}\muu_c\Bigr\|_{2 \to \infty} \lesssim d_{\max}^{1/2}
               n^{-1} D^{-\gamma/2} \log n
               \end{equation*}
               with high probability,
               which is negligible as $n, D \rightarrow \infty$. In
               summary the right hand side of
               Eq.~\eqref{eq:equivalence_PCA1} and
               Eq.~\eqref{eq:equivalence_PCA2} are the same, as when
               $D$ increases,
                $\sigma^2\mi$ is also
               negligible compared with $\mLambda^{(i)}$, and notice that $(\mLambda^{(i)})^{1/2}\muu^{(i)\top}\muu_c(\mLambda_c^{(i)})^{-1}=(\mLambda^{(i)})^{-1/2}\muu^{(i)\top}\muu_c$.
               Thus the expansion in
               Theorem~\ref{thm:PCA_UhatW-U=(I-UU^T)EVR^{-1}+xxx_previous}
                 is conceptually equivalent to that in
                 Theorem~\ref{thm:PCA_UhatW-U=(I-UU^T)EVR^{-1}+xxx},
               with the only difference being the analysis of the
               lower-order term $\mq_{\muu_c}$ due to the relationship between
               $n$ and $D$ (see Table~\ref{tb:PCA}). 


\end{sloppypar}

\end{document}